\documentclass[a4paper,12pt,frenchb]{article}

\usepackage{amsmath,amsbsy,amsfonts,amssymb,amsthm}
\usepackage[french]{babel}
\newcommand{\ass}[2]{\vskip0.3cm\noindent
{\bf {#1}}. { \sl {#2}}\vskip0.3cm\noindent
}

\begin{document}

  \title{  Repr\'esentations et quasi-caract\`eres de niveau $0$; endoscopie }
\author{J.-L. Waldspurger}
\date{30 octobre 2018}
 \maketitle
 
 {\bf Introduction}
 \bigskip
 
  Soit $F$ un corps local non archim\'edien de caract\'eristique nulle et soit $G$ un groupe r\'eductif connexe d\'efini sur $F$. On s'int\'eresse ici aux repr\'esentations admissibles et irr\'eductibles de $G(F)$ dans des espaces vectoriels complexes. On note $Irr(G)$ l'ensemble des classes d'isomorphismes de ces repr\'esentations. Dans le cas o\`u $G$ est un groupe classique, on sait associer \`a toute $\pi\in Irr(G)$ un param\`etre de Langlands, cela gr\^ace aux r\'esultats de Harris et Taylor, Henniart, Arthur et Mok. Dans le cas g\'en\'eral, l'existence de ce param\`etre reste conjecturale. 
 Moy et Prasad ont d\'efini le niveau de $\pi$. Notons $Irr(G)^0$ le sous-ensemble des repr\'esentations de niveau $0$. Dans le cas o\`u $G$ est adjoint, Lusztig a d\'efini un param\'etrage de l'ensemble $Irr(G)^0$, qui est un bon candidat pour \^etre celui de Langlands. L'hypoth\`ese que $G$ est adjoint a \'et\'e r\'ecemment lev\'ee par Solleveld. 
  Hormis le cas des groupes classiques, savoir quelles conditions caract\'erisent le param\'etrage de Langlands n'est pas clair, du moins pour l'auteur. Il y a en tout cas une condition minimale: le param\'etrage doit v\'erifier des conditions de compatibilit\'e \`a l'endoscopie. D'o\`u une question pr\'ealable que l'on formule ici en termes vagues:  le niveau $0$ se conserve-t-il par endoscopie?   Dans cet article, nous r\'epondons positivement \`a cette question, pourvu que la caract\'eristique r\'esiduelle $p$ de $F$ soit grande relativement \`a $G$ (ou que $G$ soit petit relativement \`a $p$, c'est une question de point de vue).

 Pour tout ensemble $X$, notons ${\mathbb C}[X]$ l'espace vectoriel complexe de base $X$. En associant \`a toute repr\'esentation son caract\`ere-distribution, on d\'efinit une injection  $\Theta:{\mathbb C}[Irr(G)] \to I(G)^*$, o\`u $I(G)^*$ est l'espace des distributions sur $G(F)$ invariantes par conjugaison. D\'efinissons le projecteur de Bernstein $p^0:{\mathbb C}[Irr(G)]\to {\mathbb C}[Irr(G)^0]$ qui annule toute repr\'esentation irr\'eductible qui n'est pas de niveau $0$. 
 
 Supposons d'abord que $G$ est quasi-d\'eploy\'e. On sait d\'efinir le sous-espace $SI(G)^*\subset I(G)^*$ des distributions stables. On note ${\mathbb C}[Irr(G)]^{st}$ le sous-espace des $\pi\in {\mathbb C}[Irr(G)]$ telles que $\Theta_{\pi}\in SI(G)^*$ (ici, $\pi$ n'est plus irr\'eductible, c'est une combinaison lin\'eaire \`a coefficients complexes de  repr\'esentations irr\'eductibles). Une cons\'equence de l'article  \cite{A3} est qu'il existe une projection naturelle $p^{st}:{\mathbb C}[Irr(G)]\to {\mathbb C}[Irr(G)]^{st}$. Notre premier r\'esultat est le
 
 \ass{Th\'eor\`eme 1}{ Supposons que $G$ soit quasi-d\'eploy\'e et que l'hypoth\`ese $(Hyp)_{endo}(G)$ soit v\'erifi\'ee. Alors on a l'\'egalit\'e $p^{st}\circ p^0=p^0\circ p^{st}$.}
 
Cf. le corollaire \ref{lecasquasideploye}. L'hypoth\`ese $(Hyp)_{endo}(G)$, pr\'ecis\'ement \'enonc\'ee en \ref{donneesendoscopiques}, est l'hypoth\`ese sur $p$ \'evoqu\'ee plus haut. Grosso-modo, elle suppose $p\geq c(G)val_{F}(p)$, o\`u $c(G)$ est un entier d\'ependant de $G$ et $val_{F}$ est la valuation usuelle de $F$.

Revenons au cas d'un groupe $G$ quelconque. Soit ${\bf G}'$ une donn\'ee endoscopique elliptique de $G$, cf.  \ref{donneesendoscopiques}. Une telle donn\'ee est un triplet dont l'un des termes est un groupe endoscopique $G'$, qui est quasi-d\'eploy\'e sur $F$. En g\'en\'eral, on doit fixer des donn\'ees auxiliaires pour d\'efinir un facteur de transfert. Pour simplifier l'introduction, supposons que cela ne soit pas n\'ecessaire et fixons un facteur de transfert d\'efini sur un sous-ensemble de $G'(F)\times G(F)$. Toujours en cons\'equence de \cite{A3}, il y a alors un homomorphisme de transfert spectral $transfert:{\mathbb C}[Irr(G')]^{st}\to {\mathbb C}[Irr(G)]$.

 \ass{Th\'eor\`eme 2}{ Supposons que  l'hypoth\`ese $(Hyp)_{endo}(G)$ soit v\'erifi\'ee. Alors on a l'\'egalit\'e $p^0\circ transfert=transfert\circ p^0$.}
 
 Cf. le th\'eor\`eme \ref{letheoreme}. Evidemment, le premier $p^0$ vit sur $G$ et le second sur $G'$. Remarquons que ce dernier conserve l'espace ${\mathbb C}[Irr(G')]^{st}$ d'apr\`es le premier th\'eor\`eme. 
 
 L'application $p^0$ est un projecteur de Bernstein et il s'en d\'eduit un tel projecteur sur divers objets, par exemple l'espace $C_{c}^{\infty}(G(F))$ des fonctions \`a valeurs complexes sur $G(F)$, localement constantes et \`a support compact. Les th\'eor\`emes ci-dessus ont des contreparties pour les espaces de fonctions, cf. \ref{derechef} et \ref{transfertdefonctions}.
 
 La m\'ethode utilis\'ee pour prouver ces th\'eor\`emes est de caract\'eriser les repr\'esentations de niveau $0$ par le d\'eveloppement local de leurs caract\`eres. Notons $\mathfrak{g}$ l'alg\`ebre de Lie de $G$. On sait d\'efinir l'exponentielle $exp$ qui envoie un voisinage de $0$ dans $\mathfrak{g}(F)$ sur un voisinage de $1$ dans $G(F)$. Sous l'hypoth\`ese $(Hyp)_{endo}(G)$, ces voisinages peuvent \^etre choisis les plus gros possibles. A savoir l'ensemble $\mathfrak{g}_{tn}(F)$ des \'el\'ements topologiquement nilpotents dans $\mathfrak{g}(F)$ et l'ensemble $G_{tu}(F)$ des \'el\'ements topologiquement unipotents dans $G(F)$. Il en est de m\^eme pour tout sous-groupe r\'eductif  connexe de $G$. Soit $D\in I(G)^*$ une distribution invariante et localement int\'egrable sur 
 $G(F)$, donc associ\'ee \`a une fonction $\theta_{D}$ d\'efinie presque partout sur $G(F)$, invariante par conjugaison et localement int\'egrable. Nous dirons que $D$ est un quasi-caract\`ere si et seulement si $\theta_{D}$ v\'erifie la condition suivante.  Soit $x\in G(F)$ un \'el\'ement semi-simple, notons $G_{x}$ la composante neutre de son commutant dans $G$. Alors il existe un voisinage $\mathfrak{V}_{x}$ de $0$ dans $\mathfrak{g}_{x}(F)$ de sorte que la fonction $X\mapsto \theta_{D}(x exp(X))$ d\'efinie presque partout sur $\mathfrak{V}_{x}$ soit combinaison lin\'eaire de transform\'ees de Fourier d'int\'egrales orbitales nilpotentes sur $\mathfrak{g}_{x}(F)$. Harish-Chandra a d\'emontr\'e que le caract\`ere de toute repr\'esentation irr\'eductible \'etait un quasi-caract\`ere (d'o\`u la terminologie "quasi-caract\`ere"). Appelons $p'$-\'el\'ement  un \'el\'ement $\epsilon\in G(F)$ qui est semi-simple et v\'erifie la condition suivante. Fixons une extension finie $F'$ de $F$ contenant toutes les valeurs propres de l'op\'erateur $ad(\epsilon)$ agissant dans $\mathfrak{g}$ et notons $\Sigma$ l'ensemble de ces valeurs propres. Notons aussi $\vert .\vert _{F'}$ la valeur absolue usuelle de $F'$. Alors, pour tout $\sigma\in \Sigma$ telle que $\vert \sigma\vert _{F'}=1$, $\sigma$ est une racine de l'unit\'e d'ordre premier \`a $p$. Soit $D\in I(G)^*$ une distribution invariante localement int\'egrable. Nous dirons que $D$ est un quasi-caract\`ere de niveau $0$ si et seulement si, pour tout $p'$-\'el\'ement $\epsilon\in G(F)$, la fonction $X\mapsto \theta_{D}(\epsilon exp(X))$ d\'efinie presque partout sur $\mathfrak{g}_{\epsilon,tn}(F)$ est combinaison lin\'eaire de transform\'ees de Fourier d'int\'egrales orbitales nilpotentes sur $\mathfrak{g}_{\epsilon}(F)$. Autrement dit, le d\'eveloppement est valable sur le plus gros voisinage possible dans $\mathfrak{g}_{\epsilon}(F)$. 
 
 \ass{Th\'eor\`eme 3}{Soit $\pi\in {\mathbb C}[Irr(G)]$. Alors $\pi\in {\mathbb C}[Irr(G)^0]$ si et seulement si $\Theta_{\pi}$ est un quasi-caract\`ere de niveau $0$.}
 
 Cf. \ref{caracterisation}. Le fait que le caract\`ere d'une repr\'esentation de niveau $0$ se d\'eveloppe sur de gros voisinages se trouvait d\'ej\`a dans la litt\'erature, par exemple dans des articles de S. Debacker, J.-L. Kim, F. Murnaghan. Mais je ne crois pas que l'on y trouvait la r\'eciproque.
 
 Il s'av\`ere que cette caract\'erisation des repr\'esentations de niveau $0$ "passe bien" \`a l'endoscopie et cela nous permet d'en d\'eduire les th\'eor\`emes 1 et 2.
 
 Pour d\'emonstrer le th\'eor\`eme 3, on doit d\'ecrire l'espace des quasi-caract\`eres de niveau $0$. Pour que cette introduction reste d'une longueur raisonnable, donnons un simple exemple. Soit ${\cal F}$ un sommet de l'immeuble de Bruhat-Tits  du groupe adjoint de $G$. On lui associe un groupe parahorique $K_{{\cal F}}^0\subset G(F)$. Notons $K_{{\cal F}}^+$ son plus grand sous-groupe distingu\'e pro-$p$-unipotent. On sait qu'il existe un groupe r\'eductif connexe ${\bf G}_{{\cal F}}$ d\'efini sur le corps r\'esiduel $k_{F}$ de $F$, de sorte que $K_{{\cal F}}^0/K_{{\cal F}}^+\simeq {\bf G}_{{\cal F}}(k_{F})$. Soit $f:{\bf G}_{{\cal F}}(k_{F})\to {\mathbb C}$ une fonction invariante par conjugaison et cuspidale. Par l'isomorphisme pr\'ec\'edent, on la consid\`ere comme une fonction sur $K_{{\cal F}}^0$, invariante par translations par $K_{{\cal F}}^+$, et on l'\'etend en une fonction sur $G(F)$, nulle hors de $K_{{\cal F}}^0$. Notons $A_{G}$ le plus grand sous-tore central de $G$ qui soit d\'eploy\'e sur $F$ et fixons des mesures de Haar sur $G(F)$ et $A_{G}(F)$. On d\'efinit la distribution $D_{f}$ qui, \`a une fonction $\varphi\in C_{c}^{\infty}(G(F))$, associe
 $$D_{f}(\varphi)=\int_{A_{G}(F)\backslash G(F)}\int_{G(F)}\varphi(g^{-1}xg)f(x)\,dx\,dg.$$
 L'int\'egrale n'est pas absolument convergente mais converge dans l'ordre indiqu\'e. On montre que $D_{f}$ est un quasi-caract\`ere de niveau $0$. On g\'en\'eralise cette construction dans deux directions. D'une part en consid\'erant non pas des groupes parahoriques mais des fixateurs de points de l'immeuble, qui donnent naissance \`a des groupes non connexes sur $k_{F}$. D'autre part, en induisant de telles distributions d\'efinies sur des groupes de Levi de $G$. On montre alors que tout quasi-caract\`ere de niveau $0$ est obtenue par cette construction  convenablement g\'en\'eralis\'ee, cf. \ref{undeuxiemetheoreme}. Ce qui relie les distributions ci-dessus aux repr\'esentations de niveau $0$ est la formule calculant le caract\`ere d'une telle repr\'esentation que l'on a obtenue dans un article ant\'erieur, cf. \cite{W2}. 
 
 Signalons un r\'esultat curieux. Les constructions ci-dessus conduisent naturellement \`a la d\'efinition d'un sous-espace de l'espace des quasi-caract\`eres de niveau $0$, not\'e $D^G[{\cal D}_{cusp}(G)]$ dans l'article. Ce sous-espace est reli\'e \`a celui des caract\`eres des repr\'esentations de niveau $0$ qui sont elliptiques au sens d'Arthur. Mais il ne lui est pas \'egal. On montre en   \ref{stabilite} et \ref{finale} qu'il poss\`ede n\'eanmoins les m\^emes propri\'et\'es magiques de ces espaces de caract\`eres, d\'emontr\'ees par Arthur dans \cite{A3}. A savoir que, dans le cas o\`u $G$ est quasi-d\'eploy\'e, leur stabilit\'e se lit sur les \'el\'ements elliptiques du groupe et que, dans le cas g\'en\'eral, le transfert entre tels quasi-caract\`eres se lit lui-aussi sur les \'el\'ements elliptiques des groupes en question.  
 
 L'hypoth\`ese $(Hyp)_{endo}(G)$ est utilis\'ee de deux fa\c{c}ons. D'une part, elle entra\^{\i}ne diverses propri\'et\'es concernant le groupe $G$. Par exemple, il existe une extension $F'$ de $F$ de degr\'e  fini et premier \`a $p$ telle que $G$ soit d\'eploy\'e sur $F'$. Ou bien, l'ordre du groupe des composantes connexes du centre de $G$ est premier \`a $p$. Pour ces propri\'et\'es, une hypoth\`ese plus faible serait suffisante. D'autre part, comme on l'a d\'ej\`a dit, l'hypoth\`ese $(Hyp)_{endo}(G)$ entra\^{\i}ne que l'exponentielle est d\'efinie sur les plus gros voisinages possibles. Pour certains groupes, on peut remplacer l'exponentielle par un substitut qui converge beaucoup mieux (par exemple $X\mapsto 1+X$ pour le groupe $GL(n)$).  De nouveau, on peut dans ce cas  remplacer $(Hyp)_{endo}(G)$ par une hypoth\`ese plus faible.  
 
 Les liens entre niveau, param\'etrage et endoscopie ont fait l'objet de plusieurs travaux r\'ecents. Citons \cite{Oi} et aussi \cite{Lanard1} et \cite{Lanard2}, o\`u  l'auteur \'etudie des sous-cat\'egories de celle des repr\'esentations de niveau $0$.

 \section{Les donn\'ees}
 
 \subsection{Le corps local}
 
 Soit $F$ un corps local non archim\'edien de caract\'eristique nulle. On note $\mathfrak{o}_{F}$ son anneau d'entiers, $\mathfrak{o}_{F}^{\times}$ le groupe des unit\'es, $ \mathfrak{p}_{F}$ l'id\'eal maximal, $k_{F}=\mathfrak{o}_{F}/\mathfrak{p}_{F}$ le corps r\'esiduel, $p$ la caract\'eristque de $k_{F}$, $\vert .\vert _{F}$ la valeur absolue usuelle, $val_{F}$ la valuation. 
 
 Fixons une cl\^oture alg\'ebrique $\bar{F}$ de $F$, resp. $\bar{k}_{F}$ de $k_{F}$. Notons  $F^{nr}$ la plus grande extension  non ramifi\'ee de $F$ contenue dans $\bar{F}$. On note $\Gamma_{F}$ le groupe de Galois de $\bar{F}/ F$ et $I_{F}$ son sous-groupe d'inertie, c'est-\`a-dire le groupe de Galois de $\bar{F}/F^{nr}$.  
 
 \subsection{Notations diverses}
 
 Quand un groupe abstrait $H$ agit sur un ensemble $X$, on note $X^{H}$ le sous-ensemble des points fixes. Si $Y$ est un sous-ensemble de $X$, on note $Norm_{H}(Y)$ le sous-ensemble des \'el\'ements de $H$ dont l'action conserve $Y$. 
 
 Si $X$ est un ensemble, on note ${\mathbb C}[X]$ le ${\mathbb C}$-espace vectoriel de base $X$.

\subsection{Groupes alg\'ebriques\label{groupesalgebriques}}

Posons $K=F$ ou $K=k_{F}$. Soit $H$ un groupe alg\'ebrique d\'efini sur  $K$. On note $H^0$ sa composante neutre et $Z(H)$ le centre de $H$.     Les sous-groupes alg\'ebriques de $H$ que l'on consid\'erera seront implicitement suppos\'es d\'efinis sur $K$,  sauf mention explicite du contraire. 

Supposons $H$ connexe. On appelle Levi de $H$ une composante de Levi d'un sous-groupe parabolique de $G$ (tous deux d\'efinis sur $K$ comme on vient de le dire).   On note $H_{AD}=H/Z(H)$ le groupe adjoint. Pour $x\in H$, on note  $x_{ad}$ son image dans $H_{AD}$, $x_{ss}$ la partie semi-simple de $x$, $Z_{H}(x)$ le commutant de $x$ dans $H$ et $H_{x}=Z_{H}(x)^0$.

 On utilise les notations d'Arthur concernant les Levi et paraboliques: si $M$ est un Levi, ${\cal L}(M)$ est l'ensemble des Levi contenant $M$ et ${\cal P}(M)$ est l'ensemble des sous-groupes paraboliques de $H$ de composante de Levi $M$. Si $P$ est un sous-groupe parabolique de $H$, on note $U_{P}$ son radical unipotent.

Soit $T$ un tore d\'efini sur $K$. On note $X^{*}(T)$, resp. $X_{*}(T)$,  les groupes de caract\`eres alg\'ebriques de $T$, resp. de sous-groupes \`a un param\`etre. Pour ces d\'efinitions, $T$ est vu comme un  tore sur la cl\^oture alg\'ebrique $\bar{K}$, c'est-\`a-dire que les caract\`eres ou cocaract\`eres ne sont pas forc\'ement d\'efinis sur $K$.   Si $T$ est d\'efini sur $F$, on note $T(F)_{c}$ le plus grand sous-groupe compact de $T(F)$. La d\'efinition de $X^*(T)$ s'\'etend au cas o\`u $T$ est un groupe diagonalisable,  c'est-\`a-dire un  sous-groupe alg\'ebrique d'un tore.

Soit $H$ un groupe r\'eductif connexe d\'efini sur $K$. On note $A_{H}$ le plus grand sous-tore de $Z(H)$ qui soit d\'efini et d\'eploy\'e sur $K$. On pose $a_{H}=dim(A_{H})$ et ${\cal A}_{H}=X_{*}(A_{H})\otimes_{{\mathbb Z}}{\mathbb R}$. On appelle sous-tore d\'eploy\'e maximal de $H$ un sous-tore de $H$ qui est d\'eploy\'e sur $K$ et est maximal parmi les sous-tores d\'eploy\'es sur $K$. Ces sous-tores d\'eploy\'es maximaux sont en bijection avec les groupes de Levi minimaux de $H$: \`a un sous-tore $A$ est associ\'e son commutant $Z_{H}(A)$. 

Soit $H$ un groupe r\'eductif connexe d\'efini sur $k_{F}$. On appelle espace tordu sous $H$ une vari\'et\'e alg\'ebrique $\tilde{H}$ d\'efinie sur $k_{F}$, munie de deux actions alg\'ebriques \`a droite et \`a gauche de $H$ telles que, pour chacune des actions, $\tilde{H}$ soit un espace principal homog\`ene sous $H$. On impose de plus ici que $\tilde{H}(k_{F})\not=\emptyset$. Cette terminologie a \'et\'e introduite par Labesse. Une partie de la th\'eorie habituelle des groupes r\'eductifs s'\'etend aux espaces tordus. Ainsi, on d\'efinit les notions de sous-espaces paraboliques ou d'espace de Levi. On envoie pour tout cela au premier chapitre de \cite{MW}.

  \subsection{Le groupe $G$\label{legroupeG}}

  On fixe pour tout l'article un groupe r\'eductif connexe $G$ d\'efini sur $F$.    Soit $N\geq1$ un entier tel qu'il existe un plongement de $G$ dans $GL(N)$. 
 On impose l'hypoth\`ese
  
  $(Hyp)(G)$: $p\geq (2+val_{F}(p))N$.
  
      Cela 
    entra\^{\i}ne  qu'il existe  une extension  $F'$ de $F$ de degr\'e premier \`a $p$ telle que $G$ soit d\'eploy\'e sur $F'$. Plus g\'en\'eralement, pour tout tore d\'efini sur $F$ de dimension inf\'erieure ou \'egale au rang de $G$, il existe une telle extension telle que le tore soit  d\'eploy\'e sur $F'$.  L'hypoth\`ese $(Hyp)(G)$ implique $Hyp(H)$ pour tout sous-groupe r\'eductif connexe de $G$. 
    
   A partir du paragraphe \ref{donneesendoscopiques}, nous renforcerons cette hypoth\`ese $(Hyp)(G)$ en une hypoth\`ese $(Hyp)_{endo}(G)$.

   On fixe un Levi minimal $M_{min}$ de $G$. On pose simplement $A=A_{M_{min}}$,  ${\cal A}={\cal A}_{M_{min}}$, ${\cal L}_{min}={\cal L}(M_{min})$. On d\'efinit $W^G=Norm_{G}(A)/M_{min}$.  
         
   Quand un objet a \'et\'e d\'efini relativement \`a notre groupe $G$, nous noterons souvent l'objet analogue d\'efini relativement \`a un autre groupe $H$ en ajoutant un exposant $H$ dans la notation. Par exemple, si $M\in {\cal L}_{min}$, on note ${\cal L}_{min}^M$ le sous-ensemble des $L\in {\cal L}_{min}$ qui sont contenus dans $M$. Le groupe $W^G$ agit sur ${\cal L}_{min}$.  
Pour $M\in {\cal L}_{min}$, on note $N_{W^G}(M)$ le fixateur de $M$ dans $W^G$ et $W^G(M)=N_{W^G}(M)/W^M$.

   On munit $G(F)$ d'une mesure de Haar. Plus g\'en\'eralement, pour tout sous-groupe ferm\'e $J$ de $G(F)$,  on suppose $J$ muni d'une mesure de Haar.

  \section{L'immeuble de $G$}
  
  \subsection{Facettes, groupes parahoriques\label{facettes}}
  
  On note $Imm(G_{AD})$ l'immeuble de Bruhat-Tits sur $F$ du groupe adjoint $G_{AD}$. Cet ensemble est r\'eunion d'appartements   associ\'es aux sous-tores d\'eploy\'es maximaux de $G$ ou encore aux Levi minimaux de $G$. On note $App(A_{M})$ l'appartement associ\'e \`a un  tel Levi  $M$.   L'appartement $App(A_{M})$ est un espace affine euclidien  sous l'espace vectoriel r\'eel  ${\cal A}_{M}/{\cal A}_{G}$. Pour $x,y\in App(A_{M})$, on note $x-y$ l'\'el\'ement de ${\cal A}_{M}/{\cal A}_{G}$ tel que $x=y+(x-y)$. 
  
  L'immeuble se d\'ecompose aussi en r\'eunion disjointe de facettes, chaque facette \'etant contenue dans (au moins) un appartement. On  note $Fac(G)$ l'ensemble des facettes et $Fac(G,A)$ le sous-ensemble des facettes contenues dans $App(A)$. Le groupe $G(F)$ agit sur l'immeuble. Pour ${\cal F}\in Fac(G)$, notons $K_{{\cal F}}^{\dag}$ le stabilisateur de ${\cal F}$ dans $G(F)$.
  
   Introduisons le groupe $\hat{G}$ dual de Langlands de $G$.  Modulo le choix d'une paire de Borel \'epingl\'ee de $\hat{G}$, ce groupe est muni d'une action de $\Gamma_{F}$.  Le groupe $\Gamma_{F}/I_{F}$ agit sur $X^*(Z(\hat{G})^{I_{F}})$. Notons ${\cal N}=X^{*}(Z(\hat{G})^{I_{F}})^{\Gamma_{F}/I_{F}}$ le sous-groupe des invariants. Kottwitz a d\'efini un homomorphisme surjectif $w_{G}:G(F)\to {\cal N}$. 
 
 Soit ${\cal F}\in Fac(G)$. Notons ${\cal N}({\cal F})$ l'image de $K_{{\cal F}}^{\dag}$ par $w_{G}$ et, pour tout $\nu\in {\cal N}$, notons $K_{{\cal F}}^{\nu}$ l'ensemble des $g\in K_{{\cal F}}^{\dag}$ tels que $w_{G}(g)=\nu$ (on a donc $ K_{{\cal F}}^{\nu}\not=\emptyset$ si et seulement si $\nu\in {\cal N}({\cal F})$). L'ensemble $K_{{\cal F}}^{0}$ est le sous-groupe parahorique de $G(F)$ associ\'e \`a ${\cal F}$, cf. \cite{HR} proposition 3. Notons $K_{{\cal F}}^{+}$ le plus grand sous-groupe distingu\'e et pro-$p$-unipotent dans $K_{{\cal F}}^0$. Bruhat et Tits ont  associ\'e \`a ${\cal F}$ un sch\'ema en groupes  ${\cal G}_{{\cal F}}$ d\'efini sur $\mathfrak{o}_{F}$. On a ${\cal G}_{{\cal F}}(\mathfrak{o}_{F})=K_{{\cal F}}^{0}$. Notons ${\bf G}_{{\cal F}}$ la   partie r\'eductive de la fibre sp\'eciale de ${\cal G}_{{\cal F}}$. C'est un groupe r\'eductif connexe d\'efini sur $k_{F}$ et   $K_{{\cal F}}^0/K_{{\cal F}}^+$ est isomorphe \`a ${\bf G}_{{\cal F}}(k_{F})$, cf. \cite{HV} proposition 3.7. Pour $\nu\in {\cal N}({\cal F})$, il existe un espace tordu ${\bf G}_{{\cal F}}^{\nu}$ sous ${\bf G}_{{\cal F}}$ de sorte que $K_{{\cal F}}^{\nu}/K_{{\cal F}}^{+}$ s'identifie \`a ${\bf G}_{{\cal F}}^{\nu}(k_{F})$.  
  
  Le groupe $K_{{\cal F}}^0$ fixe tout point de ${\cal F}$. Pour $\nu\in {\cal N}({\cal F})$, il existe une permutation isom\'etrique $\sigma_{{\cal F},\nu}$ de ${\cal F}$ telle que tout \'el\'ement de $K_{{\cal F}}^{\nu}$ agisse dans ${\cal F}$ par cette permutation. On note ${\cal F}^{\nu}$ le sous-ensemble des points fixes de $\sigma_{{\cal F},\nu}$ dans ${\cal F}$.  Pour tout appartement contenant ${\cal F}$,   ${\cal F}^{\nu}$ est l'intersection de ${\cal F}$ avec un sous-espace affine de  cet appartement. Soit $x\in {\cal F}$. Alors le fixateur de $x$ dans $G(F)$ est la r\'eunion des $K_{{\cal F}}^{\nu}$ sur les $\nu\in {\cal N}({\cal F})$ tels que $\sigma_{{\cal F},\nu}$ fixe $x$. 
  
   On note $Fac^*(G)$ l'ensemble des couples $({\cal F},\nu)$, o\`u ${\cal F}\in Fac(G)$ et $\nu\in {\cal N}({\cal F})$. On note $Fac^*_{max}(G)$ le sous-ensemble des $({\cal F},\nu)$ tels que ${\cal F}^{\nu}$ est r\'eduit \`a un point. On note $Fac^*(G,A)$ et $Fac^*_{max}(G,A)$ les sous-ensembles des $({\cal F},\nu)$ tels que ${\cal F}\subset App(A)$. 
   
   {\bf Remarque.} L'indice $max$ est contestable puisqu'il s'agit  de facettes de dimension minimale. Ce sont plut\^ot les groupes qui leur sont attach\'es qui sont maximaux. En tout cas, on conserve cette notation de \cite{W2} pour simplifier les r\'ef\'erences.
   \bigskip

     Du sch\'ema en groupes ${\cal G}_{{\cal F}}$  se d\'eduit une sous-$\mathfrak{o}_{F}$-alg\`ebre de Lie $\mathfrak{k}_{{\cal F}}$ de $\mathfrak{g}_{F}$. On note $\mathfrak{k}_{{\cal F}}^+$ son radical pro-$p$-nilpotent. Le quotient $\mathfrak{k}_{{\cal F}}/\mathfrak{k}_{{\cal F}}^+$ s'identifie \`a    l'espace des points sur $k_{F}$ de la partie r\'eductive de l'alg\`ebre de Lie de ${\cal G}_{{\cal F}}$.
     
      On peut aussi consid\'erer l'immeuble \'etendu $Imm(G)=Imm(G_{AD})\times {\cal A}_{G}$, qui est  muni d'une action de $G(F)$. Pour $x=(y,a)\in Imm(G)$, avec $y\in Imm(G_{AD})$ et $a\in {\cal A}_{G}$, le fixateur de $x$ dans $G(F)$ est le sous-groupe des \'el\'ements $g\in G(F)$ qui fixent $y$ et tels que $w_{G}(g)$ appartienne au sous-groupe de torsion ${\cal N}_{tors}$ de ${\cal N}$. 
      
      \subsection{Description des facettes\label{description}}
      
      Pour $a\in A(F)$, on note $a_{{\mathbb Z}}$ l'\'el\'ement de $ X_{*}(A)$ tel que $<x^*,a_{{\mathbb Z}}>=-val_{F}(x^*(a))$ pour tout $x^*\in X^*(A)$.  Le noyau de cet homomorphisme $a\mapsto a_{{\mathbb Z}}$ n'est autre que le plus grand sous-groupe compact $A(F)_{c}$ de $A(F)$.  L'action sur l'immeuble du  groupe $A(F)$ conserve $App(A)$. Pour $a\in A(F)$ et $x\in App(A)$, on a $ax= x+a_{{\mathbb Z}}$, o\`u ici $X_{*}(A)$ est vu comme un sous-groupe de ${\cal A}$. L'action du  groupe $Norm_{G}(A)(F)$ conserve aussi l'appartement $App(A)$. Notons $M_{min}(F)_{c}$ l'unique sous-groupe compact maximal de  $M_{min}(F)$. L'action de $Norm_{G}(A)(F)$ se quotiente en une action du groupe $Norm_{G}(A)(F)/M_{min}(F)_{c}$.
      
      Notons $\Sigma$ l'ensemble des racines r\'eduites de $A$ dans $G$. Fixons un sommet sp\'ecial  de $App(A)$ que l'on note $0$ qui nous permet d'identifier $App(A)$ et ${\cal A}/{\cal A}_{G}$. A tout $\alpha\in \Sigma$ est associ\'e un sous-ensemble $\Gamma_{\alpha}$ de ${\mathbb Q}$, qui est l'image r\'eciproque dans ${\mathbb Q}$ d'un sous-ensemble fini de ${\mathbb Q}/{\mathbb Z}$. Pour $c\in \Gamma_{\alpha}$, on note $c^+$ le plus petit \'el\'ement de $\Gamma_{\alpha}$ strictement sup\'erieur \`a $c$ et $c^-$ le plus grand \'el\'ement de $\Gamma_{\alpha}$ strictement inf\'erieur \`a $c$. On note $H_{\alpha,c}$ l'hyperplan affine de $App(A)$ d\'efini par l'\'equation $\alpha(x)=c$. Alors la d\'ecomposition en facettes de $App(A)$ est d\'efinie par la famille d'hyperplans $(H_{\alpha,c})_{\alpha\in \Sigma, c\in \Gamma_{\alpha}}$. A toute facette ${\cal F}\in Fac(G,A)$ sont associ\'es un sous-ensemble $\Sigma_{{\cal F}}\subset \Sigma$ et, pour tout $\alpha\in \Sigma$, un \'el\'ement $c_{\alpha,{\cal F}}$ de sorte que ${\cal F}$ soit le sous-ensemble des \'el\'ements $x\in App(A)$ qui v\'erifient les relations
      
      (1) $\alpha(x)=c_{\alpha,{\cal F}}$ pour tout $\alpha\in \Sigma_{{\cal F}}$;
      
      (2) $c_{\alpha,{\cal F}}<\alpha(x)<c_{\alpha,{\cal F}}^+$ pour tout $\alpha\in \Sigma-\Sigma_{{\cal F}}$.
      
      Pour $\alpha\in \Sigma$, notons $U_{\alpha}$ le groupe radiciel associ\'e \`a $\alpha$. A tout $c\in \Gamma_{\alpha}$ est associ\'e un sous-groupe ouvert compact $U_{\alpha,c}$ de $U_{\alpha}(F)$ de sorte que les propri\'et\'es suivantes soient v\'erifi\'ees:
      
      (3) si $c,c'\in \Gamma_{\alpha}$ et $c<c'$, alors $U_{\alpha,c}\subsetneq U_{\alpha,c'}$;
      
      (4) quelle que soit ${\cal F}\in Fac(G,A)$, on a $U_{\alpha}(F)\cap K_{{\cal F}}^0=U_{\alpha,c_{\alpha,{\cal F}}}$; on a $U_{\alpha}(F)\cap K_{{\cal F}}^+=U_{\alpha,c_{\alpha,{\cal F}}^-}$ si $\alpha\in \Sigma_{{\cal F}}$ et $U_{\alpha}(F)\cap K_{{\cal F}}^+=U_{\alpha,c_{\alpha,{\cal F}}}$ si $\alpha\not \in \Sigma_{{\cal F}}$.
      
      Notons $M_{min}(F)^{\star}=\{m\in M_{min}(F); w_{G}(m)=0\}$. C'est un sous-groupe d'indice fini de $M_{min}(F)_{c}$.  Pour toute facette ${\cal F}\in Fac(G,A)$, $K_{{\cal F}}^0$ est le sous-groupe de $G(F)$ engendr\'e par les $U_{\alpha}(F)\cap K_{{\cal F}}^0$ et par le groupe  $M_{min}(F)^{\star}$.  Les $\mathfrak{o}_{F}$-alg\`ebres $\mathfrak{k}_{{\cal F}}$ et $\mathfrak{k}_{{\cal F}}^+$ admettent une description similaire. Il en r\'esulte que 
    les groupes $K_{{\cal F}}^0$ et les $\mathfrak{o}_{F}$-alg\`ebres $\mathfrak{k}_{{\cal F}}$ et $\mathfrak{k}_{{\cal F}}^+$ caract\'erisent les facettes. C'est-\`a-dire
     
     (5) soient ${\cal F},{\cal F}'\in Fac(G)$; supposons $K_{{\cal F}}^0=K_{{\cal F}'}^0$, ou $\mathfrak{k}_{{\cal F}}=\mathfrak{k}_{{\cal F}'}$ ou $\mathfrak{k}^+_{{\cal F}}=\mathfrak{k}^+_{{\cal F}'}$; alors ${\cal F}={\cal F}'$. 
     
     Pour ${\cal F}\in Fac(G,A)$, il existe un unique Levi $M_{{\cal F}}\in {\cal L}_{min}$ tel que ${\cal A}_{M_{{\cal F}}}/{\cal A}_{G}$ soit l'espace vectoriel r\'eel engendr\'e par les $x-y$ pour $x,y\in {\cal F}$. Pour $({\cal F},\nu)\in Fac^*(G,A)$, il existe un unique Levi $M_{{\cal F},\nu}\in {\cal L}_{min}$ tel que ${\cal A}_{M_{{\cal F}}}/{\cal A}_{G}$ soit l'espace vectoriel r\'eel engendr\'e par les $x-y$ pour $x,y\in {\cal F}^{\nu}$. On a \'evidemment $M_{{\cal F}}\subset M_{{\cal F},\nu}$.
     
     \subsection{Groupes de Levi\label{groupesdelevi}}
     Soit $M\in {\cal L}_{min}$. Consid\'erons le sous-ensemble $Imm(G_{AD},M)$ de $Imm(G_{AD})$ r\'eunion des appartements $App(A_{M'})$ associ\'es aux Levi minimaux $M'$ contenus dans $M$. Le groupe ${\cal A}_{M}$ agit sur chacun de ces appartements et ces actions se recollent en une action sur la r\'eunion.  De plus, $Imm(G_{AD},M)$ est conserv\'e par l'action de $M(F)$. Le quotient de $Imm(G_{AD},M)$ par l'action de ${\cal A}_{M}$, muni de son action de $M(F)$, s'identifie canoniquement \`a l'immeuble $Imm(M_{AD})$. En particulier, $App(A)/{\cal A}_{M}$ s'identifie \`a l'appartement $App^M(A)$ associ\'e \`a $M_{min}$ dans l'immeuble $Imm(M_{AD})$. On note $p_{M}:App(A)\to App^M(A)$ cette projection. 
     
      L'action sur l'immeuble du  groupe $Norm_{G}(M)(F)$ conserve $Imm(G_{AD},M)$ et se descend en une action sur $Imm(M_{AD})$. Cette action permute les \'elements de $Fac(M)$. Pour ${\cal F}_{M}\in Fac(M)$, on note $K_{{\cal F}}^{\dag,G}$ le sous-groupe des \'el\'ements de $Norm_{G}(M)(F)$ qui conservent ${\cal F}_{M}$. Il contient $K_{{\cal F}_{M}}^{\dag}$ comme sous-groupe distingu\'e. 
      
      Notons $Norm_{G}(M,A)$ le normalisateur commun de $M$ et $A$ dans $G$. L'action du  groupe $Norm_{G}(M,A)(F)$ dans $Imm(M_{AD})$ conserve $App^M(A)$. 
     
     Pour ${\cal F}\in Fac(G,A)$, il existe une unique facette not\'ee ${\cal F}^{M}\in Fac(M,A)$ telle que $p_{M}({\cal F})\subset {\cal F}^{M}$. Si ${\cal F}$ est d\'ecrit par les relations (1) et (2) de \ref{description}, ${\cal F}^{M}$ est l'ensemble des $x\in App^M(A)$ qui v\'erifient les relations
     
     (1) $\alpha(x)=c_{\alpha,{\cal F}}$ pour tout $\alpha\in \Sigma_{{\cal F}}\cap \Sigma^M$;
     
     (2) $c_{\alpha,{\cal F}}<\alpha(x)<c_{\alpha,{\cal F}}$ pour tout $\alpha\in \Sigma^M-\Sigma_{{\cal F}}\cap \Sigma^M$.
     
    De l'inclusion $Z(\hat{G})\to Z(\hat{M})$ se d\'eduit   un homomorphisme naturel ${\cal N}^M\to {\cal N}$. Le diagramme
    $$\begin{array}{ccc}M(F)&\to&G(F)\\ w_{M}\downarrow\,\,&&w_{G}\downarrow\,\,\\ {\cal N}^M&\to &{\cal N}\\ \end{array}$$ 
    est commutatif. Posons $M_{ad}=M/Z(G)$. On a aussi un homomorphisme naturel ${\cal N}^M\to {\cal N}^{M_{ad}}$. Notons ${\cal N}_{G-comp}^M$ l'image r\'eciproque dans ${\cal N}^M$ du sous-groupe de torsion de ${\cal N}^{M_{ad}}$. Alors l'homomorphisme ${\cal N}^M\to {\cal N}$ se restreint en un homomorphisme injectif ${\cal N}^M_{G-comp}\to {\cal N}$, cf. \cite{W2} 6(1). On identifie ${\cal N}^M_{G-comp}$ \`a son image dans ${\cal N}$. On note $Fac^*_{G-comp}(M)$ le sous-ensemble des $({\cal F}_{M},\nu)\in Fac^*(M)$ tels que $\nu\in {\cal N}^M_{G-comp}$ (avec les variantes $Fac^*_{max,G-comp}(M)$ etc...).  Remarquons que, pour $({\cal F}_{M},\nu)\in Fac^*_{G-comp}(M)$ et pour $n\in K_{{\cal F}_{M}}^{\dag,G}$, la conjugaison par $n$ conserve l'ensemble $K_{{\cal F}_{M}}^{\nu}$.

  \section{Espaces de fonctions et de distributions}
  
  \subsection{Les espaces $I(G)$ et $I(G)^*$\label{lesespaces}}
  
   Le groupe $G(F)$ agit sur l'espace $C_{c}^{\infty}(G(F))$ par conjugaison: pour $g\in G(F)$ et $f\in C_{c}^{\infty}(G(F))$, $^gf$ est la fonction $x\mapsto f(g^{-1}xg)$. On note $I(G)$ le quotient de $C_{c}^{\infty}( G(F))$ par le sous-espace complexe engendr\'e par les $^gf-f$ pour $f\in C_{c}^{\infty}(G(F))$ et $g\in G(F)$. 
   
    On note $G_{reg}$ le sous-ensemble des \'el\'ements fortement r\'eguliers de $G$. Soit $f\in C_{c}^{\infty}(G(F))$.  Pour $x\in G_{reg}(F)$, on d\'efinit l'int\'egrale orbitale
  $$I^G(x,f)=D^G(x)^{1/2}\int_{A_{G_{x}}(F)\backslash G(F)}f(g^{-1}xg)\,dg,$$
  o\`u $D^G$ est le discriminant de Weyl usuel. L'espace $I(G)$ est aussi le quotient de $C_{c}^{\infty}(G(F))$ par le sous-espace des $f\in C_{c}^{\infty}(G(F))$ telles que $I^G(x,f)=0$ pour tout $x\in G_{reg}(F)$. 
  
   On appelle distribution invariante sur $G(F)$ une forme lin\'eaire sur $C_{c}^{\infty}(G(F))$ qui se quotiente en une forme lin\'eaire sur $I(G)$. On identifie une telle distribution \`a la forme lin\'eaire quotient. On note $I(G)^*$ l'espace des distributions invariantes sur $G(F)$. 
   
     Soient $M$ un Levi de $G$ et $D^M\in I(M)^*$. On d\'efinit une distribution $Ind_{M}^G(D^M)\in I(G)^*$ de la fa\c{c}on suivante. On fixe un sous-groupe parabolique $P\in {\cal P}(M)$. Comme on l'a dit, des mesures de Haar sont  fix\'ees sur $M(F)$, $U_{P}(F)$ et $G(F)$. Il s'en d\'eduit  ce que l'on peut appeler une pseudo-mesure invariante \`a droite sur $P(F)\backslash G(F)$. Pr\'ecis\'ement, c'est une forme lin\'eaire non pas sur $C_{c}^{\infty}(P(F)\backslash G(F))$ mais sur l'espace des fonctions localement constantes $\varphi:G(F)\to {\mathbb C}$  qui v\'erifient la relation $\varphi(mug)=\delta_{P}(m)\varphi(g)$ pour tous $m\in M(F)$, $u\in U_{P}(F)$ et $g\in G(F)$, o\`u $\delta_{P}$ est le module usuel. Cette pseudo-mesure est caract\'eris\'ee par l'\'egalit\'e
  $$\int_{G(F)}f(g)\,dg\,=\int_{P(F)\backslash G(F)}\int_{M(F)}\int_{U_{P}(F)}f(mug)\,du\,dm\,dg$$
    pour tout $f\in C_{c}^{\infty}(G(F))$.
 On d\'efinit la fonction $f_{U_{P}}$ sur $M(F)$ par 
 $$f_{U_{P}}(m)=\delta_{P}(m)^{1/2}\int_{U_{P}(F)}f(mu)\,du,$$
    puis la distribution $Ind_{M}^G(D^M)$ par l'\'egalit\'e
  $$Ind_{M}^G(D^M)(f)=\int_{P(F)\backslash G(F)}D^M(({^g}f)_{U_{P}})\,dg.$$
  Cela ne d\'epend pas du choix de $P$. 
  
    On dit que $x\in G_{reg}(F)$ est elliptique si le commutant $T$ de $x$ dans $G$ est un tore elliptique modulo $Z(G)$, c'est-\`a-dire si $A_{T}=A_{G}$. On note $G_{ell}(F)$ le sous-ensemble des \'el\'ements elliptiques dans $G_{reg}(F)$ (la notation \'etant un peu abusive: il n'y a pas de sous-ensemble alg\'ebrique $G_{ell}$ de $G$). On dit qu'une fonction $f\in C_{c}^{\infty}(G(F))$ est cuspidale si et seulement si les int\'egrales orbitales de $f$ sont nulles en tout point $x\in G_{reg}(F)-G_{ell}(F)$. Autrement dit si $f$ est annul\'ee par  $Ind_{M}^G(D^M)$ pour tout Levi propre $M$ et tout $D^M\in I(M)^*$. On note $C_{cusp}(G(F))$ l'espace des fonctions cuspidales et $I_{cusp}(G)$ son image dans $I(G)$.
  
  On dit que $f$ est tr\`es cuspidale si et seulement si, pour tout sous-groupe parabolique propre $P$ de $G$, de composante de Levi $M$, la fonction  $f_{U_{P}}$ est nulle. En fait, $I_{cusp}(G)$ est aussi l'image dans $I(G)$ de l'espace engendr\'e par les fonctions tr\`es cuspidales, cf. \cite{W3} lemme 2.7. Pour une fonction  $f$ tr\`es cuspidale, on d\'efinit une distribution $D_{f}\in I(G)^*$ par l'\'egalit\'e
  $$D_{f}(\varphi)=\int_{A_{G}(F)\backslash G(F)}\int_{G(F)}\varphi(x^{-1}gx)f(g)\,dg\,dx.$$
  La double int\'egrale n'est pas absolument convergente mais converge dans cet ordre, cf. \cite{W2} lemme 9. 
  
  \subsection{Filtrations\label{filtrations}}
  
  On fixe un ensemble de repr\'esentants $\underline{{\cal L}}_{min}\subset {\cal L}_{min}$ des classes de conjugaison de Levi de $G$. Pour tout entier $n\in {\mathbb Z}$, notons $\underline{{\cal L}}_{min}^n$ le sous-ensemble des $M\in \underline{{\cal L}}_{min}$ tels que $ a_{M}=n$ (on peut \'evidemment se limiter aux $n$ appartenant \`a l'intervalle $\{a_{G},...,a_{M_{min}}\}$). 
  
  On d\'efinit une filtration sur $I(G)$ de la fa\c{c}on suivante. Pour tout $n \in {\mathbb Z}$, notons $Fil^nI(G)$ l'image dans $I(G)$ du sous-espace des $f\in C_{c}^{\infty}(G(F))$ qui v\'erifient la condition: pour tout Levi $M$ tel que $a_{M}>n$ et tout $m\in G_{reg}(F)\cap M(F)$, on a $I^G(m,f)=0$.  On a
$$Fil^{a_{G}-1}I(G)=\{0\}\subset Fil^{a_{G}}I(G)=I_{cusp}(G)\subset...\subset Fil^{a_{M_{min}}}I(G)=I(G),$$
et, en posant $Gr^nI(G)=Fil^nI(G)/Fil^{n-1}I(G)$, on a
$$Gr^nI(G)\simeq \oplus_{M\in \underline{{\cal L}}_{min}^n}I_{cusp}(M)^{ W^G(M)},$$
cf. \ref{legroupeG} pour la d\'efinition du groupe $W^G(M)$, qui agit naturellement sur $I_{cusp}(M)$.

Notons $Ann^nI(G)^*$ l'annulateur de $Fil^{n-1}I(G)$ dans $I(G)^*$. On a
$$Ann^{a_{M_{min}}+1}I(G)^*=\{0\}\subset Ann^{a_{M_{min}}}I(G)^*\subset...\subset Ann^{a_{G}}I(G)^*=I(G)^*,$$
et, en posant $Gr^nI(G)^*=Ann^nI(G)^*/Ann^{n+1}I(G)^*$, on a
$$(1) \qquad Gr^nI(G)^*\simeq \oplus_{M\in \underline{{\cal L}}_{min}^n}I_{cusp}(M)^{* W^G(M)}.$$
 On v\'erifie que $Ann^nI(G)^*$ est le sous-espace de $I(G)^*$ engendr\'e par les distributions induites $Ind_{M}^G(I(M)^*)$ pour les Levi $M$ de $G$ tels que $a_{M}\geq n$. Pour $M\in \underline{{\cal L}}_{min}^n$ et $d^M\in I(M)^*$, l'image de $Ind_{M}^G(d^M)$ dans $Gr^nI(G)^*$ est nulle dans les composantes $I_{cusp}(M')^{* W^G(M')}$ de (1) pour $M'\not=M$ et est l'image naturelle de $d^M$ dans $I_{cusp}(M)^{* W^G(M)}$ (c'est-\`a-dire la restriction de $d^{M}$ \`a $I_{cusp}(M)^{ W^G(M)}$). 
  
  \subsection{El\'ements compacts et topologiquement unipotents\label{elementscompacts}}
  Soit $x\in G(F)$. On dit que $x$ est compact si et seulement s'il est contenu dans un sous-groupe compact de $G(F)$, c'est-\`a-dire si l'adh\'erence $\overline{x^{{\mathbb Z}}}$ du groupe $x^{{\mathbb Z}}$ engendr\'e par $x$ est compacte. On dit qu'il est compact mod $Z(G)$ si et seulement si l'image $x_{ad}$ de $x$ dans $G_{AD}(F)$ est compacte.  On dit que $x$ est  topologiquement unipotent si et seulement si $lim_{n\to \infty}x^{p^n}=1$
  
  Fixons un sous-tore maximal de $G$ contenant  la partie semi-simple $x_{ss}$ de $G$  et fixons une extension finie de $F$ telle que $T$ soit d\'eploy\'ee sur $F'$. Alors:
 
  $x$ est compact si et seulement si $ \chi(x_{ss})\in \mathfrak{o}^{\times}_{F'}$ pour tout $\chi\in X^*(T)$; $x$ est compact mod $Z(G)$ si et seulement si $ \chi(x_{ss})\in \mathfrak{o}^{\times}_{F'}$ pour tout $\chi\in X^*(T)$ tel que la restriction de $\chi$ \`a $Z(G)$ soit triviale;  $x$ est topologiquement unipotent si et seulement si  $\chi(x_{ss})\in 1+\mathfrak{p}_{F'}$ pour tout $\chi\in X^*(T)$.

D'autre part, $x$ est compact mod $Z(G)$ si et seulement si il existe ${\cal F}\in Fac(G)$ tel que $x\in K_{{\cal F}}^{\dag}$. L'ensemble des \'el\'ements compacts, resp. compacts mod $Z(G)$, est un sous-ensemble ouvert et ferm\'e de $G(F)$ invariant par conjugaison. Soit $M$ un Levi de $G$, soit $({\cal F}_{M},\nu)\in Fac^*(M)$ et soit $x\in K_{{\cal F}_{M}}^{\nu}$. Alors $x$ est compact mod $Z(G)$ si et seulement si $\nu\in {\cal N}^M_{G-comp}$. 

On note $G_{comp}(F)$ l'ensemble des \'el\'ements de $G(F)$ qui sont compacts mod $Z(G)$ (il serait plus correct de le noter $G_{comp\, mod \,Z(G)}$...).   On note $G_{tu}(F)$ l'ensemble des \'el\'ements topologiquement unipotents de $G(F)$. On note $I(G)^*_{comp}$ l'ensemble des distributions dont le support est contenu dans  $G_{comp}(F)$.

  \subsection{Facettes, fonctions et distributions\label{facettesfonctions}}
     Soit $({\cal F},\nu)\in Fac_{max}^*(G)$. On note $C_{cusp}({\bf G}^{\nu}_{{\cal F}})$ l'espace des fonctions sur ${\bf G}^{\nu}_{{\cal F}}(k_{F})$ qui sont invariantes par conjugaison par ${\bf G}_{{\cal F}}(k_{F})$ et qui sont cuspidales. On munit cet espace du produit hermitien d\'efini par
  $$<f,f'>=\vert {\bf G}_{{\cal F}}(k_{F})\vert ^{-1}\sum_{x\in {\bf G}^{\nu}_{{\cal F}}(k_{F})}\bar{f}(x)f'(x).$$
  Il est d\'efini positif. Soit $f\in C_{cusp}({\bf G}^{\nu}_{{\cal F}})$.  L'espace ${\bf G}_{{\cal F}}^{\nu}(k_{F})$ s'identifie \`a $K_{{\cal F}}^{\nu}/K_{{\cal F}}^+$. On identifie $f$ \`a une fonction sur ce quotient $K_{{\cal F}}^{\nu}/K_{{\cal F}}^+$, on la rel\`eve en une fonction sur $K_{{\cal F}}^{\nu}$ puis on l'\'etend en une fonction sur $G(F)$ par $0$ hors de $K_{{\cal F}}^{\nu}$. On note $f_{{\cal F}}$ la fonction ainsi obtenue. Elle est tr\`es cuspidale, cf. \cite{W2} lemme 10. On d\'eduit de $f_{{\cal F}}$ une distribution $D_{f_{{\cal F}}}$, cf. \ref{lesespaces}, que l'on note simplement $D_{f}$.  Cette distribution est \`a support compact mod $Z(G)$. 
  
  Notons $\boldsymbol{{\cal D}}_{cusp}(G)$ le sous-ensemble des familles
  $$(f_{{\cal F},\nu})_{({\cal F},\nu)\in Fac^*_{max}(G)}\in \prod_{({\cal F},\nu)\in Fac^*_{max}(G)}C_{cusp}({\bf G}^{\nu}_{{\cal F}})$$
  qui v\'erifient la condition: pour tout $\nu\in {\cal N}$, l'ensemble des ${\cal F}$ tels que $({\cal F},\nu)\in Fac^*_{max}(G)$ et $f_{{\cal F},\nu}\not=0$ est fini.
  
  Pour ${\bf f}=(f_{{\cal F},\nu})_{({\cal F},\nu)\in Fac^*_{max}(G)}\in \boldsymbol{{\cal D}}_{cusp}(G)$, la somme $\sum_{{\cal F},\nu}D_{f_{{\cal F},\nu}}$ est d\'efinie. En effet, chaque $D_{f_{{\cal F},\nu}}$ est \`a support dans l'ouvert ferm\'e $w_{G}^{-1}(\nu)$ et, pour tout $\nu$, il n'y a qu'un nombre fini de $({\cal F},\nu')$ tels que $\nu'=\nu$ et $f_{{\cal F},\nu'}\not=0$. Pour tout compact $C$ de $G(F)$, il n'y a donc qu'un nombre fini de $({\cal F},\nu)$ tels que le support de $D_{f_{{\cal F},\nu}}$ coupe $C$. On pose $D_{{\bf f}}=\sum_{{\cal F},\nu}D_{f_{{\cal F},\nu}}$. Cela d\'efinit une application lin\'eaire
  $$D: \boldsymbol{{\cal D}}_{cusp}(G)\to I(G)_{comp}^*.$$
  
  {\bf Remarque.} Quand un \'el\'ement de $\boldsymbol{{\cal D}}_{cusp}(G)$ sera not\'e ${\bf f}$, on notera son image $D_{{\bf f}}$ ou $D^G_{{\bf f}}$ comme ci-dessus. Quand l'\'el\'ement sera not\'e plus symboliquement $d$, on notera son image $D[d]$ ou $D^G[d]$.
  
  \bigskip
  Le groupe $G(F)$ agit naturellement sur $\boldsymbol{{\cal D}}_{cusp}(G)$. Il est clair que, pour ${\bf f}\in \boldsymbol{{\cal D}}_{cusp}(G)$ et $g\in G(F)$, on a $D({^g{\bf f}}-{\bf f})=0$. Notons ${\cal D}_{cusp}(G)$ l'espace des coinvariants, c'est-\`a-dire le quotient de $\boldsymbol{{\cal D}}_{cusp}(G)$ par le sous-espace engendr\'e par les ${^g{\bf f}}-{\bf f}$ pour ${\bf f}\in \boldsymbol{{\cal D}}_{cusp}(G)$ et $g\in G(F)$. Alors $D$ se quotiente en une application lin\'eaire 
  
  $$D:{\cal D}_{cusp}(G)\to I(G)_{comp}^*.$$

  Notons $\boldsymbol{{\cal D}}(G)$ la somme directe des $\boldsymbol{{\cal D}}_{cusp}(M)$ o\`u $M$ parcourt les Levi de $G$. On d\'efinit une application lin\'eaire
  $$D^G:  \boldsymbol{{\cal D}}(G)\to I(G)^*$$
  qui, \`a $\sum_{M}d^M$, associe $\sum_{M}Ind_{M}^G(D^M[d^M])$.
  Le groupe $G(F)$ agit naturellement sur l'espace $\boldsymbol{{\cal D}}(G)$: un \'el\'ement $g\in G(F)$ envoie un Levi $M$ sur le Levi $gMg^{-1}$, une facette $({\cal F}_{M},\nu)\in Fac^*_{max}(M)$ sur une facette $(g{\cal F}_{M},g\nu)\in Fac^*_{max}(gMg^{-1})$ et une fonction $f_{{\cal F}^M,\nu}\in C_{cusp}({\bf M}_{{\cal F}_{M}}^{\nu})$ sur une fonction $^g(f_{{\cal F}_{M},\nu})\in C_{cusp}(({\bf gMg^{-1}})_{g{\cal F}_{M}}^{g\nu})$. Il est clair que, pour $g\in G(F)$ et ${\bf f}\in \boldsymbol{{\cal D}}(G)$, on a $D^G[{^g{\bf f}}-{\bf f}]=0$. On note ${\cal D}(G)$ l'espace des coinvariants pour cette action de $G(F)$ dans $\boldsymbol{{\cal D}}(G)$. L'application pr\'ec\'edente se quotiente en une application lin\'eaire
  $$D^G={\cal D}(G)\to I(G)^*.$$ 
  
    Remarquons  que  l'application $\boldsymbol{{\cal D}}_{cusp}(M)\to \boldsymbol{{\cal D}}(G)\to {\cal D}(G)$ se quotiente par ${\cal D}_{cusp}(M)$. 
  
  \subsection{Variantes de ${\cal D}(G)$\label{variantes}}  
    Soit $M$ un Levi de $G$.  On note $\boldsymbol{{\cal D}}_{cusp,G-comp}(M)$ le sous-espace des \'el\'ements  ${\bf f}=(f_{{\cal F}_{M},\nu})_{({\cal F}_{M},\nu)\in Fac^*_{max}(M)}$ de $ \boldsymbol{{\cal D}}_{cusp}(M)$ tels que $f_{{\cal F}_{M},\nu}=0$ si $\nu\not\in {\cal N}^M_{G-comp}$. On note $\boldsymbol{{\cal D}}(G)$ la somme directe des $\boldsymbol{{\cal D}}_{cusp,G-comp}(M)$ o\`u $M$ parcourt les Levi de $G$. On d\'efinit comme dans le paragraphe pr\'ec\'edent l'espace de coinvariants ${\cal D}_{G-comp}(G)$, qui s'identifie \`a un sous-espace de ${\cal D}(G)$. L'application $D^G$ se restreint en une application lin\'eaire
    $$D^G:{\cal D}_{G-comp}(G)\to I(G)^*_{G-comp}.$$
    
     On note ${\cal D}_{cusp,G-comp}(M)$ l'image de $\boldsymbol{{\cal D}}_{cusp,G-comp}(M)$ dans ${\cal D}_{cusp}(M)$. Il y a une projection naturelle $\boldsymbol{{\cal D}}_{cusp}(M)\to \boldsymbol{{\cal D}}_{cusp,G-comp}(M)$ qui, \`a ${\bf f}=(f_{{\cal F}_{M},\nu})_{({\cal F}_{M},\nu)\in Fac^*_{max}(M)}\in \boldsymbol{{\cal D}}_{cusp}(M)$ associe ${\bf f}=(f_{{\cal F}_{M},\nu})_{({\cal F}_{M},\nu)\in Fac^*_{max,G-comp}(M)}\in \boldsymbol{{\cal D}}_{cusp}(M)$. Cette projection se quotiente en une projection de ${\cal D}_{cusp}(M)$ sur ${\cal D}_{cusp,G-comp}(M)$.

    Fixons $\nu\in {\cal N}$. Notons $I(G)^{*\nu}_{G-comp}$ le sous-espace des distributions invariantes \`a support   contenu dans $G_{comp}(F)\cap w_{G}^{-1}(\nu)$. Pour un Levi $M$ de $G$, notons $\boldsymbol{{\cal D}}^{\nu}_{cusp,G-comp}(M)$  le sous-espace des ${\bf f}=(f_{{\cal F}_{M},\nu'})_{({\cal F}_{M},\nu')\in Fac^*_{max}(M)}\in \boldsymbol{{\cal D}}_{cusp,G-comp}(M)$ tels que $f_{{\cal F}_{M},\nu'}=0$ si $\nu'\not=\nu$.  A l'aide de ces espaces, on d\'efinit comme ci-dessus  un espace ${\cal D}^{\nu}_{G-comp}(G)$. On a des isomorphismes naturels
    $${\cal D}_{G-comp}(G)\simeq \prod_{\nu\in {\cal N}}{\cal D}_{G-comp}^{\nu}(G),$$
    $$I(G)^*_{G-comp}\simeq \prod_{\nu\in {\cal N}}I(G)^{*\nu}_{G-comp},$$
    et l'application $D^G$ s'identifie au produit de ses restrictions
    $$D^G:{\cal D}_{G-comp}^{\nu}(G)\to I(G)^{*\nu}_{G-comp}.$$

   \subsection{Calcul de $D_{f'}^G(f_{{\cal F}})$\label{calcul}}
   
    Soit $({\cal F},\nu)\in Fac^*(G,A)$. On a d\'efini le Levi $M_{{\cal F},\nu}\in {\cal L}_{min}$ au paragraphe \ref{description}. Posons simplement $M=M_{{\cal F},\nu}$. On associe \`a ${\cal F}$ la facette ${\cal F}^{M}\in Fac(M,A)$. L'\'el\'ement $\nu$ appartient \`a ${\cal N}_{G-comp}^M({\cal F}^{M})$, $({\cal F}^M,\nu)$ appartient \`a $Fac_{max}(M,A)$, ${\bf G}_{{\cal F}}$ s'identifie \`a ${\bf M}_{{\cal F}}$ et ${\bf G}_{{\cal F}}^{\nu}$ s'identifie \`a ${\bf M}_{{\cal F}^M}^{\nu}$, cf. \cite{W2} lemme 6. Soit $f\in C_{cusp}({\bf M}^{\nu}_{{\cal F}^M})$. On identifie $f$  \`a une fonction sur ${\bf G}_{{\cal F}}^{\nu}(k_{F})$, puis \`a une fonction sur $K_{{\cal F}}^{\nu}/K_{{\cal F}}^+$, que l'on rel\`eve en une fonction sur $K_{{\cal F}}^{\nu}$. On prolonge cette fonction en une fonction sur $G(F)$, nulle hors de $K_{{\cal F}}^{\nu}$. On note $f_{{\cal F}}$ la fonction sur $G(F)$ obtenue ainsi. Elle est par construction invariante par conjugaison par $K_{{\cal F}}^0$. 
    
  Soient $L\in {\cal L}_{min}$, $({\cal F}',\nu')\in Fac^*_{max,G-comp}(L,A)$ et $f'\in C_{cusp}({\bf L}_{{\cal F}'}^{\nu'})^{K_{{\cal F}'}^{\dag,G}}$. On a d\'efini   la distribution $D^G_{f'}=Ind_{L}^G(D^L_{f'})$ sur $G(F)$.  Notons $N(M,{\cal F}^{M},L,{\cal F}')$ l'ensemble des $n\in Norm_{G}(A)(F)$ tels que $nMn^{-1}=L$ et $n{\cal F}^{M}={\cal F}'$. Cet ensemble (qui peut \^etre vide) est invariant \`a gauche par $A_{L}(F)(K_{{\cal F}'}^0\cap Norm_{G}(A)(F))$. On fixe un ensemble de repr\'esentants $\underline{N}(M,{\cal F}^{M},L,{\cal F}')$ du quotient
  $$A_{L}(F)(K_{{\cal F}'}^0\cap Norm_{G}(A)(F))\backslash N(M,{\cal F}^{M},L,{\cal F}').$$
    
  \ass{Proposition}{(i) $D^G_{f'}(f_{{\cal F}})\not=0$ seulement si $\nu=\nu'$ et  $N(M,{\cal F}^{M},L,{\cal F}')\not=\emptyset$. 
    
  (ii) Supposons $\nu=\nu'$ et $N(M,{\cal F}^{M},L,{\cal F}')\not=\emptyset$. Alors
  $$D^G_{f'}(f_{{\cal F}})= mes(K_{{\cal F}}^0)mes(K_{{\cal F}'}^0) mes(A_{L}(F)_{c})^{-1}\sum_{n\in \underline{N}(M,{\cal F}^{M},L,{\cal F}')}<\overline{^nf},f'>.$$ }
  
  {\bf Remarques.} (1) La fonction $^nf$ est celle d\'efinie en \ref{facettesfonctions}.
  
  (2) Si $f'$ est invariante par $K_{{\cal F}'}^{\dag,G}$, la formule du (ii) se simplifie et devient
  $$D_{f'}^G(f_{{\cal F}})=mes(K_{{\cal F}}^0)mes(K_{{\cal F}'}^0)mes(A_{L}(F)_{c})^{-1}[K_{{\cal F}'}^{\dag,G}:A_{L}(F)K_{{\cal F}'}^0]<\bar{f'},f>.$$
  \bigskip
    
  Preuve. Puisque $f$ est \`a support dans $w_{G}^{-1}(\nu)$ et que $D^G_{f'}$ est \`a support dans $w_{G}^{-1}(\nu')$ il est clair que $D^G_{f'}(f)=0$ si $\nu\not=\nu'$. Supposons d\'esormais $\nu=\nu'$.

 Fixons un sous-groupe parabolique $Q\in {\cal P}(L)$. Commen\c{c}ons par calculer $f_{{\cal F},U_{Q}}$. On sait d\'efinir la facette ${\cal F}^L\in Imm(L_{AD})$.  Si $M\subset L$, on a \'evidemment $({\cal F}^L)^M={\cal F}^M$ et la fonction $f_{{\cal F}^L}$ sur $L(F)$ est bien d\'efinie. Montrons que
 
 (3) si $M\not\subset L$, $ f_{{\cal F},U_{Q}}=0$; si $M\subset L$, $f_{{\cal F},U_{Q}}=mes(U_{Q}(F)\cap K_{{\cal F}}^+)f_{{\cal F}^L}$.
 
 Rappelons que, pour $l\in L(F)$, on a $f_{{\cal F},U_{Q}}(l)=\int_{U_{Q}(F)}f_{{\cal F}}(lu)\,du$. Le groupe $Q$, resp. $L$,  d\'etermine un sous-groupe parabolique ${\bf Q}$ de ${\bf G}_{{\cal F}}$, resp. une composante de Levi ${\bf L}$ de ${\bf Q}$, de sorte que l'image de $Q(F)\cap K_{{\cal F}}^0$, resp. $L(F)\cap K_{{\cal F}}^0$,  dans ${\bf G}_{{\cal F}}(k_{F})$ soit ${\bf Q}(k_{F})$, resp. ${\bf L}(k_{F})$. Supposons que $f_{{\cal F},U_{Q}}\not=0$. Alors $Q(F)\cap K_{{\cal F}}^{\nu}\not=\emptyset$ et il existe  un sous-espace parabolique ${\bf Q}^{\nu}$ de ${\bf G}^{\nu}_{{\cal F}}$, associ\'e au parabolique ${\bf Q}$, de sorte que l'image dans ${\bf G}^{\nu}_{{\cal F}}(k_{F})$ de $Q(F)\cap K_{{\cal F}}^{\nu}$ soit ${\bf Q}^{\nu}(k_{F})$. Notons ${\bf L}^{\nu}$ le normalisateur de ${\bf L}$ dans ${\bf Q}^{\nu}$. Posons $V=U_{Q}(F)\cap K_{{\cal F}}^0$. Fixons $l\in L(F)$ tel que $f_{{\cal F},U_{Q}}(l)\not=0$.  On \'ecrit
 $$f_{{\cal F},U_{Q}}(l)=\int_{U_{Q}(F)/V}\int_{V}f_{{\cal F}}(luv)\,dv\,du.$$
 Fixons $u\in U_{Q}(F)$ tel que l'int\'egrale int\'erieure soit non nulle. L'\'el\'ement $lu$ appartient \`a $Q(F)\cap K_{{\cal F}}^{\nu}$. D\'ecomposons son image dans ${\bf Q}^{\nu}(k_{F})$ en ${\bf l}{\bf u}$ avec ${\bf l}\in {\bf L}^{\nu}(k_{F})$ et ${\bf u}\in {\bf U}_{{\bf Q}}(k_{F})$. On voit alors que
 $$\int_{V}f_{{\cal F}}(luv)\,dv=c\sum_{{\bf v}\in {\bf U}_{{\bf Q}}(k_{F})}f({\bf l}{\bf v}),$$
 o\`u $c>0$ est une constante provenant des mesures et o\`u on a identifi\'e $f$ \`a une fonction sur ${\bf G}^{\nu}_{{\cal F}}(k_{F})$. Cette fonction est cuspidale. La somme ci-dessus est donc  nulle  si ${\bf Q}^{\nu}\not={\bf G}_{{\cal F}}^{\nu}$. Notre hypoth\`ese de non-nullit\'e implique donc ${\bf Q}^{\nu}={\bf G}_{{\cal F}}^{\nu}$, d'o\`u aussi ${\bf Q}={\bf L}={\bf G}_{{\cal F}}$. Le tore $A$ d\'etermine un sous-tore d\'eploy\'e maximal ${\bf A}$ de ${\bf G}_{{\cal F}}$. On a $X_{*}(A)\simeq X_{*}({\bf A})$.  Les sous-tores $A_{M_{{\cal F}}}$ et $A_{L}$ d\'eterminent des sous-tores ${\bf A}_{M_{{\cal F}}}$ et ${\bf A}_{L}$. Par d\'efinition de $M_{{\cal F}}$, ${\bf A}_{M_{{\cal F}}}$ est le plus grand sous-tore d\'eploy\'e contenu dans le centre de ${\bf G}_{{\cal F}}$. Par d\'efinition de ${\bf L}$, ce groupe est le commutant de ${\bf A}_{L}$ dans ${\bf G}_{{\cal F}}$. L'\'egalit\'e ${\bf L}={\bf G}_{{\cal F}}$ entra\^{\i}ne donc ${\bf A}_{L}\subset {\bf A}_{M_{\cal F}}$. D'o\`u aussi $A_{L}\subset A_{M_{{\cal F}}}$.   Puisque ${\bf Q}^{\nu}={\bf G}_{{\cal F}}^{\nu}$, on a $K_{{\cal F}}^{\nu}=(Q(F)\cap K_{{\cal F}}^{\nu})K_{{\cal F}}^+$. On a aussi  $K_{{\cal F}}^+=(K_{{\cal F}}^+\cap Q(F))(K_{{\cal F}}^+\cap U_{\bar{Q}}(F))$, o\`u $\bar{Q}$ est le parabolique de composante de Levi $L$ oppos\'e \`a $Q$. Donc $K_{{\cal F}}^{\nu}=(Q(F)\cap K_{{\cal F}}^{\nu})(K_{{\cal F}}^+\cap U_{\bar{Q}}(F))$. On sait que l'ensemble $K_{{\cal F}}^{\nu}\cap Norm_{G}(A)(F)$ est non vide. Mais un \'el\'ement de $Norm_{G}(A)(F)$ qui appartient \`a $Q(F)U_{\bar{Q}}(F)$ appartient forc\'ement \`a $Norm_{L}(A)(F)$. Soit donc $w\in K_{{\cal F}}^{\nu}\cap Norm_{L}(A)(F)$. Il agit naturellement dans $A$ et conserve $A_{M_{{\cal F}}}$. Par d\'efinition de $M=M_{{\cal F},\nu}$, $A_{M}$ est le plus grand sous-tore de $A_{M_{{\cal F}}}$ contenu dans l'ensemble des points fixes de l'action de $w$. Puisque $w\in L(F)$, ce sous-tore contient $A_{L}$. Donc $A_{L}\subset A_{M}$ et $M\subset L$. Cela d\'emontre la premi\`ere assertion de (3). 
 
 Supposons maintenant $M\subset L$. Alors, pour $l\in L(F)$ et $u\in U_{Q}(F)$, on a $lu\in K_{{\cal F}}^{\nu}$ si et seulement si $l\in K_{{\cal F}^L}^{\nu}$ et $u\in U_{Q}(F)\cap K_{{\cal F}}^+$ (cela r\'esulte de \cite{W2} lemme 6). On voit  que la fonction $f_{{\cal F},U_{Q}}$ est \`a support dans $K_{{\cal F}^L}^{\nu}$ et que, pour un \'el\'ement $l$ de ce groupe, on a $f_{{\cal F},U_{Q}}(l)=mes(U_{Q}(F)\cap K_{{\cal F}}^+)f(\bar{l})$, o\`u $\bar{l}$ est l'image de $l$ dans ${\bf L}_{{\cal F}^L}^{\nu}(k_{F})\simeq {\bf M}_{{\cal F}^M}^{\nu}(k_{F})$. Par d\'efinition de la fonction $f_{{\cal F}^L}$, on obtient la deuxi\`eme assertion de (3), ce qui ach\`eve la preuve de cette assertion. 
 
 Supposons maintenant $M\subset L$. Posons
 $$I_{{\cal F}^L,{\cal F}'}(f,f')=\int_{L(F)}f_{{\cal F}^L}(l)f'_{{\cal F}'}(l)\, dl.$$
 Montrons que
 
 (4) si ${\cal F}^L\not={\cal F}'$, $I_{{\cal F}^L,{\cal F}'}(f,f')=0$; si ${\cal F}^L={\cal F}'$, alors $M=L$, ${\cal F}^M={\cal F}'$ et $I_{{\cal F}^L,{\cal F}'}(f,f')=mes(K_{{\cal F}'}^0)<\bar{f},f'>$.
 
 Ici, tout se passe dans $L$, on peut aussi bien supposer $L=G$ pour simplifier les notations. On a donc ${\cal F}^L={\cal F}$. Supposons ${\cal F}\not={\cal F}'$. Dans l'appartement $App(A)$, fixons
un segment $[x,x']$ joignant un point $x\in {\cal F}^{\nu}$ \`a un point $x'\in {\cal F}^{'\nu}$. L'hypoth\`ese ${\cal F}\not={\cal F}'$ entra\^{\i}ne que $x\not=x'$. Il y a une facette ${\cal F}''\in Fac(G,A)$ et un \'el\'ement $x''\in ]x,x'[$ tel que le segment $[x'',x'[$ soit contenu dans ${\cal F}''$.   D'apr\`es ce que sont les supports de $f_{{\cal F}}$ et $f'_{{\cal F}'}$, on a
$$ I_{{\cal F},{\cal F}'}(f,f')=\int_{K_{{\cal F}}^{\nu}\cap K_{{\cal F}'}^{\nu}}f_{{\cal F}}(g)f'_{{\cal F}'}(g)\, dg.$$
Si  $K_{{\cal F}}^{\nu}\cap K_{{\cal F}'}^{\nu}=\emptyset$, cette int\'egrale est nulle et on a d\'emontr\'e la premi\`ere assertion de (4). Sinon, consid\'erons un \'el\'ement $g\in K_{{\cal F}}^{\nu}\cap K_{{\cal F}'}^{\nu}$. L'action de $g$ sur l'immeuble fixe $x$ et $x'$, donc aussi tout le segment $[x,x']$. En particulier, il fixe $[x'',x'[$, donc $g\in K_{{\cal F}''}^{\dag}$. On a $w_{G}(g)=\nu$, donc $g\in K_{{\cal F}''}^{\nu}$. Alors $[x'',x'[\subset {\cal F}^{''\nu}$. Cela entra\^{\i}ne ${\cal F}''\not={\cal F}'$ car l'hypoth\`ese $({\cal F}',\nu)\in Fac^*_{max}(G,A)$ (rappelons que l'on a suppos\'e $L=G$) signifie que ${\cal F}^{'\nu}$ est r\'eduit \`a un unique point, qui est donc $x'$. Le point $x'$ est adh\'erent \`a ${\cal F}''$, donc toute la facette ${\cal F}'$ est contenue dans l'adh\'erence de ${\cal F}''$. On sait qu'alors il existe un sous-groupe parabolique propre ${\bf P}$ de ${\bf G}_{{\cal F}'}$ et un sous-espace parabolique ${\bf P}^{\nu}$ de ${\bf G}_{{\cal F}'}^{\nu}$ associ\'e \`a ${\bf P}$, de sorte que l'image de $K_{{\cal F}''}^{0}\cap K_{{\cal F}'}^0$ dans ${\bf G}_{{\cal F}'}(k_{F})$ soit ${\bf P}(k_{F})$ et que l'image de $K_{{\cal F}''}^{\nu}\cap K_{{\cal F}'}^{\nu}$ dans ${\bf G}_{{\cal F}'}^{\nu}(k_{F})$ soit ${\bf P}^{\nu}(k_{F})$. Il est facile d'identifier l'ensemble ${\bf U}_{{\bf P}}(k_{F})$ en d\'ecrivant les facettes ${\cal F}'$ et ${\cal F}''$ comme en \ref{description}. On a forc\'ement $\Sigma_{{\cal F}''}\subset \Sigma_{{\cal F}'}$, $c_{\alpha,{\cal F}''}=c_{\alpha,{\cal F}'}$ si $\alpha\in \Sigma_{{\cal F}''}$ et $c_{\alpha,{\cal F}''}=c_{\alpha,{\cal F}'}$ ou $c_{\alpha,{\cal F}'}^-$ si $\alpha\in \Sigma_{{\cal F}'}-\Sigma_{{\cal F}''}$. Notons $\Xi$ l'ensemble des $\alpha\in \Sigma_{{\cal F}'}-\Sigma_{{\cal F}''}$ tels que $c_{\alpha,{\cal F}''}=c_{\alpha,{\cal F}'}$. Alors ${\bf U}_{{\bf P}}(k_{F})$ est l'image dans ${\bf G}_{{\cal F}'}(k_{F})$ de $\prod_{\alpha\in \Xi}U_{\alpha,c_{\alpha,{\cal F}'}}$. Notons $V$ ce dernier groupe. Il est inclus dans $K_{{\cal F}'}^0$. Montrons qu'il est inclus dans $K_{{\cal F}}^+$. En effet, soit $\alpha\in \Xi$. Les d\'efinitions entra\^{\i}nent que $\alpha(x'')> \alpha(x')=c_{\alpha,{\cal F}'}$. A fortiori, $\alpha(x)> \alpha(x')> c_{\alpha,{\cal F}'}$. Si $\alpha\in \Sigma_{{\cal F}}$, on a $c_{\alpha,{\cal F}}=\alpha(x)> c_{\alpha,{\cal F}'}$ donc $c_{\alpha,{\cal F}}^-\geq c_{\alpha,{\cal F}'}$ et $K_{{\cal F}}^+\cap U_{\alpha}(F)=U_{\alpha,c_{\alpha,{\cal F}}^-}\supset U_{\alpha,c_{\alpha,{\cal F}'}}$. Si $\alpha\not\in \Sigma_{{\cal F}}$, $c_{\alpha,{\cal F}}$ est le plus grand \'el\'ement de $\Gamma_{\alpha}$ qui soit strictement inf\'erieur \`a $\alpha(x)$, donc $c_{\alpha,{\cal F}}\geq c_{\alpha,{\cal F}'}$. On a encore $K_{{\cal F}}^+\cap U_{\alpha}(F)=U_{\alpha,c_{\alpha,{\cal F}}}\supset U_{\alpha,c_{\alpha,{\cal F}'}}$. Cela d\'emontre l'assertion. On peut alors \'ecrire
$$ I_{{\cal F},{\cal F}'}(f,f')=\int_{(K_{{\cal F}}^{\nu}\cap K_{{\cal F}'}^{\nu})/V}\int_{V}f_{{\cal F}}(gv)f'_{{\cal F}'}(gv)\, dv\, dg.$$
On vient de voir que $V\subset K_{{\cal F}}^+$ et $f_{{\cal F}}$ est invariante par ce groupe donc l'int\'egrale se r\'ecrit
$$ I_{{\cal F},{\cal F}'}(f,f')=\int_{(K_{{\cal F}}^{\nu}\cap K_{{\cal F}'}^{\nu})/V}f_{{\cal F}}(g)\int_{V}f'_{{\cal F}'}(gv)\, dv\, dg.$$
Pour $g$ intervenant dans cette int\'egrale, notons $\bar{g}$ son image dans ${\bf G}^{\nu}_{{\cal F}'}(k_{F})$. On a en fait $\bar{g}\in {\bf P}^{\nu}(k_{F})$. A une constante positive pr\`es provenant des mesures, l'int\'egrale int\'erieure n'est autre que
$$\sum_{\bar{u}\in {\bf U}_{{\bf P}}(k_{F})}f'(\bar{g}\bar{ u}).$$
Ceci est nul car $f'$ est cuspidale. Donc $I_{{\cal F},{\cal F}'}(f,f')=0$, ce qui d\'emontre la premi\`ere assertion de (4).

Supposons maintenant ${\cal F}={\cal F}'$. Alors $M=M_{{\cal F},\nu}=M_{{\cal F}',\nu}=G$ puisque $({\cal F}',\nu)\in Fac^*_{max}(G,A)$. D'o\`u ${\cal F}^M={\cal F}={\cal F}'$. Un calcul imm\'ediat donne  
la derni\`ere formule de (4), ce qui d\'emontre cette assertion.

Apr\`es ces pr\'eliminaires, d\'emontrons la proposition. Notons $N$ un ensemble de repr\'esentants du quotient
  $$Norm_{L}(A)(F)\backslash Norm_{G}(A)(F)/(Norm_{G}(A)(F)\cap K_{{\cal F}}^0).$$
  On sait que $G(F)$ est union disjointe des ensembles $L(F)nK_{{\cal F}}^0$ quand $n$ d\'ecrit $N$.  Fixons un sous-groupe parabolique $Q \in {\cal P}(L)$. On d\'efinit comme en \ref{lesespaces} une pseudo-mesure sur $Q(F)\backslash G(F)$. Pour tout $n\in N$, introduisons la fonction $\varphi_{n}$ sur $G(F)$ \`a support dans $Q(F)nK_{{\cal F}}^0$ telle que $\varphi(lunk)=\delta_{Q}(l)$ pour tous $l\in L(F)$, $u\in U_{Q}(F)$, $k\in K_{{\cal F}}^0$. On note $m_{n}$ la valeur de son int\'egrale contre la pseudo-mesure sur $Q(F)\backslash G(F)$. 
  Puisque $f_{{\cal F}}$ est invariante par $K_{{\cal F}}^0$, la d\'efinition de $D^G_{f'}(f_{{\cal F}})$ se r\'ecrit
  $$(5) \qquad D^G_{f'}(f_{{\cal F}})= \sum_{n\in N}m_{n}D_{f'}^L((^n(f_{{\cal F}}))_{U_{Q}}).$$
  Pour $n\in N^0$, l'action de $n$ transporte ${\cal F}$ en une facette $n{\cal F}$, le Levi $M=M_{{\cal F},\nu}$ en le Levi $nMn^{-1}=M_{n{\cal F},\nu}$, la facette ${\cal F}^M$ en une facette $n{\cal F}^M\in Fac(nMn^{-1},A)$ et $f$ en une fonction $^nf\in C_{cusp}(({\bf nMn^{-1}})_{n{\cal F}^M}^{\nu})$. On a $^n(f_{{\cal F}})=(^nf)_{n{\cal F}}$. On applique (3) en y rempla\c{c}ant ${\cal F}$ par $n{\cal F}$ et $f_{{\cal F}}$ par $(^nf)_{n{\cal F}}$. Cette assertion nous dit que la contribution $D_{f'}^L((^nf_{{\cal F}})_{U_{Q}})$ de $n$ \`a la formule (5) est nulle si $nMn^{-1}\not\subset L$. On note $N^0$ le sous-ensemble des $n\in N$ tels que $nMn^{-1}\subset L$. Pour $n\in N^0$,   cette contribution  est $mes(U_{Q}(F)\cap K_{n{\cal F}}^+)D_{f'}^L((^nf)_{(n{\cal F})^L})$. Par d\'efinition, on a
  $$(6) \qquad D_{f'}^L((^nf)_{(n{\cal F})^L})=\int_{A_{L}(F)\backslash L(F)} \int_{L(F)}(^nf)_{(n{\cal F})^L}(x^{-1}lx)f'_{{\cal F}'}(l)\,dl\,dx.$$
 Notons $ N'_{n}$ un  ensemble de repr\'esentants du quotient
 $$A_{L}(F)(K_{{\cal F}'}^0\cap Norm_{L}(A)(F))\backslash Norm_{L}(A)(F)/(K_{(n{\cal F})^L}^0\cap Norm_{L}(A)(F)).$$
 On utilise maintenant la d\'ecomposition en union disjointe
  $$L(F)=\sqcup_{n'\in N'_{n}} A_{L}(F)K_{{\cal F}'}^0n'K_{(n{\cal F})^L}^0.$$
  Pour $n'\in  N'_{n}$, notons $m'_{n'}$ la mesure de 
  $$A_{L}(F)\backslash  A_{L}(F)K_{{\cal F}'}^0n'K_{(n{\cal F})^L}^0.$$
  Puisque les fonctions $(^nf)_{(n{\cal F})^L}$, resp. $f'_{{\cal F}'}$, sont invariantes par conjugaison par $K_{(n{\cal F})^L}^0$, resp. $K_{{\cal F}'}^0$, la formule (6) se r\'ecrit
  $$ D_{f'}^L((^nf)_{(n{\cal F})^L})=\sum_{n'\in  N'_{n}}m'_{n'}\int_{L(F)}{^{n'}(^nf)_{(n{\cal F})^L}}(l)f'_{{\cal F}'}(l)\,dl.$$
  Pour $n'$ apparaissant ci-dessus, on a $^{n'}(^nf)_{(n{\cal F})^L}=(^{n'n}f)_{(n'n{\cal F})^L}$ avec des d\'efinitions similaires aux pr\'ec\'edentes. Remarquons que $M_{n'n{\cal F},\nu}=n'nM(n'n)^{-1}$.  La formule ci-dessus se r\'ecrit
  $$ D_{f'}^L((^nf)_{(n{\cal F})^L})=\sum_{n'\in  N'_{n}}m'_{n'}I_{(n'n{\cal F})^L,{\cal F}'}(^{n'n}f,f').$$
  Pour tout $n'$, on applique (4) en y rempla\c{c}ant ${\cal F}$ par $n'n{\cal F}$ et $f$ par $^{n'n}f$. Cette assertion nous permet de nous limiter aux $n'$ tels que $n'nM(n'n)^{-1}=L$, $(n'n{\cal F})^L=n'n{\cal F}^M={\cal F}'$. Notons $N_{n}^{'0}$ l'ensemble des $n'\in N'_{n}$ v\'erifiant ces conditions. En utilisant la derni\`ere assertion de (4), on obtient
  $$ D_{f'}^L((^nf)_{(n{\cal F})^L})=\sum_{n'\in  N^{'0}_{n}}m'_{n'}mes(K_{{\cal F}'}^0)<\overline{^{n'n}f},f'>.$$
  En rassemblant ces calculs, on obtient
  $$(7) \qquad D^G_{f'}(f_{{\cal F}})= \sum_{n\in N^0}\sum_{n'\in N_{n}^{'0}}m_{n}mes(U_{Q}(F)\cap K_{n{\cal F}}^+) m'_{n'}mes(K_{{\cal F}'}^0)<\overline{^{n'n}f},f'>.$$
  On v\'erifie que l'application
  $$\begin{array}{ccc}\{(n,n'); n\in N^0, n'\in N_{n}^{'0}\}&\to &Norm_{G}(A)(F)\\ (n,n')&\mapsto &n'n\\ \end{array}$$
  est une bijection de l'ensemble de d\'epart sur le quotient
  $$A_{L}(F)(K_{{\cal F}'}^0\cap Norm_{G}(A)(F))\backslash N(M,{\cal F}^{M},L,{\cal F}') .$$
  Si $N(M,{\cal F}^{M},L,{\cal F}')$ est vide, on a donc $D^G_{f'}(f_{{\cal F}})=0$, ce qui d\'emontre le (i) de l'\'enonc\'e. Supposons $N(M,{\cal F}^M,L,{\cal F}')\not=\emptyset$. On peut supposer que $\underline{N}(M,{\cal F}^M,L,{\cal F}')$ est \'egal \`a l'ensemble des $n'n$ pour $(n,n')$  comme ci-dessus. 
Pour un tel couple, on v\'erifie les propri\'et\'es suivantes

$K_{n{\cal F}}^0\cap Q(F)=(K_{n{\cal F}}^0\cap L(F))(K_{n{\cal F}}^0\cap U_{Q}(F))$, $K_{n{\cal F}}^0\cap L(F)=K_{(n{\cal F})^L}^0$, $K_{n{\cal F}}^0\cap U_{Q}(F)=K_{n{\cal F}}^+\cap U_{Q}(F)$, cf. \cite{W2} lemme 6;

$m_{n}=mes(K_{{\cal F}}^0)mes(K_{(n{\cal F})^L}^0)^{-1}mes(K_{n{\cal F}}^+\cap U_{Q}(F))^{-1}$;

$mes(K_{(n{\cal F})^L}^0)=mes(K_{{\cal F}'}^0)$;

$m'_{n'}=mes(A_{L}(F)\backslash A_{L}(F)K_{{\cal F}'}^0)=mes(K_{{\cal F}'}^0)mes(A_{L}(F)_{c})^{-1}$.

En utilisant ces propri\'et\'es, la formule (7) devient celle du (ii) de l'\'enonc\'e. $\square$

\subsection{Un espace de fonctions\label{unespacedefonctions}}

  Pour $({\cal F},\nu)\in Fac^*(G)$, notons $C({\bf G}_{{\cal F}}^{\nu})$ l'espace des fonctions \`a valeurs complexes sur ${\bf G}_{{\cal F}}^{\nu}(k_{F})$ et notons ${\cal E}_{{\cal F}}^{\nu}$ l'espace des fonctions sur $G(F)$, \`a support dans $K_{{\cal F}}^{\nu}$ et invariantes par multiplication \`a droite ou \`a gauche par  $K_{{\cal F}}^+$. Ces espaces s'identifient comme en \ref{facettesfonctions}: pour une fonction $f\in  C({\bf G}_{{\cal F}}^{\nu})$, on identifie $f$ \`a une fonction sur $K_{{\cal F}}^{\nu}/K_{{\cal F}}^+$, on la rel\`eve en une fonction sur $K_{{\cal F}}^{\nu}$ puis on l'\'etend en une fonction sur $G(F)$ par $0$ hors de $K_{{\cal F}}^{\nu}$. On note $f_{{\cal F}}$ la fonction obtenue, qui appartient \`a ${\cal E}_{{\cal F}}^{\nu}$. On note ${\cal E}(G)$ le sous-espace de $C_{c}^{\infty}(G(F))$ engendr\'e par les ${\cal E}_{{\cal F}}^{\nu}$ quand $({\cal F},\nu)$ d\'ecrit $Fac^*(G)$. On note $I{\cal E}(G)$ l'image de ${\cal E}(G)$ dans $I(G)$.

Pour $M\in \underline{{\cal L}}_{min}$, le groupe $Norm_{G}(A)(F)\cap Norm_{G}(M)(F)$ agit sur l'ensemble $Fac^*_{max,G-comp}(M,A)$. Fixons un ensemble de repr\'esentants  $ \underline{Fac}^*_{max,G-comp}(M,A)$ des orbites. Pour chaque $({\cal F}_{M},\nu)\in \underline{Fac}^*_{max,G-comp}(M,A)$, le groupe $K_{{\cal F}_{M}}^{\dag,G}$ agit naturellement sur ${\bf M}^{\nu}_{{\cal F}_{M}}(k_{F})$ et sur l'espace de fonctions sur ce groupe. Cette action conserve l'espace $C_{cusp}({\bf M}^{\nu}_{{\cal F}_{M}})$ et on note $C_{cusp}({\bf M}^{\nu}_{{\cal F}_{M}})^{K_{{\cal F}_{M}}^{\dag,G}}$ le sous-espace des invariants. 
On pose
$${\cal E}(G,M)=\sum_{({\cal F}_{M},\nu)\in \underline{Fac}^*_{max,G-comp}(M,A)}C_{cusp}({\bf M}^{\nu}_{{\cal F}_{M}})^{K_{{\cal F}_{M}}^{\dag,G}}.$$
Pour tout $({\cal F}_{M},\nu)\in \underline{Fac}^*_{max,G-comp}(M,A)$, fixons une facette ${\cal F}_{M}^G\in Fac(G,A)$ telle que $({\cal F}_{M}^G,\nu)\in Fac^*(G,A)$, $M_{{\cal F}_{M}^G,\nu}=M$ et $({\cal F}_{M}^G)^M={\cal F}_{M}$. C'est possible d'apr\`es \cite{W2} lemme 7. Pour $f\in C_{cusp}({\bf M}^{\nu}_{{\cal F}_{M}})^{K_{{\cal F}_{M}}^{\dag,G}}$, on a d\'efini la fonction $f_{{\cal F}_{M}^G}$ sur $G(F)$. L'application $f\mapsto f_{{\cal F}_{M}^G}$ se prolonge par lin\'earit\'e en une application de ${\cal E}(G,M)$ dans ${\cal E}(G)$. 
 On note $e(G,M):{\cal E}(G,M)\to I{\cal E}(G)$ la compos\'ee de cette application et de la projection ${\cal E}(G)\to I{\cal E}(G)$. On note $ Im(e(G,M))$ l'image de $e(G,M)$. 

\ass{Lemme}{On a l'\'egalit\'e $I{\cal E}=\sum_{M\in \underline{{\cal L}}_{min}}Im(e(G,M))$.}

Preuve. Soit $({\cal F},\nu)\in Fac^*(G)$ et $f\in C({\bf G}_{{\cal F}}^{\nu})$. On veut prouver que l'image de $f_{{\cal F}}$ dans $I(G)$ appartient \`a la somme des $Im(e(G,M))$. On ne perd rien \`a supposer $f$ invariante par conjugaison par ${\bf G}_{{\cal F}}(k_{F})$.  On sait que l'on peut \'ecrire $f$ comme somme finie de fonctions "induites" $Ind_{{\bf M}^{\nu}}^{{\bf G}^{\nu}}(f_{{\bf M}^{\nu}})$ (la d\'efinition est rappel\'ee en (1) ci-dessous), o\`u ${\bf M}^{\nu}$ est un espace de Levi de ${\bf G}^{\nu}$ et $f_{{\bf M}^{\nu}}$ est une fonction cuspidale sur ${\bf M}^{\nu}(k_{F})$ invariante par conjugaison par ${\bf M}(k_{F})$. On peut aussi bien supposer que $f$ est l'une de ces fonctions induites, disons $f=Ind_{{\bf M}^{\nu}}^{{\bf G}^{\nu}}(f')$. Fixons une telle fonction et un espace parabolique ${\bf P}^{\nu}$ de composante de Levi  ${\bf M}^{\nu}$. Notons $f''$ la fonction sur ${\bf G}^{\nu}(k_{F})$ qui est \`a support dans ${\bf P}^{\nu}(k_{F})$ et v\'erifie $f''(mu)=f'(m)$ pour tous $m\in {\bf M}^{\nu}(k_{F})$ et $u\in {\bf U}_{{\bf P}}(k_{F})$. Par d\'efinition, on a
$$(1) \qquad f(g)=\vert {\bf P}(k_{F})\vert ^{-1}\sum_{x\in {\bf G}(k_{F})}f''(x^{-1}gx)$$
pour tout $g\in {\bf G}^{\nu}(k_{F})$. Il en r\'esulte que l'image de $f_{{\cal F}}$ dans $I(G)$ est proportionnelle \`a celle de $f''_{{\cal F}}$. Au sous-groupe parabolique ${\bf P}$ correspond une facette ${\cal F}'\in Fac(G)$ dont l'adh\'erence contient ${\cal F}$. On a encore $({\cal F}',\nu)\in Fac^*(G)$ d'apr\`es l'existence d'un espace ${\bf P}^{\nu}$ de parabolique associ\'e ${\bf P}$. On v\'erifie que $f''_{{\cal F}}=f'_{{\cal F}'}$. Cela nous ram\`ene au probl\`eme de d\'epart o\`u l'on a remplac\'e ${\cal F}$ et $f$ par ${\cal F}'$ et $f'$. En oubliant cette construction, on peut supposer $f$ cuspidale. 

 On peut conjuguer ${\cal F}$ et $f$ par un \'el\'ement de $G(F)$, cela ne change pas l'image de $f_{{\cal F}}$ dans $I(G)$. On peut donc supposer ${\cal F}\in App(A)$. Posons $M=M_{{\cal F},\nu}$. On a d\'efini la facette ${\cal F}^M\in Imm(M_{AD})$ associ\'ee \`a ${\cal F}$. On a $\nu\in {\cal N}_{G-comp}^M({\cal F}^M)$ d'apr\`es le (i) du lemme 6 de \cite{W2}. La derni\`ere assertion de ce lemme  et l'\'egalit\'e $M=M_{{\cal F},\nu}$ entra\^{\i}nent que ${\cal F}^{M,\nu}$ est r\'eduit \`a un point, c'est-\`a-dire que $({\cal F}^M,\nu)\in Fac^*_{max}(M,A)$. De nouveau, par conjugaison, on peut supposer $M\in \underline{{\cal L}}_{min}$ et $({\cal F}^M,\nu)\in \underline{Fac}^*_{max,G-comp}(M,A)$. Les espaces ${\bf G}_{{\cal F}}^{\nu}$ et ${\bf M}_{{\cal F}^M}^{\nu}$ s'identifient, on peut consid\'erer $f$ comme un \'el\'ement de $C_{cusp}({\bf M}_{{\cal F}^M}^{\nu})$.
 On peut moyenner $f$ par l'action par conjugaison de $K_{{\cal F}^M}^{\dag,G}$, cela ne change pas l'image de $f_{{\cal F}}$ dans $I(G)$. Alors $f\in {\cal E}(G,M)$. A $f$, on a associ\'e ci-dessus une fonction $f_{({\cal F}^M)^G}\in {\cal E}(G)$. Pour achever la preuve, il suffit de prouver que $f_{{\cal F}}$ a m\^eme image dans $I(G)$ que $f_{({\cal F}^M)^G}$. Notons simplement ${\cal F}'=({\cal F}^M)^G$.  La facette ${\cal F}'$ v\'erifie exactement les m\^emes propri\'et\'es que ${\cal F}$, \`a savoir  $({\cal F}',\nu)\in Fac^*(G,A)$, $M_{{\cal F}',\nu}=M$ et ${\cal F}^{'M}={\cal F}^M$.   Si ${\cal F}'={\cal F}$, on a termin\'e: $f_{{\cal F}}=f_{{\cal F}'}$. Supposons ${\cal F}\not={\cal F}'$. Notons $X$ l'image r\'eciproque de ${\cal F}^{M,\nu}$ dans $Imm(G_{AD})$. C'est un espace affine sous ${\cal A}_{M}/{\cal A}_{G}$ puisque ${\cal F}^{M,\nu}$ est r\'eduit \`a un point. Les ensembles ${\cal F}^{\nu}\cap X$ et ${\cal F}^{'\nu}\cap X$ sont ouverts dans $X$, cf. \cite{W2} lemme 7. Fixons des points $x\in {\cal F}^{\nu}\cap X$ et $x'\in {\cal F}^{'\nu}\cap X$. Le d\'ecoupage de $App(A)$ en facettes induit un d\'ecoupage du segment $[x,x']$ en points et segments ouverts. C'est-\`a-dire que l'on a des points $x_{i}$ pour $i=1,...,n$, des facettes ${\cal F}_{i}$ pour $i=1,...,n-1$ et des facettes ${\cal F}''_{i}$ pour $i=1,...,n$ de sorte que
 
 $[x,x_{1}[=[x,x']\cap {\cal F}$, $]x_{n},x']=[x,x']\cap {\cal F}'$; 
 
 $]x_{i},x_{i+1}[=[x,x']\cap {\cal F}_{i}$ pour $i=1,...,n-1$;
 
 $\{x_{i}\}=[x,x']\cap {\cal F}''_{i}$ pour $i=1,...,n$. 
 
 Un \'el\'ement $k\in K_{{\cal F}^M}^{\nu}$ appartient \`a la fois \`a $K_{{\cal F}}^{\nu}$ et $K_{{\cal F}'}^{\nu}$. Donc son action sur $Imm(G)$ fixe tout le segment $[x,x']$. Il en r\'esulte que $\nu$ appartient \`a ${\cal N}({\cal F}_{i})$ et \`a ${\cal N}({\cal F}''_{i})$ pour tout $i$ et que l'on peut aussi bien ajouter des indices $\nu$ dans les \'egalit\'es ci-desssus, par exemple $]x_{i},x_{i+1}[=[x,x']\cap {\cal F}_{i}^{\nu}$. On d\'ecrit les facettes comme en \ref{description}. Parce que $M=M_{{\cal F},\nu}\supset M_{{\cal F}}$, on a $\Sigma_{{\cal F}}=\Sigma_{{\cal F}^M}$ et $c_{\alpha,{\cal F}}=c_{\alpha,{\cal F}^M}$ pour tout $\alpha\in \Sigma_{{\cal F}^M}$. Il en est de m\^eme en rempla\c{c}ant ${\cal F}$ par ${\cal F}'$.  Fixons $i=1,...,n-1$. L'ensemble $\Sigma_{{\cal F}_{i}}$ est celui des $\alpha$ tels que, pour un point $y\in {\cal F}_{i}$, ou pour tout point $y\in {\cal F}_{i}$, $\alpha(y)$ appartient \`a $\Gamma_{\alpha}$. Puisque $[x,x']\cap {\cal F}_{i}$ est ouvert dans $[x,x']$, il revient au m\^eme de dire que $\alpha\in \Sigma_{{\cal F}}\cap \Sigma_{{\cal F}'}$ et $c_{\alpha,{\cal F}}=c_{\alpha,{\cal F}'}$. La description ci-dessus entra\^{\i}ne alors $\Sigma_{{\cal F}_{i}}=\Sigma_{{\cal F}}=\Sigma_{{\cal F}'}=\Sigma_{{\cal F}^M}$. Donc $M_{{\cal F}_{i}}=M_{{\cal F}}$. L'espace ${\cal A}_{M_{{\cal F}_{i},\nu}}$ est le sous-espace  des points fixes de l'action de l'\'el\'ement $k$ ci-dessus dans ${\cal A}_{M_{{\cal F}_{i}}}$. Il en est de m\^eme en rempla\c{c}ant ${\cal F}_{i}$ par ${\cal F}$. Puisque $M_{{\cal F}_{i}}=M_{{\cal F}}$, cela entra\^{\i}ne $M_{{\cal F}_{i},\nu}=M_{{\cal F},\nu}=M$. Enfin, on a \'evidemment ${\cal F}_{i}^M={\cal F}^M$. Cela d\'emontre que chaque facette ${\cal F}_{i}$ v\'erifie les m\^emes conditions que ${\cal F}$ et ${\cal F}'$. On peut d\'efinir une fonction $f_{{\cal F}_{i}}$ pour tout $i$ et il nous suffit de d\'emontrer successivement que les fonctions $f_{{\cal F}}$, $f_{{\cal F}_{1}}$,..., $f_{{\cal F}_{n-1}}$ et $f_{{\cal F}'}$ ont m\^eme image dans $I(G)$. On est ramen\'e au cas de facettes cons\'ecutives, autrement dit au cas $n=1$ dans notre construction.  Notons simplement ${\cal F}''={\cal F}''_{1}$. Les facettes ${\cal F}$ et ${\cal F}'$ d\'eterminent des espaces paraboliques ${\bf P}^{\nu}$ et ${\bf P}^{'\nu}$ et ${\cal G}_{{\cal F}''}^{\nu}$. Ils ont un espace de Levi commun ${\bf M}^{\nu}$, qui est isomorphe \`a ${\bf G}_{{\cal F}}^{\nu}$ et \`a ${\bf G}_{{\cal F}'}^{\nu}$:  ${\bf M}^{\nu}(k_{F})$ est l'image de $K_{{\cal F}^M}^{\nu}\cap K_{{\cal F}''}^{\nu}$ dans ${\bf G}_{{\cal F}''}^{\nu}(k_{F})$. Par la m\^eme construction que ci-dessus, les fonctions $f_{{\cal F}}$ et $f_{{\cal F}'}$ ont m\^eme image dans $I(G)$ que la fonction  $f''_{{\cal F}''}$, o\`u $f''=Ind_{{\bf M}^{\nu}}^{{\bf G}_{{\cal F}''}^{\nu}}(f)$, le point \'etant que cette induction ne d\'epend pas de l'espace parabolique choisi. Cela ach\`eve la d\'emonstration. $\square$

 \subsection{Deux espaces en dualit\'e\label{deuxespaces}}

  Pour $M\in \underline{{\cal L}}_{min}$, notons ${\cal D}_{G-comp}(G,M)$  l'image dans ${\cal D}_{G-comp}(G)$ de
  l'espace  $\boldsymbol{{\cal D}}_{cusp,G-comp}(M)$. Si $M=G$, ${\cal D}_{G-comp}(G;G)$ est simplement l'espace $D_{cusp}(G)$ d\'ej\`a d\'efini. On a ${\cal D}_{G-comp}(G)=\sum_{M\in \underline{{\cal L}}_{min}}{\cal D}_{G-comp}(G,M)$. Pour tout $M\in {\cal L}_{min}$,
    on v\'erifie que le sous-espace 
   $$(1) \qquad \prod_{({\cal F}_{M},\nu)\in \underline{Fac}^*_{max,G-comp}(M,A)}C_{cusp}({\bf M}^{\nu}_{{\cal F}_{M}})^{K_{{\cal F}_{M}}^{\dag,G}}\subset \boldsymbol{{\cal D}}_{G-comp}(G)$$
   s'envoie bijectivement sur ${\cal D}_{G-comp}(G,M)$.

\ass{Proposition}{(i) Pour tout $M\in \underline{{\cal L}}_{min}$, l'application $e(G,M)$ est injective et la restriction de l'application $D^G$ \`a ${\cal D}_{G-comp}(G,M)$ est injective.

(ii)   $I{\cal E}$, resp. $D^G({\cal D}_{G-comp}(G))$,  est somme directe des $Im(e(G,M))$, resp.   $D^G({\cal D}_{G-comp}(G,M))$, quand $M$ d\'ecrit $\underline{{\cal L}}_{min}$. 

(iii) Soient $M,L\in \underline{{\cal L}}_{min}$. L'application
$$\begin{array}{ccc}{\cal D}_{G-comp}(G,M)&\to&Im(e(G;L))^*\\ {\bf f}&\mapsto& (D^G_{{\bf f}})_{\vert Im(e(G;L))}\\ \end{array}$$
est nulle si $M\not=L$ et est un isomorphisme si $M=L$.}

\bigskip

Preuve. Tous ces espaces sont sommes ou produits d'espaces index\'es par des couples $({\cal F}_{M},\nu)$. Pour des raisons de support, on peut fixer $\nu\in {\cal N}$ et remplacer tous les espaces par les sous-espaces analogues o\`u on se limite au couples $({\cal F}_{M},\nu')$ tels que $\nu'=\nu$, cf. \ref{variantes}. On gagne que ces espaces sont de dimension finie. En fait, d'apr\`es (1), les espaces ${\cal D}_{G-comp}(G,M)$ et ${\cal E}(G,M)$ deviennent isomorphes. Fixons une base $({\bf f}_{i,M})_{i=1,...,n_{M}}$ de  ${\cal E}(G,M)$,  r\'eunion de bases orthogonales des espaces $C_{cusp}({\bf M}^{\nu}_{{\cal F}_{M}})^{K_{{\cal F}_{M}}^{\dag,G}}$ intervenant. En appliquant la proposition \ref{calcul} et la remarque (2) qui la suit, on voit
  que la matrice
$$\big(D^G_{{\bf f}_{i,M}}(e(G,L)(\overline{{\bf f}_{j,L}}))\big)_{M,L\in \underline{{\cal L}}_{min},i=1,...,n_{M},j=1,...,n_{L}}$$
est diagonale, de coefficients diagonaux non nuls. Compte tenu du lemme pr\'ec\'edent, cela entra\^{\i}ne toutes les assertions de l'\'enonc\'e. $\square$

\subsection{Un corollaire\label{uncorollaire}}
\ass{Corollaire}{(i) L'application lin\'eaire $D^G:{\cal D}_{G-comp}(G)\to I(G)^*$ est injective.

(ii) L'application lin\'eaire compos\'ee ${\cal D}_{cusp}(G)\stackrel{D^G}{\to}I(G)^*\to I_{cusp}(G)^*$ est injective.}

Cela r\'esulte de la proposition pr\'ec\'edente, en se rappelant que ${\cal E}(G;G)$ est form\'e d'apr\`es sa d\'efinition de fonctions cuspidales.  
$\square$ 

\subsection{Injectivit\'e de $D^G$\label{injectivite}}
 Pour $M\in \underline{{\cal L}}_{min}$, notons ${\cal D}(G,M)$  l'image dans ${\cal D}(G)$ de
  l'espace  $\boldsymbol{{\cal D}}_{cusp}(M)$.  On a ${\cal D}(G)=\sum_{M\in \underline{{\cal L}}_{min}}{\cal D}(G,M)$. 

\ass{Proposition}{(i) L'application $D^G:{\cal D}(G)\to I(G)^*$ est injective.

(ii) Pour tout $n\in {\mathbb Z}$, l'espace $D^G({\cal D}(G))\cap Ann^nI(G)^*$ est \'egal \`a la somme des $D^G({\cal D}(G,M))$ o\`u $M$ parcourt les \'el\'ements de $\underline{{\cal L}}_{min}$ tels que $a_{M}\geq n$. }

Preuve. Notons ${\cal D}^n(G)$ la somme des ${\cal D}(G,M)$ o\`u $M$ parcourt l'ensemble de Levi indiqu\'e ci-dessus. Par construction, on a
$${\cal D}^n(G)={\cal D}^{n+1}(G)+\sum_{M\in \underline{{\cal L}}^n_{min}}{\cal D}(G,M).$$
Il r\'esulte  des d\'efinitions que, pour tout $M\in \underline{{\cal L}}^n_{min}$, ${\cal D}(G,M)$ est l'image naturelle dans ${\cal D}(G)$ de ${\cal D}_{cusp}(M)^{W^G(M)}$. 
On a aussi $D^G({\cal D}^n(G))\subset D^G({\cal D}(G))\cap Ann^nI(G)^*$. On en d\'eduit une suite d'applications
$$\oplus_{M\in \underline{{\cal L}}^n_{min}} {\cal D}_{cusp}(M)^{W^G(M)}\to \oplus_{M\in \underline{{\cal L}}^n_{min}}{\cal D}(G,M)\to {\cal D}^n(G)/{\cal D}^{n+1}(G)\stackrel{\delta^n}{\to}$$
$$ (D^G({\cal D}(G))\cap Ann^nI(G)^*)/(D^G({\cal D}(G))\cap Ann^{n+1}I(G)^*)\to Gr^nI(G)^*\simeq \oplus_{M\in \underline{{\cal L}}^n_{min}}I_{cusp}(M)^{*W^G(M)}.$$
 Le corollaire pr\'ec\'edent, appliqu\'e aux Levi $M\in\underline{{\cal L}}_{min}^n$, implique que l'application compos\'ee est injective. Les deux premi\`eres applications  de la suite sont surjectives. Il en r\'esulte que $\delta_{n}$ est injective. Pour $n=a_{G}$, on a $D^G({\cal D}(G))\cap Ann^{a_{G}}I(G)^*=D^G({\cal D}(G))=D^G({\cal D}^n(G))$. L'injectivit\'e de $\delta_{n}$ pour tout $n$ entra\^{\i}ne alors par r\'ecurrence que, pour tout $n$, on a $D^G({\cal D}(G))\cap Ann^nI(G)^*=D^G({\cal D}^n(G))$ et $Ker(D^G)\subset {\cal D}^n(G)$. Pour $n=a_{M_{min}}+1$, cette derni\`ere relation implique l'injectivit\'e de $D^G$. $\square$

\subsection{Variantes avec caract\`ere central\label{varianteaveccaracterecentral}}
Pour un groupe topologique ab\'elien et localement compact $X$, nous appelons caract\`ere de $X$ un homomorphisme continu de $X$ dans ${\mathbb C}^{\times}$. 
Soit $\xi$ un caract\`ere de $A_{G}(F)$. On note $C_{c,\xi}^{\infty}(G(F))$ l'espace des fonctions sur $G(F)$, \`a valeurs complexes, localement constantes, telles que $f(ag)=\xi(a)^{-1}f(g)$ pour tous $a\in A_{G}(F)$ et tout $g\in G(F)$ et telles que l'image dans $A_{G}(F)\backslash G(F)$ du support de $f$ soit compacte. Pour $f\in C_{c}^{\infty}(G(F))$, notons $f_{\xi}$ la fonction d\'efinie par $f_{\xi}(g)=\int_{A_{G}(F)}f(ag)\xi(a)\,da$. L'application lin\'eaire $f\mapsto f_{\xi}$ est une surjection de $C_{c}^{\infty}(G(F))$ sur $ C_{c,\xi}^{\infty}(G(F))$. De m\^eme que l'on a d\'efini $I(G)$, on  d\'efinit l'espace $I_{\xi}(G)$. L'action de $A_{G}(F)$ sur $C_{c}^{\infty}(G(F))$ se descend en une action sur $I(G)$.  Le groupe $A_{G}(F)$ agit dualement sur $I(G)^*$. L'espace  $I_{\xi}(G)^*$ dual de $I_{\xi}(G)$ s'identifie au sous-espace des \'el\'ements de $I(G)^*$ qui se transforment selon $\xi$. Pr\'ecis\'ement, un \'el\'ement $d$ de ce sous-espace s'identifie \`a l'\'el\'ement de $I_{\xi}(G)^*$ qui envoie $f_{\xi}$ sur $d(f)$ pour tout $f\in C_{c}^{\infty}(G(F))$.  

Soit ${\cal F}\in Fac(G)$. Le groupe $A_{G}(F)$ est contenu dans $K_{{\cal F}}^{\dag}$. Il agit par multiplication sur ce groupe. La multiplication par $a\in A_{G}(F)$ envoie un sous-ensemble $K_{{\cal F}}^{\nu}$ sur $K_{{\cal F}}^{\nu+\nu_{a}}$, o\`u on a pos\'e $\nu_{a}=w_{G}(a)$. Si $({\cal F},\nu)\in Fac^*_{max}(G)$, on a aussi $({\cal F},\nu+\nu_{a})\in Fac^*_{max}(G)$. Dans ce cas, la multiplication par $a$ se descend en un isomorphisme encore not\'e $a:C_{cusp}({\bf G}_{{\cal F}}^{\nu})\to C_{cusp}({\bf G}_{{\cal F}}^{\nu+\nu_{a}})$. Ces isomorphismes d\'efinissent une action de $A(F)$ sur l'espace $\boldsymbol{{\cal D}}_{cusp}(G)$ (on prendra soin de la distinguer de l'action par conjugaison, qui est triviale). L'action du sous-groupe $A_{G,tu}(F)$ est triviale. L'action se descend en une action sur ${\cal D}_{cusp}(G)$ que l'on note $(a,{\bf f})\mapsto {\bf f}^{a}$. Soit $\xi$ un caract\`ere mod\'er\'ement ramifi\'e de $A_{G}(F)$, c'est-\`a-dire trivial sur $A_{G,tu}(F)$. On note ${\cal D}_{cusp,\xi}(G)$ le sous-espace des \'el\'ements ${\bf f}\in {\cal D}_{cusp}(G)$ tels que ${\bf f}^{a}=\xi(a){\bf f}$ pour tout $a\in A_{G}(F)$. 

Plus g\'en\'eralement, on d\'efinit de m\^eme les variantes "\`a caract\`ere central" de beaucoup d'objets d\'ej\`a d\'efinis (au sens ci-dessus: il s'agit d'un caract\`ere $\xi$ de $A_{G}(F)$). On les note en ajoutant un indice $\xi$ dans les notations. Notons en particulier l'\'egalit\'e

(1) $D^G({\cal D}(G))\cap I(G)^*_{\xi}=D^G({\cal D}_{\xi}(G))$, 

\noindent qui r\'esulte ais\'ement de l'injectivit\'e de $D^G$.
Plus g\'en\'eralement, si $\xi_{1},...,\xi_{n}$ sont des caract\`eres  de $A_{G}(F)$, 

(2) $D^G({\cal D}(G))\cap (\sum_{i=1,...,n}I(G)^*_{\xi_{i}})=D^G(\sum_{i=1,...,n}{\cal D}_{\xi_{i}}(G))$. 

Preuve. On peut supposer les $\xi_{i}$ distincts. Par interpolation, si $d$ appartient au membre de gauche, les composantes $d_{i}$ de $d$ dans chaque $I(G)^*_{\xi_{i}}$ sont combinaisons lin\'eaires finies de translat\'es $d^{a}$ pour des $a\in A_{G}(F)$. Or ces $d^{a}$ appartiennent tous \`a $D^G({\cal D}(G))$. Donc $d_{i}\in D^G({\cal D}(G))\cap I(G)^*_{\xi_{i}}$ et il reste \`a appliquer l'\'egalit\'e (1). $\square$

 \section{El\'ements compacts, $p'$-\'el\'ements}
 
 \subsection{Retour sur les \'el\'ements topologiquement unipotents\label{retour}}
  Pour $x\in G_{tu}(F)$, l'homomorphisme $n\mapsto x^n$ de ${\mathbb Z}$ dans $G(F)$ se prolonge  en un homomorphisme continu $z\mapsto x^z$ de ${\mathbb Z}_{p}$ dans $G(F)$. Il en r\'esulte que, 
  
  (1) si $J\subset G(F)$ est un sous-groupe ferm\'e et s'il existe un entier $c\geq1$ premier \`a $p$ tel que $x^c\in J$, alors $x\in J$.
  
  On a
 
 (2) l'application naturelle $G_{tu}(F)\to G_{AD,tu}(F)$ est surjective; ses fibres sont les orbites de l'action par multiplication de $(Z(G)^0)_{tu}(F)$ dans $G_{tu}(F)$. 
 
 Preuve.   
 Notons $\pi:G\to G_{AD}$ l'homomorphisme naturel. L'hypoth\`ese $(Hyp)(G)$ entra\^{\i}ne que $G_{AD}(F)/\pi(G(F))$ est d'ordre premier \`a $p$. Pour $y\in G_{AD,tu}(F)$, il y a donc un entier $c\geq1$ premier \`a  $p$ tel que $y^c\in \pi(G(F))$.  D'apr\`es (1), $y$ appartient \`a $\pi(G(F))$. Soit $x\in G(F)$ tel que $\pi(x)=y$. On d\'ecompose $x=x_{ss}x_{u}$, o\`u $x_{ss}$ est semi-simple, $x_{u}$ est unipotent et $x_{ss}$ et $x_{u}$ commutent. Fixons un sous-tore maximal $T$ de $G$ tel que $x_{ss}\in T(F)$ et fixons une extension galoisienne finie $F'$ de $F$,  de degr\'e $d$ premier \`a $p$, de sorte que $T$ soit d\'eploy\'e sur $F'$. Posons $T_{ad}=T/Z(G)$. L'\'el\'em\'ent $\pi(x_{ss})$ est topologiquement unipotent. La projection $T_{tu}(F')\to T_{ad,tu}(F')$ s'identifie \`a $X_{*}(T)\otimes_{{\mathbb Z}}(1+\mathfrak{p}_{F'})\to X_{*}(T_{ad})\otimes_{{\mathbb Z}}(1+\mathfrak{p}_{F'})$. Elle est surjective car l'hypoth\`ese $(Hyp)(G)$ implique que l'image de $X_{*}(T)$ dans $X_{*}(T_{ad})$ est un sous-groupe d'indice fini premier \`a $p$. 
  On peut donc choisir $x_{1}\in T_{tu}(F')$ tel que $\pi(x_{1})=\pi(x_{ss})$. L'\'el\'ement $x_{2}=Norme_{F'/F}(x_{1})$ appartient \`a $T_{tu}(F)$ et on a $\pi(x_{2})=\pi(x_{ss}^d)$.  D'apr\`es (1), cela implique qu'il existe $x_{3}\in T_{tu}(F)$ tel que $\pi(x_{3})=\pi(x_{ss})$. L'\'el\'ement $x_{3}x_{u}$ appartient \`a $G_{tu}(F)$ et v\'erifie $\pi(x_{3}x_{u})=y$. Cela d\'emontre la premi\`ere assertion.

  L'hypoth\`ese $(Hyp)(G)$ entra\^{\i}ne que $Z(G)/Z(G)^0$ est d'ordre premier \`a $p$.  La seconde assertion s'en d\'eduit par le m\^eme argument qui prouvait ci-dessus que $y$ appartenait \`a $\pi(G(F))$. $\square$
 
 On a
 
 (3) pour $x\in G_{tu}(F)$, l'image r\'eciproque de $Z_{G_{AD}}(x_{ad})$ dans $G$ est \'egale \`a $Z_{G}(x)$.
 
 Preuve. Notons simplement $Z_{G}(x_{ad})$ cette image r\'eciproque. Il est bien connu que $Z_{G}(x)$ est un sous-groupe distingu\'e d'indice fini dans $Z_{G}(x_{ad})$ et l'hypoth\`ese $(Hyp)(G)$ implique que cet indice est premier \`a $p$. Soit $g\in Z_{G}(x_{ad})$. Il existe donc un entier $c\geq1$ premier \`a $p$ tel que $g^c$ commute \`a $x$. Puisque $g\in Z_{G}(x_{ad})$ il existe $z\in Z(G)$  tel que $ g^{-1}xg=zx$. On  a alors $z^c=1$. Or, puisque $x$ et $g^{-1}xg$ sont topologiquement unipotents, $z$ est lui-aussi topologiquement unipotent. L'\'egalit\'e $z^c=1$ entra\^{\i}ne alors que $z=1$, donc $g\in Z_{G}(x)$. $\square$ 
 
   \subsection{El\'ements topologiquement nilpotents et exponentielle\label{elementstopologiquement}}
  On note $\mathfrak{g}$ l'alg\`ebre de Lie de $G$. On appelle "conjugaison" par $G$ l'action adjointe de $G$ dans $\mathfrak{g}$. Pour $g\in G$, on  note cette action de $g$, soit $ad(g)$, soit $X\mapsto gXg^{-1}$.  On note $\mathfrak{g}_{reg}$ l'ensemble des \'el\'ements semi-simples r\'eguliers de $\mathfrak{g}$.
  
   Soit  $X\in \mathfrak{g}(F)$, fixons un sous-tore maximal $T$ de $G$ tel que $X_{ss}\in \mathfrak{t}(F)$ et fixons une extension finie $F'$ de $F$ telle que $T$ soit d\'eploy\'e sur $F'$. L'\'el\'ement $X$  est dit topologiquement nilpotent si et seulement si $ \chi(X_{ss})\in \mathfrak{p}_{F'}$ pour tout $\chi\in X^*(T)$ (cela ne d\'epend pas du choix de $T$). Une autre caract\'erisation est la suivante. Fixons une sous-alg\`ebre d'Iwahori $\mathfrak{b}$ de $\mathfrak{g}$, c'est-\`a-dire $\mathfrak{b}=\mathfrak{k}_{{\cal F}}$ pour une facette ${\cal F}\in Fac(G)$ de dimension maximale. Notons $\mathfrak{u}$ son radical pro-$p$-nilpotent, c'est-\`a-dire $\mathfrak{u}=\mathfrak{k}_{{\cal F}}^+$ pour la m\^eme facette. Alors $X$ est topologiquement nilpotent si et seulement si $X$ est conjugu\'e par un \'el\'ement de $G(F)$ \`a un \'el\'ement de $\mathfrak{u}(F)$. On note $\mathfrak{g}_{tn}(F)$ l'ensemble des \'el\'ements topologiquement nilpotents de $\mathfrak{g}(F)$. 
   
   On peut d\'efinir une application exponentielle $exp$ qui envoie un voisinage de $0$ dans $\mathfrak{g}(F)$ sur un voisinage de $1$ dans $G(F)$ et qui est \'equivariante pour les actions par conjugaison de $G(F)$. Sous l'hypoth\`ese $(Hyp)(G)$, ces voisinages sont les plus grands possibles, c'est-\`a-dire qu'on dispose de l'application exponentielle
   $$exp:\mathfrak{g}_{tn}(F)\to G_{tu}(F)$$
   \'equivariante pour les actions par conjugaison de $G(F)$ et qui
    un hom\'eomorphisme entre les ensembles de d\'epart et d'arriv\'ee, cf. \cite{DR} appendice B pour cette propri\'et\'e et celles ci-dessous. 
    
    Pour ${\cal F}\in Fac(G)$, l'exponentielle se restreint en des hom\'eomorphismes de $\mathfrak{g}_{tn}(F)\cap \mathfrak{k}_{{\cal F}}$ sur $G_{tu}(F)\cap K_{{\cal F}}^0$
    et de $\mathfrak{k}_{{\cal F}}^+$ sur $K_{{\cal F}}^+$. Il s'en d\'eduit une bijection
    $$(\mathfrak{g}_{tn}(F)\cap \mathfrak{k}_{{\cal F}})/\mathfrak{k}_{{\cal F}}^+\to (G_{tu}(F)\cap K_{{\cal F}}^0)/K_{{\cal F}}^+.$$
    
    On peut prolonger les paires $K_{{\cal F}}^0\supset K_{{\cal F}}^+$ et $\mathfrak{k}_{{\cal F}}\supset \mathfrak{k}_{{\cal F}}^+$
en des suites $(K_{{\cal F},n})_{n\in {\mathbb N}}$ et $(\mathfrak{k}_{{\cal F},n})_{n\in {\mathbb N}}$ de sorte que

$K_{{\cal F}}^0=K_{{\cal F},0}$, $K_{{\cal F}}^+ =K_{{\cal F},1}$, $K_{{\cal F},n}\supset K_{{\cal F},n+1}$, $\cap_{n\in {\mathbb N}}K_{{\cal F},n}=\{1\}$;

$\mathfrak{k}_{{\cal F}}=\mathfrak{k}_{{\cal F},0}$, $\mathfrak{k}_{{\cal F}}^+ =\mathfrak{k}_{{\cal F},1}$, $\mathfrak{k}_{{\cal F},n}\supset \mathfrak{k}_{{\cal F},n+1}$, $\cap_{n\in {\mathbb N}}\mathfrak{k}_{{\cal F},n}=\{0\}$;

pour $n\geq1$, $K_{{\cal F},n}$ est un sous-groupe distingu\'e de $K_{{\cal F}}^{\dag}$ et $\mathfrak{k}_{{\cal F},n}$ est un $\mathfrak{o}_{F}$-id\'eal de $\mathfrak{k}_{{\cal F}}$;

pour $n\geq1$, $K_{{\cal F},n}=exp(\mathfrak{k}_{{\cal F},n})$ et l'exponentielle se r\'eduit en un isomorphisme de groupes de $\mathfrak{k}_{{\cal F},n}/\mathfrak{k}_{{\cal F},n+1}$ sur $K_{{\cal F},n}/K_{{\cal F},n+1}$.

 \subsection{ El\'ements $p'$-compacts\label{elementsp}}
 
 Soit $x\in G(F)$.     On dit que $x$ est 
 
  $p'$-compact si et seulement s'il existe un entier $c\geq1$ premier \`a $p$ tel que $x^c=1$;
  
 $p'$-compact mod $Z(G)$ si et seulement si l'image  $x_{ad}$ de $x$ dans $G_{AD}(F)$ est   $p'$-compacte.

  Pour $x\in G(F)$, fixons un sous-tore maximal de $G$ contenant   $x_{ss}$  et fixons une extension finie de $F$ telle que $T$ soit d\'eploy\'ee sur $F'$. Alors:

  $x$ est $p'$-compact si et seulement $x=x_{ss}$ et il existe un entier $c\geq1$ premier \`a $p$ tel que $\chi(x)^c=1$ pour tout $\chi\in X^*(T)$;
  
  $x$ est $p'$-compact mod $Z(G)$ si et seulement $x=x_{ss}$ et il existe un entier $c\geq1$ premier \`a $p$ tel que $\chi(x)^c=1$ pour tout $\chi\in X^*(T)$ dont la restriction \`a $Z(G)$ est triviale.

  Puisque les classes de conjugaison de sous-tores maximaux de $G$ sont en nombre fini, on peut ci-dessus choisir $F'$ ind\'ependant de $x$. Puisque $F^{'\times}$ ne contient qu'un nombre fini de racines de l'unit\'e, on voit qu'il n'y a qu'un nombre fini de classes de conjugaison par $G(F)$ d'\'el\'ements $p'$-compacts. A fortiori,  il existe un entier $c\geq1$ premier \`a $p$ tel que $x^c=1$ pour tout \'el\'ement $p'$-compact de $G(F)$. On a
  
 (1)  tout \'el\'ement compact $x\in G(F)$ s'\'ecrit de fa\c{c}on unique $x=x_{p'}x_{tu}$, o\`u $x_{p'}$ est $p'$-compact, $x_{tu}\in G_{tu}(F)$ et $x_{p'}$ et $x_{tu}$ commutent; de plus, $x_{p'}$ et $x_{tu}$ appartiennent \`a $\overline{x^{{\mathbb Z}}}$.
  
  Preuve. Pour $c$ comme ci-dessus, il r\'esulte de ce qui pr\'ec\`ede que $x^c$ est topologiquement unipotent. Soit $c'$ l'inverse de $c$ dans ${\mathbb Z}_{p}$. On pose $x_{tu}=(x^c)^{c'}$ et $x_{p'}=xx_{tu}^{-1}$. Ces termes v\'erifient les conditions requises et on voit que ce sont les seules solutions possibles. $\square$
  
  On d\'eduit de (1) et de \ref{retour} (2) et (3) que
  
  (2) tout \'el\'ement $x\in G(F)$ compact mod $Z(G)$ s'\'ecrit $x=x_{p'}x_{tu}$, o\`u $x_{p'}$ est $p'$-compact, $x_{tu}\in G_{tu}(F)$ et $x_{p'}$ et $x_{tu}$ commutent;  le groupe $(Z(G)^0)_{tu}(F)$ agit sur l'ensemble des solutions:  un \'el\'ement  $z\in (Z(G)^0)_{tu}(F)$ envoie le couple $(x_{p'},x_{tu})$ sur $(x_{p'}z^{-1},x_{tu}z)$; les solutions forment une unique orbite pour cette action; de plus, pour toute solution $(x_{p'},x_{tu})$, les deux \'el\'ements $x_{p'}$ et $x_{tu}$ appartiennent \`a l'adh\'erence du groupe $Z(G)(F)x^{\mathbb Z}$.
  
 \subsection{$p'$-\'el\'ements\label{pelements}}

Soit $x\in G$. On lui associe un sous-groupe parabolique $Q[x]$ de $G$ et une composante de Levi $L[x]$ de $Q[x]$ de la fa\c{c}on suivante, cf. \cite{C}. L'\'el\'ement $x_{ss}$ agit par conjugaison sur $\mathfrak{g}$. On fixe une extension galoisienne finie $F'$ de $F$ telle que toutes les valeurs propres appartiennent \`a $F^{'\times}$. On note $\Sigma$ l'ensemble de ces valeurs propres et, pour $\sigma\in \Sigma$, $\mathfrak{g}_{\sigma}$ l'espace propre. Alors l'alg\`ebre de Lie $\mathfrak{q}[x]$ est la somme des $\mathfrak{g}_{\sigma}$ sur les $\sigma\in \Sigma$ telles que $\vert \sigma\vert _{F'}\leq 1$ et $\mathfrak{l}[x]$ est la somme des $\mathfrak{g}_{\sigma}$ sur les $\sigma\in \Sigma$ telles que $\vert \sigma\vert _{F'}= 1$.
Le couple $(Q[x], L[x])$ \'etant uniquement d\'efini et ne d\'ependant que de $x_{ss}$, il est d\'efini sur $F$ et conserv\'e par $Z_{G}(x_{ss})$. Cela entra\^{\i}ne

(1) $Z_{G}(x_{ss})\subset L[x]$, a fortiori $ Z_{G}(x)\subset  L[x]$ et $x\in L[x]$.

Par construction, $x$, vu comme \'el\'ement de $L[x]$, est compact mod $Z(L[x])$ (on dira que $x$ est compact dans $L[x]$ mod $Z(L[x])$). Le Levi $L[x]$ est le plus grand Levi $L$ de $G$ tel que $x$ appartienne \`a $L$ et que $x$ soit compact dans $L$ mod $Z(L)$.
  
 Nous dirons que $x$ est un $p'$-\'el\'ement si et seulement si $x$ est $p'$-compact dans $L[x]$ mod $Z(L[x])$.  Avec les notations ci-dessus, cela \'equivaut \`a ce que $x=x_{ss}$ soit semi-simple et, pour tout $\sigma\in \Sigma$ telle que $\vert \sigma\vert _{F}=1$, $\sigma$ soit une racine de l'unit\'e d'ordre premier \`a $p$.  On note $G(F)_{p'}$  l'ensemble des $p'$-\'el\'ements de $G(F)$.  
 
 \ass{Lemme}{(i) Tout \'el\'ement $x\in G(F)$ s'\'ecrit $x=x_{p'}x_{tu}$ o\`u $x_{p'}$ est un $p'$-\'el\'ement, $x_{tu}$ est topologiquement unipotent et $x_{p'}$ et $x_{tu}$ commutent.
 
 (ii) Pour une telle d\'ecomposition, les \'el\'ements $x_{p'}$ et $x_{tu}$ appartiennent \`a $L[x](F)$ et on a $L[x_{p'}]=L[x]$.
 
 (iii) Pour une telle d\'ecomposition, on a $Z_{G}(x)=Z_{G}(x_{p'})\cap Z_{G}(x_{tu})$ et $G_{x}=(G_{x_{p'}})_{x_{tu}}$.
 
 (iv)  
  La d\'ecomposition est unique modulo l'action de $(Z(L[x])^0)_{tu}(F)$ similaire \`a celle de \ref{elementsp}(2).}
  
  Preuve. Consid\'erons une d\'ecomposition $x=yx_{tu}$ o\`u $x_{tu}$ est topologiquement unipotent et $y$ et $x_{tu}$ commutent. Montrons que
 
 (2) $y,x_{tu}\in L[x]$ et $L[x]=L[y]$.
 
 Les \'el\'ements $y$ et $x_{tu}$ commutent \`a $x$ et la premi\`ere assertion r\'esulte de (1). Notons ici $x_{ad}$, $y_{ad}$ et $x_{tu,ad}$ les images de $x$, $y$ et $x_{tu}$ dans $L[x]_{AD}(F)$. L'adh\'erence $\overline{y_{ad}^{{\mathbb Z}}}$ du  groupe engendr\'e par $y_{ad}$ est contenue dans le produit des deux groupes $\overline{x_{ad}^{{\mathbb Z}}}$ et $\overline{x_{tu,ad}^{{\mathbb Z}}}$. Ce dernier groupe est compact puisque $x_{tu}$ est topologiquement unipotent. Le premier est    compact   par d\'efinition de $L[x]$. Donc   $\overline{y_{ad}^{{\mathbb Z}}}$ est compacte, c'est-\`a-dire que $y$ est compact mod $Z(L[x])$. Or $L[y]$ est le plus grand Levi $L$ tel que $y\in L$ et $y$ soit compact mod $Z(L)$. Donc $L[x]\subset L[y]$. On obtient l'inclusion oppos\'ee en \'echangeant les r\^oles de $x$ et $y$ (on a $y=xx_{tu}^{-1}$). Cela prouve (2).

  On a vu que $x$ appartenait \`a $ L[x]$ et \'etait compact mod $Z(L[x])$. En appliquant \ref{elementsp}(2) dans le groupe $L[x]$, on obtient une d\'ecomposition $x=x_{p'}x_{tu}$ o\`u $x_{p'}, x_{tu}\in L[x](F)$, $x_{p'}$ est $p'$-compact mod $Z(L[x])$, $x_{tu}$ est topologiquement unipotent et $x_{p'}$ et $x_{tu}$ commutent. D'apr\`es (2), on a $L[x_{p'}]=L[x]$, donc $x_{p'}$ est $p'$-compact mod $Z(L[x_{p'}])$, donc c'est un $p'$-\'el\'ement. Cela d\'emontre (i).

  Pour une d\'ecomposition $x=x_{p'}x_{tu}$ comme en (i), les \'el\'ements $x_{p'}$ et $x_{tu}$ commutent \`a $x$ donc appartiennent \`a $L[x]$ d'apr\`es (1). Cela d\'emontre la premi\`ere assertion de (ii) et la seconde r\'esulte de (2). En d\'efinitive, les d\'ecompositions $x=x_{p'}x_{tu}$ v\'erifiant (i) sont exactement les d\'ecompositions dans $L[x](F)$ o\`u l'on impose que $x_{p'}$ est $p'$-compact mod $Z(L[x])$. Le (iv) r\'esulte donc de \ref{elementsp}(2) appliqu\'e dans le groupe $L[x]$.
  
 Pour (iii), les commutants $Z_{G}(x)$ et $Z_{G}(x_{p'})$ sont contenus dans $L[x]$ d'apr\`es (1) et la derni\`ere assertion de (ii). On ne perd rien \`a supposer $L[x]=G$. Il est clair que $Z_{G}(x_{p'})\cap Z_{G}(x_{tu})\subset Z_{G}(x)$.   L'image $x_{tu,ad}$ de $x_{tu}$ dans $G_{AD}$ appartient \`a l'adh\'erence du groupe engendr\'e par $x_{ad}$. L'image $g_{ad}$ d'un \'el\'ement  $g\in Z_{G}(x)$ commute donc \`a $x_{tu,ad}$. D'apr\`es \ref{retour}(2), cela implique que $g\in Z_{G}(x_{tu})$. Alors $g$ appartient forc\'ement  aussi \`a $Z_{G}(x_{p'})$. Cela d\'emontre la premi\`ere \'egalit\'e de (iii). La composante neutre de $Z_{G}(x_{p'})\cap Z_{G}(x_{tu})$ est clairement $(G_{x_{p'}})_{x_{tu}}$, d'o\`u la seconde \'egalit\'e. $\square$
 
 Pour $x\in G(F)$,  appelons $p'$-d\'ecomposition de $x$ une d\'ecomposition $x=x_{p'}x_{tu}$ v\'erifiant les conditions du (i) de l'\'enonc\'e.

 \subsection{D\'ecomposition de $G(F)$ associ\'ee aux $p'$-\'el\'ements\label{decomposition}}
 
   Pour tout $\epsilon\in G(F)_{p'}$, posons $C(\epsilon)=\epsilon G_{\epsilon,tu}(F)=\epsilon\,exp(\mathfrak{g}_{\epsilon,tn}(F))$. Evidemment, l'action de $G(F)$ par conjugaison conserve $G(F)_{p'}$ et, pour $g\in G(F)$, on a $gC(\epsilon)g^{-1}=C(g\epsilon g^{-1})$. Posons
   $$C_{G}(\epsilon)=\cup_{g\in G(F)}C(g^{-1}\epsilon g)=\cup_{g\in G(F)}g^{-1}C(\epsilon )g.$$
On a d\'ej\`a utilis\'e le discriminant de Weyl $D^G(x)$ d'un \'el\'ement $x\in G(F)$. Il y a de m\^eme un  discriminant de Weyl $D^G(X)$ pour un \'el\'ement $X\in \mathfrak{g}(F)$.

\ass{Lemme}{ 
(i)  Pour $\epsilon\in G(F)_{p'}$, l'ensemble $C_{G}(\epsilon) $ est ouvert et ferm\'e.

(ii)  Pour $\epsilon\in G(F)_{p'}$, il existe un nombre r\'eel $d^G(\epsilon)>0$ tel que, pour tout  $X\in \mathfrak{g}_{\epsilon,tn}(F)$, on ait
  $$D^G(\epsilon exp(X))=d^G(\epsilon)D^{G_{\epsilon}}(X).$$
  
  (iii) On  a les \'egalit\'es
$$G(F)=\bigcup_{\epsilon\in G(F)_{p'}}C(\epsilon)=\bigcup_{\epsilon\in G(F)_{p'}}C_{G}(\epsilon).$$

(iv)  Soient $\epsilon ,\epsilon '\in G(F)_{p'}$. Alors $C(\epsilon )=C(\epsilon ')$ ou $C(\epsilon )\cap C(\epsilon ')=\emptyset$. L'\'egalit\'e a lieu si et seulement s'il existe $z\in  (Z_(L[\epsilon ])^0)_{tu}(F)$ tel que $\epsilon '=\epsilon z$. Dans ce cas, on a $L[\epsilon ]=L[\epsilon ']$, $G_{\epsilon}=G_{\epsilon'}$ et $z\in G_{\epsilon,tu}(F)$.
}

Preuve. La propri\'et\'e suivante est bien connue: soit $x\in G(F)$ un \'el\'ement semi-simple; alors il existe un voisinage $V$ de $1$ dans $G_{x}(F)$ tel que l'ensemble $\{ g^{-1} xyg; y\in V, g\in G(F)\}$ soit ouvert et ferm\'e. On voit que, pour $\epsilon\in G(F)_{p'}$,   $C(\epsilon)$ est r\'eunion finie de tels ensembles. D'o\`u (i). 

 Les deux parties de la formule (ii) sont insensibles au remplacement de $X$ par sa partie semi-simple. On peut donc supposer $X$ semi-simple. On fixe un sous-tore maximal $T$ de $G_{\epsilon}$ tel que $X\in \mathfrak{t}(F)$. On a aussi $\epsilon\in T(F)$ puisque $T$ commute \`a $\epsilon$. Fixons une extension finie $F'$ de $F$ telle que $T$ soit d\'eploy\'e sur $F'$. On note $\Sigma$ l'ensemble des racines de $T$ dans $G$, $\Sigma^{L[\epsilon]}$, resp. $\Sigma^{G_{\epsilon}}$, le sous-ensemble des racines dans $L[\epsilon]$, resp. $G_{\epsilon}$. On a $\Sigma^{G_{\epsilon}}\subset \Sigma^{L[\epsilon]}$. Par d\'efinition,
  $$D^G(\epsilon)=\vert \prod_{\alpha\in \Sigma}(1-\alpha(\epsilon exp(X)))\vert _{F}$$
  $$D^{G_{\epsilon}}(X)=\vert \prod_{\alpha\in \Sigma^{G_{\epsilon}}}\alpha(X)\vert _{F}.$$
  Puisque $exp(X)$ est topologiquement unipotent, $\alpha(exp(X))$ est un \'el\'ement  $1+\mathfrak{p}_{F'}$  pour tout $\alpha\in \Sigma$. A fortiori $\vert \alpha(exp(X))\vert _{F'}=1$. 
  Si $\alpha\in \Sigma-\Sigma^{L[\epsilon]}$, on a $\vert \alpha(\epsilon)\vert _{F'}\not=1$ par d\'efinition de $L[\epsilon]$. Puisque  $\vert \alpha(exp(X))\vert _{F'}=1$, on a $\vert 1-\alpha(\epsilon exp(X))\vert _{F'}=\vert 1-\alpha(\epsilon)\vert _{F'}$ qui est non nul. Si $\alpha\in \Sigma^{L[\epsilon]}$, $\alpha(\epsilon)$ est une racine de l'unit\'e d'ordre premier \`a $p$ puisque $\epsilon$ est $p'$-compact mod $Z(L[\epsilon])$. Cette racine est \'egale \`a $1$ si et seulement si $\alpha\in \Sigma^{G_{\epsilon}}$. Si $\alpha\in \Sigma^{L[\epsilon]}-\Sigma^{G_{\epsilon}}$, on a $ \alpha(\epsilon)\in \mathfrak{o}_{F'}^{\times}-(1+\mathfrak{p}_{F'})$ et $\alpha(exp(X))\in 1+\mathfrak{p}_{F'}$ donc encore   $\vert 1-\alpha(\epsilon exp(X))\vert _{F'}=\vert 1-\alpha(\epsilon)\vert _{F'}=1$. Enfin, si $\alpha\in \Sigma^{G_{\epsilon}}$, on a $\alpha(\epsilon)=1$ donc $\vert 1-\alpha (\epsilon exp(X))\vert _{F'}=\vert 1-\alpha(exp(X))\vert _{F'}=\vert 1-exp(\alpha(X))\vert _{F'}$ avec $\alpha(X)\in \mathfrak{p}_{F'}$, donc $\vert 1-\alpha (\epsilon exp(X))\vert _{F'}=\vert \alpha(X)\vert _{F'}$. L'assertion (ii) en r\'esulte. $\square$

    Le (i) du lemme \ref{pelements} implique la premi\`ere \'egalit\'e du (iii) d'o\`u trivialement la deuxi\`eme.
  
  Soient $\epsilon,\epsilon'\in G(F)_{p'}$.  Supposons $C(\epsilon )\cap C(\epsilon ')\not=\emptyset$ et fixons $x$ dans cette intersection. On peut \'ecrire $x=\epsilon u=\epsilon 'u'$ avec $u\in G_{\epsilon ,tu}(F)$ et $u'\in G_{\epsilon ',tu}(F)$. Ces deux d\'ecompositions sont des $p'$-d\'ecompositions.  D'apr\`es le (iv) du lemme \ref{pelements}, il existe   $z\in  (Z(L[\epsilon ])^0)_{tu}(F)$ tel que $\epsilon '=\epsilon z$. Inversement, s'il existe un tel $z$, l'assertion (2) de \ref{pelements} appliqu\'ee \`a $x=\epsilon'=\epsilon z$ nous dit que $L[\epsilon ]=L[\epsilon ']$. D'apr\`es l'assertion (1) de \ref{pelements}, on a $G_{\epsilon }=L[\epsilon ]_{\epsilon }$ et $G_{\epsilon '}=L[\epsilon ']_{\epsilon '}=L[\epsilon ]_{\epsilon '}$. Puisque $z\in Z(L[\epsilon ])^0(F)$, l'\'egalit\'e $\epsilon '=\epsilon z$ entra\^{\i}ne que $L[\epsilon ]_{\epsilon '}=L[\epsilon ]_{\epsilon }$ donc aussi $G_{\epsilon '}=G_{\epsilon }$. On a aussi $z\in G_{\epsilon ,tu}(F)$ et on en d\'eduit que $C(\epsilon ')=\epsilon 'G_{\epsilon ',tu}(F)=\epsilon zG_{\epsilon ,tu}(F)=\epsilon G_{\epsilon ,tu}(F)=C(\epsilon )$. Cela d\'emontre (iv). $\square$
  
  \subsection{Le cas d'un Levi\label{lecasdunlevi}}
  
    Soit $M$ un Levi de $G$. On a dans $M$ la m\^eme propri\'et\'e que dans $G$, \`a savoir
  $$M(F)=\bigcup_{\epsilon\in M(F)_{p'}}C^M(\epsilon),$$
  o\`u $C^M(\epsilon)=\epsilon exp(\mathfrak{m}_{\epsilon,tn}(F))$. Le lemme ci-dessous \'enonce  une propri\'et\'e plus fine.
  
  \ass{Lemme}{Soit $M$ un Levi de $G$. Alors
  
  (i) $G(F)_{p'}\cap M(F)\subset M(F)_{p'}$;
  
  (ii) $M(F)=\bigcup_{\epsilon\in G(F)_{p'}\cap M(F)}C^M(\epsilon)$.}
  
  Preuve.  Soit $\epsilon\in G(F)_{p'}\cap M(F)$.  Pour d\'emontrer que $\epsilon\in M(F)_{p'}$, on doit prouver que toute valeur propre de $ad(\epsilon)$ dans $\mathfrak{m}(F)$ qui est de valeur absolue $1$ dans une extension convenable de $F$, est une racine de l'unit\'e d'ordre premier \`a $p$. Mais une valeur propre de $ad(\epsilon)$ dans $\mathfrak{m}(F)$ est aussi une valeur propre de $ad(\epsilon)$ dans $\mathfrak{g}(F)$. La propri\'et\'e voulue r\'esulte du fait que $\epsilon\in G(F)_{p'}$.

 Soit $x\in M(F)$.     
    Ecrivons une $p'$-d\'ecomposition de $x$ dans $G(F)$: $x=\epsilon exp(X)$, o\`u $X\in \mathfrak{g}_{\epsilon,tn}(F)$. On a $L[x]=L[\epsilon]$ d'apr\`es le (ii) du lemme \ref{pelements}. On a $A_{M}\subset G_{x} $ puisque $x\in M(F)$ et $G_{x}\subset L[x]$ d'apr\`es \ref{pelements} (1). Donc $Z(L[\epsilon])=Z(L[x])$ est contenu dans le commutant de $A_{M}$, qui est $M$. 
 On fixe un entier $c\geq1$ premier \`a $p$ tel que  $\epsilon^{c}\in Z(L[\epsilon])(F)$.  Alors $exp(cX)=x^c\epsilon^{-c}$ appartient \`a $M(F)$. Cela entra\^{\i}ne comme toujours $ X\in \mathfrak{m}_{tn}(F)$ donc aussi $\epsilon=x exp(-X)\in M(F)$ et $X\in \mathfrak{m}_{\epsilon,tn}(F)$.  La d\'ecomposition $x=\epsilon exp(X)$ montre que $x\in C^M(\epsilon)$.   $\square$

  \subsection{Un lemme sur les classes de conjugaison  et les \'el\'ements $p'$-compacts mod $Z(G)$ \label{commutant}}
 
 \ass{Lemme}{
 Soit $\epsilon\in G(F)$ un \'el\'ement $p'$-compact mod $Z(G)$ et soit ${\cal F}\in Fac(G)$. Supposons $\epsilon\in K_{{\cal F}}^{\dag}$. Soit $X\in \mathfrak{g}_{\epsilon,tn}(F)\cap \mathfrak{k}_{{\cal F}}$. Alors, pour tout \'el\'ement $x\in \epsilon exp(X) K_{{\cal F}}^+$, il existe $Y\in \mathfrak{g}_{\epsilon}(F)\cap \mathfrak{k}_{{\cal F}}^+$ tel que $x$ soit conjugu\'e \`a $\epsilon exp(X+Y)$ par un \'el\'ement de $K_{{\cal F}}^+$.}
 
 La preuve est standard, on la rappelle pour \^etre complet.  
 
 Preuve. Fixons un entier $c\geq1$ premier \`a $p$ tel que $\epsilon^c\in Z(G)(F)$. On introduit les deux polyn\^omes 
 $$P(T)=c^{-1}(T^{c-1}+...+T+1),\,\, Q(T)=c^{-1}(1-T)(c-1+(c-2)T+...+T^{c-2})$$
 \`a coefficients dans ${\mathbb Z}_{p}$. On a $P(T)+Q(T)=1$. 
  L'op\'erateur $ad(\epsilon)$  dans $\mathfrak{g}(F)$ est semi-simple et ses valeurs propres  dans $\bar{F}$ sont des racines $c$-i\`emes de l'unit\'e. Alors l'espace $\mathfrak{g}$ se d\'ecompose en somme de l'espace propre associ\'e \`a la racine $1$, qui est $\mathfrak{g}_{\epsilon}$, et de l'espace somme des espaces propres associ\'es aux racines diff\'erentes de $1$, notons-le $\mathfrak{g}_{\not=1}$. L'op\'erateur $P(ad(\epsilon))$, resp. $Q(ad(\epsilon))$, est le projecteur sur l'espace $\mathfrak{g}_{\epsilon}$, resp. $\mathfrak{g}_{\not=1}$, relativement \`a cette d\'ecomposition.
 On introduit des suites $(K_{{\cal F},n})_{n\in {\mathbb N}}$ et $(\mathfrak{k}_{{\cal F},n})_{n\in {\mathbb N}}$ comme en \ref{elementstopologiquement}. L'action $ad(\epsilon)$ par conjugaison conserve $\mathfrak{k}_{{\cal F},n}$ pour tout $n$. On a donc
 $$\mathfrak{k}_{{\cal F},n}=\mathfrak{k}_{{\cal F},n,\epsilon}\oplus \mathfrak{k}_{{\cal F},n,\not=1},$$
 o\`u $\mathfrak{k}_{{\cal F},n,\epsilon}=\mathfrak{k}_{{\cal F},n}\cap \mathfrak{g}_{\epsilon}(F)=P(ad(\epsilon))(\mathfrak{k}_{{\cal F}_{n}})$ et $\mathfrak{k}_{{\cal F},n,\not=1}=\mathfrak{k}_{{\cal F},n}\cap \mathfrak{g}_{\not=1}(F)=Q(ad(\epsilon))(\mathfrak{k}_{{\cal F}_{n}})$. Posons $\bar{\mathfrak{k}}_{{\cal F},n}=\mathfrak{k}_{{\cal F},n}/\mathfrak{k}_{{\cal F},n+1}$. L'op\'erateur $ad(\epsilon)$ se r\'eduit en un op\'erateur de cet espace, que l'on note encore $ad(\epsilon)$. On a encore
 $$(1) \qquad \bar{\mathfrak{k}}_{{\cal F},n}=\bar{\mathfrak{k}}_{{\cal F},n,\epsilon}\oplus \bar{\mathfrak{k}}_{{\cal F},n,\not=1},$$
 o\`u $\bar{\mathfrak{k}}_{{\cal F},n,\epsilon}=\mathfrak{k}_{{\cal F},n,\epsilon}/\mathfrak{k}_{{\cal F},n+1,\epsilon}=P(ad(\epsilon))(\bar{\mathfrak{k}}_{{\cal F},n})$ et $\bar{\mathfrak{k}}_{{\cal F},n,\not=1}=\mathfrak{k}_{{\cal F},n,\not=1}/\mathfrak{k}_{{\cal F},n+1,\not=1}=Q(ad(\epsilon))(\bar{\mathfrak{k}}_{{\cal F},n})$.   
 
   On va prouver par r\'ecurrence sur $n\geq1$ qu'il existe $Y_{n}\in  \mathfrak{g}_{\epsilon}(F)\cap \mathfrak{k}_{{\cal F}}^+$, $Z_{n}\in \mathfrak{k}_{{\cal F},n}$ et $k_{n}\in K_{{\cal F}}^+$ tels que
 $k_{n}^{-1}xk_{n}=\epsilon exp(X+Y_{n})exp(Z_{n})$. Pour $n=1$, il suffit de prendre $Y_{1}=0$, $k_{1}=1$ et pour $Z_{1}$ l'\'el\'ement de $\mathfrak{k}_{{\cal F}}^+$ tel que $x=\epsilon exp(X)exp(Z_{1})$. Supposons ces termes d\'efinis au rang $n$. Posons $x_{n}=\epsilon exp(X+Y_{n})$. Notons $\bar{Z}_{n}$ la r\'eduction de $Z_{n}$ dans $\bar{\mathfrak{k}}_{{\cal F},n}$ et d\'ecomposons $\bar{Z}_{n}$ en $\bar{Z}_{n,\epsilon}+\bar{Z}_{n,\not=1}$ conform\'ement \`a la d\'ecomposition (1).  Relevons $\bar{Z}_{n,\epsilon}$ en un \'el\'ement  $Y'_{n}\in \mathfrak{k}_{{\cal F},n,\epsilon}$. Les op\'erateurs $ad(exp(X+Y_{n}))$ et $ad(x_{n})$ conservent eux-aussi les espaces $\mathfrak{k}_{{\cal F},n'}$ pour tout $n'$ et ils se descendent en des op\'erateurs sur $\bar{\mathfrak{k}}_{{\cal F},n}$. Parce que $exp(X+Y_{n})$ commute \`a $\epsilon$, ces op\'erateurs commutent \`a $ad(\epsilon)$ donc pr\'eservent la d\'ecomposition (1). Sur l'espace $\bar{\mathfrak{k}}_{{\cal F},\not=1}$, les valeurs propres de $ad(\epsilon)$ sont toutes diff\'erentes de $1$ tandis que l'op\'erateur $ad(X+Y_{n})$ est nilpotent. Il en r\'esulte que $1-ad(x_{n}^{-1})$ est inversible sur cet espace. On peut donc fixer un \'el\'ement $Y''_{n}\in \mathfrak{k}_{{\cal F},n}$ tel que $\bar{Z}_{n,\not=1}$ soit la r\'eduction de $(1-ad(x_{n}^{-1})(Y''_{n})$. Posons $h_{n}=exp(Y''_{n})$. C'est un \'el\'ement de $K_{{\cal F},n}$.On a 
 $$h_{n}x_{n}exp(Z_{n})h_{n}^{-1}=x_{n}exp(ad(x_{n}^{-1})(Y''_{n}))exp(Z_{n})exp(-Y''_{n}).$$
 Mais on voit que $exp(ad(x_{n}^{-1})(Y''_{n}))exp(Z_{n})exp(-Y''_{n})\in exp(Y'_{n})K_{{\cal F},n+1}$.
 Il existe $Y_{n+1}\in  \mathfrak{g}_{\epsilon}(F)\cap \mathfrak{k}_{{\cal F}}^+ $ tel que $exp(X+Y_{n})exp(Y'_{n})=exp(X+Y_{n+1})$. En posant $k_{n+1}=h_{n}k_{n}$, on obtient 
 $k_{n+1}^{-1}xk_{n+1}=\epsilon exp(X+Y_{n+1})z_{n+1}$ avec $z_{n+1}\in K_{{\cal F},n+1}$ et il reste \`a prendre pour $Z_{n+1}$ l'\'el\'ement de $\mathfrak{k}_{{\cal F},n+1}$ tel que $z_{n+1}=exp(Z_{n+1})$ pour obtenir les \'el\'ements cherch\'es au rang $n+1$. 
 
 On peut extraire des suites $(Y_{n})_{n\geq1}$ et $(k_{n})_{n\geq1}$ des sous-suites convergentes (en fait, on n'en a m\^eme pas besoin, les suites d\'efinies ci-dessus convergent). A la limite, on obtient des \'el\'ements $Y\in \mathfrak{g}_{\epsilon}(F)\cap \mathfrak{k}_{{\cal F}}^+$ et $k\in K_{{\cal F}}^+$ tels que $k^{-1}xk=\epsilon exp(X+Y)$. Cela ach\`eve la preuve. $\square$
 
  \subsection{Un lemme sur les \'el\'ements $p'$-compacts et les espaces de Levi\label{unlemme}}
  \ass{Lemme}{
 Soit $\epsilon\in G(F)$ un \'el\'ement $p'$-compact mod $Z(G)$ et soit $({\cal F},\nu)\in Fac^*(G)$. Supposons $\epsilon\in K_{{\cal F}}^{\nu}$. Notons $\bar{\epsilon}$ la r\'eduction de $\epsilon$ dans ${\bf G}_{{\cal F}}^{\nu}(k_{F})$. Soit ${\bf P}^{\nu}$ un espace parabolique de ${\bf G}_{{\cal F}}^{\nu}$. Supposons $\bar{\epsilon}\in {\bf P}^{\nu}(k_{F})$. Alors il existe une composante de Levi ${\bf M}^{\nu}$ de ${\bf P}^{\nu}$ telle que $\bar{\epsilon}\in {\bf M}^{\nu}(k_{F})$.}
 
 La preuve s'inspire de celle du lemme 9 de \cite{P}.
 
 Preuve. On note ${\bf P}$ le sous-groupe parabolique de ${\bf G}_{{\cal F}}$ associ\'e \`a ${\bf P}^{\nu}$. L'espace $\boldsymbol{\mathfrak{u}}_{{\bf P}}$ poss\`ede une filtration finie $(\boldsymbol{\mathfrak{u}}_{i})_{i=1,...,n}$  d\'efinie par r\'ecurrence par $\boldsymbol{\mathfrak{u}}_{1}=\boldsymbol{\mathfrak{u}}_{{\bf P}}$ et, pour $i\geq1$, $\boldsymbol{\mathfrak{u}}_{i+1}$ est l'espace engendr\'e par les $[U_{1},U_{i}]$ pour $U_{1}\in \boldsymbol{\mathfrak{u}}_{1}$ et $U_{i}\in \boldsymbol{\mathfrak{u}}_{i}$. La filtration se termine par $\boldsymbol{\mathfrak{u}}_{n}=\{0\}$. On note ${\bf U}_{i}=exp(\boldsymbol{\mathfrak{u}}_{i})$. Les ${\bf U}_{i}$ sont des sous-groupes distingu\'es de ${\bf U}_{{\bf P}}$, les quotients ${\bf U}_{i}/{\bf U}_{i+1}$ sont ab\'eliens et isomorphes aux $\boldsymbol{\mathfrak{u}}_{i}/\boldsymbol{\mathfrak{u}}_{i+1}$. On a ${\bf U}_{n}=\{1\}$. Nous allons prouver par r\'ecurrence sur $i$ que
 
 (1) il existe une composante de Levi ${\bf M}_{i}$  et un \'el\'ement $u_{i}\in {\bf U}_{i}(k_{F})$ tels que $ad(\bar{\epsilon})({\bf M}_{i})=ad(u_{i})({\bf M}_{i})$.
 
 Pour $i=1$, on choisit pour ${\bf M}_{1}$ n'importe quelle composante de Levi de ${\bf P}$. Puisque $\bar{\epsilon}\in {\bf P}^{\nu}(k_{F})$, $ad(\bar{\epsilon})({\bf M}_{1})$ est encore une telle composante de Levi de ${\bf P}$ et on sait que deux telles composantes sont conjugu\'es par un \'el\'ement de ${\bf U}_{{\bf P}}(k_{F})={\bf U}_{1}(k_{F})$. On choisit pour $u_{1}$ l'\'el\'ement qui conjugue ${\bf M}_{1}$ en $ad(\bar{\epsilon})({\bf M}_{1})$. Supposons le probl\`eme r\'esolu au rang $i\geq1$. 
 Fixons un entier $c\geq1$ premier \`a $p$ tel que $\epsilon^c\in Z(G)$. On voit par r\'ecurrence sur un entier $m\geq1$ que 
 $$ad(\bar{\epsilon}^m)({\bf M}_{i})=ad( u_{i,m})({\bf M}_{i}),$$
 o\`u $u_{i,m}=ad(\bar{\epsilon}^{m-1})(u_{i})...ad(\bar{\epsilon})(u_{i})u_{i}$. Pour $m=c$, on obtient ${\bf M}_{i}=ad(u_{i,c})({\bf M}_{i})$, donc $u_{i,c}=1$. Remarquons que $ad(\bar{\epsilon})$ conserve la filtration $({\bf U}_{j})_{j=1,...,n}$. En notant $X_{i}$ la r\'eduction de $u_{i}$ dans ${\bf U}_{i}(k_{F})/{\bf U}_{i+1}(k_{F})$ et en \'ecrivant ce groupe additivement, on obtient
 $$ad(\bar{\epsilon}^{c-1})(X_{i})+...+ad(\bar{\epsilon})(X_{i})+X_{i}=0.$$
 Par le m\^eme argument que dans la preuve pr\'ec\'edente, cela entra\^{\i}ne l'existence d'un \'el\'ement $Y_{i}\in  {\bf U}_{i}(k_{F})/{\bf U}_{i+1}(k_{F})$ tel que  $X_{i}=Y_{i}-ad(\bar{\epsilon})(Y_{i})$. On rel\`eve $Y_{i}$ en $v_{i}\in {\bf U}_{i}(k_{F})$ et on pose ${\bf M}_{i+1}=ad(v_{i})({\bf M}_{i})$. On calcule
 $$ad(\bar{\epsilon})({\bf M}_{i+1})=ad(u_{i+1})({\bf M}_{i+1}),$$
 o\`u $u_{i+1}=ad(\bar{\epsilon})(v_{i})u_{i}v_{i}^{-1}$. Par construction de $v_{i}$, on voit que $u_{i+1}\in {\bf U}_{i+1}(k_{F})$, ce qui r\'esout le probl\`eme en $i+1$.
 
 Pour $i=n$, le Levi ${\bf M}={\bf M}_{n}$ est conserv\'e par $ad(\bar{\epsilon})$. Alors ${\bf M}^{\nu}=\bar{\epsilon}{\bf M}$ est une composante de Levi de ${\bf P}^{\nu}$ contenant  $\bar{\epsilon}$. $\square$

  \subsection{Points fixes dans $Imm(G_{AD})$ d'un \'el\'ement $p'$-compact mod $Z(G)$\label{pointsfixes}}
Soit $\epsilon \in G(F)$ un \'el\'ement $p'$-compact mod $Z(G)$. Notons $M$ le commutant dans $G$ du tore d\'eploy\'e $A_{G_{\epsilon}}$. C'est un Levi de $G$. On a $\epsilon\in M$ et $G_{\epsilon}\subset M$ ($\epsilon$ et $G_{\epsilon}$ commutent \`a $A_{G_{\epsilon}}$). On a aussi $A_{M}=A_{G_{\epsilon}}$. En effet, puisque $G_{\epsilon}\subset M$, on a $ A_{G_{\epsilon}}\subset M$ et, puisque $M$ commute \`a $A_{G_{\epsilon}}$ par construction de $M$, on a $A_{G_{\epsilon}}\subset A_{M}$. Inversement, puisque $\epsilon\in M$, on a $A_{M}\subset G_{\epsilon}$ et, puisque $A_{M}$ commute \`a $M$, donc aussi au sous-ensemble $G_{\epsilon}$, on a $A_{M}\subset A_{G_{\epsilon}}$. Notons $J$ la r\'eunion dans $Imm(G_{AD})$ des appartements $App(A_{M'})$  associ\'es aux Levi minimaux $M'$  de $M$.  L'action de $M(F)$  sur $Imm(G_{AD})$ conserve le sous-ensemble $J$, il en est donc de m\^eme des actions de $\epsilon$ et de $G_{\epsilon}(F)$.  L'espace vectoriel ${\cal A}_{M} /{\cal A}_{G}$ agit naturellement sur $J$ (il agit sur chaque $App(A_{M'})$ et ces actions se recollent). Cette action commute \`a celles de $\epsilon$ et de $G_{\epsilon}(F)$. Notons $Imm(G_{AD})^{\epsilon}$ et $J^{\epsilon}$ les sous-ensembles de points fixes de l'action de $\epsilon$ dans $Imm(G_{AD})$, resp. $J$. L'ensemble $J^{\epsilon}$ est stable par l'action de ${\cal A}_{M}/{\cal A}_{G}$. On note $J^{\epsilon}/({\cal A}_{M}/{\cal A}_{G})$ l'ensemble quotient. L'action de $G_{\epsilon}(F)$ se descend en une action sur ce quotient. 

\ass{Proposition}{(i) On a l'\'egalit\'e $Imm(G_{AD})^{\epsilon}=J^{\epsilon}$.

(ii) L'ensemble $J^{\epsilon}/({\cal A}_{M}/{\cal A}_{G})$ muni de son action de $G_{\epsilon}(F)$ s'identifie canoniquement \`a $Imm(G_{\epsilon,AD})$.

(iii) Soit ${\cal F}\in Fac(G)$ telle que ${\cal F}\cap Imm(G_{AD})^{\epsilon}\not=\emptyset$. Alors l'image de ${\cal F}\cap Imm(G_{AD})^{\epsilon}$ est contenue dans une facette ${\cal F}'\in Fac(G_{\epsilon})$. L'action naturelle de $\epsilon$ sur $K_{{\cal F}}^0$ se descend en une action alg\'ebrique sur ${\bf G}_{{\cal F}}$ et le groupe ${\bf G}_{\epsilon,{\cal F}'}$ s'identifie \`a la composante neutre ${\bf G}_{{\cal F}}^{\epsilon,0}$ du sous-groupe des points fixes  ${\bf G}_{{\cal F}}^{\epsilon}$ par cette action. Le groupe $K_{{\cal F}'}^0$ est le groupe des $g\in G_{\epsilon}(F)\cap K_{{\cal F}}^0$ tels que $w_{G_{\epsilon}}(g)=0$. On a $K_{{\cal F}'}^+=G_{\epsilon}(F)\cap K_{{\cal F}}^+$, $\mathfrak{k}_{{\cal F}'}=\mathfrak{g}_{\epsilon}(F)\cap \mathfrak{k}_{{\cal F}}$, $\mathfrak{k}_{{\cal F}'}^+=\mathfrak{g}_{\epsilon}(F)\cap \mathfrak{k}_{{\cal F}}^+$.}

Preuve. Dans la seconde preuve du th\'eor\`eme 1.9 de \cite{PY}, G. Prasad et J.-K. Yu d\'emontrent que $Imm(G_{AD})^{\epsilon}=J^{\epsilon}$ et que cet ensemble, muni de son action de $G_{\epsilon}(F)$,  s'identifie \`a l'immeuble \'etendu du groupe $G_{\epsilon,ad}=G_{\epsilon}/Z(G)$.  L'identification est canonique \`a translations pr\`es par les \'el\'ements de ${\cal A}_{M}/{\cal A}_{G}$. Cela \'equivaut aux assertions (i) et (ii). Soit ${\cal F}\in Fac(G)$ telle que ${\cal F}\cap Imm(G_{AD})^{\epsilon}\not=\emptyset$. Soit $x\in {\cal F}\cap Imm(G_{AD})^{\epsilon}\not=\emptyset$, notons $y$ son image dans $Imm(G_{\epsilon,AD})$. Soit ${\cal F}'\in Fac(G_{\epsilon})$ la facette \`a laquelle appartient $y$. Plongeons les immeubles pour les groupes adjoints dans les immeubles \'etendus.  Comme on le sait, on peut d\'efinir un sch\'ema en groupes ${\cal G}_{x}$ d\'efini sur $\mathfrak{o}_{F}$ v\'erifiant entre autres que ${\cal G}(\mathfrak{o}_{F})$ est le sous-groupe des \'el\'ements de $G(F)$ dont l'action sur l'immeuble \'etendu $Imm(G)$ fixe $x$. On sait que la partie r\'eductive de la composante neutre de sa fibre sp\'eciale n'est autre que ${\bf G}_{{\cal F}}$. Dans le paragraphe 3 de \cite{P}, cf. en particulier les d\'efinitions de 3.1 de cette r\'ef\'erence, Prasad montre que  $\epsilon$ agit naturellement sur ${\cal G}_{x}$ et il d\'efinit la composante  neutre ${\cal G}_{x}^{\epsilon,0}$  du sous-sch\'ema des points fixes ${\cal G}_{x}^{\epsilon}$.  Il prouve que  la composante neutre ${\cal G}_{\epsilon,y}^{0}$ de ${\cal G}_{\epsilon,y}$ s'identifie \`a ce sch\'ema   ${\cal G}_{x}^{\epsilon,0}$. En passant aux parties r\'eductives des  fibres sp\'eciales, on obtient que ${\bf G}_{\epsilon,{\cal F}'}$ s'identifie \`a ${\bf G}_{{\cal F}}^{\epsilon,0}$.  Cela entra\^{\i}ne

(1) $K_{{\cal F}'}^0\subset K_{{\cal F}}^0\cap G_{\epsilon}(F)$.

En passant aux alg\`ebres de Lie, l'\'egalit\'e ${\cal G}_{\epsilon,y}^{0}={\cal G}_{x}^{\epsilon,0}$ entra\^{\i}ne aussi

(2)   $\mathfrak{k}_{{\cal F}'}=\mathfrak{k}_{{\cal F}}\cap \mathfrak{g}_{\epsilon}(F)$, $\mathfrak{k}^+_{{\cal F}'}=\mathfrak{k}^+_{{\cal F}}\cap \mathfrak{g}_{\epsilon}(F)$.

\noindent L'\'egalit\'e

(3) $K_{{\cal F}'}^+=K_{{\cal F}}^+\cap G_{\epsilon}(F)$

 \noindent se d\'eduit par l'exponentielle de la seconde \'egalit\'e de (2). 

 On peut pr\'eciser (1). D'apr\`es la d\'efinition de $K_{{\cal F}'}^0$,  $K_{{\cal F}'}^0$ est contenu dans l'ensemble des $g\in K_{{\cal F}}^0\cap G_{\epsilon}(F) $ tels que $w_{G_{\epsilon}}(g)=0$. Inversement, un tel $g$ fixe $x$ donc aussi $y$ (la projection de $Imm(G_{AD})^{\epsilon}$ sur $Imm(G_{\epsilon,AD})$ \'etant compatible avec l'action de $G_{\epsilon}(F)$). Donc $g\in K_{{\cal F}'}^{\dag}$. Puisque $w_{G_{\epsilon}}(g)=0$, on a $g\in K_{{\cal F}'}^0$. Cela d\'emontre l'assertion de l'\'enonc\'e concernant $K_{{\cal F}'}^0$. Mais alors, ce groupe ne d\'epend pas de $x$ et, d'apr\`es \ref{description} (5),  cela entra\^{\i}ne que ${\cal F}'$ est uniquement d\'etermin\'ee, ce qui est la premi\`ere assertion de (iii). $\square$

 \section{Quasi-caract\`eres}
 \subsection{Transform\'ees de Fourier et int\'egrales orbitales\label{transformees}}
 Plusieurs notions que l'on a introduites sur le groupe $G(F)$  ont des analogues sur l'alg\`ebre de Lie $\mathfrak{g}(F)$.   Pour $f\in C_{c}^{\infty}(\mathfrak{g}(F))$ et $X\in \mathfrak{g}_{reg}(F)$, on d\'efinit l'int\'egrale orbitale $I^G(X,f)$. On d\'efinit aussi l'espace $I(\mathfrak{g})$ quotient de $C_{c}^{\infty}(\mathfrak{g}(F))$ par le sous-espace des fonctions dont toutes les int\'egrales orbitales sont nulles.

 On fixe un caract\`ere continu $\psi$ de $F$ de conducteur $\mathfrak{p}_{F}$ et une forme bilin\'eaire sym\'etrique non d\'eg\'en\'er\'ee $<.,.>$ sur $\mathfrak{g}(F)$, invariante par conjugaison par $G(F)$. On sait que l'on peut la choisir telle que, pour toute ${\cal F}\in Fac(G)$, $\mathfrak{k}_{{\cal F}}^+$ soit le "dual" de $\mathfrak{k}_{{\cal F}}$, c'est-\`a-dire que $\mathfrak{k}_{{\cal F}}^+$ est l'ensemble des $X\in \mathfrak{g}(F)$ tels que $<X,Y>\in \mathfrak{p}_{F}$ pour tout $Y\in \mathfrak{k}_{{\cal F}}$. On suppose qu'il en est ainsi. On d\'efinit la transformation de Fourier $f\mapsto \hat{f}$ dans $C_{c}^{\infty}(\mathfrak{g}(F))$ par
 $$\hat{f}(X)=\int_{\mathfrak{g}(F)}f(Y)\psi(<X,Y>)\,dY,$$
 o\`u $dY$ est la mesure auto-duale. La transformation  $f\mapsto \hat{f}$ de $C_{c}^{\infty}(\mathfrak{g}(F))$ se descend en une transformation de $I(\mathfrak{g})$. 
 
 Notons $Nil(\mathfrak{g})$ l'ensemble des orbites nilpotentes dans $\mathfrak{g}(F)$. Pour tout ${\cal O}\in Nil(\mathfrak{g})$, on fixe une mesure sur ${\cal O}$ invariante par conjugaison. Cela permet de d\'efinir l'int\'egrale orbitale $I_{{\cal O}}$ sur $C_{c}^{\infty}(\mathfrak{g}(F))$, puis sa transform\'ee de Fourier $f\mapsto I_{{\cal O}}(\hat{f})$. Notons $\mathfrak{g}_{reg}$ l'ensemble des \'el\'ements semi-simples et r\'eguliers de $\mathfrak{g}$. D'apr\`es Harish-Chandra, il existe une fonction $\hat{j}({\cal O})$ localement int\'egrable sur $\mathfrak{g}(F)$ et localement constante sur $\mathfrak{g}_{reg}(F)$, de sorte que
 $$I_{{\cal O}}(\hat{f})=\int_{\mathfrak{g}(F)}f(Y)\hat{j}({\cal O},Y)\,dY$$
 pour toute $f\in C_{c}^{\infty}(\mathfrak{g}(F))$. 
 
 Notons $E$ l'espace des fonctions sur $\mathfrak{o}_{F}$ engendr\'e par les fonctions $\lambda\mapsto \vert \lambda\vert ^n$ pour $n\in {\mathbb Z}$. Il est imm\'ediat que, si deux \'el\'ements de $E$ co\"{\i}ncident sur $\mathfrak{p}_{F}^n$ pour un entier $n\geq0$, alors, ils sont \'egaux. 
 
 On sait que, pour tout ${\cal O}\in Nil(\mathfrak{g})$ et pour tout $Y\in \mathfrak{g}_{reg}(F)$, la fonction $\lambda\mapsto \hat{j}({\cal O},\lambda^2 Y)$, d\'efinie sur $\mathfrak{o}_{F}$, appartient \`a $E$. 
 
 \subsection{Quasi-caract\`eres et quasi-caract\`eres de niveau $0$\label{quasicaracteres}}
 
 Si $\theta$ est une fonction localement int\'egrable sur $G(F)$ et invariante par conjugaison, il lui est associ\'ee la distribution invariante $D$ sur $G(F)$ d\'efinie par
 $$D(f)=\int_{G(F)}f(g)\theta(g)\,dg.$$
 
 On appelle quasi-caract\`ere de $G(F)$ une distribution invariante $D$ sur $G(F)$ associ\'ee \`a une fonction $\theta_{D}$ sur $G(F)$ localement int\'egrable  et invariante par conjugaison, qui  v\'erifie la propri\'et\'e suivante:
 
 (1) pour tout $x\in G(F)$ semi-simple, il existe un voisinage $\mathfrak{V}$ de $0$ dans $\mathfrak{g}_{x}(F)$ et, pour tout ${\cal O}\in Nil(\mathfrak{g}_{x})$, il existe un nombre complexe $c_{D,{\cal O}}$ de sorte que, pour presque tout $Y\in \mathfrak{V}$, on ait l'\'egalit\'e
 $$\theta_{D}(x exp(Y))=\sum_{{\cal O}\in Nil(\mathfrak{g}_{x})}c_{D,{\cal O}}\hat{j}({\cal O},Y).$$
 
  {\bf Remarques.} (2) Le voisinage $\mathfrak{V}$ n'est \'evidemment pas uniquement d\'etermin\'e, par contre les constantes $c_{D,{\cal O}}$ le sont car les fonctions $\hat{j}({\cal O})$ sont lin\'eairement ind\'ependantes dans tout voisinage de $0$.

(3)  Appliqu\'ee \`a $x$ fortement r\'egulier, cette propri\'et\'e implique que, quitte \`a modifier $\theta_{D}$ sur un ensemble de mesure nulle, on peut supposer $\theta_{D}$ d\'efinie et localement constante sur $G_{reg}(F)$.  

(4) La condition que $\theta_{D}$ est localement int\'egrable est redondante car une fonction  invariante par conjugaison et v\'erifiant la condition (1) est automatiquement localement int\'egrable. 

(5)  Soit $D\in I(G)^*$ et soit $\alpha$ une fonction   sur $G(F)$, localement constante et invariante par conjugaison. On d\'efinit la distribution $\alpha D$ par $(\alpha D)(f)=D(\alpha f)$. Si $D$ est un quasi-caract\`ere, alors $\alpha D$ l'est aussi.
 
\bigskip

Soit $f\in C_{c}^{\infty}(G(F))$ une fonction tr\`es cuspidale. On a d\'efini la distribution $D_{f}$ en \ref{lesespaces}. On a

(6) $D_{f}$ est un quasi-caract\`ere.

Cf. \cite{W1} 5.9. La fonction $\theta_{D_{f}}$ se calcule de la fa\c{c}on suivante. Soit $x\in G_{reg}(F)$. Notons $M^x$ le commutant de $A_{G_{x}}$ dans $G$. C'est le plus petit Levi $M'$ tel que $x\in M'(F)$. De plus $x$ est elliptique dans $M^{x}(F)$. On sait d\'efinir l'int\'egrale orbitale pond\'er\'ee
$$J^G_{M^{x}}(x,f)=D^G(x)^{1/2}\int_{A_{M^{x}}(F)\backslash G(F)}f(g^{-1}xg)v_{M^{x}}(g)\,dg.$$
Le poids $v_{M^{x}}$  est calcul\'e relativement \`a un sous-groupe compact sp\'ecial de $G(F)$  fix\'e, cf. \cite{A1}, mais, parce que $f$ est tr\`es cuspidale, l'int\'egrale ci-dessus ne d\'epend pas de ce  choix, cf. \cite{W1} lemme 5.2. Pour tout tore d\'eploy\'e $A'$, posons $m(A')=mes(A'(F)_{c})$. D'apr\`es \cite{W2} 9(1), on a alors
$$(7) \qquad \theta_{D_{f}}(x)=(-1)^{a_{M^{x}}-a_{G}}D^G(x)^{-1/2}m(A_{M^{x}})m(A_{G})^{-1}J^G_{M^{x}}(x,f).$$
 
 \bigskip
 
On appelle quasi-caract\`ere de niveau $0$ une distribution invariante $D$ associ\'ee \`a une fonction $\theta_{D}$ sur $G(F)$ localement int\'egrable et invariante par conjugaison, qui v\'erifie la condition

 (8) pour tout $\epsilon\in G(F)_{p'}$, pour tout ${\cal O}\in Nil(\mathfrak{g}_{\epsilon})$, il existe un nombre complexe $c_{D,{\cal O}}$ de sorte que, pour presque tout $Y\in \mathfrak{g}_{\epsilon,tn}(F)$, on ait l'\'egalit\'e
 $$\theta_{D}(\epsilon exp(Y))=\sum_{{\cal O}\in Nil(\mathfrak{g}_{\epsilon})}c_{D,{\cal O}}\hat{j}({\cal O},Y).$$

 On a des variantes   de ces d\'efinitions pour l'alg\`ebre de Lie. Soit $D$ une distribution invariante sur $\mathfrak{g}(F)$ associ\'ee \`a une  fonction $\theta_{D}$ sur $\mathfrak{g}(F)$ localement int\'egrable et invariante par conjugaison. On dit que $D$ est un quasi-caract\`ere si $\theta_{D}$ v\'erifie la condition
 
 (9)  pour tout $X\in \mathfrak{g}(F)$ semi-simple, il existe un voisinage $\mathfrak{V}$ de $0$ dans $\mathfrak{g}_{X}(F)$ et, pour tout ${\cal O}\in Nil(\mathfrak{g}_{X})$, il existe un nombre complexe $c_{D,{\cal O}}$ de sorte que, pour presque tout $Y\in \mathfrak{V}$, on ait l'\'egalit\'e
 $$\theta_{D}(X+Y)=\sum_{{\cal O}\in Nil(\mathfrak{g}_{X})}c_{D,{\cal O}}\hat{j}({\cal O},Y).$$
 
 On dit que c'est un quasi-caract\`ere de niveau $0$ si $\theta_{D}$ est \`a support dans $\mathfrak{g}_{tn}(F)$ et que,  pour tout ${\cal O}\in Nil(\mathfrak{g})$, il existe un nombre complexe $c_{D,{\cal O}}$ de sorte que, pour presque tout $Y\in \mathfrak{g}_{tn}(F)$, on ait l'\'egalit\'e
 $$\theta_{D}(Y)=\sum_{{\cal O}\in Nil(\mathfrak{g})}c_{D,{\cal O}}\hat{j}({\cal O},Y).$$
 
 Si $f\in C_{c}^{\infty}(\mathfrak{g}(F))$ est une fonction tr\`es cuspidale, on d\'efinit la distribution $D_{f}$ sur $\mathfrak{g}(F)$ comme on l'a fait sur le groupe. Cette distribution est un quasi-caract\`ere.

\ass{Lemme}{(i) Tout quasi-caract\`ere sur $G(F)$, resp. $\mathfrak{g}(F)$, de niveau $0$ est un quasi-caract\`ere.

(ii) soit $D$ un quasi-caract\`ere sur $G(F)$; alors $D$ est de niveau $0$ si et seulement si, pour tout $\epsilon\in G(F)_{p'}$ et presque tout $X\in \mathfrak{g}_{\epsilon,tn}(F)$, la fonction $\lambda\mapsto \theta_{D}(\epsilon exp(\lambda^2X))$ d\'efinie sur $\mathfrak{o}_{F}$ appartient \`a $E$, cf. \ref{transformees} .

 (iii) Fixons une d\'ecomposition $G(F)=\bigcup_{\epsilon\in B}C_{G}(\epsilon)$, o\`u $B$ est un sous-ensemble de $G(F)_{p'}$. Soit $D$ une distribution invariante associ\'ee \`a une fonction $\theta_{D}$ qui v\'erifie la   condition (8) restreinte aux $\epsilon\in B$. Alors $D$ est un quasi-caract\`ere. }

Preuve.   Il r\'esulte de \cite{W1} lemme 6.3(iii)  que, pour tout quasi-caract\`ere $D$ de niveau $0$ sur $\mathfrak{g}(F)$, il existe    une fonction tr\`es cuspidale $f$ telle que $D$ co\"{\i}ncide avec $D_{f}$ sur $\mathfrak{g}_{tn}(F)$.  Puisque $D_{f}$ est un quasi-caract\`ere,  $D$ l'est aussi, cf. remarque (6).  

Soit maintenant $D$ un quasi-caract\`ere de niveau $0$ sur $G(F)$. Soit $x\in G(F)$ un \'el\'ement semi-simple. Ecrivons    une $p'$-d\'ecomposition $x=\epsilon exp(X)$  avec $X\in \mathfrak{g}_{\epsilon,tn}(F)$. Soit $Y\in \mathfrak{g}_{x}(F)$. D'apr\`es \ref{pelements}(1), on a $Y\in \mathfrak{g}_{\epsilon}(F)$ et $Y$ commute \`a $X$. Donc $xexp(Y)=\epsilon exp(X+Y)$. Si $Y$ est assez petit, $X+Y$ est topologiquement nilpotent. Par hypoth\`ese sur $D$, on a donc
$$\theta_{D}(x exp(Y))=\theta_{D}(\epsilon exp(X+Y))=\sum_{{\cal O}\in Nil(\mathfrak{g}_{\epsilon})}c(D,{\cal O})\hat{j}({\cal O},X+Y).$$
 Notons $D'$ le quasi-caract\`ere de niveau $0$ sur $\mathfrak{g}_{\epsilon}(F)$ dont la fonction $\theta_{D'}$ associ\'ee est d\'efinie par
 $$ \theta_{D'}(Z)= \sum_{{\cal O}\in Nil(\mathfrak{g}_{\epsilon})}c(D,{\cal O})\hat{j}({\cal O},Z)$$
 pour $Z\in \mathfrak{g}_{\epsilon,tn}(F)$. On vient de voir que c'est un quasi-caract\`ere. L'hypoth\`ese que $x$ est semi-simple implique que $X$ l'est aussi. Il existe donc un voisinage $\mathfrak{V}$ de $0$ dans $(\mathfrak{g}_{\epsilon})_{X}(F)$ tel que, pour presque tout $Y\in \mathfrak{V}$, on ait l'\'egalit\'e
 $$\theta_{D'}(X+Y)=\sum_{{\cal O}\in Nil((\mathfrak{g}_{\epsilon})_{X})}c_{D',{\cal O}}\hat{j}({\cal O},Y).$$
 Autrement dit, $\theta_{D}(x exp(Y))$ est \'egal \`a l'expression de droite pour presque tout $Y\in \mathfrak{V}$. Puisque $(\mathfrak{g}_{\epsilon})_{X}=\mathfrak{g}_{x}$, c'est exactement la condition (1) requise. Donc $D$ est un quasi-caract\`ere, ce qui d\'emontre (i).

Soit $D$ un quasi-caract\`ere sur $G(F)$. Si  $D$ est de niveau $0$, la propri\'et\'e \'enonc\'ee au (ii)  r\'esulte de ce que les fonctions $\lambda\mapsto \hat{j}({\cal O},\lambda^2X)$ appartiennent \`a $E$. Inversement, supposons la propri\'et\'e  en question v\'erifi\'ee. Soit $\epsilon\in G(F)_{p'}$, fixons  un voisinage $\mathfrak{V}$ et des constantes $c_{D,{\cal O}}$, de sorte que (1) soit v\'erifi\'ee. Soit $X\in \mathfrak{g}_{\epsilon,tn}(F)$  et supposons $\epsilon exp(\lambda X)\in G_{reg}$ pour tout $\lambda\in \mathfrak{o}_{F}$ (ceci est v\'erifi\'e pour presque tout $X$). Les fonctions 
$$\lambda\mapsto \theta(\epsilon exp(\lambda^2X))$$
 et 
 $$\lambda\mapsto \sum_{{\cal O}\in Nil(\mathfrak{g}_{\epsilon})}c_{D,{\cal O}}\hat{j}({\cal O},\lambda^2X)$$
 appartiennent toutes deux  \`a $E$ (la premi\`ere par hypoth\`ese). Pour un entier $n$ assez grand, $\lambda^2 X$ appartient \`a $\mathfrak{V}$ pour tout $\lambda\in \mathfrak{p}_{F}^n$, donc les deux fonctions co\"{\i}ncident sur $\mathfrak{p}_{F}^n$. Elles sont alors \'egales. Pour $\lambda=1$, cela d\'emontre l'\'egalit\'e
 $$\theta_{D}(\epsilon exp(X))=\sum_{{\cal O}\in Nil(\mathfrak{g}_{\epsilon})}c_{D,{\cal O}}\hat{j}({\cal O},X)$$
 et cela pour presque tout $X\in \mathfrak{g}_{\epsilon,tn}(F)$. Donc $D$ est un quasi-caract\`ere de niveau $0$, ce qui d\'emontre (ii).
 
  Pour (iii), il s'agit de voir que (8) est aussi v\'erifi\'ee pour $\epsilon\in G(F)_{p'}-B$. Ce probl\`eme \'etant insensible \`a la conjugaison par $G(F)$, on peut supposer que $\epsilon\in C(\epsilon')$ pour un $\epsilon'\in B$. Ecrivons $\epsilon=\epsilon'exp(Z)$, avec $Z\in \mathfrak{g}_{\epsilon',tn}(F)$.  D'apr\`es le (iv) du lemme \ref{decomposition}, on a  $G_{\epsilon}=G_{\epsilon'}$ et $Z$ est un \'el\'ement central dans $\mathfrak{g}_{\epsilon,tn}(F)$.  Pour $Y\in \mathfrak{g}_{\epsilon,tn}(F)$, on a alors $\epsilon exp(Y)=\epsilon' exp(Z+Y)$ et, en appliquant l'hypoth\`ese (8) pour $\epsilon'$, on a
 $$\theta_{D}(\epsilon exp(Y))=\sum_{{\cal O}\in Nil(\mathfrak{g}_{\epsilon})}c_{D,{\cal O}}\hat{j}({\cal O},Z+Y).$$
 Mais les fonctions $\hat{j}({\cal O})$ sont invariantes par translations par tout \'el\'ement central dans $\mathfrak{g}_{\epsilon}(F)$. Le $Z$ dispara\^{\i}t de l'expression ci-dessus et on aboutit \`a une expression comme en (8) de $\theta_{D}(\epsilon exp(Y))$. $\square$
 \bigskip
 
 {\bf Remarque.} Cette preuve et la remarque (4) entra\^{\i}nent qu'une fonction $\theta_{D}$ invariante par conjugaison par $G(F)$ et v\'erifiant la condition (8) est forc\'ement localement int\'egrable.

 \subsection{Induction de quasi-caract\`eres\label{induction}}
 
\ass{Lemme}{ Soient $M$ un Levi de $G$ et $D^M$ un quasi-caract\`ere de $M(F)$. Posons $D= Ind_{M}^G(D^M)$ .

(i) La distribution $D$ est un quasi-caract\`ere.

(ii) Si $D^M$ est de niveau $0$, $D$ est de niveau $0$.}

 Preuve. Le (i) est d\'emontr\'e dans \cite{W3} lemme 2.3. 

On note $\theta_{D^M}$ la fonction associ\'ee \`a $D^M$.  Pour $x\in G_{reg}(F)$,  l'ensemble $\{g^{-1}xg; g\in G(F)\cap M(F)$ se d\'ecompose en un nombre fini de classes de conjugaison par $M(F)$. Fixons un ensemble $X_{M}(x)$ de repr\'esentants de ces classes.  Par un calcul d'int\'egration facile, $D$ est associ\'e \`a la fonction $\theta_{D}$ d\'efinie sur $G_{reg}(F)$ par la formule
$$(1) \qquad \theta_{D}(x)=\sum_{y\in X_{M}(x)}D^G(x)^{-1/2}D^M(y)^{1/2}\theta_{D^M}(y).$$
Soit $\epsilon\in G(F)_{p'}$ et $X\in \mathfrak{g}_{\epsilon,tn}(F)$. Posons $x=\epsilon exp(X)$ et supposons $x\in G_{reg}(F)$. Pour tout $y\in X_{M}(x)$, soit $g_{y}\in G(F)$ tel que $g_{y}^{-1}xg_{y}=y$. Posons $\epsilon_{y}=g_{y}^{-1}\epsilon g_{y}$ et $X_{y}=g_{y}^{-1}\epsilon g_{y}$. On a $y=\epsilon_{y}exp(X_{y})$ et ceci est une $p'$-d\'ecomposition de $y$. Le groupe $A_{M}$ commute \`a $y$ puisque $y\in M(F)$. Il commute \`a $\epsilon_{y}$ et $X_{y}$ d'apr\`es \ref{pelements} (1). Donc $\epsilon_{y}\in M(F)$ et $X_{y}\in \mathfrak{m}_{\epsilon_{y}}(F)$. De plus $\epsilon_{y}$ est un $p'$-\'el\'ement dans $M(F)$ d'apr\`es le (i) du lemme \ref{lecasdunlevi}. Soit $\lambda\in F^{\times}$ tel que $\lambda X$ soit encore topologiquement nilpotent. Posons   $x'=\epsilon exp(\lambda X)$. L'\'el\'ement $x'$ appartient encore \`a $G_{reg}(F)$.  Les \'el\'ements $\epsilon_{y}exp(\lambda X_{y})=g_{y}^{-1}x'g_{y}$ sont tous des \'el\'ements de $M(F)$ conjugu\'es \`a $x'$ par un \'el\'ement de $G(F)$. Montrons que

(2) si $y,y'\in X_{M}(x)$, avec $y\not=y'$, alors $\epsilon_{y}exp(\lambda X_{y})$ et $\epsilon_{y'}exp(\lambda X_{y'})$ ne sont pas conjugu\'es par un \'el\'ement de $M(F)$. 

Supposons qu'il existe $m\in M(F)$ tel que $m^{-1}\epsilon_{y}exp(\lambda X_{y})m=\epsilon_{y'}exp(\lambda X_{y'})$. En posant $h=g_{y}mg_{y'}^{-1}$, on a alors $h^{-1}x'h=x'$. Puisque $x'\in G_{reg}(F)$, on a $h\in G_{x'}(F)$, d'o\`u $h\in G_{\epsilon}(F)$ d'apr\`es le (iii) 
  du lemme \ref{pelements}. On voit alors que $m^{-1}\epsilon_{y}m=\epsilon_{y'}$, d'o\`u aussi $m^{-1}\lambda X_{y}m=\lambda X_{y'}$. Il en r\'esulte que $m^{-1}X_{y}m =X_{y'}$, puis que $m^{-1}ym=y'$. Mais cela est contradictoire avec la d\'efinition de $X_{M}(x)$, ce qui d\'emontre (2). 
  
  En cons\'equence de (2), l'ensemble $X_{M}(x')$ a au moins autant d'\'el\'ements que $X_{M}(x)$. La situation \'etant sym\'etrique en $x$ et $x'$, ces deux ensembles ont m\^eme nombre d'\'el\'ements. Alors (2) nous dit que l'on peut choisir pour $X_{M}(x')$ l'ensemble $\{\epsilon_{y}exp(\lambda X_{y}); y\in X_{M}(x)\}$. On remplace maintenant $\lambda$ par $\lambda^2$, avec $\lambda\in \mathfrak{o}_{F}$. L'\'egalit\'e (1) pour $\epsilon exp(\lambda^2X)$ devient
  $$\theta_{D}(\epsilon exp(\lambda^2X))=\sum_{y\in X_{M}(x)}D^G(\epsilon exp(\lambda^2X))^{-1/2}D^M(\epsilon_{y}exp(\lambda^2X_{y}))^{1/2}\theta_{D^M}(\epsilon_{y}exp(\lambda^2X_{y})).$$
  D'apr\`es le (ii) du lemme \ref{decomposition}, cela se r\'ecrit
 $$\theta_{D}(\epsilon exp(\lambda^2X))=\sum_{y\in X_{M}(x)}d^G(\epsilon_{y})^{-1/2}D^{G_{\epsilon_{y}}}(\lambda^2X)^{-1/2} d^M(\epsilon_{y})^{1/2}D^{M_{\epsilon_{y}}}(\lambda^2X_{y})^{1/2}\theta_{D^M}(\epsilon_{y}exp(\lambda^2X_{y})).$$
 Comme fonction de $\lambda$, le terme  $D^{G_{\epsilon_{y}}}(\lambda^2X)^{-1/2}$, resp.  $D^{M_{\epsilon_{y}}}(\lambda^2X_{y})^{1/2}$, est produit d'une constante et d'une puissance enti\`ere de $\vert \lambda\vert _{F}$. Si $D^M$ est de niveau $0$, le terme $ \theta_{D^M}(\epsilon_{y}exp(\lambda^2X_{y}))$ appartient \`a $E$ d'apr\`es le (ii) du lemme \ref{quasicaracteres}. Donc la fonction $\lambda\mapsto \theta_{D}(\epsilon exp(\lambda^2X))$ appartient \`a $E$, ce qui prouve que $D$ est un quasi-caract\`ere de niveau $0$ d'apr\`es le m\^eme lemme. $\square$

Un cas particulier de la formule (1) nous servira plus tard. Supposons $x\in G_{reg}(F)\cap M_{ell}(F)$. Fixons un ensemble de repr\'esentants $\underline{N}(M)$ du quotient $Norm_{G}(M)(F)/ M(F)$. On v\'erifie que l'on peut choisir pour ensemble $X_{M}(x)$ l'ensemble $X_{M}(x)=\{n^{-1}xn; n\in \underline{N}(M)\}$. On a aussi $D^M(n^{-1}xn)=D^M(x)$ pour tout $n$. La formule devient
$$(3) \qquad \theta_{D}(x)=D^G(x)^{-1/2}D^M(x)^{1/2}\sum_{n\in \underline{N}(M)}\theta_{D^M}(n^{-1}xn).$$

 \subsection{L'espace ${\cal D}(\mathfrak{g})$\label{lespace}}
 
 Soit ${\cal F}\in Fac(G)$. On note $\boldsymbol{\mathfrak{g}}_{{\cal F}}$ l'alg\`ebre de Lie du groupe ${\bf G}_{{\cal F}}$ et $\boldsymbol{\mathfrak{g}}_{{\cal F},nil}$ son sous-ensemble des \'el\'ements nilpotents. On note $C_{cusp}(\boldsymbol{\mathfrak{g}}_{{\cal F},nil})$ l'espace des fonctions sur $\boldsymbol{\mathfrak{g}}_{{\cal F}}(k_{F})$, \`a valeurs complexes, \`a support nilpotent, qui sont invariantes par conjugaison par ${\bf G}_{{\cal F}}(k_{F})$ et cuspidales. On note $Fac_{max}(G)$ l'ensemble des facettes r\'eduites \`a un point (c'est-\`a-dire les sommets de l'immeuble). 
 
 A l'aide de ces objets, on d\'efinit des espaces $\boldsymbol{{\cal D}}_{cusp}(\mathfrak{g})$, ${\cal D}_{cusp}(\mathfrak{g})$, $\boldsymbol{{\cal D}}(\mathfrak{g})$ et ${\cal D}(\mathfrak{g})$ de fa\c{c}on similaire \`a ceux de \ref{facettes}. On pose
 $$\boldsymbol{{\cal D}}_{cusp}(\mathfrak{g})=\sum_{{\cal F}\in Fac_{max}(G)}C_{cusp}(\boldsymbol{\mathfrak{g}}_{{\cal F},nil}),$$
$$ \boldsymbol{{\cal D}}(\mathfrak{g})=\oplus_{M}\boldsymbol{{\cal D}}_{cusp}(\mathfrak{m}),$$
o\`u $M$ parcourt les Levi de $G$.  Le groupe $G(F)$ agit naturellement sur ces espaces et on note ${\cal D}_{cusp}(\mathfrak{g})$ et ${\cal D}(\mathfrak{g})$ les espaces de coinvariants. On d\'efinit une application lin\'eaire $D^G:\boldsymbol{{\cal D}}(\mathfrak{g})\to I(\mathfrak{g})^*$: pour un Levi $M$, une facette ${\cal F}_{M}\in Fac_{max}(M)$ et une fonction $f\in C_{cusp}(\boldsymbol{\mathfrak{m}}_{{\cal F}_{M},nil})$, on rel\`eve $f$ en une fonction $f_{{\cal F}_{M}}$ sur $M(F)$ qui est tr\`es cuspidale et on pose $D^G_{f}=Ind_{M}^G(D^M_{f_{{\cal F}_{M}}})$. L'application $D^G$ se quotiente en une application d\'efinie sur ${\cal D}(\mathfrak{g})$. 

\ass{Proposition}{(i) L'application $D^G:{\cal D}(\mathfrak{g})\to I(\mathfrak{g})^*$ est injective.

(ii) Son image est l'espace des quasi-caract\`eres de niveau $0$ sur $\mathfrak{g}(F)$.}

Preuve. Pour ${\cal F}\in Fac(G)$, notons  $ {\cal E}_{{\cal F},nil}$ l'espace des fonctions sur $\mathfrak{g}(F)$  \`a valeurs complexes \`a support dans $\mathfrak{k}_{{\cal F}}\cap \mathfrak{g}_{tn}(F)$ et  invariantes par  translations par $\mathfrak{k}_{{\cal F}}^+$. Notons ${\cal E}(\mathfrak{g})$ le sous-espace de $C_{c}^{\infty}(\mathfrak{g}(F))$ engendr\'e par les ${\cal E}_{{\cal F},nil}$ pour ${\cal F}\in Fac(G)$. On note $I{\cal E}(\mathfrak{g})$ son image dans $I(\mathfrak{g})$. Les d\'emonstrations des paragraphes \ref{calcul} \`a \ref{uncorollaire} qui concernent des fonctions sur $G(F)$ s'adaptent \`a l'alg\`ebre de Lie et conduisent aux m\^emes conclusions: 

(1) la compos\'ee de $D^G$ et de la restriction $I(\mathfrak{g})^*\to I{\cal E}(\mathfrak{g})^*$ est un isomorphisme de ${\cal D}(\mathfrak{g})$ sur cet espace;

(2) l'application lin\'eaire compos\'ee ${\cal D}_{cusp}(\mathfrak{g})\stackrel{D^G}{\to}I(\mathfrak{g})^*\to I_{cusp}(\mathfrak{g})^*$ est injective. 

L'assertion (1) entra\^{\i}ne que $D^G$ est injective.

Notons ${\cal H}$ le sous-espace de $ C_{c}^{\infty}(\mathfrak{g}(F)$ engendr\'e par les fonctions $f$ pour lesquelles il existe une sous-alg\`ebre d'Iwahori $\mathfrak{b}$ de $\mathfrak{g}$ de sorte que $f$ soit invariante par $\mathfrak{b}$. En notant $\mathfrak{u}$ le radical pro-$p$-nilpotent de $\mathfrak{b}$, cette invariance \'equivaut \`a ce que $\hat{f}$ soit \`a support dans $\mathfrak{u}$. En se rappelant que $\mathfrak{g}_{tn}(F)$ est l'ensemble des $X\in \mathfrak{g}(F)$ tels qu'il existe une telle alg\`ebre $\mathfrak{u}$ de sorte que $X\in \mathfrak{u}$, on voit facilement que ${\cal H}$ est l'ensemble des fonctions $f$ tels que le support de $\hat{f}$ soit contenu dans $\mathfrak{g}_{tn}(F)$. Remarquons que,  pour $f\in {\cal E}(\mathfrak{g})$, le support de $f$ est contenu dans $\mathfrak{g}_{tn}(F)$, donc $\hat{f}\in {\cal H}$. 
Notons $I{\cal H}$ l'image de ${\cal H}$ dans $I(\mathfrak{g})$ et $\hat{I}{\cal E}(\mathfrak{g})$ l'image de $I{\cal E}(\mathfrak{g})$ par transformation de Fourier. 

 Rappelons qu'un \'el\'ement $X\in \mathfrak{g}(F)$ est dit entier s'il existe une facette ${\cal F}\in Fac(G)$ de sorte que $X\in \mathfrak{k}_{{\cal F}}$.  Notons $I(\mathfrak{g})^*_{ent}$ l'espace des distributions invariantes sur $\mathfrak{g}(F)$ \`a support entier. 
 
Soient   $M$ un Levi de $G$, ${\cal F}_{M}\in Fac_{max}(M)$ et $f\in C_{cusp}(\boldsymbol{\mathfrak{m}}_{{\cal F}_{M},nil})$. On a $D^G_{f}=Ind_{M}^G(D^M_{f_{{\cal F}_{M}}})$. On v\'erifie que la transform\'ee de Fourier $\hat{D}^G_{f}$ de la distribution $D^G_{f}$ est \'egale \`a $ Ind_{M}^G(D^M_{\hat{f}_{{\cal F}_{M}}})$.  La fonction $\hat{f}_{{\cal F}_{M}}$ \'etant \`a support dans $k_{{\cal F}_{M}}$ donc entier, on a $\hat{D}^G_{f}\in  I(\mathfrak{g})^*_{ent}$. On note $\hat{D}^G({\cal D}(\mathfrak{g}))$ l'image de  $D^G({\cal D}(\mathfrak{g}))$ par transform\'ee de Fourier. 

Notons $i:{\mathbb C}[Nil(\mathfrak{g}]\to I(\mathfrak{g})^*$ l'application qui, \`a une orbite nilpotente ${\cal O}$, associe l'int\'egrale orbitale $I_{{\cal O}}$. Son image est bien s\^ur contenue dans $I(\mathfrak{g})^*_{ent}$. On a un diagramme commutatif
$$\begin{array}{ccccccc}&&&&I{\cal H}^*&&\\ &&&res\nearrow\,\,&&\,\,\searrow r&\\

\hat{D}^G({\cal D}(\mathfrak{g}))&\stackrel{j}{\to}&I(\mathfrak{g})^*_{ent}&&p_{1}\uparrow\,\,&&\hat{I}{\cal E}(\mathfrak{g})\\ &&&i\nwarrow\,\,&&\,\,\nearrow p_{2}&\\ &&&&{\mathbb C}[Nil(\mathfrak{g})]&&\\ \end{array}$$
L'application $j$ est l'injection naturelle. Les applications $res$ et $r$ sont les restrictions. Les applications $p_{1}$ et $p_{2}$ sont celles qui rendent le diagramme commutatif. Dans \cite{D} th\'eor\`eme 2.1.5, Debacker a prouv\'e les assertions suivantes:

(3) les applications $res$ et $p_{1}$ ont m\^eme image;

(4) $p_{1}$ et $p_{2}$ sont injectives.

Soit ${\bf f}\in {\cal D}(\mathfrak{g})$. Par construction, $D^G_{{\bf f}}$ est \`a support dans $\mathfrak{g}_{tn}(F)$. Pour prouver que $D^G_{{\bf f}}$ est un quasi-caract\`ere de niveau $0$, il suffit de prouver que $D^G_{{\bf f}}$ co\"{\i}ncide sur cet ensemble avec une combinaison lin\'eaire de transform\'ees de Fourier d'int\'egrales orbitales nilpotentes. D'apr\`es ce que l'on a dit ci-dessus, il suffit de prouver que $\hat{D}^G_{{\bf f}}$ co\"{\i}ncide sur ${\cal H}$ avec une combinaison lin\'eaire d'int\'egrales orbitales nilpotentes. Autrement dit, il suffit de prouver que $res\circ j(\hat{D}^G_{{\bf f}})$ appartient \`a l'image de $p_{1}$. Cela r\'esulte de (3). 

Inversement, soit $d$ un quasi-caract\`ere de niveau $0$. Notons $\hat{d}_{1}\in I{\cal H}^*$ la restriction de $\hat{d}$ \`a $I{\cal H}$.  Le quasi-caract\`ere $d$ co\"{\i}ncide sur $\mathfrak{g}_{tn}(F)$ avec une combinaison lin\'eaire de transform\'ees de Fourier d'int\'egrales orbitales nilpotentes.  Donc $\hat{d}$ co\"{\i}ncide sur ${\cal H}$ avec une combinaison lin\'eaire d'int\'egrales orbitales nilpotentes.  C'est-\`a-dire que $\hat{d}_{1}$ appartient \`a l'image de $p_{1}$. Soit $N_{1}\in {\mathbb C}[Nil(\mathfrak{g})]$ tel que $\hat{d}_{1}=p_{1}(N_{1})$. Par transformation de Fourier, l'assertion (1) dit que l'application $r\circ res\circ j$ est un isomorphisme. Il existe donc ${\bf f}\in {\cal D}(\mathfrak{g})$ tel que $r(\hat{d}_{1})=r\circ res\circ j(\hat{D}^G_{{\bf f}})$. D'apr\`es (3), il existe $N_{2}\in {\mathbb C}[Nil(\mathfrak{g})]$ tel que $res\circ j(\hat{D}^G_{{\bf f}})=p_{1}(N_{2})$. On a alors
$$p_{2}(N_{2})=r\circ p_{1}(N_{2})=r\circ res\circ j(\hat{D}^G_{{\bf f}})=r(\hat{d}_{1})=r\circ p_{1}(N_{1})=p_{2}(N_{1}).$$
D'apr\`es (4), cela entra\^{\i}ne $N_{1}=N_{2}$. D'o\`u
$$\hat{d}_{1}=p_{1}(N_{1})=p_{1}(N_{2})=res\circ j(\hat{D}^G_{{\bf f}}).$$
Autrement dit, les distributions $\hat{d}$ et $\hat{D}^G_{{\bf f}}$ co\"{\i}ncident sur ${\cal H}$. Donc $d$ et $D^G_{{\bf f}}$ co\"{\i}ncident sur $\mathfrak{g}_{tn}(F)$. Mais elles sont toutes deux \`a support dans cet ensemble. Donc $d=D^G_{{\bf f}}$, ce qui ach\`eve la d\'emonstration. $\square$

\subsection{Filtrations sur l'alg\`ebre de Lie\label{filtrationssurlalgebredelie}}
Pour tout $n \in {\mathbb Z}$, notons $Fil^nI(\mathfrak{g})$ l'image dans $I(\mathfrak{g})$ du sous-espace des $f\in C_{c}^{\infty}(\mathfrak{g}(F))$ qui v\'erifient la condition: pour tout Levi $M$ tel que $dim(A_{M})>n$ et tout $X\in \mathfrak{g}_{reg}(F)\cap \mathfrak{m}(F)$, on a $I^G(X,f)=0$.   On a
$$Fil^{a_{G}-1}I(\mathfrak{g})=\{0\}\subset Fil^{a_{G}}I(\mathfrak{g})=I_{cusp}(\mathfrak{g})\subset...\subset Fil^{a_{M_{min}}}I(\mathfrak{g})=I(\mathfrak{g}),$$
et, en posant $Gr^nI(\mathfrak{g})=Fil^nI(\mathfrak{g})/Fil^{n-1}I(\mathfrak{g})$, on a
$$Gr^nI(\mathfrak{g})\simeq \oplus_{M\in \underline{{\cal L}}^n}I_{cusp}(\mathfrak{m})^{W^G(M)}.$$

Notons $Ann^nI(\mathfrak{g})^*$ l'annulateur de $Fil^{n-1}I(\mathfrak{g})$ dans $I(\mathfrak{g})^*$. On a
$$Ann^{a_{M_{min}}+1}I(\mathfrak{g})^*=\{0\}\subset Ann^{a_{M_{min}}}I(\mathfrak{g})^*\subset...\subset Ann^{a_{G}}I(\mathfrak{g})^*=I(\mathfrak{g})^*,$$
et, en posant $Gr^nI(\mathfrak{g})^*=Ann^nI(\mathfrak{g})^*/Ann^{n+1}I(\mathfrak{g})^*$,
$$(1) \qquad Gr^nI(\mathfrak{g})^*\simeq \oplus_{M\in \underline{{\cal L}}^n}I_{cusp}(\mathfrak{m})^{*W^G(M)}.$$
Notons $Qc_{0}(\mathfrak{g})$ l'espace des quasi-caract\`eres de niveau $0$ sur $\mathfrak{g}(F)$. Notons $Qc_{0}^n(\mathfrak{g})$ la somme des $Ind_{M}^G(Qc_{0}(\mathfrak{m}))$ sur les Levi $M$ tels que $dim(A_{M})\geq n$. Notons aussi ${\cal D}^n(\mathfrak{g})$ la somme des images dans ${\cal D}(\mathfrak{g})$ des  ${\cal D}_{cusp}(\mathfrak{m})$ sur les m\^emes Levi $M$.

\ass{Lemme}{Pour tout $n\in {\mathbb Z}$, on a l'\'egalit\'e
$$Qc_{0}^n(\mathfrak{g})=Qc_{0}(\mathfrak{g})\cap Ann^nI(\mathfrak{g})^*=D^G({\cal D}^n(\mathfrak{g})),$$
et l'isomorphisme
$$Qc_{0}^n(\mathfrak{g})/Q_{0}^{n+1}(\mathfrak{g})^*\simeq \oplus_{M\in \underline{{\cal L}}^n}D^M({\cal D}_{cusp}(\mathfrak{m})^{W^G(M)}).$$}

Preuve. L'\'egalit\'e $Qc_{0}^n(\mathfrak{g})=D^G({\cal D}^n(\mathfrak{g}))$ r\'esulte des d\'efinitions et de la proposition \ref{lespace} appliqu\'ee aux Levi $M$ tels que $dim(A_{M})\geq n$. Notons $I(\mathfrak{g})^*_{nil}$ le sous-espace des distributions sur $\mathfrak{g}(F)$ \`a support nilpotent. Autrement dit $I(\mathfrak{g})^*_{nil}=i({\mathbb C}[Nil])$ avec les notations de la preuve pr\'ec\'edente.  On a prouv\'e en \cite{MW} proposition I.5.3 que $I(\mathfrak{g})^*_{nil}\cap Ann^nI(\mathfrak{g})^*$ \'etait la somme des $Ind_{M}^G(I(\mathfrak{m})^*_{nil})$ o\`u $M$ parcourt les Levi tels que $dim(A_{M})\geq n$. Les filtrations sont invariantes par transform\'ees de Fourier. Cela r\'esulte de ce que, pour toute $f\in C_{c}^{\infty}(\mathfrak{g}(F))$, tout Levi $M$ et tout $X\in \mathfrak{m}(F)\cap \mathfrak{g}_{reg}(F)$, $I^G(X,\hat{f})$ ne d\'epend que des $I^G(Y,f)$ pour $Y\in \mathfrak{m}(F)\cap \mathfrak{g}_{reg}(F)$. Notons $\hat{I}(\mathfrak{g})^*_{nil}$ l'image de $I(\mathfrak{g})^*_{nil}$ par transform\'ee de Fourier. On obtient que $\hat{I}(\mathfrak{g})^*_{nil}\cap Ann^nI(\mathfrak{g})^*$ est la somme des $Ind_{M}^G(\hat{I}(\mathfrak{m})^*_{nil})$ o\`u $M$ parcourt les Levi tels que $dim(A_{M})\geq n$. Soit $d\in Qc_{0}(\mathfrak{g})\cap Ann^nI(\mathfrak{g})^* $. Par d\'efinition d'un quasi-caract\`ere de niveau $0$, $d$ est la restriction \`a $\mathfrak{g}_{tn}(F)$ d'un \'el\'ement $d_{0}\in \hat{I}(\mathfrak{g})^*_{nil}$. Puisque $d\in Ann^nI(\mathfrak{g})^*$, $d_{0}$ annule tout \'el\'ement de $Fil^{n-1}I(\mathfrak{g})$  \`a support topologiquement nilpotent. L'espace $Fil^{n-1}I(\mathfrak{g})$  est invariant par l'action de $F^{\times}$ par homoth\'etie. Puisque $\hat{I}(\mathfrak{g})^*_{nil}$ est engendr\'e par des \'el\'ements homog\`enes, sinon par l'action de $F^{\times}$, du moins par l'action du sous-groupe des carr\'es, dire que $d_{0}$ annule tout \'el\'ement de $Fil^{n-1}I(\mathfrak{g})$  \`a support topologiquement nilpotent \'equivaut \`a dire qu'il annule tout $Fil^{n-1}I(\mathfrak{g})$, c'est-\`a-dire que $d_{0}\in \hat{I}(\mathfrak{g})^*_{nil}\cap Ann^nI(\mathfrak{g})^*$. Donc $d_{0}=\sum_{M}Ind_{M}^G(d_{M,0})$,   o\`u l'on somme sur les Levi $M$  tels que $dim(A_{M})\geq n$ et o\`u $d_{M,0}$ est un certain \'el\'ement de  $\hat{I}(\mathfrak{m})^*_{nil}$. Donc $d= \sum_{M}Ind_{M}^G(d_{M})$, o\`u $d_{M}$ est la restriction de $d_{M,0}$ \`a $\mathfrak{m}_{tn}(F)$. Cet \'el\'ement $d_{M}$ appartient \`a $Qc_{0}(\mathfrak{m})$. Donc $d\in Qc_{0}^n(\mathfrak{g})$. Cela d\'emontre l'inclusion $Qc_{0}(\mathfrak{g})\cap Ann^nI(\mathfrak{g})^*\subset Qc_{0}^n(\mathfrak{g})$ et l'inclusion oppos\'ee est \'evidente.

On a
$$Qc_{0}^n(\mathfrak{g})/Qc_{0}^{n+1}(\mathfrak{g})^*=D^G({\cal D}^n(\mathfrak{g}))/D^G({\cal D}^{n+1}(\mathfrak{g}))$$
et, d'apr\`es ce que l'on vient de prouver, cet espace s'envoie injectivement dans $Gr^nI(\mathfrak{g})^*$.
Par d\'efinition, $D^G({\cal D}^n(\mathfrak{g}))$ est la somme de $D^G({\cal D}^{n+1}(\mathfrak{g}))$ et de $D^G(\oplus_{M\in \underline{{\cal L}}^n}{\cal D}_{cusp}(\mathfrak{m}))$. Ce dernier espace a donc m\^eme image dans $Gr^nI(\mathfrak{g})^*$ que $Qc_{0}^n(\mathfrak{g})/Qc_{0}^{n+1}(\mathfrak{g})^*$. Soit $M\in \underline{{\cal L}}^n$ et ${\bf f}\in {\cal D}_{cusp}(\mathfrak{m})$. L'image de $D^G_{{\bf f}}$ dans $Gr^nI(\mathfrak{g})^*$  est nulle dans  les composantes $I_{cusp}(\mathfrak{m}')^{*W^G(M')}$ de la d\'ecomposition (1) pour $M'\not=M$ et   est   l'image naturelle de $D^M_{{\bf f}}$ dans $I_{cusp}(\mathfrak{m})^{* W^G(M)}$ pour $M'=M$ (c'est-\`a-dire l'image de  $D^M_{{\bf f}}$ dans $I_{cusp}(\mathfrak{m})^*$ que l'on restreint au sous-espace $I_{cusp}(\mathfrak{m})^{W^G(M)}$). En  moyennant sur le groupe $ W^G(M)$, on obtient que $Qc_{0}^n(\mathfrak{g})/Qc_{0}^{n+1}(\mathfrak{g})^*$ s'identifie \`a la somme des images de ${\cal D}_{cusp}(\mathfrak{m})^{W^G(M)}$ dans $I_{cusp}(\mathfrak{m})^{*W^G(M)}$. 
 D'apr\`es le (2) de \ref{lespace}, ${\cal D}_{cusp}(\mathfrak{m})^{W^G(M)}$ s'envoie injectivement dans $I_{cusp}(\mathfrak{m})^{*W^G(M)}$, ce qui ach\`eve la d\'emonstration. $\square$

\subsection{ Un premier th\'eor\`eme\label{unpremiertheoreme}}

\ass{Th\'eor\`eme}{ L'image de l'application $D^G:{\cal D}(G)\to I(G)^*$ est form\'ee de  quasi-caract\`eres de niveau $0$ sur $G(F)$.}

Preuve. En vertu du (ii) du lemme \ref{induction}, il suffit de prouver que, pour tout \'el\'ement  ${\bf f}\in {\cal D}_{cusp}(G)$, $D^G_{{\bf f}}$ est un quasi-caract\`ere de niveau $0$ sur $G(F)$.  La famille ${\bf f}$ peut \^etre infinie mais la d\'efinition des quasi-caract\`eres   de niveau $0$ est locale et, localement, seuls un nombre fini de fonctions de la famille ${\bf f}$ interviennent. On peut donc aussi bien supposer qu'il y a une seule fonction. C'est-\`a-dire que l'on peut fixer $({\cal F},\nu)\in Fac^*_{max}(G)$ et $f\in C_{cusp}({\bf G}_{{\cal F}}^{\nu})$ et supposer que ${\bf f}$ est r\'eduit \`a l'unique fonction $f$. On a alors $D^G_{{\bf f}}=D^G_{f_{{\cal F}}}$. C'est un quasi-caract\`ere d'apr\`es \ref{quasicaracteres} (6). On note $\theta$ sa fonction localement int\'egrable associ\'ee. Soit $\epsilon$ un \'el\'ement $p'$-compact de $G(F)$.  On doit calculer $\theta(\epsilon exp(X))$ pour tout $X\in \mathfrak{g}_{\epsilon,tn}(F)$ tel que $\epsilon exp(X)$ soit fortement r\'egulier.  On peut \'evidemment supposer $\nu=w_{G}(\epsilon)$, sinon cette fonction est nulle. Fixons $X$. Notons $T$ le commutant de $X$ dans $G_{\epsilon}$, qui est aussi le commutant de $\epsilon exp(X)$ dans $G$. Notons $M$ le commutant de $A_{T}$ dans $G$. C'est un Levi de $G$ contenant $\epsilon$ et $exp(X)$ d'apr\`es le (iii) du lemme \ref{pelements}. On pose $M_{\epsilon}=G_{\epsilon}\cap M$. C'est un Levi de $G_{\epsilon}$ et on a $X\in \mathfrak{m}_{\epsilon}(F)$. L'\'el\'ement  $X$ est elliptique dans $\mathfrak{m}_{\epsilon}(F)$ et $\epsilon exp(X)$ est elliptique dans $M(F)$. Notons que $A_{T}=A_{M_{\epsilon}}=A_{M}$.  D'apr\`es \ref{quasicaracteres} (7), on a
$$(1) \qquad \theta(\epsilon exp(X))=(-1)^{a_{M}-a_{G}}D^G(\epsilon exp(X))^{-1/2}m(A_{M})m(A_{G})^{-1}J^G_{M}(\epsilon exp(X),f_{{\cal F}})$$
 $$=(-1)^{a_{M}-a_{G}}m(A_{M})m(A_{G})^{-1} \int_{A_{M}(F)\backslash G(F)}f_{{\cal F}}(g^{-1} \epsilon exp(X)g)v_{M}(g)\, dg.$$
 Pour $g\in G(F)$, on a

(2) $g^{-1}\epsilon exp(X)g\in K_{{\cal F}}^{\nu}$ si et seulement si $g^{-1}\epsilon g\in K_{{\cal F}}^{\nu}$ et $g^{-1} exp(X)g\in K_{{\cal F}}^0$.

En effet, soit $c\geq1$ un entier premier \`a $p$ tel que $\epsilon^c\in Z(G)(F)$. Supposons $g^{-1}\epsilon exp(X)g\in K_{{\cal F}}^{\nu}$. Alors $g^{-1} \epsilon^cexp(cX)g\in K_{{\cal F}}^{\dag}$. On a aussi $g^{-1}\epsilon^cg\in Z(G)(F)\subset K_{{\cal F}}^{\dag}$. Donc $g^{-1}exp(cX)g\in K_{{\cal F}}^{\dag}$. Il en r\'esulte que $g^{-1}exp(X)g\in K_{{\cal F}}^{\dag}$ et, puisque c'est un \'el\'ement topologiquement unipotent, on a forc\'ement $g^{-1}exp(X)g\in K_{{\cal F}}^0$. Puisque $g^{-1}\epsilon exp(X)g\in K_{{\cal F}}^{\nu}$, cela entra\^{\i}ne $g^{-1}\epsilon g\in K_{{\cal F}}^{\nu}$. La r\'eciproque est \'evidente. D'o\`u (2). 

On a aussi:

(3) la classe de conjugaison par $G(F)$ de $\epsilon$ coupe $K_{{\cal F}}^{\nu}$ en un nombre fini de classes de conjugaison par $K_{{\cal F}}^0$. 

En effet, notons $Cl(\epsilon)$ cette classe de conjugaison par $G(F)$. C'est un sous-ensemble ferm\'e de $G(F)$, donc son intersection avec $K_{{\cal F}}^{\nu}$ est compacte. On sait que l'application 
$$\begin{array}{ccc}Z_{G}(\epsilon)(F)\backslash G(F)&\to&Cl(\epsilon)\\ g&\mapsto& g^{-1}\epsilon g\\ \end{array}$$
est un hom\'eomorphisme. Puisque les orbites de l'action de $K_{{\cal F}}^{\nu}$ sur l'ensemble de d\'epart sont ouvertes, il en est de m\^eme de celles sur l'ensemble d'arriv\'ee. Puisque $Cl(\epsilon)\cap K_{{\cal F}}^{\nu}$ est compact, il n'y a donc qu'un nombre fini de telles orbites, d'o\`u (3).

  Notons $\Gamma_{\epsilon}$ l'ensemble des $g\in G(F)$ tels que $g^{-1}\epsilon g\in K_{{\cal F}}^{\nu}$. En cons\'equence de (3), l'ensemble $G_{\epsilon}(F)\backslash \Gamma_{\epsilon}/K_{{\cal F}}^0$ est fini. Fixons-en un ensemble de repr\'esentants $\boldsymbol{\gamma}$. Pour tout $\gamma\in \boldsymbol{\gamma}$, d\'efinissons comme en \ref{facettesfonctions} la facette $\gamma{\cal F}$ et posons
  $$m_{\gamma}=mes(G_{\epsilon}(F)\cap K_{\gamma{\cal F}}^0)^{-1},$$
  la mesure \'etant calcul\'ee relativement \`a la mesure implicitement fix\'ee sur $G_{\epsilon}(F)$.  On v\'erifie la formule d'int\'egration
  $$\int_{A_{M}(F)\backslash\Gamma_{\epsilon}}\varphi(g)\,dg=\sum_{\gamma\in \boldsymbol{\gamma}}m_{\gamma}\int_{A_{M}(F)\backslash G_{\epsilon}(F)}\int_{K_{{\cal F}}^0}\varphi(g\gamma k)\,dk\,dg$$
  pour toute fonction int\'egrable $\varphi$ sur $A_{M}(F)\backslash \Gamma_{\epsilon}$. On applique cela \`a la fonction $\varphi(g)=f_{{\cal F}}(g^{-1}\epsilon exp(X) g)v_{M}(g)$. 
  La fonction $f_{{\cal F}} $ est invariante par conjugaison par $K_{{\cal F}}^0$. On voit que 
  $f_{{\cal F}}(k^{-1}\gamma^{-1}g^{-1}\epsilon exp(X)g \gamma k)=(^{\gamma}f)_{\gamma{\cal F}}(\epsilon exp(g^{-1}Xg))$ pour tous $\gamma\in \boldsymbol{\gamma}$, $g\in G_{\epsilon}(F)$ et  $k\in K_{{\cal F}}^0$. Pour de m\^emes
$\gamma$ et $g$, posons 
$$v_{M,\gamma}(g)=\int_{K_{{\cal F}}^0}v_{M}(g\gamma k)\,dk.$$  
  Gr\^ace \`a (1) et (2) et \`a l'\'egalit\'e $A_{M_{\epsilon}}=A_{M}$, on obtient
$$(4) \qquad \theta(\epsilon exp(X))=(-1)^{a_{M_{\epsilon}}-a_{G}}m(A_{M_{\epsilon}})m(A_{G})^{-1}\sum_{\gamma\in \boldsymbol{\gamma}}m_{\gamma}$$
$$\int_{A_{M_{\epsilon}}(F)\backslash G_{\epsilon}(F)}(^{\gamma}f)_{\gamma{\cal F}}(\epsilon exp(g^{-1}X g))v_{M,\gamma} (g)\,dg.$$

  Fixons $\gamma\in \boldsymbol{\gamma}$ et simplifions provisoirement la notation en supposant $\gamma=1$.  Par  d\'efinition de $\boldsymbol{\gamma}$, on a $\epsilon\in K_{{\cal F}}^{\nu}$.   Introduisons la facette ${\cal F}'\in Imm(G_{\epsilon,AD})$ associ\'ee \`a ${\cal F}$, cf.  proposition \ref{pointsfixes} (iii). Le groupe ${\bf G}_{\epsilon,{\cal F}'}$ est la composante neutre de $({\bf G}_{{\cal F}})^{\epsilon}$.
Notons  $\bar{\epsilon}$ la r\'eduction de $\epsilon$ dans ${\bf G}^{\nu}_{{\cal F}}(k_{F})$ et $f'$ la fonction sur ${\bf G}_{\epsilon,{\cal F}'}(k_{F})$, \`a support unipotent,  d\'efinie par $f'(x)=f(\bar{\epsilon}x)$, pour tout \'el\'ement unipotent $x\in {\bf G}_{\epsilon,{\cal F}'}(k_{F})$.  De $f'$ se d\'eduit une fonction $f'_{{\cal F}'}$ sur $G_{\epsilon}(F)$.
Pour  $y\in K_{{\cal F}}^0\cap G_{\epsilon}(F)$, la r\'eduction $\bar{y}$ de $y$ dans ${\bf G}_{{\cal F}}(k_{F})$ appartient \`a     $({\bf G}_{{\cal F}})^{\epsilon}(k_{F})$. Mais, si $y$ est topologiquement unipotent, on peut remplacer dans cette relation le groupe $({\bf G}_{{\cal F}})^{\epsilon}$ par sa composante neutre ${\bf G}_{\epsilon,{\cal F}'}$. De plus cette r\'eduction $\bar{y}$ est unipotente. On voit alors que, pour $g\in G_{\epsilon}(F)$, on a l'\'egalit\'e 

(5) $f_{{\cal F}}(\epsilon exp(g^{-1} Xg))=f'_{{\cal F}'}(exp(g^{-1}Xg))$. 

Montrons que

(6) si $dim(A_{G_{\epsilon}})> dim(A_{G})$ ou si ${\cal F}'$ n'est pas r\'eduit \`a un point, $ f'=0$.

Notons ${\bf A}_{{\cal F}'}$ le plus grand sous-tore central et d\'eploy\'e de ${\bf G}_{\epsilon,{\cal F}'}$. Notons ${\bf M}'$ son commutant dans ${\bf G}_{{\cal F}}$, qui est un Levi de ce groupe. Fixons un \'el\'ement $x_{*}\in X_{*}({\bf A}_{{\cal F}'})$ en position g\'en\'erale, notons ${\bf P}'$ le sous-groupe de ${\bf G}_{{\cal F}}$ engendr\'e par ${\bf M}'$ et par les groupes radiciels relatifs \`a l'action de ${\bf A}_{{\cal F}'}$ associ\'es aux racines $\alpha$ telles que $<\alpha,x_{*}>>0$. C'est un \'el\'ement de ${\cal P}({\bf M}')$. 
Le tore ${\bf A}_{{\cal F}'}$ est contenu dans l'ensemble des points fixes de l'action de $\bar{\epsilon}$ sur ${\bf G}_{{\cal F}}$. Il en r\'esulte que cette action conserve ${\bf M}'$ et  ${\bf P}'$.  Alors ${\bf P}^{'\nu}=\epsilon{\bf P}'$ est un espace parabolique de ${\bf G}^{\nu}_{{\cal F}}$ d'espace de Levi ${\bf M}^{'\nu}=\epsilon{\bf M}'$. Sur la cl\^oture alg\'ebrique $\bar{k}_{F}$, les valeurs propres de l'action de $\bar{\epsilon}$ dans $\boldsymbol{\mathfrak{u}}_{{\bf P}'}$ sont des racines de l'unit\'e, puisque $\epsilon$ est $p'$-compact mod $Z(G)$. L'espace des points fixes de l'action de $\bar{\epsilon}$ dans $\boldsymbol{\mathfrak{g}}_{{\cal F}}$ est $\boldsymbol{\mathfrak{g}}_{\epsilon,{\cal F}'}$, qui est inclus dans $\boldsymbol{\mathfrak{m}}'$ par d\'efinition de ${\bf M}'$. Donc les valeurs propres de l'action de $\bar{\epsilon}$ dans $\boldsymbol{\mathfrak{u}}_{{\bf P}'}$ sont toutes diff\'erentes de $1$. Il en r\'esulte ais\'ement que, pour tout \'el\'ement unipotent $\bar{n}\in {\bf M}'(k_{F})$, l'action  $1-ad(\bar{\epsilon}\bar{n})$  dans $\boldsymbol{\mathfrak{u}}_{{\bf P}'}(k_{F})$ est un isomorphisme. Pour un tel \'el\'ement $\bar{n}$, il s'en d\'eduit l'\'egalit\'e
$$\sum_{\bar{u}\in {\bf U}_{{\bf P}'}(k_{F})}f(\bar{\epsilon}\bar{n}\bar{u})=\sum_{\bar{u}\in {\bf U}_{{\bf P}'}(k_{F})}f(\bar{u}^{-1}\bar{\epsilon}\bar{n}\bar{u}).$$
Puisque $f$ est invariante par conjugaison, le membre de droite n'est autre que $\vert  {\bf U}_{{\bf P}'}(k_{F})\vert  f(\bar{\epsilon}\bar{n})$. Supposons que ${\bf P}^{'\nu}$ soit un  espace parabolique propre de ${\bf G}_{{\cal F}}^{\nu}$. Alors le membre de gauche est nul puisque $f$ est cuspidale et invariante par conjugaison. Dans ce cas $f(\bar{\epsilon}\bar{n})=0$. A fortiori, cela est vrai pour tout \'el\'ement unipotent $\bar{n}\in {\bf G}_{\epsilon,{\cal F}'}(k_{F})$, auquel cas $f(\bar{\epsilon}\bar{n})=f'(\bar{n})$. D'o\`u $f'=0$.   Pour achever la preuve de (6), il reste \`a prouver que, sous les hypoth\`eses de cette assertion, l'espace parabolique ${\bf P}^{'\nu}$ est propre, ou encore que ${\bf M}'\not={\bf G}_{{\cal F}}$. On d\'emontre la contrapos\'ee de cette assertion. Supposons donc ${\bf M}'={\bf G}_{{\cal F}}$. Notons ${\bf A}_{{\cal F}}$ le plus grand sous-tore central d\'eploy\'e dans ${\bf G}_{{\cal F}}$ et notons ${\bf A}_{{\cal F}}^{\nu}$ le plus grand sous-tore contenu dans l'ensemble des points fixes de l'action de $\bar{\epsilon}$ dans ${\bf A}_{{\cal F}}$. L'hypoth\`ese  ${\bf M}'={\bf G}_{{\cal F}}$ signifie que ${\bf A}_{{\cal F}'}\subset {\bf A}_{{\cal F}}$. Puisque $\bar{\epsilon}$ agit trivialement sur ${\bf A}_{{\cal F}'}$, cela entra\^{\i}ne ${\bf A}_{{\cal F}'}\subset {\bf A}_{{\cal F}}^{\nu}$ d'o\`u l'\'egalit\'e ${\bf A}_{{\cal F}'}= {\bf A}_{{\cal F}}^{\nu}$ car l'inclusion oppos\'ee est imm\'ediate. Puisque $({\cal F},\nu)\in Fac_{max}(G)$, on a $dim( {\bf A}_{{\cal F}}^{\nu})=dim(A_{G})$. D'o\`u $dim({\bf A}_{{\cal F}'})=dim(A_{G})$. Or on a \'evidemment $dim({\bf A}_{{\cal F}'})\geq dim(A_{G_{\epsilon}})\geq dim(A_{G})$. Ces trois dimensions sont donc \'egales. L'\'egalit\'e $dim({\bf A}_{{\cal F}'})= dim(A_{G_{\epsilon}})$ \'equivaut \`a ce que ${\cal F}'$ soit r\'eduit \`a un point. L'\'egalit\'e des trois dimensions pr\'ec\'edentes est donc le contraire de l'hypoth\`ese de (6), ce qui ach\`eve la d\'emonstration de cette assertion.

Supposons maintenant que  $dim(A_{G_{\epsilon}})=dim(A_{G})$ et que ${\cal F}'$ soit r\'eduit \`a un point, c'est-\`a-dire ${\cal F}'\in Fac_{max}(G_{\epsilon})$. Montrons que

(7) $f'$ est cuspidale.

Fixons un sous-groupe parabolique propre ${\bf P}'_{\epsilon}$ de ${\bf G}_{\epsilon,{\cal F}'}$ de composante de Levi ${\bf M}'_{\epsilon}$. Notons ${\bf A}'$ le plus grand sous-tore central d\'eploy\'e de ${\bf M}'_{\epsilon}$ et ${\bf M}'$ le  commutant de ${\bf A}'$ dans ${\bf G}_{{\cal F}}$. Fixons un \'el\'ement $x_{*}\in X_{*}({\bf A}')$ tel que $<\alpha,x_{*}>>0$ pour toute racine $\alpha$ de ${\bf A}'$ dans $\boldsymbol{\mathfrak{u}}_{{\bf P}'_{\epsilon}}$.  On d\'efinit un sous-groupe parabolique ${\bf P}'$ de ${\bf G}_{{\cal F}}$ par la condition: $\boldsymbol{\mathfrak{u}}_{{\bf P}'}$ est la somme des espaces radiciels associ\'es aux racines $\alpha$ de ${\bf A}'$ dans $\boldsymbol{\mathfrak{g}}_{{\cal F}}$ telles que $<\alpha,x_{*}>>0$.  En supposant $x_{*}$ en position g\'en\'erale, ${\bf P}'$ a pour composante de Levi ${\bf M}'$ et on a ${\bf P}'\cap {\bf G}_{\epsilon,{\cal F}'}={\bf P}'_{\epsilon}$. De plus ${\bf P}'$ est propre puisque ${\bf P}'_{\epsilon}$ l'est. L'argument est maintenant similaire \`a celui de la preuve de (6). L'ensemble ${\bf P}^{'\nu}=\bar{\epsilon}{\bf P}'$ est un espace parabolique de ${\bf G}_{{\cal F}}^{\nu}$ d'espace de Levi ${\bf M}^{'\nu}=\bar{\epsilon}{\bf M}'$. Pour un \'el\'ement unipotent $\bar{n}\in {\bf M}'_{\epsilon}(k_{F})$, on a les \'egalit\'es
$$\sum_{\bar{u}\in {\bf U}_{{\bf P}'_{\epsilon}}(k_{F})}f'(\bar{n}\bar{u})=\sum_{\bar{u}\in {\bf U}_{{\bf P}'_{\epsilon}}(k_{F})}f(\bar{\epsilon}\bar{n}\bar{u})$$
$$=\vert {\bf U}_{{\bf P}'}(k_{F})\vert ^{-1}\sum_{\bar{x}\in  {\bf U}_{{\bf P}'}(k_{F})}\sum_{\bar{u}\in {\bf U}_{{\bf P}'_{\epsilon}}(k_{F})}f(\bar{x}\bar{\epsilon}\bar{n}\bar{u}\bar{x}^{-1})=\vert  {\bf U}_{{\bf P}'}(k_{F})\vert ^{-1}\vert {\bf U}_{{\bf P}'_{\epsilon}}(k_{F}\vert \sum_{\bar{x}\in {\bf U}_{{\bf P}'}(k_{F})}f(\bar{\epsilon}\bar{n}\bar{x})$$
et cette derni\`ere expression est nulle car $f$ est cuspidale. La nullit\'e de la premi\`ere expression signifie que $f'$ l'est. D'o\`u (7).

 R\'etablissons le $\gamma$ que l'on avait suppos\'e \'egal \`a $1$ pour simplifier la notation. On note  ${\cal F}'_{\gamma}$ la facette de $Imm(G_{\epsilon,AD})$ que l'on avait not\'ee ${\cal F}'$ et $(^{\gamma}f)'$ l'analogue de $f'$ pour la fonction $^{\gamma}f$. On note $\boldsymbol{\gamma}_{max}$ le sous-ensemble des $\gamma\in \boldsymbol{\gamma}$ tels que ${\cal F}'_{\gamma}$ soit r\'eduit \`a un point. Si $dim(A_{G_{\epsilon}})> dim(A_{G})$, l'assertion (6) et la formule (4) impliquent que $\theta(\epsilon exp(X))=0$ et on a termin\'e. Supposons d\'esormais $dim(A_{G_{\epsilon}})=dim(A_{G})$. Alors les assertions (5) et (6) entra\^{\i}nent  que la formule (4) se transforme en 
 $$(8) \qquad \theta(\epsilon exp(X))=(-1)^{a_{M_{\epsilon}}-a_{G}}m(A_{M_{\epsilon}})m(A_{G})^{-1}\sum_{\gamma\in \boldsymbol{\gamma}_{max}}m_{\gamma}$$
 $$\int_{A_{M_{\epsilon}}(F)\backslash G_{\epsilon}(F)}(^{\gamma}f)'_{{\cal F}'_{\gamma}}(exp(g^{-1}X g))v_{M,\gamma} (g)\,dg.$$
  Fixons $\gamma\in \boldsymbol{\gamma}_{max}$ et calculons $v_{M,\gamma}(g)$ pour $g\in G_{\epsilon}(F)$. Le poids $v_{M}$ est calcul\'e gr\^ace au choix d'un sous-groupe compact sp\'ecial de $G(F)$. Fixons aussi un tel sous-groupe compact  $K_{\epsilon}$ du groupe $G_{\epsilon}(F)$, qui permet de d\'efinir un poids $v_{M_{\epsilon}}$. Une formule d'Arthur, que l'on a reprise en \cite{W1} lemme 3.3, dit que, pour $g\in G_{\epsilon}(F)$, $v_{M,\gamma}(g)-v_{M_{\epsilon}}(g)$ est une somme finie de termes $v'_{Q_{\epsilon}}(g)$, o\`u $Q_{\epsilon}$ est un sous-groupe parabolique propre de $G_{\epsilon}$ contenant  $M_{\epsilon}$ et $v'_{Q_{\epsilon}}$ est une fonction localement int\'egrable sur $G_{\epsilon}(F)$ invariante \`a gauche par $M_{\epsilon}(F)$ et par $U_{Q_{\epsilon}}(F)$. Pour de telles donn\'ees, notons $L_{\epsilon}$ la composante de Levi de $Q_{\epsilon}$ contenant $M_{\epsilon}$.  En utilisant la d\'ecomposition d'Iwasawa, on calcule
$$\int_{A_{M_{\epsilon}}(F)\backslash G_{\epsilon}(F)}(^{\gamma}f)'_{{\cal F}'_{\gamma}}(exp(g^{-1}X g))v_{Q_{\epsilon}}' (g)\,dg=$$
$$c\int_{K_{\epsilon}}\int_{A_{M_{\epsilon}}(F)\backslash L_{\epsilon}(F)}\int_{U_{Q_{\epsilon}}(F)}(^{\gamma}f)'_{{\cal F}'_{\gamma}}(k^{-1}u^{-1}l^{-1}exp(X)luk)v'_{Q_{\epsilon}}(lk)\,du\,dl\,dk,$$
o\`u $c>0$ ne d\'epend que des mesures de Haar. Puisque $X$ est r\'egulier, l'int\'egrale int\'erieure en $u$ se transforme en 
$$c(X)\int_{U_{Q_{\epsilon}}(F)}(^{\gamma}f)'_{{\cal F}'_{\gamma}}(k^{-1}l^{-1}exp(X)luk)\, du,$$
o\`u $c(X)>0$ est un certain d\'eterminant. Or cette int\'egrale est nulle car $(^{\gamma}f)'_{{\cal F}'_{\gamma}}$ est tr\`es cuspidale. Les fonctions $v'_{Q_{\epsilon}}$ ne contribuent donc pas \`a la formule (8) et  on peut r\'ecrire cette formule
$$\theta(\epsilon exp(X))=(-1)^{a_{M_{\epsilon}}-a_{G}}m(A_{M_{\epsilon}})m(A_{G})^{-1}\sum_{\gamma\in \boldsymbol{\gamma}_{max}}m_{\gamma}$$
$$\int_{A_{M_{\epsilon}}(F)\backslash G_{\epsilon}(F)}(^{\gamma}f)'_{{\cal F}'_{\gamma}}(exp(g^{-1}X g))v_{M_{\epsilon}} (g)\,dg.$$
 En comparant avec \ref{quasicaracteres} (7), on obtient
$$\theta(\epsilon exp(X))=(-1)^{a_{G_{\epsilon}}-a_{G}}m(A_{G_{\epsilon}})m(A_{G})^{-1}\sum_{\gamma\in \boldsymbol{\gamma}_{max}}m_{\gamma}\theta_{\gamma}(exp(X)),$$
o\`u on a not\'e $\theta_{\gamma}$ la fonction associ\'ee au quasi-caract\`ere $D^{G_{\epsilon}}_{(^{\gamma}f)'}$. Pour tout $\gamma\in \boldsymbol{\gamma}_{max}$, notons $\varphi_{\gamma}$ la fonction sur $\boldsymbol{\mathfrak{g}}_{\epsilon,{\cal F}'_{\gamma}}(k_{F})$, \`a support unipotent, telle que $\varphi_{\gamma}(Y)=(^{\gamma}f)'(exp(Y))$ pour tout \'el\'ement unipotent $Y\in \boldsymbol{\mathfrak{g}}_{\epsilon,{\cal F}'_{\gamma}}(k_{F})$. Il est clair que la fonction $Y\mapsto \theta_{\gamma}(exp(Y))$, d\'efinie sur $\mathfrak{g}_{\epsilon,tn}(F)$,  n'est autre que la fonction associ\'ee \`a la distribution $D^{G_{\epsilon}}_{\varphi_{\gamma}}$ sur $\mathfrak{g}_{\epsilon}(F)$. Notons $\theta_{\varphi_{\gamma}}$ cette fonction. On obtient la formule finale
$$(9) \qquad \theta(\epsilon exp(X))=(-1)^{a_{G_{\epsilon}}-a_{G}}m(A_{G_{\epsilon}})m(A_{G})^{-1}\sum_{\gamma\in \boldsymbol{\gamma}_{max}}m_{\gamma}\theta_{\varphi_{\gamma}}(X).$$
On applique le (ii) de la  proposition \ref{lespace}. Il entra\^{\i}ne que chaque fonction  $\theta_{\varphi_{\gamma}}$ est combinaison lin\'eaire de fonctions $Y\mapsto \hat{j}({\cal O},Y)$ o\`u ${\cal O}$ parcourt les orbites nilpotentes de $\mathfrak{g}_{\epsilon}(F)$. Il en est donc de m\^eme
de la fonction $X\mapsto \theta(\epsilon exp(X))$. Cela signifie que $D^G_{f}$ est un quasi-caract\`ere de niveau $0$ sur $G(F)$. $\square$

\subsection{Restriction d'un quasi-caract\`ere  sur $G(F)$ de niveau $0$ aux \'el\'ements elliptiques\label{restriction}}

Soit $D$ un quasi-caract\`ere sur $G(F)$, notons $\theta_{D}$ sa fonction associ\'ee. Disons que $D$ est de niveau $0$ sur les elliptiques si et seulement si $\theta_{D}$ v\'erifie la condition suivante:

(1) pour tout \'el\'ement $\epsilon\in G(F)$ qui est $p'$-compact mod $Z(G)$ et pour tout $X\in \mathfrak{g}_{\epsilon,tn}(F)$ tel que $\epsilon exp(X)\in G_{ell}(F)$, la fonction $\lambda\mapsto \theta_{D}(\epsilon exp(\lambda^2X))$ sur $\mathfrak{o}_{F}$ appartient \`a $E$.

{\bf Remarques} (2) Pour $\epsilon\in G(F)_{p'}$ et $X\in \mathfrak{g}_{\epsilon,tn}(F)$, $\epsilon exp(X)$ ne peut \^etre elliptique que si $A_{G_{\epsilon}}=A_{G}$. A fortiori, $\epsilon$ est $p'$-compact.

(3) La m\^eme preuve qu'au (ii) du lemme \ref{quasicaracteres} montre que (1) \'equivaut \`a la condition:

 pour tout \'el\'ement $\epsilon\in G(F)$ qui est $p'$-compact mod $Z(G)$ et pour tout ${\cal O}\in Nil(\mathfrak{g}_{\epsilon})$, il existe $c_{D,{\cal O}}\in {\mathbb C}$ de sorte que, pour tout $X\in \mathfrak{g}_{\epsilon,tn}(F)$ tel que $\epsilon exp(X)\in G_{ell}(F)$, on ait l'\'egalit\'e
$$\theta_{D}(\epsilon exp(X))=\sum_{{\cal O}\in Nil(\mathfrak{g}_{\epsilon})}c_{D,{\cal O}}\hat{j}({\cal O},X).$$

(4) Compte tenu de la remarque pr\'ec\'edente, la m\^eme preuve qu'au (iii) du lemme \ref{quasicaracteres} montre que, si $B$ est un sous-ensemble de $G(F)_{p'}$ tel que $G(F)=\cup_{\epsilon\in B} C_{G}(\epsilon)$, alors $D$ est un quasi-caract\`ere de niveau $0$ sur les elliptiques si et seulement si la condition (1) est v\'erifi\'ee pour tout $\epsilon\in B$. 

 \bigskip
\ass{Proposition}{Soit $D$ un quasi-caract\`ere sur $G(F)$ de niveau $0$ sur les elliptiques. Alors il existe ${\bf f}\in {\cal D}_{cusp}(G)$ tel que $D$ co\"{\i}ncide avec $D^G_{{\bf f}}$ sur $G_{ell}(F)$.  }

Preuve. Comme on l'a vu dans \ref{variantes}, l'espace ${\cal D}_{cusp}(G)$ se d\'ecompose en produit d'espaces index\'es par ${\cal N}$. Il en est de m\^eme de l'espace des quasi-caract\`eres de niveau $0$ sur les elliptiques. On peut donc fixer $\nu\in {\cal N}$ et supposer que
 $D$ est \`a support dans $w_{G}^{-1}(\nu)$. L'intersection $w_{G}^{-1}(\nu)\cap G_{ell}(F)$ est contenue dans une r\'eunion finie d'ensembles $ C_{G}(\epsilon )$, o\`u $\epsilon$ est $p'$-compact mod $Z(G)$.  Ces ensembles sont disjoints ou confondus et ils sont ouverts et ferm\'es. On peut fixer un \'el\'ement $p'$-compact mod $Z(G)$, supposer que $D$ est \`a support dans $ C_{G}(\epsilon)$ et prouver:

(5) il existe ${\bf f}\in C_{cusp}(G)$ tel que $D^G_{{\bf f}}$ soit \`a support dans $C_{G}(\epsilon)$ et co\"{\i}ncide avec $D$ sur $C_{G}(\epsilon)\cap G_{ell}(F)$. 

 Comme on l'a dit dans la remarque (2), on peut supposer $A_{G_{\epsilon}}=A_{G}$, sinon la solution de (5) est triviale. 
 
 Ecrivons comme dans  la remarque (3)
 $$\theta_{D}(\epsilon exp(X))=\sum_{{\cal O}\in Nil(\mathfrak{g}_{\epsilon})}c_{D,{\cal O}}\hat{j}({\cal O},X)$$
 pour tout $X\in \mathfrak{g}_{\epsilon,tn}(F)$ tel que $\epsilon exp(X)\in G_{ell}(F)$.
 Notons $\tau$ la fonction  sur $\mathfrak{g}_{\epsilon}(F)$, \`a support topologiquement nilpotent, qui est \'egale sur cet ensemble  au membre de droite ci-dessus.  En la moyennant sur le groupe $Z_{G}(\epsilon)(F)/G_{\epsilon}(F)$, on peut la supposer invariante par $Z_{G}(\epsilon)(F)$. La fonction $\tau$ est 
  la fonction associ\'ee \`a un quasi-caract\`ere de niveau $0$ sur $\mathfrak{g}_{\epsilon}(F)$. Celui-ci est de la forme $D^{G_{\epsilon}}_{\boldsymbol{\varphi}}$ pour un \'el\'ement $\boldsymbol{\varphi}\in {\cal D}(\mathfrak{g}_{\epsilon})$ d'apr\`es le (ii) de la proposition \ref{lespace}. Si on se restreint aux \'el\'ements elliptiques, les distributions induites disparaissent. Donc il existe $\boldsymbol{\varphi}\in {\cal D}_{cusp}(\mathfrak{g}_{\epsilon})$ tel que $\tau$ co\"{\i}ncide avec $\theta_{D^{G_{\epsilon}}_{\boldsymbol{\varphi}}}$ sur $\mathfrak{g}_{\epsilon,ell}(F)$. Fixons un tel   $\boldsymbol{\varphi}=(\varphi_{b})_{b\in  B}$, o\`u $B$ est un sous-ensemble fini de $Fac_{max}(G_{\epsilon})$.  Puisque $A_{G_{\epsilon}}=A_{G}$, $Imm(G_{\epsilon,AD})$ s'identifie \`a $Imm(G_{AD})^{\epsilon}$, cf. \ref{pointsfixes}. Pour tout $b\in B$, il y a donc une unique  facette ${\cal F}_{b}\in Imm(G_{AD})$ telle que  ${\cal F}_{b}\cap Imm(G_{AD})^{\epsilon}=b$.  On a  $\epsilon\in K_{{\cal F}_{b}}^{\nu}$ et ${\cal F}_{b}^{\nu}$ est r\'eduit \`a un point, autrement dit $({\cal F}_{b},\nu)\in Fac^*_{max}(G)$. De la fonction $\varphi_{b}$, on d\'eduit une fonction $f_{b}$ sur ${\bf G}_{{\cal F}_{b}}^{\nu}(k_{F})$ de la fa\c{c}on suivante. Par l'exponentielle, on identifie $\varphi_{b}$ \`a une fonction $f_{b,1}$ sur ${\bf G}_{\epsilon,b}(k_{F})$ \`a support nilpotent, telle que $f_{b,1}(exp(X))=\varphi_{b}(X)$ pour tout \'el\'ement nilpotent $X\in \boldsymbol{\mathfrak{g}}_{\epsilon,b}(k_{F})$. En notant $\bar{\epsilon}$ la r\'eduction de $\epsilon$ dans ${\bf G}_{{\cal F}_{b}}^{\nu}(k_{F})$, on d\'efinit une fonction $f_{b,2}$ sur ${\bf G}_{{\cal F}_{b}}^{\nu}(k_{F})$, \`a support dans $\bar{\epsilon}{\bf G}_{\epsilon,b}(k_{F})$, telle que $f_{b,2}(\bar{\epsilon}x)=f_{b,1}(x)$ pour tout $x\in {\bf G}_{\epsilon,b}(k_{F})$. Enfin, on pose
$$f_{b}(x)=\vert {\bf G}_{{\cal F}_{b}}(k_{F})\vert ^{-1}\sum_{y\in  {\bf G}_{{\cal F}_{b}}(k_{F})}f_{b,2}(y^{-1}xy)$$
pour tout $x\in {\bf G}_{{\cal F}_{b}}^{\nu}(k_{F})$. Montrons que

(6) $f_{b}\in C_{cusp}({\bf G}_{{\cal F}_{b}}^{\nu})$.

  On doit prouver que, pour tout  espace parabolique propre  ${\bf P}^{\nu}$ de ${\bf G}_{{\cal F}_{b}}^{\nu}$ et pour tout $x\in {\bf P}^{\nu}(k_{F})$, on a 
$$\sum_{u\in {\bf U}_{{\bf P}}(k_{F})}f_{b}(xu)=0.$$
Il suffit de prouver que, pour tout  espace parabolique propre  ${\bf P}^{\nu}$ de ${\bf G}_{{\cal F}_{b}}^{\nu}$ et pour tout $x\in {\bf P}^{\nu}(k_{F})$, on a 
$$(7) \qquad \sum_{u\in {\bf U}_{{\bf P}}(k_{F})}f_{b,2}(xu)=0.$$
En effet, la premi\`ere assertion pour ${\bf P}^{\nu}$ r\'esulte de la seconde appliqu\'ee \`a chaque espace $y^{-1}{\bf P}^{\nu}y$ pour $y\in {\bf G}_{{\cal F}_{b}}(k_{F})$. Fixons donc ${\bf P}^{\nu}$. Si ${\bf P}^{\nu}(k_{F})$ ne coupe pas le support de $f_{b,2}$, l'assertion (7) est claire. Supposons que ${\bf P}^{\nu}(k_{F})$ coupe ce support. Il existe donc un \'el\'ement nilpotent $X\in \boldsymbol{\mathfrak{g}}_{\epsilon,b}(k_{F})$ tel que $\bar{\epsilon}exp(X)\in {\bf P}^{\nu}(k_{F})$. Puisque $\bar{\epsilon}$ et $exp(X)$ commutent, puisque $ad(\bar{\epsilon})$ est d'ordre premier \`a $p$ tandis que $ad(exp(X))$ est d'ordre une puissance de $p$, $ad(\bar{\epsilon})$ appartient au groupe engendr\'e par $ad(\bar{\epsilon}exp(X))$. Ce dernier op\'erateur conserve ${\bf P}$, donc $ad(\bar{\epsilon})$ conserve lui-aussi ${\bf P}$, c'est-\`a-dire $\bar{\epsilon}\in {\bf P}^{\nu}(k_{F})$. D'apr\`es le lemme \ref{unlemme}, on peut fixer une composante de Levi ${\bf M}^{\nu}$ de ${\bf P}^{\nu}$ qui contient $\bar{\epsilon}$.   Notons  ${\bf M}_{\epsilon}$ et ${\bf P}_{\epsilon}$ les composantes neutres des groupes des points fixes de l'op\'erateur $ad(\bar{\epsilon})$ dans ${\bf M}$ et ${\bf P}$. Le groupe ${\bf P}_{\epsilon}$ est un sous-groupe parabolique de ${\bf G}_{\epsilon,b}$ et ${\bf M}_{\epsilon}$ en est une composante de Levi. Montrons que

(8) ${\bf P}_{\epsilon}\not={\bf G}_{\epsilon,b}$.

Notons ${\bf A}_{{\bf M}}$, resp. ${\bf A}_{{\cal F}_{b}}$ le plus grand sous-tore central d\'eploy\'e de ${\bf M}$, resp. ${\bf G}_{{\bf F}_{b}}$. Notons ${\bf A}_{{\bf M},\epsilon}$, resp. ${\bf A}_{{\cal F}_{b},\epsilon}$, la composante neutre du sous-groupe des points fixes de l'op\'erateur $ad(\bar{\epsilon})$ dans ${\bf A}_{{\bf M}}$, resp. ${\bf A}_{{\cal F}_{b}}$. Puisque ${\bf P}^{\nu}$ est propre, ${\bf A}_{{\bf M},\epsilon}$ contient strictement ${\bf A}_{{\cal F}_{b},\epsilon}$.  Donc $dim({\bf A}_{{\bf M},\epsilon})>dim({\bf A}_{{\cal F}_{b},\epsilon})\geq dim(A_{G})$. Le tore ${\bf A}_{{\bf M},\epsilon}$ est contenu dans le plus grand sous-tore central d\'eploy\'e ${\bf A}_{{\bf M}_{\epsilon}}$ de ${\bf M}_{\epsilon}$. Donc $dim({\bf A}_{{\bf M}_{\epsilon}})> dim(A_{G})$. On a suppos\'e $dim(A_{G})=dim(A_{G_{\epsilon}})$ et on a $dim(A_{G_{\epsilon}})=dim({\bf A}_{{\bf G}_{\epsilon,b}})$ puisque $b$ est un sommet de $Imm(G_{\epsilon,AD})$.  Donc $dim({\bf A}_{{\bf M}_{\epsilon}})> dim({\bf A}_{{\bf G}_{\epsilon,b}})$, ce qui signifie que ${\bf M}_{\epsilon}$ est un Levi propre de ${\bf G}_{\epsilon,b}$. Cela d\'emontre (8).

  On peut dans (7) se restreindre au cas o\`u $x\in {\bf M}^{\nu}(k_{F})$, c'est-\`a-dire $x=\bar{\epsilon}m$ avec $m\in {\bf M}(k_{F})$. Montrons que

(9) si $m\not\in exp(\boldsymbol{\mathfrak{m}}_{\epsilon,nil})(k_{F})$, l'ensemble des $u\in {\bf U}_{{\bf P}}(k_{F})$ tels que $xu$ appartient au support de $f_{b,2}$ est vide; si $m=exp(Y)$ avec $Y\in \boldsymbol{\mathfrak{m}}_{\epsilon,nil}(k_{F})$, alors cet ensemble est contenu dans celui  des $exp(N)$ pour $N\in \boldsymbol{\mathfrak{u}}_{{\bf P}_{\epsilon}}(k_{F})$.

Dire que $xu$ appartient au support de $f_{b,2}$ implique que $xu=\epsilon exp(X)$ avec $X\in \boldsymbol{\mathfrak{g}}_{\epsilon,b,nil}(k_{F})$, ou encore $mu=exp(X)$. Si on applique $ad(\bar{\epsilon})$ \`a cette \'egalit\'e, puisque $X$ commute \`a $\bar{\epsilon}$ et que $ad(\bar{\epsilon})$ conserve ${\bf M}$ et ${\bf U}_{{\bf P}}$, on obtient $mu=ad(\bar{\epsilon})(m)ad(\bar{\epsilon})(u)$, ce qui entra\^{\i}ne $m=ad(\bar{\epsilon})(m)$ et $u=ad(\bar{\epsilon})(u)$. Cette derni\`ere relation implique que $u=exp(N)$ pour un $N\in \boldsymbol{\mathfrak{u}}_{{\bf P}_{\epsilon}}(k_{F})$. La premi\`ere, jointe au fait que $m$ est forc\'ement unipotent (puisque $exp(X)$ l'est), entra\^{\i}ne que $m=exp(Y)$ pour un $Y\in \boldsymbol{\mathfrak{m}}_{\epsilon,nil}(k_{F})$. Cela d\'emontre (9). 

On peut donc supposer $m=exp(Y)$ avec $Y\in \boldsymbol{\mathfrak{m}}_{\epsilon,nil}(k_{F})$. Il r\'esulte de (9) et de la d\'efinition de $f_{b,2}$ que

$$\sum_{u\in {\bf U}_{{\bf P}}(k_{F})}f_{b,2}(xu)=\sum_{N\in  \boldsymbol{\mathfrak{u}}_{{\bf P}_{\epsilon}}(k_{F}) }\varphi_{b}(Y+N).$$
 Puisque ${\bf P}_{\epsilon}$ est un sous-groupe parabolique propre de ${\bf G}_{\epsilon,b}$,  ceci est nul puisque $\varphi_{b}$ est cuspidale. Cela ach\`eve la preuve de (6).

 On d\'efinit la distribution $D^G_{f_{b}}$. Notons $\theta_{b}$ sa fonction associ\'ee. Il est clair que son support est contenu dans la r\'eunion des conjugu\'es du support de $f_{b}$. D'apr\`es la construction de cette fonction et le lemme \ref{commutant}, ce support est contenu dans l'ensemble des $k^{-1}C(\epsilon)k$ pour $k\in K_{{\cal F}_{b}}^0$. Donc le support de $\theta_{b}$ est contenu dans $C_{G}(\epsilon)$ et il nous suffit de calculer $\theta_{b}(\epsilon exp(X))$ pour $X\in \mathfrak{g}_{\epsilon,tn}(F)$. 
C'est le calcul que l'on a fait dans le paragraphe pr\'ec\'edent.  L'ensemble $\Gamma_{\epsilon}$ contient  $Z_{G}(\epsilon)(F)K_{{\cal F}_{b}}^0$. Montrons que 

(10) on peut remplacer l'ensemble $\Gamma_{\epsilon}$ par $Z_{G}(\epsilon)(F)K_{{\cal F}_{b}}^0$.

Il suffit de montrer que, si $\gamma\in \Gamma_{\epsilon}-Z_{G}(\epsilon)(F)K_{{\cal F}_{b}}^0$ et si $g\in G_{\epsilon}(F)\gamma K_{{\cal F}_{b}}^0$, $g^{-1}\epsilon exp(X) g$ n'appartient pas au support de $f_{b}$. Supposons que $g^{-1}\epsilon exp(X) g$ appartienne \`a ce support. Comme on l'a dit ci-dessus, il existe donc $k\in K_{{\cal F}_{b}}^0$ et $Y\in \mathfrak{g}_{\epsilon,tn}(F)$ tels que $g^{-1}\epsilon exp(X)g=k^{-1}\epsilon exp(Y)k$. En fixant un entier $c\geq1$ premier \`a $p$ tel que $\epsilon^c\in Z(G)(F)$, et en \'elevant l'\'egalit\'e pr\'ec\'edente \`a la puissance $c$, on obtient
$exp(cg^{-1}Xg)=exp(ck^{-1}Yk)$, d'o\`u  $exp(g^{-1}Xg)=exp(k^{-1}Yk)$, d'o\`u aussi $g^{-1}\epsilon g=k^{-1} \epsilon k$. Mais alors $gk^{-1}\in Z_{G}(\epsilon)(F)$ et $\gamma\in Z_{G}(\epsilon)(F)K_{{\cal F}_{b}}^0 $ contrairement \`a l'hypoth\`ese. Cela d\'emontre (10).

 On peut donc supposer $\boldsymbol{\gamma}\subset Z_{G}(\epsilon)(F)$. On modifie l\'eg\`erement les calculs du paragraphe pr\'ec\'edent en remarquant que la contribution d'un $\gamma\in Z_{G}(\epsilon)(F)$ \`a $\theta_{b}(\epsilon exp(X))$ est \'egale \`a celle de $\gamma=1$, o\`u l'on  remplace $X$ par $\gamma^{-1} X\gamma$. L'assertion (9) de \ref{unpremiertheoreme} devient
 $$\theta_{b}(\epsilon exp(X))=c_{b} \vert \boldsymbol{\gamma}\vert ^{-1}\sum_{\gamma\in \boldsymbol{\gamma}}
 \theta_{\varphi_{b}}(\gamma^{-1}X\gamma),$$
 o\`u $c_{b}$ est une certaine constante non nulle. 
On pose ${\bf f}=(c_{b}^{-1}f_{b})_{b\in B}$.  C'est un \'el\'ement de ${\cal D}_{cusp}(G)$. 
D'apr\`es les calculs ci-dessus, $D^G_{{\bf f}}$ est \`a support dans $C_{G}(\epsilon)$ et, en notant $\theta'$ sa fonction associ\'ee, on a
$$\theta'(\epsilon exp(X))= \vert \boldsymbol{\gamma}\vert ^{-1}\sum_{\gamma\in \boldsymbol{\gamma}}\sum_{b\in B}\theta_{\varphi_{b}}(\gamma^{-1}X\gamma)$$
$$= \vert \boldsymbol{\gamma}\vert ^{-1}\sum_{\gamma\in \boldsymbol{\gamma}} \tau(\gamma^{-1} X\gamma)$$
$$=\tau(X),$$
puisque   la fonction $\tau$ est invariante par $Z_{G}(\epsilon)(F)$. Donc $\theta'(\epsilon exp(X))=\theta_{D}(\epsilon exp(X))$ si $X\in \mathfrak{g}_{\epsilon,tn}(F)$ et $\epsilon exp(X)\in G_{ell}(F)$. Cela  prouve (5) et la proposition. $\square$

 \subsection{Filtration sur le groupe et quasi-caract\`eres\label{filtrationsurlegroupe}}
 
  Notons $Qc(G)$ l'espace des quasi-caract\`eres sur $G(F)$. Notons $Qc^n(G)$ la somme des $Ind_{M}^G(Qc(M))$ sur les Levi $M$ tels que $dim(A_{M})\geq n$.  

\ass{Proposition}{Pour tout $n\in {\mathbb Z}$, on a l'\'egalit\'e
 $Qc^n(G)=Qc(G)\cap Ann^nI(G)^*$.}

Preuve. L'inclusion $Qc^n(G)\subset Qc(G)\cap Ann^nI(G)^*$ est \'evidente. Soit $D\in Qc(G)\cap Ann^nI(G)^*$, notons $\theta$ sa fonction associ\'ee. L'hypoth\`ese que $D$ appartient \`a $Ann^nI(G)^*$ signifie que, pour tout Levi $M$ tel que $dim(A_{M})<n$ et tout $m\in M_{ell}(F)\cap G_{reg}(F)$, on a $\theta(m)=0$. Pour $M\in \underline{{\cal L}}^n$, notons $\theta_{M}$ la fonction d\'efinie presque partout sur $M(F)$ par $\theta_{M}(x)=D^G(x)^{1/2}D^M(x)^{-1/2}\theta(x)$ pour $x\in M(F)$.  Notons  $D_{M}$ la distribution qui s'en d\'eduit. 
Par l'isomorphisme (1) de \ref{filtrations}, $D$ s'envoie sur la somme des restrictions des $D_{M}$ \`a $I_{cusp}(M)^{W^G(M)}$.  Montrons que

(2) il existe un quasi-caract\`ere $D'_{M}\in Qc(M)$ tel que $D_{M}$ et $D'_{M}$ aient m\^eme restriction \`a $I_{cusp}(M)$. 

Commen\c{c}ons par fixer un \'el\'ement semi-simple $x\in M(F)$ et prouvons que

(3) il existe un voisinage $\mathfrak{ V}$ de $0$ dans $\mathfrak{m}_{x}(F)$ et, pour tout ${\cal O}\in Nil(\mathfrak{m}_{x})$, il existe un nombre complexe $c_{D,{\cal O}}$ de sorte que, pour presque tout $Y\in \mathfrak{ V}$ tel que $x exp(Y)$ soit elliptique dans $M(F)$, on ait l'\'egalit\'e
$$\theta_{M}(x exp(Y))=\sum_{{\cal O}\in Nil(\mathfrak{m}_{x})}c_{D,{\cal O}}\hat{j}({\cal O},Y).$$

Si $dim(A_{M_{x}})> dim(A_{M})$,  un tel \'el\'ement $x exp(Y)$ n'est jamais elliptique dans $M(F)$ et l'assertion est triviale. Supposons $dim(A_{M_{x}})=dim(A_{M})=n$. Puisque $D$ est un quasi-caract\`ere, on peut en tout cas fixer un voisinage $\mathfrak{ V}'$  de $0$ dans $\mathfrak{g}_{x}(F)$ et, pour tout ${\cal O}\in Nil(\mathfrak{g}_{x})$,   un nombre complexe $c'_{D,{\cal O}}$ de sorte que, pour presque tout $Y\in \mathfrak{V}'$, on ait l'\'egalit\'e
$$(4) \qquad \theta(x exp(Y))=\sum_{{\cal O}\in Nil(\mathfrak{g}_{x})}c'_{D,{\cal O}}\hat{j}({\cal O},Y).$$
Notons $\tau$ la fonction du membre de droite, qui est d\'efinie presque partout sur $\mathfrak{g}_{x}(F)$. 
Soit $L^x$ un Levi de $G_{x}$ tel que $dim(A_{L^x})<n$ et soit $Y\in \mathfrak{l}^x_{ell}(F)\cap \mathfrak{V}'$. On note $L$ le commutant de $A_{L^x}$ dans $G$. Alors $x\in L(F)$, $L_{x}=L^x$,  $A_{L}=A_{L^x}$ et $x exp(Y)$ est elliptique dans $L(F)$. Puisque $dim(A_{L})<n$, $\theta(x exp(Y))=0$. Donc $\tau(Y)=0$. Puisque $\tau$ est  combinaison lin\'eaire de fonctions homog\`enes pour l'action de $F^{\times,2}$ par homoth\'eties, on peut supprimer ci-dessus la restriction: $Y\in \mathfrak{V}'$. La restriction de $\tau$ \`a $\mathfrak{g}_{x,tn}(F)$ est donc la fonction associ\'ee \`a un \'el\'ement de $Qc_{0}(\mathfrak{g}_{x})\cap Ann^{n}I(\mathfrak{g}_{x})^*$, avec les notations de \ref{filtrationssurlalgebredelie}. D'apr\`es le lemme de ce paragraphe, on peut \'ecrire cet \'el\'ement comme somme sur les $L^x\in \underline{{\cal L}}^{G_{x},n}$ de distributions $Ind_{L^x}^{G_{x}}(d_{L^x})$, o\`u $d_{L^x}\in Qc(\mathfrak{l}^x)$. On peut supposer que $M_{x}\in \underline{{\cal L}}^{G_{x},n}$ et que $d_{M_{x}}$ est invariante par $ W^{G_{x}}(M_{x})$. Notons $\tau_{x}$ la fonction associ\'ee \`a $d_{M_{x}}$. En utilisant \ref{induction} (3), on  v\'erifie que l'on a l'\'egalit\'e
$$\tau(Y)=cD^{G_{Y}}(Y)^{-1/2}D^{M_{Y}}(Y)^{1/2}\tau_{x}(Y)$$
pour  presque tout $Y\in 
  \mathfrak{m}_{x,ell}(F)\cap \mathfrak{m}_{x,tn}(F)$, o\`u $c=\vert W^{G_{x}}(M_{x})\vert  $. D'autre part, si $\mathfrak{V}'$ est assez petit, il existe $c'>0$ tel que  $D^G(x exp(Y))^{1/2}D^M(x exp(Y))^{-1/2}=c'D^{G_{Y}}(Y)^{1/2}D^{M_{Y}}(Y)^{-1/2}$ pour tout $Y\in \mathfrak{V}'\cap \mathfrak{m}_{x}(F)$ tel que $x exp(Y)\in G_{reg}(F)$. L'\'egalit\'e pr\'ec\'edente devient
  $$\theta_{M}(x exp(Y))=cc'\tau_{x}(Y)$$
 pour  presque tout $Y\in 
  \mathfrak{V}'\cap\mathfrak{m}_{x,ell}(F)\cap \mathfrak{m}_{x,tn}(F)$.  
  Cette fonction est de la forme du membre de droite de l'\'egalit\'e de l'assertion (3). En prenant $\mathfrak{V}=\mathfrak{V}'\cap \mathfrak{m}_{x}(F)$, cela r\'esout cette assertion (3). 

 On peut  supposer que le membre de droite de l'\'egalit\'e de cette assertion  est invariant par l'action du groupe $Norm_{G_{x}}(M_{x})(F) $. Posons $V_{x}=\{m^{-1}x exp(Y)m; m\in M(F), Y\in \mathfrak{V}\}$. Quitte \`a restreindre $\mathfrak{V}$, cet ensemble est ouvert et ferm\'e dans $M(F)$. On d\'efinit un quasi-caract\`ere $d_{x}$ de $M(F)$, \`a support dans $V_{x}$, dont la fonction associ\'ee v\'erifie
$$\theta_{d_{x}}(m^{-1}x exp(Y)m)=\sum_{{\cal O}\in Nil(\mathfrak{m}_{x})}c_{D,{\cal O}}\hat{j}({\cal O},Y)$$
pour tout $m\in M(F)$ et $Y\in \mathfrak{V}$.  
Cette fonction co\"{\i}ncide avec $\theta_{M}$ sur $V_{x}\cap M_{ell}(F)$. On construit ces ensembles $V_{x}$ et ces quasi-caract\`eres $d_{x}$ pour tout \'el\'ement semi-simple $x\in M(F)$. La r\'eunion des $V_{x}$ est $M(F)$ tout entier, on peut en extraire un recouvrement localement fini $M(F)=\cup_{j\in J}V_{x_{j}}$. Il existe une famille $(\alpha_{j})_{j\in J}$, o\`u $\alpha_{j}$ est une fonction localement constante sur $M(F)$, invariante par conjugaison et \`a support dans $V_{x_{j}}$, de sorte que $\sum_{j\in J}\alpha_{j}(m)=1$ pour tout $m\in M(F)$. En utilisant la remarque (5) de \ref{quasicaracteres}, la distribution 
$$D'_{M}=\sum_{j\in J}\alpha_{j}d_{x_{j}}$$
est bien d\'efinie puisque le recouvrement est localement fini, c'est un quasi-caract\`ere et elle co\"{\i}ncide avec $D_{M}$ sur les \'el\'ements elliptiques de $M(F)$, donc elle a m\^eme restriction que $D_{M}$ \`a $I_{cusp}(M)$. Cela d\'emontre (2). 

 La distribution
 $$\sum_{M\in \underline{{\cal L}}^n} \vert W^G(M)\vert ^{-1}Ind_{M}^G(D'_{M})$$
 appartient \`a $Qc^n(G)$ et, d'apr\`es \ref{filtrations},  son image dans  $Gr^nI(G)^*$ est  \'egale \`a celle de $D$. Cela d\'emontre l'inclusion
 $$Qc(G)\cap Ann^{n}I(G)^*\subset (Qc(G)\cap Ann^{n+1}I(G)^*)+Qc^n(G).$$
 L'inclusion $Qc(G)\cap Ann^{n}I(G)^*\subset Qc^n(G)$ s'en d\'eduit par r\'ecurrence descendante sur $n$. Cela d\'emontre la proposition.$\square$ 
 
 \subsection{Un deuxi\`eme th\'eor\`eme\label{undeuxiemetheoreme}}

   Notons $Qc_{0}(G)$ l'espace des quasi-caract\`eres de niveau $0$ sur $G(F)$. 
   
   \ass{Th\'eor\`eme}{L'application $D^G:{\cal D}(G)\to I(G)^*$  est un isomorphisme de ${\cal D}(G)$ sur $Qc_{0}(G)$.}
   
   Preuve. Compte tenu du th\'eor\`eme \ref{unpremiertheoreme}, il suffit de prouver que $Qc_{0}(G)$ est contenu dans l'image de l'application $D^G$. Soit $n\in {\mathbb Z}$. Prouvons que
   
   (1) $Qc_{0}(G)\cap Ann^nI(G)^*$ est contenu dans la somme de $D^G({\cal D}^n(G))$ et de $Qc_{0}(G)\cap Ann^{n+1}I(G)^*$.
   
 Remarquons que $D^G({\cal D}^n(G))$ est contenu dans $Qc_{0}(G)\cap Ann^nI(G)^*$  d'apr\`es \ref{unpremiertheoreme}. Soit $D\in Qc_{0}(G)\cap Ann^nI(G)^*$. Pour $M\in \underline{{\cal L}}_{min}^n$, on a d\'efini dans le paragraphe pr\'ec\'edent la distribution $D_{M}$ sur $M(F)$, de fonction associ\'ee $\theta_{M}$, et on a prouv\'e l'existence  d'un quasi-caract\`ere $D'_{M}$ sur $M(F)$ qui   co\"{\i}ncidaient avec $D_{M}$ sur  $M_{ell}(F)$. Montrons que 
 
 (2) $D'_{M}$ est de niveau $0$ sur les elliptiques.
 
  D'apr\`es les remarques (2) et (4) de \ref{restriction} et le lemme \ref{lecasdunlevi}, on peut fixer $\epsilon\in G(F)_{p'}\cap M(F)$  et d\'emontrer que, pour tout $X\in \mathfrak{m}_{\epsilon,tn}(F)$ tel que $\epsilon exp(X)\in M_{ell}(F)$, la fonction $\lambda\mapsto \theta_{D'_{M}}(\epsilon exp(\lambda^2X))$ appartient \`a $E$. Pour un tel $X$, on a par d\'efinition
 $$\theta_{D'_{M}}(\epsilon exp(\lambda^2X))=\theta_{M}(\epsilon exp(\lambda^2X))=D^G(\epsilon exp(\lambda^2X))^{1/2}D^M(\epsilon exp(\lambda^2X))^{-1/2}\theta_{D}(\epsilon exp(\lambda^2X)).$$
 D'apr\`es le (ii) du lemme \ref{decomposition}, la fonction $\lambda\mapsto D^G(\epsilon exp(\lambda^2X))^{1/2}D^M(\epsilon exp(\lambda^2X))^{-1/2}$ est produit d'une constante et d'une puissance enti\`ere de $\vert \lambda\vert $. Puisque $D\in Qc_{0}(G)$, la fonction $\lambda\mapsto \theta_{D}(\epsilon exp(\lambda^2X))$ appartient \`a $E$. Il en r\'esulte que la fonction $\lambda\mapsto \theta_{D'_{M}}(\epsilon exp(\lambda^2X))$ appartient \`a $E$, d'o\`u (2). 
 
 D'apr\`es la proposition \ref{restriction}, il existe donc $D''_{M}\in D^M({\cal D}_{cusp}(M))$ tel que $D''_{M}$ co\"{\i}ncide avec $D'_{M}$ sur $M_{ell}(F)$, donc aussi avec $D_{M}$.  La distribution
 $$\sum_{M\in \underline{{\cal L}}^n} \vert W^G(M)\vert ^{-1}Ind_{M}^G(D''_{M})$$
 appartient alors \`a $D^G({\cal D}^n(G))$ et son image dans $Gr^nI(G)^*$ est \'egale \`a celle de $D$. Cela prouve (1).
 
 Pour $n=a_{M_{min}}+1$, on a $Qc_{0}(G)\cap Ann^nI(G)^*=\{0\}=D^G({\cal D}^n(G))$. Gr\^ace \`a (1), on d\'emontre par r\'ecurrence descendante que $Qc_{0}(G)\cap Ann^nI(G)^*=D^G({\cal D}^n(G))$ pour tout $n$. Pour $n=a_{G}$, cela d\'emontre l'\'egalit\'e $Qc_{0}(G)=D^G({\cal D}(G))$. $\square$

  \section{Repr\'esentations et repr\'esentations de niveau $0$\label{representations}}

 \subsection{Repr\'esentations de niveau $0$}
Notons $Irr(G)$ l'ensemble des classes d'isomorphismes de repr\'esentations lisses irr\'eductibles de $G(F)$ (dans des espaces complexes). Une telle repr\'esentation $\pi$, dans un espace $V$, est dite de niveau $0$ si et seulement s'il existe une facette ${\cal F}\in Fac(G)$ telle que l'espace des invariants $V^{K_{{\cal F}}^+}$ soit non nul. On note $Irr(G)^0$ le sous-ensemble des (classes d'isomorphismes de) repr\'esentations de niveau $0$ et $Irr(G)^{>0}$ celui  des repr\'esentations qui ne sont pas de niveau $0$. On note $p^0 $ la projection  ${\mathbb C}[Irr(G)]\to {\mathbb C}[Irr(G)^0]$ qui  annule ${\mathbb C}[Irr(G)^{>0}]$. On sait que $p^0$ est associ\'ee \`a l'action d'un idempotent du centre de Bernstein.  Il se d\'eduit de cet idempotent une action sur divers objets, par exemple $C_{c}^{\infty}(G(F))$, que l'on note encore $p^0$. On sait que $p^0$ commute, en un sens compr\'ehensible, aux op\'erations d'induction et de passage au module de Jacquet: par exemple, si $P$ est un sous-groupe parabolique de $G$ de composante de Levi $M$ et si $\pi^M\in Irr(M)$, alors $p^0(Ind_{P}^G(\pi^M))=Ind_{P}^G(p^0(\pi^M))$. 

On note $Temp(G)$ le sous-ensemble de $Irr(G)$ des repr\'esentations temp\'er\'ees et on pose  $Temp(G)^0=Temp(G)\cap Irr(G)^0$.

Si $\xi$ est un caract\`ere de $A_{G}(F)$, on note $Irr(G)_{\xi}$  le sous-ensemble de $Irr(G)$ des repr\'esentations  dont le caract\`ere central co\"{\i}ncide avec $\xi$ sur $A_{G}(F)$. On d\'efinit de m\^eme $Irr(G)_{\xi}^0$, $Temp(G)_{\xi}$ etc... Evidemment, $Irr(G)_{\xi}^0$ est vide si $\xi$ n'est pas mod\'er\'ement ramifi\'e et $Temp(G)_{\xi}$ est vide si $\xi$ n'est pas unitaire.  

Pour toute repr\'esentation irr\'eductible $\pi$, on note $\Theta_{\pi}$
   son caract\`ere-distribution.  D'apr\`es Harish-Chandra, c'est un quasi-caract\`ere. On note $\theta_{\pi}$ sa fonction associ\'ee. Par lin\'earit\'e, l'application $\pi\mapsto \Theta_{\pi}$ se prolonge en une application
  $$\Theta:{\mathbb C}[Irr(G)]\to I(G)^*$$
  qui est injective.

\subsection{Repr\'esentations elliptiques\label{representationselliptiques}}
Arthur a d\'efini un ensemble de "repr\'esentations elliptiques", cf. \cite{A2} paragraphe 3.  Une telle repr\'esentation est une combinaison lin\'eaire finie, \`a coefficients dans ${\mathbb C}$, de repr\'esentations irr\'eductibles temp\'er\'ees.  Sa d\'efinition d\'epend de divers choix,   en particulier de  normalisations d'op\'erateurs d'entrelacement. Nous supposons ces choix fix\'es et nous notons $Ell(G)$ l'ensemble des repr\'esentations elliptiques. C'est un ensemble lin\'eairement ind\'ependant. L'espace ${\mathbb C}[Ell(G)]$ ne d\'epend pas des choix. Pour un Levi $M\in {\cal L}_{min}$, le sous-espace ${\mathbb C}[Ell(M)]\subset {\mathbb C}[Temp(M)]$ est invariant par l'action naturelle du groupe $ W^G(M)$  sur ${\mathbb C}[Temp(M)]$. On note ${\mathbb C}[Ell(M)]^{W^G(M)}$ le sous-espace des invariants. On a alors la d\'ecomposition en somme directe
$${\mathbb C}[Temp(G)]=\oplus_{M\in\underline{{\cal L}}_{min}}Ind_{M}^G({\mathbb C}[Ell(M)]^{W^G(M)}).$$

Introduisons l'application d'Harish-Chandra $H_{G}:G(F)\to {\cal A}_{G}$, notons ${\cal A}_{G}^*$ l'espace vectoriel r\'eel dual de ${\cal A}_{G}$ et  ${\cal A}_{G,{\mathbb C}}^*$ son complexifi\'e. Un \'el\'ement $\chi\in {\cal A}_{G,{\mathbb C}}^*$ d\'efinit un caract\`ere de $G(F)$, que l'on note encore $\chi$, d\'efini par $\chi(g)=exp(<\chi,H_{G}(g)>)$. Si $\pi\in Irr(G)$, on note $\pi_{\chi}$ la repr\'esentation $\pi_{\chi}(g)=\chi(g)\pi(g)$.  Si $\chi\in i{\cal A}_{G}^*$, cette op\'eration $\pi\mapsto \pi_{\chi}$ conserve l'espace ${\mathbb C}[Ell(G)]$. On note $Ell(G)_{{\mathbb R}}$ l'ensemble des $\pi_{\chi}$ pour $\pi\in Ell(G)$ et $\chi\in {\cal A}_{G}^*$. On a alors la d\'ecomposition plus g\'en\'erale
$$(1) \qquad {\mathbb C}[Irr(G)]=\oplus_{M\in\underline{{\cal L}}_{min}}Ind_{M}^G({\mathbb C}[Ell(M)_{{\mathbb R}}]^{W^G(M)}).$$

Via l'injection $\Theta$, cette d\'ecomposition est compatible avec la filtration de $I(G)^*$ d\'efinie en \ref{filtrations}.  C'est-\`a-dire que, pour $n\in {\mathbb Z}$, l'intersection de $Ann^nI(G)^*$ et de l'image de $\Theta$ est \'egale \`a l'image par cette injection de la sous-somme de l'expression ci-dessus, o\`u on limite la sommation aux $M$ tels que $dim(A_{M})\geq n$. En particulier, l'application compos\'ee
$$(2) \qquad {\mathbb C}[Ell(G)_{{\mathbb R}}]\stackrel{\Theta}{\to}I(G)^*\to I_{cusp}(G)^*$$
est injective.

En fait, chaque repr\'esentation elliptique est combinaison lin\'eaire de sous-repr\'esentations d'une m\^eme induite d'une s\'erie discr\`ete d'un Levi de $G$.  Il en r\'esulte d'une part que toute repr\'esentation elliptique admet un caract\`ere central. D'autre part que ou bien toutes ces composantes sont de niveau $0$, ou bien aucune ne l'est. On note $Ell(G)^0$, resp. $Ell(G)^{>0}$, les repr\'esentations elliptiques du premier cas, resp. du second. Puisque $p^0$ "commute" \`a l'induction, on obtient des variantes des \'egalit\'es pr\'ec\'edentes telles que
$$(3) \qquad {\mathbb C}[Temp(G)^0]=\oplus_{M\in\underline{{\cal L}}_{min}}Ind_{M}^G({\mathbb C}[Ell(M)^0]^{W^G(M)});$$
$${\mathbb C}[Irr(G)^0]=\oplus_{M\in\underline{{\cal L}}_{min}}Ind_{M}^G({\mathbb C}[Ell(M)^0_{{\mathbb R}}]^{W^G(M)});$$
$${\mathbb C}[Irr(G)^{>0}]=\oplus_{M\in\underline{{\cal L}}_{min}}Ind_{M}^G({\mathbb C}[Ell(M)_{{\mathbb R}}^{>0}]^{W^G(M)}).$$
Il y a aussi des variantes de ces \'egalit\'es  en imposant que toutes les repr\'esentations se transforment selon un caract\`ere fix\'e $\xi$ de $A_{G}(F)$.

\subsection{Le r\'esultat de \cite{W2}\label{leresultat}}

Soit $\pi\in Irr(G)^0$. En \cite{W2}, on a associ\'e \`a $\pi$ un \'el\'ement de $  {\cal D}_{cusp}(G)$ que l'on avait not\'e $\Theta^G_{\pi,cusp}$. On modifie l\'eg\`erement la d\'efinition en y incluant les mesures $mes(A_{G}(F)\backslash K_{{\cal F}}^{\dag})^{-1}$ qui intervenaient dans les formules de cette r\'ef\'erence. On note   $\Delta^G_{\pi,cusp}$ l'\'el\'ement de ${\cal D}_{cusp}(G)$ ainsi obtenu.

Soit $M$ un Levi de $G$. Choisissons un sous-groupe parabolique $P\in {\cal P}(M)$. On associe au module de Jacquet $\pi_{P}$ un \'el\'ement $\Delta^M_{\pi_{P},cusp}$. Notons $\Delta^M_{\pi_{M},cusp,G-comp}$ sa projection dans ${\cal D}_{cusp,G-comp}(M)$, autrement dit sa restriction au sous-ensemble des \'el\'ements de $M(F)$ qui sont compacts mod $Z(G)$. Elle ne d\'epend pas du choix de $P$.

Le r\'esultat principal de \cite{W2} est que

(1)  $\Theta_{\pi}$ co\"{\i}ncide sur $G_{comp}(F)$ avec la distribution
$$\sum_{M\in {\cal L}_{min}}\vert W^M\vert \vert W^G\vert ^{-1}D^G[\Delta^M_{\pi_{M},cusp,G-comp}].$$
\bigskip
Autrement dit, notons $\Delta_{\pi}$ l'image naturelle de 
$$\sum_{M\in {\cal L}_{min}}\vert W^M\vert \vert W^G\vert ^{-1}\Delta^M_{\pi_{M},cusp,G-comp}$$
dans ${\cal D}_{G-comp}(G)$. Alors $\Theta_{\pi}$ co\"{\i}ncide avec $D^G[\Delta_{\pi}]$ sur $G_{comp}(F)$. On prolonge $\pi\mapsto \Delta_{\pi}$ par lin\'earit\'e  en une application 
$$\Delta:{\mathbb C}[Irr(G)^0]\to  {\cal D}_{G-comp}(G).$$

{\bf Remarque.} La somme intervenant dans (1) peut se r\'ecrire

$$(2) \qquad\sum_{M\in \underline{{\cal L}}_{min}}\vert W^G(M)\vert ^{-1}D^G[\Delta^M_{\pi_{M},cusp,G-comp}].$$

\subsection{Surjectivit\'e\label{surjectivite}}

Soit $\xi$  un caract\`ere unitaire et   mod\'er\'ement ramifi\'e  de $A_{G}(F)$. L'application $\Delta$ se restreint en une application lin\'eaire ${\mathbb C}[Temp(G)_{\xi}^0]\to  {\cal D}_{G-comp,\xi}(G)$.

\ass{Lemme}{(i) L'application $\Delta:{\mathbb C}[Temp(G)_{\xi}^0]\to  {\cal D}_{G-comp,\xi}(G)$ est surjective.

(ii) L'application $\Delta_{cusp}:{\mathbb C}[Ell(G)_{\xi}^0]\to  {\cal D}_{cusp,\xi}(G)$ est bijective.
 }

Preuve. Notons ici $q_{\xi}$ la projection de $C_{c}^{\infty}(G(F))$ dans $C_{c,\xi}^{\infty}(G(F))$ ou de $I(G)$ dans $I_{\xi}(G)$. On a le diagramme naturel
$$\begin{array}{ccccc}&&{\cal D}_{G-comp,\xi}(G)&&\\

&\Delta^G\nearrow\,\,& &p_{1}\searrow\,\,&\\

{\mathbb C}[Temp(G)_{\xi}^0]&&D^G\downarrow\,\,&& (q_{\xi}(I{\cal E}(G)))^*\\

&\Theta\searrow\,\,&&p_{2}\nearrow\,\,&\\

&& I_{\xi}(G)^*&&\\ \end{array}$$
L'application $p_{2}$ est la restriction et $p_{1}=p_{2}\circ D^G$. Le triangle de droite est donc commutatif. Celui de gauche l'est d'apr\`es \ref{leresultat}(1).    Montrons que

(1) $p_{2}\circ \Theta$ est surjective. 

Notons que l'espace $(q_{\xi}(I{\cal E}(G)))^*$ est de dimension finie.   Chaque \'el\'ement de ${\cal E}(G)$ est invariant par un groupe $K_{{\cal F}}^+$ pour une certaine facette ${\cal F}\in Fac(G)$. Donc ${\cal E}(G)$ est annul\'e par toute repr\'esentation qui n'est pas de niveau $0$. D'apr\`es le th\'eor\`eme de Paley-Wiener (avec caract\`ere central), l'application ${\mathbb C}[Temp(G)_{\xi}]\to (q_{\xi}(I{\cal E}(G)))^*$ similaire \`a $p_{2}\circ \Theta$ est surjective. D'apr\`es la propri\'et\'e pr\'ec\'edente, on peut aussi bien remplacer $Temp(G)_{\xi}$ par $Temp(G)_{\xi}^0$, d'o\`u (1).

Une variante avec caract\`ere central de la proposition \ref{deuxespaces} nous dit que l'application $p_{1}$ est un isomorphisme.  Alors le (i) de l'\'enonc\'e se d\'eduit de (1). 

En utilisant \ref{leresultat} (1), l'injectivit\'e de  l'application du (ii)  se d\'eduit de celle de l'application  (2) de \ref{representationselliptiques}
D'apr\`es le (i) d\'ej\`a d\'emontr\'e et la relation (3) de \ref{representationselliptiques}, pour d\'emontrer la surjectivit\'e de l'application du (ii), il reste \`a prouver que,
 si $M$ est un Levi propre de $G$, si  $\pi^M\in {\mathbb C}[Ell(M)^0_{\xi}]$ et si l'on note $\pi=Ind_{M}^G(\pi^M)$, alors $\Delta^G_{\pi,cusp}=0$. La projection $I(G)^*\to I_{cusp}(G)^*$ annule toutes les distributions dont le support ne coupe pas $G(F)_{comp}$, car tout \'el\'ement elliptique r\'egulier est compact mod $Z(G)$. La projection de $\Theta_{\pi}$ est donc \'egale \`a celle de $D^G[\Delta^G_{\pi}]$. Cette projection $I(G)^*\to I_{cusp}(G)^*$ annule aussi  toutes les distributions induites \`a partir d'un Levi propre. Il en r\'esulte que la projection de $\Theta_{\pi}$ est nulle et que la projection de $D^G[\Delta^G_{\pi}]$ est \'egale \`a celle de $D^G[\Delta^G_{\pi,cusp}]$. Donc la projection dans $I_{cusp}(G)^*$ de $D^G[\Delta^G_{\pi,cusp}]$ est nulle. D'apr\`es le (ii) du corollaire  \ref{uncorollaire}, on a donc $\Delta^G_{\pi,cusp}=0$, ce qui ach\`eve la d\'emonstration. $\square$

\subsection{Produit scalaire elliptique\label{produitscalaireelliptique}}
Fixons un caract\`ere unitaire $\xi$ de $A_{G}(F)$. L'espace 
$${\mathbb C}[Temp(G)_{\xi}]=\oplus_{M\in\underline{{\cal L}}_{min}}Ind_{M}^G({\mathbb C}[Ell(M)_{\xi}]^{W^G(M)})$$
 est muni d'un produit scalaire elliptique not\'e $<.,>_{ell}$, cf. \cite{A2} paragraphe 6.   Les composantes ci-dessus  index\'ees par un Levi propre sont dans le noyau de ce produit. Par contre, le produit se restreint en un produit hermitien d\'efini positif sur ${\mathbb C}[Ell(G)_{\xi}]$, pour lequel $Ell(G)_{\xi}$ est une base orthonormale. Rappelons la d\'efinition du produit scalaire. Fixons un ensemble ${\cal T}_{ell}$ de repr\'esentants des classes de conjugaison de sous-tores elliptiques maximaux de $G$. Pour $T\in {\cal T}_{ell}$, on pose $W(T)=Norm_{G}(T)(F)/T(F)$ et on munit $A_{G}(F)\backslash T(F)$ de la mesure de Haar de masse totale $1$. Pour $\pi,\pi'\in {\mathbb C}[Temp(G)_{\xi}]$, on a alors 
$$(1) \qquad <\pi,\pi'>_{ell}=\sum_{T\in {\cal T}_{ell}}\vert W(T)\vert ^{-1}\int_{A_{G}(F)\backslash T(F)}\overline{\theta_{\pi}(t)}\theta_{\pi'}(t)D^G(t)\,dt.$$
Remarquons que ceci ne d\'epend d'aucune mesure. 

On va munir ${\cal D}_{cusp,\xi}(G)$ d'un produit scalaire que l'on appelle aussi elliptique. Pour ${\cal F},{\cal F}'\in Fac(G)$, 
    notons $N({\cal F},{\cal F}')$ l'ensemble des $g\in G(F)$ tels que $g{\cal F}={\cal F}'$ et fixons un ensemble de repr\'esentants $\underline{N}({\cal F},{\cal F}')$ du quotient
$$A_{G}(F)K_{{\cal F}'}^0\backslash N({\cal F},{\cal F}').$$
On a \'etabli en \ref{variantes} une d\'ecomposition
$$(2) \qquad  \boldsymbol{{\cal D}}_{cusp}(G)=\prod_{\nu\in {\cal N}}
  \boldsymbol{{\cal D}}_{cusp}(G)^{\nu}.$$
  On \'ecrit ${\bf f}=\prod_{\nu\in {\cal F}}{\bf f}^{\nu}$ la d\'ecomposition d'un \'el\'ement ${\bf f}\in \boldsymbol{{\cal D}}_{cusp}(G)$. 
  Fixons  $\nu\in {\cal N}$.  Soient $({\cal F},\nu),({\cal F}',\nu)\in Fac^*_{max}(G)$ et soient
  $f\in C_{cusp}({\bf G}_{{\cal F}}^{\nu})$ et $f'\in C_{cusp}({\bf G}_{{\cal F}'}^{\nu})$. On pose
$$<f,f'>_{ell}= mes(A_{G}(F)\backslash A_{G}(F)K_{{\cal F}}^0)^2\sum_{g\in \underline{N}({\cal F},{\cal F}')}<{^gf},f'>.$$
Ce produit se prolonge par lin\'earit\'e en un produit sur $\boldsymbol{{\cal D}}_{cusp}(G)^{\nu}$.
 Soient ${\bf f}, {\bf f}'\in \boldsymbol{{\cal D}}_{cusp,\xi}(G)$. On voit que $<{\bf f}^{\nu+\nu'},{\bf f}^{'\nu+\nu'}>=<{\bf f}^{\nu},{\bf f}^{'\nu}>$ pour tout $\nu\in {\cal N}$ et tout $\nu'\in w_{G}(A_{G}(F))$. On pose
$$<{\bf f},{\bf f}'>_{ell}=\sum_{\nu\in {\cal N}/w_{G}(A_{G}(F))}<{\bf f}^{\nu},{\bf f}^{'\nu}>_{ell}.$$
On v\'erifie que, pour $g\in G(F)$, $<{^g{\bf f}},{\bf f}'>_{ell}=<{\bf f},{\bf f}'>_{ell}$. C'est-\`a-dire que $^g{\bf f}-{\bf f}$ appartient au noyau du produit scalaire elliptique. Celui-ci se descend donc en un produit sur ${\cal D}_{cusp,\xi}(G)$. 

\ass{Proposition}{L'application
$$\begin{array}{ccc}{\mathbb C}[Ell(G)_{\xi}^0]&\to &{\cal D}_{cusp,\xi}(G)\\ \pi&\mapsto &\Delta^G_{\pi,cusp}\\ \end{array}$$
est un isomorphisme  isom\'etrique pour les produits scalaires elliptiques.
}

Preuve. On peut  supposer que l'ensemble de repr\'esentants $\underline{Fax}^*_{max}(G,A)$ est stable par l'op\'eration $({\cal F},\nu)\mapsto ({\cal F},\nu+\nu')$ pour tout $\nu'\in w_{G}(A_{G}(F))$. Repr\'esentons $\Delta^G_{\pi,cusp}$ comme un \'el\'ement de $\prod_{({\cal F},\nu)\in \underline{Fac}^*_{max}(G,A)}C_{cusp}({\bf G}_{{\cal F}}^{\nu})^{K_{{\cal F}}^{\dag}}$, cf. \ref{deuxespaces}. On note $(f_{\pi,{\cal F},\nu})_{({\cal F},\nu)\in \underline{Fac}^*_{max}(G,A)}$ cet \'el\'ement. 
  Soit $T\in {\cal T}_{ell}$ et $t\in T(F)\cap G_{reg}(F)$. D'apr\`es \cite{W2} 18(3), on a l'\'egalit\'e
$$\theta_{\pi}(t)=D^G(t)^{-1/2}\sum_{({\cal F},\nu)\in \underline{Fac}^*_{max}(G,A)} I^G(t,f_{\pi,{\cal F},\nu}).$$
 On a identifi\'e ici l'\'el\'ement $f_{\pi,{\cal F},\nu}$ \`a la fonction $(f_{\pi,{\cal F},\nu})_{{\cal F}}$. La somme  est finie car les termes ne sont non nuls que si $\nu=w_{G}(t)$. 

{\bf Remarque.} La formule de la r\'ef\'erence contient un terme $mes(A_{G}(F)\backslash K_{{\cal F}}^{\dag})^{-1}$ que l'on a incorpor\'e \`a la d\'efinition de $f_{\pi,{\cal F},\nu}$.
\bigskip

   Fixons un ensemble de repr\'esentants $\underline{{\cal N}}$ des orbites dans ${\cal N}$ pour l'action de $w_{G}(A_{G}(F))$ par addition. Notons $ \underline{Fac}^*_{max}(G,A;\underline{\cal N})$ le sous-ensemble des \'el\'ements $({\cal F},\nu)\in \underline{Fac}^*_{max}(G,A)$ tels que $\nu\in \underline{{\cal N}}$.
La formule ci-dessus se r\'ecrit
$$\theta_{\pi}(t)=D^G(t)^{-1/2}\sum_{({\cal F},\nu)\in \underline{Fac}^*_{max}(G,A;\underline{{\cal N}})} \sum_{\nu'\in w_{G}(A_{G}(F))}I^G(t,f_{\pi,{\cal F},\nu+\nu'}).$$
Parce que la restriction \`a $A_{G}(F)$ du caract\`ere central de $\pi$ est $\xi$, on voit que la somme int\'erieure n'est autre que
$$mes(A_{G}(F)_{c})^{-1}\int_{A_{G}(F)}I^G(ta,f_{\pi,{\cal F},\nu})\xi(a)^{-1}\,da.$$
D'o\`u
$$\theta_{\pi}(t)=mes(A_{G}(F)_{c})^{-1}\int_{A_{G}(F)} D^G(ta)^{-1/2}\sum_{({\cal F},\nu)\in \underline{Fac}^*_{max}(G,A;\underline{{\cal N}})}I^G(ta,f_{\pi,{\cal F},\nu})\xi(a)^{-1}\,da.$$
La formule (1) devient
$$<\pi,\pi'>_{ell}=mes(A_{G}(F)_{c})^{-1}\sum_{({\cal F},\nu)\in \underline{Fac}^*_{max}(G,A;\underline{{\cal N}})} $$
$$\sum_{T\in {\cal T}_{ell}}\vert W(T)\vert ^{-1}\int_{T(F)}D^G(t)^{1/2}I^G(t,\overline{f_{\pi,{\cal F},\nu}})\theta_{\pi'}(t)\,dt.$$
Les fonctions $f_{\pi,{\cal F},\nu}$ sont cuspidales, leurs int\'egrales orbitales sont nulles hors des \'el\'ements elliptiques. La formule de Weyl conduit donc \`a l'\'egalit\'e
$$<\pi,\pi'>_{ell}=mes(A_{G}(F)_{c})^{-1}\sum_{({\cal F},\nu)\in \underline{Fac}^*_{max}(G,A;\underline{{\cal N}})} \Theta_{\pi'}(\overline{f_{\pi,{\cal F},\nu}}).$$
Fixons $({\cal F},\nu)$ intervenant ci-dessus. Puisque la fonction $f_{\pi,{\cal F},\nu}$ est \`a support compact mod $Z(G)$, on peut exprimer $\Theta_{\pi'}(\overline{f_{\pi,{\cal F},\nu}})$ \`a l'aide de  $\Delta_{\pi'}$. Puisque la fonction $f_{\pi,{\cal F},\nu}$  est cuspidale, les distributions induites ne comptent pas, d'o\`u
$$\Theta_{\pi'}(\overline{f_{\pi,{\cal F},\nu}})=D^G[\Delta^G_{\pi',cusp}](\overline{f_{\pi,{\cal F},\nu}}).$$
Ceci est calcul\'e par la proposition \ref{calcul} et la remarque (2) qui la suit.  On obtient
$$\Theta_{\pi'}(\overline{f_{\pi,{\cal F},\nu}})=c_{{\cal F}}<f_{\pi,{\cal F},\nu},f_{\pi',{\cal F},\nu}>,$$
o\`u 
$$c_{{\cal F}}=mes(K_{{\cal F}}^0)^2mes(A_{G}(F)_{c})^{-1}[K_{{\cal F}}^{\dag}:A_{G}(F)K_{{\cal F}}^0].$$
De la m\^eme fa\c{c}on, on a
$$<f_{\pi,{\cal F},\nu},f_{\pi',{\cal F},\nu}>_{ell}=d_{{\cal F}}<f_{\pi,{\cal F},\nu},f_{\pi',{\cal F},\nu}>,$$
o\`u 
$$d_{{\cal F}}=mes(A_{G}(F)\backslash A_{G}(F)K_{{\cal F}}^0)^2[K_{{\cal F}}^{\dag}:A_{G}(F)K_{{\cal F}}^0].$$
  D'o\`u
$$\Theta_{\pi'}(\overline{f_{\pi,{\cal F},\nu}})= c_{{\cal F}}d_{{\cal F}}^{-1}<f_{\pi,{\cal F},\nu},f_{\pi',{\cal F},\nu}>_{ell},$$
puis
$$<\pi,\pi'>_{ell}=mes(A_{G}(F)_{c})^{-1}\sum_{({\cal F},\nu)\in \underline{Fac}^*_{max}(G,A;\underline{{\cal N}})} c_{{\cal F}}d_{{\cal F}}^{-1}<f_{\pi,{\cal F},\nu},f_{\pi',{\cal F},\nu}>_{ell}.$$
On v\'erifie que
$$mes(A_{G}(F)_{c})^{-1} c_{{\cal F}}d_{{\cal F}}=1.$$
D'o\`u
$$<\pi,\pi'>_{ell}=\sum_{({\cal F},\nu)\in \underline{Fac}^*_{max}(G,A;\underline{{\cal N}})}<f_{\pi,{\cal F},\nu},f_{\pi',{\cal F},\nu}>_{ell}.$$
Par d\'efinition, la membre de droite n'est autre que $<\Delta^G_{\pi,cusp},\Delta^G_{\pi',cusp}>$. Cela d\'emontre que l'application de l'\'enonc\'e est isom\'etrique. Elle est  bijective d'apr\`es le (ii) du lemme \ref{surjectivite}. Cela ach\`eve la d\'emonstration. $\square$

\subsection{Un lemme pr\'eliminaire\label{unlemme}}

Rappelons que l'on a d\'efini en \ref{restriction} la notion de quasi-caract\`ere sur $G(F)$ de niveau $0$ sur les elliptiques. 

\ass{Lemme}{Soit $\pi\in {\mathbb C}[Ell(G)_{{\mathbb R}}]$. Supposons que $\Theta_{\pi}$ soit un quasi-caract\`ere de niveau $0$ sur les elliptiques. Alors $\pi$ appartient \`a ${\mathbb C}[Ell(G)^0_{{\mathbb R}}]$.}

Preuve. On peut d\'ecomposer la repr\'esentation $\pi$ en $\pi=\sum_{i=1,...,h}\pi_{i}$ o\`u   chaque  $\pi_{i}$ est combinaison lin\'eaire de repr\'esentations irr\'eductibles dont le caract\`ere central se restreint \`a   $A_{G}(F)$ en un m\^eme caract\`ere $\xi_{i}$. On suppose les $\xi_{i}$ tous distincts. Par interpolation, on peut trouver une famille de nombres complexes $(c_{j})_{j=1,...,h}$ et une famille $(a_{j})_{j=1,...,h}$ d'\'el\'ements de $A_{G}(F)$ de sorte que
 $$\Theta_{\pi_{i}}(g)=\sum_{j=1,...,h}c_{j}\Theta_{\pi}(a_{j}g)$$
 pour tout $g\in G(F)$. Puisque $\Theta_{\pi}$ est un quasi-caract\`ere de niveau $0$ sur les elliptiques, cette formule entra\^{\i}ne que $\Theta_{\pi_{i}}$ l'est aussi. Cela nous ram\`ene au cas o\`u $\pi$ elle-m\^eme est combinaison lin\'eaire de repr\'esentations irr\'eductibles dont le caract\`ere central se restreint \`a   $A_{G}(F)$ en un m\^eme caract\`ere que l'on note $\xi$. On  peut aussi tordre le probl\`eme par un \'el\'ement $\chi\in {\cal A}_{G}^*$: $\Theta_{\pi_{\chi}}$ est encore un quasi-caract\`ere de niveau $0$ sur les elliptiques et, si l'on prouve que $\pi_{\chi}$ appartient \`a 
 ${\mathbb C}[Ell(G)^0_{{\mathbb R}}]$, la m\^eme assertion s'ensuit pour $\pi$. On peut donc supposer $\xi$ unitaire. D'apr\`es la proposition \ref{restriction}, il existe $d'\in {\cal D}_{cusp}(G)$ tel que $\Theta_{\pi}$ co\"{\i}ncide avec $D^G[d']$ sur $G_{ell}(F)$. Rappelons que $A_{G}(F)$ agit naturellement sur ${\cal D}_{cusp}(G)$. On note $(a,\delta)\mapsto \delta^{a}$ cette action. 
  L'\'el\'ement $d'$ n'a pas forc\'ement de caract\`ere central pour l'action de $A_{G}(F)$ mais on peut le remplacer par un \'el\'ement qui admet $\xi$ pour tel caract\`ere central. En effet, on peut d'abord remplacer $d'$ par l'int\'egrale $\int_{A_{G}(F)_{c}}\xi(a)^{-1}d^{'a} \,da$. Cela ne modifie pas la propri\'et\'e ci-dessus puisque $\pi$ se transforme par $A_{G}(F)$ selon $\xi$.   Mais cela assure que $A_{G}(F)_{c}$ agit sur $d'$  selon le caract\`ere $\xi$. 
  Utilisons   la d\'ecomposition  (2) de \ref{produitscalaireelliptique} et \'ecrivons $d'=\prod_{\nu\in {\cal N}}d^{'\nu}$. Pour tout $\nu\in {\cal N}$, $D^G(d^{'\nu})$ est \`a support dans $w_{G}^{-1}(\nu)$ et co\"{\i}ncide avec $\Theta_{\pi}$ sur $G_{ell}(F)\cap w_{G}^{-1}(\nu)$.  
   Fixons un ensemble de repr\'esentants $\underline{{\cal N}}$ des orbites de l'action par addition de $w_{G}(A_{G}(F))$ dans ${\cal N}$. On voit qu'il existe un unique $ d\in {\cal D}_{cusp,\xi}(G)$ tel que $ d^{\nu}=d^{'\nu}$ pour tout $\nu\in \underline{{\cal N}}$. Puisque $\pi$ se transforme par $A_{G}(F)$ selon $\xi$, $\Theta_{\pi}$ co\"{\i}ncide avec $D^G(d)$ sur $G_{ell}(F)$. D'apr\`es le (ii) du lemme \ref{surjectivite}, il existe $\pi^0\in {\mathbb C}[Ell(G)_{\xi}^0]$ tel que $\Delta_{cusp}(\pi^0)=d$. Alors, d'apr\`es l'assertion (1) de \ref{leresultat}, $\Theta_{\pi^0}$ co\"{\i}ncide sur $G_{ell}(F)$ avec $D^G(d)$, donc aussi avec $\Theta_{\pi}$. D'apr\`es l'injectivit\'e de l'application (2) de \ref{representationselliptiques}, on a $\pi=\pi^0$, ce qui ach\`eve la d\'emonstration. $\square$

 \subsection{Caract\'erisation des repr\'esentations de niveau $0$\label{caracterisation}}
 \ass{Th\'eor\`eme}{Soit $\pi\in {\mathbb C}[Irr(G)]$. Alors $\pi\in {\mathbb C}[Irr(G)^0]$ si et seulement si $\Theta_{\pi}$ est un quasi-caract\`ere de niveau $0$ sur $G(F)$.}
 
 Preuve. Supposons $\pi\in {\mathbb C}[Irr(G)^0]$. D'apr\`es \ref{leresultat} (1) et le th\'eor\`eme \ref{unpremiertheoreme}, $\Theta_{\pi}$ co\"{\i}ncide avec un quasi-caract\`ere de niveau $0$ sur les \'el\'ements de $G(F)$ qui sont compacts mod $Z(G)$. Soit $\epsilon\in G(F)_{p'}$, supposons que $\epsilon$ n'est pas compact mod $Z(G)$. Alors $L[\epsilon]$ est un Levi propre et on a $L[\epsilon exp(X)]=L[\epsilon]$ et $Q[\epsilon exp(X)]=Q[\epsilon]$ pour tout $X\in \mathfrak{g}_{\epsilon,tn}(F)$. D'apr\`es le   th\'eor\`eme 5.2 de \cite{C}, on a $\theta_{\pi}(\epsilon exp(X))=\delta_{Q[\epsilon]}(\epsilon exp(X))^{1/2}\theta_{\pi_{Q[\epsilon]}}(\epsilon exp(X))$ pour tout $X$ comme ci-dessus. Un m\^eme calcul qu'au (ii) du lemme \ref{decomposition} montre que $\delta_{Q[\epsilon]}(\epsilon exp(X))=\delta_{Q[\epsilon]}(\epsilon)$. La repr\'esentation $\pi_{Q[\epsilon]}$ est de niveau $0$. Puisque $L[\epsilon]\subsetneq G$, on peut raisonner par r\'ecurrence sur le rang semi-simple de $G$ et suppos\'e prouv\'e que $\Theta_{\pi_{Q[\epsilon]}}$ est un quasi-caract\`ere de niveau $0$ sur $L[\epsilon](F)$. Autrement dit que, pour des constantes $c_{{\cal O}}$ convenables, on a 
 $$\theta_{\pi_{Q[\epsilon]}}(\epsilon exp(X))=\sum_{{\cal O}\in Nil(\mathfrak{l}[\epsilon]_{\epsilon})}c_{{\cal O}}\hat{j}({\cal O},X)$$
 pour tout $X$ comme ci-dessus. Puisque $\mathfrak{l}[\epsilon]_{\epsilon}=\mathfrak{g}_{\epsilon}$, on en d\'eduit un d\'eveloppement analogue de $\theta_{\pi}(\epsilon exp(X))$. Donc $\Theta_{\pi}$ co\"{\i}ncide avec un quasi-caract\`ere de niveau $0$ sur $C_{G}(\epsilon)$.  Ceci est donc vrai pour tout $\epsilon\in G(F)_{p'}$ et cela signifie que $\Theta_{\pi}$ est un quasi-caract\`ere de niveau $0$.
 
 Pour $n\in {\mathbb Z}$, notons ${\mathbb C}[Irr(G)]^n$ la sous-somme de l'expression \ref{representationselliptiques} (1) o\`u l'on somme sur les Levi $M\in \underline{{\cal L}}^m_{min}$ pour $m\geq n$. On va prouver
 
 (1) soit $\pi\in {\mathbb C}[Irr(G)]^n$; supposons que $\Theta_{\pi}$ est un quasi-caract\`ere de niveau $0$ sur $G(F)$; alors $\pi$ appartient \`a la somme de ${\mathbb C}[Irr(G)]^{n+1}$ et de
 $$\oplus_{M\in \underline{{\cal L}}^n_{min}}Ind_{M}^G({\mathbb C}[Ell(M)_{{\mathbb R}}^0]^{W^G(M)}.$$
 
 Ecrivons $\pi=\sum_{M\in \underline{{\cal L}}_{min}}Ind_{M}^G(\pi^M)$, o\`u $\pi^M\in {\mathbb C}[Ell(M)_{{\mathbb R}}]^{W^G(M)}$ et $\pi^M=0$ si $a_{M}<n$. Le caract\`ere $\Theta_{\pi}$ appartient \`a $Qc_{0}(G)$ par hypoth\`ese donc \`a $Qc_{0}(G)\cap Ann^nI(G)^*$. D'apr\`es le (ii) de la proposition  \ref{injectivite} et le th\'eor\`eme \ref{undeuxiemetheoreme}, cet espace a m\^eme image dans $Gr^nI(G)^*$ que 
 $$\oplus_{M\in \underline{{\cal L}}^n_{min}}Ind_{M}^G(D^M({\cal D}_{cusp}(M)^{W^G(M)})).$$
 Autrement dit, pour tout $M\in \underline{{\cal L}}^n_{min}$, on peut fixer $d^M\in {\cal D}_{cusp}(M)^{W^G(M)}$ tel que $\Theta_{\pi^M}$ co\"{\i}ncide avec $D^M(d^M)$ sur les \'el\'ements elliptiques de $M(F)$. A fortiori, $\pi^M$ est un quasi-caract\`ere sur $M(F)$ de niveau $0$ sur les elliptiques. D'apr\`es le lemme \ref{unlemme}, $\pi^M$ appartient \`a ${\mathbb C}[Ell(M)_{{\mathbb R}}^0]$ et donc aussi \`a ${\mathbb C}[Ell(M)_{{\mathbb R}}^0]^{W^G(M)}$. Posons 
 $\sigma=\sum_{M\in \underline{{\cal L}}^n_{min}}Ind_{M}^G(\pi^M)$. Alors $\Theta_{\pi}$ et $\Theta_{\sigma}$ ont m\^eme image dans $Gr^nI(G)^*$. Donc $\pi-\sigma$ appartient \`a ${\mathbb C}[Irr(G)]^{n+1}$. Puisque $\sigma$ appartient au dernier espace indiqu\'e dans l'assertion (1), cela d\'emontre cette assertion.
  
  D'apr\`es le sens "seulement si" d\'ej\`a d\'emontr\'e du th\'eor\`eme, l'assertion (1) implique le sens "si" par r\'ecurrence descendante sur $n$. $\square$
  
  \section{Endoscopie}
  
  \subsection{Donn\'ees endoscopiques\label{donneesendoscopiques}}
  On note $W_{F}$ le groupe de Weil de $F$. Il contient le groupe d'inertie $I_{F}$. On note  $I_{F}^s$ le groupe d'inertie sauvage de $F$, c'est-\`a-dire le $p$-Sylow du groupe  $I_{F}$. 
  
 La th\'eorie de l'endoscopie a \'et\'e d\'evelopp\'ee par Langlands et ses successeurs. Nous ne la reprendrons pas et nous adopterons sa formulation telle qu'elle figure dans \cite{MW}.   Une donn\'ee endoscopique de $G$ est un triplet ${\bf G}'=(G',{\cal G}',s)$, o\`u $G'$ est un groupe r\'eductif connexe d\'efini et quasi-d\'eploy\'e sur $F$, ${\cal G}'$ est un sous-groupe du $L$-groupe $^LG=\hat{G}\rtimes W_{F}$ et $s$ est un \'el\'ement semi-simple du groupe complexe $\hat{G}$. Le groupe $\hat{G}'$ s'identifie \`a la composante neutre $Z_{\hat{G}}(s)^0$ du commutant $Z_{\hat{G}}(s)$.   On dit que la donn\'ee est elliptique si et seulement si $Z(\hat{G})^{\Gamma}$ est un sous-groupe d'indice fini de $Z(\hat{G}')^{\Gamma}$.
 
 L'hypoth\`ese $(Hyp)(G)$ a les cons\'equences suivantes: 
 
 (1) il existe une extension de $F$ de degr\'e premier \`a $p$ telle que le groupe $G'$ soit  d\'eploy\'e sur cette extension; 
  
  (2)   $W_{F}\cap {\cal G}'$  est un sous-groupe de $W_{F}$ d'indice fini premier \`a $p$.
 
 Preuve. Le (1) r\'esulte de ce que le rang de $G'$ est le m\^eme que celui de $G$ et que l'on a d\'ej\`a remarqu\'e que, pour tout tore sur $F$  de dimension inf\'erieure ou \'egale \`a ce rang, il existe   une extension $F'$ de $F$  de degr\'e premier \`a $p$ telle que ce tore soit d\'eploy\'e sur $F'$. Pour (2), fixons une telle extension galoisienne $F'$ de $F$ de degr\'e premier \`a $p$, telle que $G'$ soit d\'eploy\'e sur $F$. L'action  de $\Gamma_{F'}$ sur $\hat{G}'$ est triviale, a fortiori l'action de $\Gamma_{F'}$ sur $\hat{G}$ est triviale.  Le noyau de l'homomorphisme naturel $W_{F}\to Gal(F'/F)$ est $W_{F'}$.
Par d\'efinition, la projection ${\cal G}'\to W_{F}$ est scind\'ee, c'est-\`a-dire que l'on peut fixer un homomorphisme continu $i:W_{F}\to {\cal G}'$ qui est une section de cette projection.  Pour $w\in W_{F}$, \'ecrivons $i(w)=(j(w),w)$. La restriction de $j$ \`a $W_{F'}$ est un caract\`ere de ce groupe \`a valeurs dans $Z_{\hat{G}}(s)$. Notons $W^1$ le sous-groupe des $w\in W_{F'}$ tels que $j(w)\in Z_{\hat{G}}(s)^0$. L'hypoth\`ese $(Hyp)(G)$ implique que $Z_{\hat{G}}(s)^0$ est un sous-groupe de $Z_{\hat{G}}(s)$ d'indice premier \`a $p$. Donc $W^1$ est d'indice fini premier \`a $p$ dans $W_{F'}$, donc aussi dans $W_{F}$. Pour $w\in W^1$, on a $(1,w)=j(w)^{-1}i(w)$ et ce terme appartient \`a $ {\cal G}'$ puisque $j(w)\in Z_{\hat{G}}(s)^0=\hat{G}'\subset {\cal G}'$. Donc $W^1\subset W_{F}\cap {\cal G}'$ et (2) est d\'emontr\'e. $\square$

 On a besoin de munir toute donn\'ee endoscopique  ${\bf G}'$ de donn\'ees auxiliaires $G'_{1},C_{1},\hat{\xi}_{1}$, cf. \cite{MW} I.2.1.  Le tore $C_{1}(F)$ est naturellement muni d'un caract\`ere $\lambda_{1}:C_{1}(F)\to {\mathbb C}^{\times}$.  Il y a une notion d'unitarit\'e pour de telles donn\'ees, cf. loc. cit. I.7.1, qui implique que $\lambda_{1}$ est unitaire.
 Nous dirons que les donn\'ees auxiliaires sont mod\'er\'ement ramifi\'ees si et seulement si elles v\'erifient les conditions suivantes:
 
 (3) $G'_{1}$ et $C_{1}$ sont d\'eploy\'es sur une extension de $F$ de degr\'e premier \`a $p$; 
 
 (4) on a $\hat{\xi}_{1}(1,w)=(1,w)$ pour tout $w\in I_{F}^s$.
 
 \ass{Lemme}{Sous notre hypoth\`ese $(Hyp)(G)$, il existe des donn\'ees auxiliaires mod\'er\'ement ramifi\'ees et unitaires.}
 
 Preuve. Choisir une suite
 $$1\to C_{1}\to G'_{1}\to G'\to 1$$
 \'equivaut \`a choisir une suite duale
 $$1\to \hat{G'}\to \hat{G}'_{1}\to \hat{C}_{1}\to 1$$
 Pour celle-ci, il y a le choix standard suivant. Fixons une paire de Borel $(\hat{B}',\hat{T})$ de $\hat{G}'$ pr\'eserv\'ee par l'action galoisienne. On pose $\hat{G}'_{1}=(\hat{G}'\times \hat{T})/diag(Z(\hat{G}'))$, o\`u $diag $ est le plongement diagonal. Le centre de $\hat{G}'_{1}$ est isomorphe \`a $\hat{T}$, donc est connexe. Le tore $\hat{C}_{1}$ est $\hat{T}/Z(\hat{G})$ et on sait bien que c'est un tore induit. Si $F'$ est une extension finie de $F$ de degr\'e \`a $p$ sur laquelle $G'$ est d\'eploy\'e, les actions galoisiennes sur $\hat{G}'_{1}$ et $\hat{C}_{1}$ sont triviales sur $\Gamma_{F'}$. Dualement, $C_{1}$ est d\'eploy\'e sur $F'$ et, puisque $G'$ l'est aussi, $G'_{1}$ l'est \'egalement.

 Pour construire $\hat{\xi}_{1}$ v\'erifiant (4), le mieux est de reprendre la preuve de Langlands \cite{L} lemme 4. On s'aper\c{c}oit qu'il y a un unique argument \`a pr\'eciser. La cl\'e de la preuve est que, quand $S_{2}\to S_{1}$ est une injection de tores d\'efinis sur $F$, tout caract\`ere  $\chi_{2}:S_{2}(F)\to {\mathbb C}^{\times}$    se prolonge en un caract\`ere $\chi_{1}:S_{1}(F)\to {\mathbb C}^{\times}$ (rappelons que, pour nous, caract\`ere signifie homomorphisme continu). On doit pr\'eciser que tout caract\`ere  mod\'er\'ement ramifi\'e $\chi_{2}$ se prolonge en un caract\`ere mod\'er\'ement ramifi\'e $\chi_{1}$. Le sous-groupe $S_{1,tu}(F)$ de $S_{1}(F)$ est compact, donc le produit $S_{2}(F)S_{1,tu}(F)$ est ferm\'e dans $S_{1}(F)$. Puisque $\chi_{2}$ est mod\'er\'ement ramifi\'e, il est trivial sur $S_{2,tu}(F)=S_{2}(F)\cap S_{1,tu}(F)$. Donc on peut prolonger $\chi_{2}$ en un caract\`ere $\chi':S_{2}(F)S_{1,tu}(F)\to {\mathbb C}^{\times}$ qui est trivial sur $S_{1,tu}(F)$. Puisque $S_{2}(F)S_{1,tu}(F)$ est ferm\'e dans $S_{1}(F)$, on peut ensuite prolonger $\chi'$ en un caract\`ere $\chi_{1}:S_{1}(F)\to {\mathbb C}^{\times}$. Il est encore trivial sur $S_{1,tu}(F)$, c'est-\`a-dire qu'il est mod\'er\'ement ramifi\'e. 
 
 Avec ce compl\'ement,la preuve de \cite{L} permet de construire $\hat{\xi}_{1}$ v\'erifiant (4). Le m\^eme  argument qu'en \cite{MW} I.7.1(3) permet de le modifier afin d'assurer que les donn\'ees sont unitaires. $\square$
 
   On fixe un ensemble de repr\'esentants ${\cal E}(G)$ des classes d'\'equivalence de donn\'ees endoscopiques elliptiques de $G$. Pour tout ${\bf G}'=(G',{\cal G}',s)\in {\cal E}(G)$, on fixe des donn\'ees auxiliaires $G'_{1},C_{1},\hat{\xi}_{1}$.
 On suppose que  ces donn\'ees auxiliaires sont mod\'er\'ement ramifi\'ees et unitaires. On peut alors  fixer un facteur de transfert unitaire $\Delta_{1}$. Signalons que, pour nous, ce facteur ne contient pas le terme $\Delta_{IV}$ de \cite{LS}, que nous avons incorpor\'e \`a la d\'efinition des int\'egrales orbitales en suivant Arthur. Il n'est pas clair que l'hypoth\`ese $(Hyp)(G)$ que nous avons pos\'ee  entra\^{\i}ne $(Hyp)(G')$ pour    tout groupe endoscopique $G'$ de $G$, et encore moins qu'elle entra\^{\i}ne $(Hyp)(G'_{1})$ pour un groupe auxiliaire $G'_{1}$. Nous renfor\c{c}ons l'hypoth\`ese $(Hyp)(G)$  en l'hypoth\`ese suivante:
 
 $(Hyp)_{endo}(G)$: l'hypoth\`ese $(Hyp)(G)$ est v\'erifi\'ee; pour tout ${\bf G}'\in {\cal E}(G)$, les hypoth\`eses $(Hyp)(G')$ et $(Hyp)(G'_{1})$  sont v\'erifi\'ees.

   Toute donn\'ee endoscopique ${\bf G}'$ de $G$ appara\^{\i}t comme une "donn\'ee de Levi" d'une donn\'ee endoscopique elliptique ${\bf G}''$. Les donn\'ees auxiliaires que l'on a fix\'ees pour cette derni\`ere se restreignent en des donn\'ees auxiliaires pour ${\bf G}'$. Les hypoth\`eses $(Hyp)(G'')$ et $(Hyp)(G''_{1})$ entra\^{\i}nent $(Hyp)(G')$ et $(Hyp)(G'_{1})$ puisque $G'$, resp. $G'_{1}$, est un groupe de Levi de $G''$, resp. $G''_{1}$. En cons\'equence, l'hypoth\`ese $(Hyp)_{endo}(G)$ implique que, pour toute donn\'ee endoscopique (pas forc\'ement elliptique), on peut supposer fix\'ees des donn\'ees auxiliaires  $G'_{1},C_{1},\hat{\xi}_{1}$ mod\'er\'ement ramifi\'ees et unitaires, de sorte que $(Hyp)(G')$ et $(Hyp)(G'_{1})$ soient v\'erifi\'ees. On suppose aussi fix\'e un facteur de transfert unitaire  $\Delta_{1}$.

 \subsection{Repr\'esentations de niveau $0$ et donn\'ees auxiliaires\label{representationsetdonnees}}
 
  Soit ${\bf G}'=(G',{\cal G}',s)$ une donn\'ee endoscopique de $G$.  On s'int\'eresse exclusivement aux repr\'esentations irr\'eductibles de $G'_{1}(F)$ dont le caract\`ere central co\"{\i}ncide avec $\lambda_{1}$ sur $C_{1}(F)$. On note $Irr (G'_{1})_{\lambda_{1}}$ l'ensemble de ces repr\'esentations, avec les variantes  $Temp(G'_{1})_{\lambda_{1}}$,  $Irr(G'_{1})^0_{\lambda_{1}}$ etc...  Ces ensembles d\'ependent du choix des donn\'ees auxiliaires. Pour que les r\'esultats qui suivent aient vraiment un sens, il vaut mieux qu'ils n'en d\'ependent pas trop. Selon le formalisme que l'on a d\'evelopp\'e en \cite{MW}, pour deux choix de donn\'ees auxiliaires, il y a une bijection canonique entre les ensembles  associ\'es \`a  chacune des donn\'ees. Expliquons cela. 
 Choisissons d'autres donn\'ees auxiliaires mod\'er\'ement ramifi\'ees $(G'_{2},C_{2},\hat{\xi}_{2})$. Notons $G'_{12}$ le produit fibr\'e de $G'_{1}$ et $G'_{2}$ au-dessus de $G'$.   On a deux suites exactes
  $$1\to C_{2}\to G'_{12}\to G'_{1}\to 1,\,\,1\to C_{1}\to G'_{12}\to G'_{2}\to 1.$$
  On a introduit en \cite{MW} I.2.5 un homomorphisme $\lambda_{12}:G'_{12}(F)\to {\mathbb C}^{\times}$. Sa restriction \`a $C_{1}(F)\times C_{2}(F)$ est le produit des caract\`eres $\lambda_{1}$ et $\lambda_{2}^{-1}$.  Soit $\pi_{1}\in Irr(G'_{1})_{\lambda_{1}}$, r\'ealis\'ee dans un espace complexe $V$. Par la  premi\`ere suite ci-dessus, $\pi_{1}$ se rel\`eve en une repr\'esentation $\pi_{12}$ de $G'_{12}(F)$. La restriction \`a $C_{1}(F)$, resp. $C_{2}(F)$, de son caract\`ere central  est $\lambda_{1}$, resp. $1$. Posons $\pi_{21}=\lambda_{12}^{-1}\otimes \pi_{12}$. Alors la restriction \`a $C_{1}(F)$, resp. $C_{2}(F)$, du caract\`ere central de $\pi_{21}$  est $1$, resp. $\lambda_{2}$. La repr\'esentation $\pi_{21}$ se descend par la deuxi\`eme suite ci-dessus en une repr\'esentation  $\pi_{2}$ de $G'_{2}(F)$ dans $V$ dont le caract\`ere central co\"{\i}ncide avec $\lambda_{2}$ sur $C_{2}(F)$. L'application $\pi_{1}\mapsto \pi_{2}$ est une bijection de $Irr (G'_{1})_{\lambda_{1}}$ sur $Irr (G'_{2})_{\lambda_{2}}$. Parce que   $\lambda_{1}$ et $\lambda_{2}$ sont unitaires, le caract\`ere $\lambda_{12}$ l'est aussi et la bijection envoie   $Temp(G'_{1})_{\lambda_{1}}$ sur $Temp(G'_{2})_{\lambda_{2}}$.  Parce que les donn\'ees auxiliaires sont mod\'er\'ement ramifi\'ees,  le caract\`ere $\lambda_{12}$ est "mod\'er\'ement ramifi\'e" au sens que sa restriction \`a l'ensemble $G'_{12,tu}(F)$ est triviale. Remarquons de plus que $Imm(G'_{1,AD})=Imm(G'_{AD})=Imm(G'_{2,AD})$. Une facette ${\cal F}\in Fac(G')$ d\'etermine un sous-groupe $K_{{\cal F}}^+$ de $G'(F)$ et des sous-groupes similaires dans $G'_{1}(F)$ et $G'_{2}(F)$, que l'on note  $K_{1,{\cal F}}^+$ et $K_{2,{\cal F}}^+$. On voit  que si $\pi_{1}$ admet des invariants non nuls par le sous-groupe  $K_{1,{\cal F}}^{+}$, alors $\pi_{2}$ admet   des invariants non nuls par le sous-groupe  $K_{2,{\cal F}}^{+}$. La bijection envoie donc $Irr(G'_{1})^0_{\lambda_{1}}$ sur $Irr(G'_{2})^0_{\lambda_{2}}$.
  
  \subsection{Correspondance entre \'el\'ements semi-simples\label{correspondance}}
  Soit ${\bf G}'=(G',{\cal G}',s)$ une donn\'ee endoscopique de $G$. On d\'efinit comme en \cite{MW} I.1.10 une correspondance entre classes de conjugaison stable d'\'el\'ements semi-simples dans $G'(F)$ et $G(F)$ (on parlera aussi bien d'une correspondance entre \'el\'ements semi-simples dans $G'(F)$ et $G(F)$). Soient $x'\in G'_{reg}(F)$ et $x\in G_{reg}(F)$ deux \'el\'ements qui se correspondent. Notons $T'$ et $T$ les commutants de $x'$ dans $G'$ et de $x$ dans $G$. Ce sont des tores et il y a un unique  isomorphisme $\xi_{T,T'}:T\to T'$ d\'efini sur $F$ de sorte que $\xi_{T,T'}(x)=x'$, cf. loc. cit..  D'o\`u un isomorphisme encore not\'e $\xi_{T,T'}:\mathfrak{t}\to \mathfrak{t'}$.
  
   \ass{Lemme}{(i) Si $x=\epsilon exp(X)$ est une $p'$-d\'ecomposition de $x$, alors $x'=\xi_{T,T'}(\epsilon)exp(\xi_{T,T'}(X))$ est une $p'$-d\'ecomposition de $x'$. 
   
   (ii) Soit $x=\epsilon exp(X)$  une $p'$-d\'ecomposition de $x$. Pour tout $\lambda\in F^{\times}$ tel que $ \lambda X\in \mathfrak{g}_{\epsilon,tn}(F)$, $\xi_{T,T'}(\epsilon)exp(\xi_{T,T'}(\lambda X))$ est un \'el\'ement de $G'_{reg}(F)$ qui correspond \`a  $\epsilon exp(\lambda X)$.
   
   (iii) Supposons ${\bf G}'$ elliptique. Alors $x$ est compact mod $Z(G)$ si et seulement si $x'$ est compact mod $Z(G')$. Si ces conditions sont r\'ealis\'ees, l'application $(\epsilon,X)\mapsto (\epsilon',X')=(\xi_{T,T'}(\epsilon),\xi_{T,T'}(X))$ r\'ealise une bijection entre les $p'$-d\'ecompositions $x=\epsilon exp(X)$ de $x$ et les $p'$-d\'ecompositions $x'=\epsilon'exp(X')$ de $x'$.}
  
   Preuve. Puisque $\xi_{T,T'}$ est un isomorphisme, l'\'el\'ement $\xi_{T,T'}(X)$ est topologiquement nilpotent. Posons $\epsilon'=\xi_{T,T'}(\epsilon)$, introduisons le Levi $L[\epsilon]$ de $G$ et le Levi similaire $L'[\epsilon']$ de $G'$. Notons $\Sigma$ l'ensemble de racines de $T$ dans $G$, $\Sigma'$ celui des racines de $T'$ dans $G'$ et $\Sigma^{L[\epsilon]}$ et $\Sigma^{L'[\epsilon']}$ les sous-ensembles des racines dans $L[\epsilon]$, resp. $L'[\epsilon']$. L'isomorphisme $\xi_{T,T'}$ transporte $\Sigma$ en un sous-ensemble de $X^*(T')$ et on sait que $\Sigma'\subset \xi_{T,T'}(\Sigma)$. Il r\'esulte alors des d\'efinitions que $\Sigma^{L'[\epsilon']}=\Sigma'\cap \xi_{T,T'}(\Sigma^{L[\epsilon]})$. D'o\`u $\xi_{T,T'}(Z(L[\epsilon]))\subset Z(L'[\epsilon'])$. Il existe un entier $c\geq1$ premier \`a $p$ tel que $\epsilon^c\in Z(L[\epsilon])$. Cela entra\^{\i}ne $(\epsilon')^c\in Z(L'[\epsilon'])$. Donc $\epsilon'$ est un $p'$-\'el\'ement dans $G'(F)$. Cela d\'emontre le (i).
   
   Le (ii) est imm\'ediat.
   
   Supposons ${\bf G}'$ elliptique. Alors $Z(G)(F)$ est un sous-groupe d'indice fini de $Z(G')(F)$ et la premi\`ere assertion du (iii) en r\'esulte. Compte tenu de (i), il reste \`a voir que, si $x'$ est compact mod $Z(G')$ et si $x'=\epsilon'exp(X')$ est une $p'$-d\'ecomposition de $x'$, alors $x=\xi_{T,T'}^{-1}(\epsilon')exp(\xi_{T,T'}^{-1}(X'))$ est une $p'$-d\'ecomposition de $x$. Comme ci-dessus, $\xi_{T,T'}^{-1}(X')$ est topologiquement nilpotent. Puisque $x'$ est compact mod $Z(G')$, $\epsilon'$ est $p'$-compact. Mais l'hypoth\`ese $(Hyp)(G)$ entra\^{\i}ne que l'indice de $Z(G)(F)$ dans $Z(G')(F)$ est d'ordre premier \`a $p$. Il en r\'esulte que $\xi_{T,T'}^{-1}(\epsilon')$ est $p'$-compact mod $Z(G)$, ce qui ach\`eve la preuve. 
   $\square$

  \subsection{Facteur de transfert\label{facteurdetransfert}}
  Soit ${\bf G}'=(G',{\cal G}',s)$ une donn\'ee endoscopique de $G$.    Rappelons que le facteur de transfert est une fonction d\'efinie sur l'ensemble  ${\bf D}_{1}$ des couples $(x'_{1},x)$ tels que: $x'_{1}\in G'_{1,reg}(F)$, $x\in G_{reg}(F)$ et, en notant $x'$ l'image de $x'_{1}$ dans $G'(F)$, les classes de conjugaison stable de  $x'$ et de  $x$ se correspondent. Consid\'erons un tel couple et fixons une
   une $p'$-d\'ecomposition $x=\epsilon exp(X)$. Par le (i) du 
   lemme \ref{correspondance}, on en d\'eduit une $p'$-d\'ecomposition $x'=\epsilon' exp(X')$. L'application $\mathfrak{g}'_{1,tn}(F)\to \mathfrak{g}'_{tn}(F)$ est surjective. Fixons $X'_{1}\in \mathfrak{g}'_{1,tn}(F)$ au-dessus de $X'$. D\'efinissons $\epsilon'_{1}\in G'_{1}(F)$ par l'\'egalit\'e $x'_{1}=\epsilon'_{1}exp(X'_{1})$. On voit que $\epsilon'_{1}$ est un $p'$-\'el\'ement, donc que l'\'egalit\'e $x'_{1}=\epsilon'_{1}exp(X'_{1})$ est une $p'$-d\'ecomposition.   Pour $\lambda\in \mathfrak{o}_{F}-\{0\}$, le couple $(\epsilon'_{1}exp(\lambda X'_{1}),\epsilon exp(\lambda X))$ appartient encore \`a l'ensemble de d\'efinition  ${\bf D}_{1}$ du facteur de transfert. 
 
 \ass{Lemme}{Pour tout $\lambda\in \mathfrak{o}_{F}-\{0\}$, on a l'\'egalit\'e
 $$\Delta_{1}(\epsilon'_{1}exp(\lambda^2 X'_{1}),\epsilon exp(\lambda^2 X))=\Delta_{1}(\epsilon'_{1}exp( X'_{1}),\epsilon exp( X)).$$}
 
 Preuve. Notons $T$, $T'$ et $T'_{1}$ les commutants de $x$, $x'$ et $x'_{1}$ dans $G$, $G'$ et $G'_{1}$. Pour calculer le facteur de transfert, on peut choisir des $\chi$-data mod\'er\'ement ramifi\'ees. Les termes $\Delta_{I}$ et $\Delta_{III,1}$ de \cite{LS} ne d\'ependent que des tores $T$ et $T'_{1}$. Le terme $\Delta_{III,2}$ est la valeur en $x'_{1}$ d'un caract\`ere de $T'_{1}(F)$. Celui-ci est construit \`a l'aide des $\chi$-data et est mod\'er\'ement ramifi\'e. Donc $\Delta_{III,2}(\epsilon'_{1}exp( X'_{1}),\epsilon exp( X))$ ne d\'epend pas du couple $(X'_{1},X)$. Il reste le terme $\Delta_{II}$. Un calcul similaire \`a celui de \cite{W3} montre que, comme fonction de $(X'_{1},X)$, $\Delta_{II}(\epsilon'_{1}exp( X'_{1}),\epsilon exp( X))$ est proportionnel  au terme analogue $\Delta_{II}(X',X)$ d\'efini sur les alg\`ebres de Lie, cf. \cite{W4}  2.3. D'apr\`es le lemme 3.2.1 de \cite{F}, la fonction $\lambda\mapsto \Delta_{II}(\lambda X',\lambda X)$ d\'efinie sur $F^{\times}$ est un caract\`ere quadratique. Elle est donc constante sur les carr\'es. Le lemme en r\'esulte. $\square$

  \subsection{Actions des centres\label{actiondescentres}}
  Soit ${\bf G}'=(G',{\cal G}',s)\in {\cal E}(G)$. Puisque cette donn\'ee est elliptique, le groupe $Z(G)(F)$ s'identifie \`a un sous-groupe d'indice fini de $Z(G')(F)$. Gr\^ace \`a  l'hypoth\`ese $(Hyp)(G)$, cet indice est premier \`a $p$. Il en r\'esulte une \'egalit\'e $Z(G)_{tu}(F)=Z(G')_{tu}(F)$. D'autre part, le sous-groupe $A_{G}(F)\subset Z(G)(F)$ est \'egal au sous-groupe $A_{G'}(F)\subset Z(G')(F)$. Le groupe $A_{G'_{1}}(F)$ s'envoie surjectivement sur $A_{G'}(F)$ (parce que $C_{1}$ est un tore induit). Notons ${\cal Z}_{1}$ le groupe des $z'_{1}\in Z(G'_{1})(F)$ tels que l'image de $z'_{1}$ dans $Z(G')(F)$ appartienne \`a $Z(G)(F)$.  Ce groupe  contient le groupe engendr\'e par $C_{1}(F)$, $A_{G'_{1}}(F)$ et $Z(G'_{1})_{tu}(F)$.  Pour tout couple $(x'_{1},x)\in {\bf D}_{1}$  et pour tout $z'_{1}\in {\cal Z}_{1}$, le couple $(z'_{1}x'_{1},zx)$  appartient aussi \`a ${\bf D}_{1}$, o\`u $z$ est l'image de $z'_{1}$ dans $Z(G)(F)$. D'apr\`es \cite{LS} lemme 4.4.A, il existe un caract\`ere $\zeta_{1}$ de ${\cal Z}_{1}$ de sorte que, pour $(x'_{1},x)\in {\bf D}_{1}$, on ait l'\'egalit\'e
  $$\Delta_{1}(z'_{1}x'_{1},zx)=\zeta(z'_{1})\Delta_{1}(x'_{1},x).$$
  Ce caract\`ere co\"{\i}ncide avec $\lambda_{1}^{-1}$ sur $C_{1}(F)$.
  Gr\^ace \`a l'hypoth\`ese $(Hyp)_{endo}(G)$, on a
  
  (1) le caract\`ere $\zeta_{1}$ est trivial sur le sous-ensemble  $Z(G'_{1})_{tu}(F)$ de ${\cal Z}_{1}$. 
  
  La preuve est similaire \`a celle du lemme pr\'ec\'edent.
  \subsection{D\'ecomposition de ${\mathbb C}[Ell(G)]$\label{decomposition}}
  
  Il y a une d\'ecomposition
  $$I_{cusp}(G)=\oplus_{{\bf G}'\in{\cal E}(G)}I_{cusp}(G)_{{\bf G}'}$$
  o\`u chaque facteur $I_{cusp}(G)_{{\bf G}'}$ est le sous-espace des $f\in I_{cusp}(G)$ dont les transferts \`a tout ${\bf G}''\not={\bf G}'$ sont nuls, cf. \cite{MW} proposition I.4.11. Fixons un caract\`ere unitaire $\xi$ de $A_{G}(F)$. Il y a une variante de la d\'ecomposition ci-dessus o\`u l'on impose que les fonctions se transforment selon ce caract\`ere $\xi$.  Dualement, il y a une d\'ecomposition
  $$(1)\qquad {\mathbb C}[Ell(G)_{\xi}]=\oplus_{{\bf G}'\in{\cal E}(G)}{\mathbb C}[Ell(G)_{\xi}]_{{\bf G}'}.$$
  Elle est orthogonale pour le produit scalaire elliptique. Pour $\pi\in {\mathbb C}[Ell(G)_{\xi}]$, notons $\pi=\sum_{{\bf G}'\in {\cal E}(G)}\pi_{{\bf G}'}$ l'\'ecriture de $\pi$ selon cette d\'ecomposition.
  
  \ass{Proposition}{(i) Si $\pi\in {\mathbb C}[Ell(G)^0_{\xi}]$, alors $\pi_{{\bf G}'}\in {\mathbb C}[Ell(G)^0_{\xi}]$ pour tout ${\bf G}'\in {\cal E}(G)$. 
 
 (ii) Si $\pi\in {\mathbb C}[Ell(G)^{>0}_{\xi}]$, alors $\pi_{{\bf G}'}\in {\mathbb C}[Ell(G)^{>0}_{\xi}]$ pour tout ${\bf G}'\in {\cal E}(G)$. }
 
 Preuve. Soient $\pi\in {\mathbb C}[Ell(G)_{\xi}]$ et ${\bf G}'\in {\cal E}(G)$. Soit $x\in G_{ell}(F)$. On sait que les classes de conjugaison par $G(F)$ contenues dans la classe de conjugaison stable de $x$ forment un espace principal homog\`ene sous l'action d'un groupe ab\'elien fini $K(x)$.  Pour $k\in K(x)$, notons $kx$ un repr\'esentant de l'image par $k$ de la classe de conjugaison de $x$ (on choisit $1x=x$).  Notons $K(x)^{\vee}$ le groupe des caract\`eres de $K(x)$. Chaque ${\bf G}'\in {\cal E}(G)$ d\'etermine un sous-ensemble $K(x)^{\vee}_{{\bf G}'}$ de $K(x)^{\vee}$ et $K(x)^{\vee}$ est union disjointe des  $K(x)^{\vee}_{{\bf G}'}$ quand ${\bf G}'$ d\'ecrit ${\cal E}(G)$. Pour tout ${\bf G}'\in {\cal E}(G)$, on a alors l'\'egalit\'e
 $$(2) \qquad \theta_{\pi_{{\bf G}'}}(x)=\vert K(x)\vert ^{-1}\sum_{\kappa\in K(x)^{\vee}_{{\bf G}'}}\sum_{k\in K(x)}\kappa(k)\theta_{\pi}(kx).$$
 Soient $\epsilon\in G(F)$ un \'el\'ement $p'$-compact mod $Z(G)$ et $X\in \mathfrak{g}_{\epsilon,tn}(F)$. Posons $x=\epsilon exp(X)$ et supposons $x\in G_{ell}(F)$. Pour tout $k\in K(x)$, fixons $g_{k}\in G(\bar{F})$ tel que $kx=g_{k}^{-1}xg_{k}$. Posons $\epsilon_{k}=g_{k}^{-1}\epsilon g_{k}$ et $X_{k}=g_{k}^{-1} X g_{k}$. Alors $x_{k}=\epsilon_{k}exp(X_{k})$ est une $p'$-d\'ecomposition de $x_{k}$. Soit $\lambda\in \mathfrak{o}_{F}-\{0\}$, posons $y=\epsilon exp(\lambda^2 X)$. On a $K(x)=K(y)$ et on voit comme dans la preuve du lemme \ref{induction}  que l'on peut choisir $y_{k}=\epsilon_{k}exp(\lambda^2 X_{k})$ pour tout $k\in K(y)$. L'\'egalit\'e (2) pour $y$ devient
 $$\theta_{\pi_{{\bf G}'}}(\epsilon exp(\lambda^2 X))=\vert K(x)\vert ^{-1}\sum_{\kappa\in K(x)^{\vee}_{{\bf G}'}}\sum_{k\in K(x)}\kappa(k)\theta_{\pi}(\epsilon_{k}exp(\lambda^2 X)).$$
 Supposons $\pi\in {\mathbb C}[Ell(G)^0_{\xi}]$. Comme fonctions de $\lambda$, les termes du membre de droite appartiennent \`a $E$. Donc aussi la fonction $\lambda\mapsto \theta_{\pi_{{\bf G}'}}(\epsilon exp(\lambda^2 X))$. C'est-\`a-dire que $\theta_{\pi_{{\bf G}'}}$ est un quasi-caract\`ere de niveau $0$ sur les elliptiques. D'apr\`es le lemme \ref{unlemme}, $\pi_{{\bf G}'}$ appartient \`a ${\mathbb C}[Ell(G)_{\xi}^0]$. Cela d\'emontre le (i) de l'\'enonc\'e.

 Supposons $\pi\in {\mathbb C}[Ell(G)^{>0}_{\xi}]$, soit  $\sigma\in {\mathbb C}[Ell(G)^0_{\xi}]$. Parce que la d\'ecomposition (1) est orthogonale, on a
 $$<\pi_{{\bf G}'},\sigma>_{ell}=<\pi_{{\bf G}'},\sigma_{{\bf G}'}>_{ell}=<\pi,\sigma_{{\bf G}'}>_{ell}.$$
 On vient de prouver que $\sigma_{{\bf G}'}$ \'etait de niveau $0$. Donc le produit de droite est nul. Il en estde m\^eme de celui de gauche. Cela \'etant vrai pour tout $\sigma\in {\mathbb C}[Ell(G)^0_{\xi}]$, cela entra\^{\i}ne $\pi_{{\bf G}'}\in {\mathbb C}[Ell(G)^{>0}_{\xi}]$. $\square$
 
 \subsection{Le cas quasi-d\'eploy\'e\label{lecasquasideploye}}
 Supposons $G$ quasi-d\'eploy\'e. On note $SI(G)$ le quotient de $C_{c}^{\infty}(G(F))$ par le sous-espace des fonctions dont toutes les int\'egrales orbitales stables sont nulles. En identifiant les int\'egrales orbitales \`a des formes lin\'eaires sur $I(G)$, c'est aussi le quotient de $I(G)$  par le sous-espace des \'el\'ements dont toutes les int\'egrales orbitales stables sont nulles. Notons
 $SI(G)^*$ l'espace des distributions stables sur $G(F)$, c'est-\`a-dire le dual de $SI(G)$. On note ${\mathbb C}[Irr(G)]^{st}$ le sous-espace des $\pi\in {\mathbb C}[Irr(G)]$ telles que $\Theta_{\pi}\in SI(G)^*$. On d\'efinit de m\^eme des espaces ${\mathbb C}[Ell(G)]^{st}$ etc...   
  
 On d\'efinit une projection lin\'eaire $p^{st}:{\mathbb C}[Irr(G)]\to {\mathbb C}[Irr(G)]^{st}$ de la fa\c{c}on suivante.
 Il y a une donn\'ee endoscopique elliptique maximale ${\bf G}=(G,{^LG},1)$ que l'on suppose appartenir \`a ${\cal E}(G)$. Dans la d\'ecomposition (1) de \ref{decomposition}, on a l'\'egalit\'e ${\mathbb C}[Ell(G)_{\xi}]^{st} ={\mathbb C}[Ell(G)_{\xi}]_{{\bf G}}$.  Alors $p^{st}:{\mathbb C}[Ell(G)_{\xi}]\to {\mathbb C}[Ell(G)_{\xi}]^{st}$ est la projection associ\'ee \`a cette d\'ecomposition (1) de   \ref{decomposition}, c'est-\`a-dire qu'elle annule les composantes ${\mathbb C}[Ell(G)_{\xi}]_{{\bf G}'}$ pour ${\bf G}'\not={\bf G}$.  
Plus g\'en\'eralement, le sous-espace ${\mathbb C}[Ell(G)_{{\mathbb R}}]^{st}$ est engendr\'e par les $\pi_{\chi}$, o\`u $\pi$ parcourt ${\mathbb C}[Ell(G)_{\xi}]^{st}$, $\xi$ parcourt le groupe des carcat\`eres unitaires de $A_{G}(F)$ et $\chi$ parcourt ${\cal A}_{G}^*$. On d\'efinit $p^{st}:{\mathbb C}[Ell(G)_{{\mathbb R}}]\to {\mathbb C}[Ell(G)_{{\mathbb R}}]^{st}$ de sorte que, si $\pi\in {\mathbb C}[Ell(G)_{\xi}]$ et $\chi\in {\cal A}_{G}^*$, on ait $p^{st}(\pi_{\chi})=p^{st}(\pi)_{\chi}$. 
  On a l'\'egalit\'e suivante, similaire \`a   \ref{representationselliptiques} (1):
 $$(1) \qquad{\mathbb C}[Irr(G)]^{st}=\oplus_{M\in \underline{{\cal L}}_{min}}Ind_{M}^G({\mathbb C}[Ell(M)_{{\mathbb R}}]^{st,W^G(M)}).$$
 Des  projections  que l'on vient de d\'efinir pour chaque Levi se d\'eduit alors la projection  $p^{st}$ cherch\'ee.  
 
 \ass{Corollaire} {(i) Pour $\pi\in {\mathbb C}[Irr(G)^0]$, on a $p^{st}(\pi)\in {\mathbb C}[Irr(G)^0]^{st}$. Pour $\pi\in {\mathbb C}[Irr(G)^{>0}]$, on a $p^{st}(\pi)\in {\mathbb C}[Irr(G)^{>0}]^{st}$.
 
 (ii) On a l'\'egalit\'e $p^0\circ p^{st}=p^{st}\circ p^0$. }
 
 Preuve. Le (i) r\'esulte imm\'ediatement de la proposition \ref{decomposition} et le (ii) est \'equivalent au (i). $\square$
 
 \subsection{Un premier r\'esultat de transfert\label{unpremierresultat}}
 Soit ${\bf G}'=(G',{\cal G}',s)\in {\cal E}(G)$.   On note $SI_{\lambda_{1}}(G_{1})$ la variante de $SI(G_{1})$ o\`u l'on impose que les fonctions se transforment par translations par $C_{1}(F)$ selon le caract\`ere $\lambda_{1}^{-1}$. On d\'efinit de fa\c{c}on similaire $SI_{\lambda_{1}}(G_{1})^*$ (c'est le dual de  $SI_{\lambda_{1}}(G_{1})$). On note $Aut({\bf G}')$ le groupe d'automorphismes de ${\bf G}'$, cf. \cite{MW} I.1.5. Ce groupe agit naturellement sur $SI_{\lambda_{1}}(G_{1})$ et $SI_{\lambda_{1}}(G_{1})^*$, cf. \cite{MW} I.2.6. 
 
 {\bf Remarque.} La d\'efinition adopt\'ee ici et dans \cite{MW} n'est pas celle que l'on trouve dans d'autres r\'ef\'erences.  Elle diff\`ere de la d\'efinition adopt\'ee notamment par Arthur par torsion par un caract\`ere de sorte que le facteur $\Delta_{1}$ soit invariant par cette action. Cela explique la disparition dans nos formules des caract\`eres que l'on trouve dans celles d'Arthur.
 \bigskip
 
 On sait d\'efinir un transfert $I(G)\to SI_{\lambda_{1}}(G_{1})$, dont l'image est contenue dans le sous-espace des invariants par $Aut({\bf G}')$. Une cons\'equence de \cite{A3}  est qu'il existe dualement une application lin\'eaire de transfert spectral
 $$transfert:{\mathbb C}[Irr(G_{1})_{\lambda_{1}}]^{st}\to{\mathbb C}[Irr(G)].$$
 Elle est invariante par l'action de $Aut({\bf G}')$, c'est-\`a-dire qu'elle se factorise en
 $$transfert:{\mathbb C}[Irr(G_{1})_{\lambda_{1}}]^{st}\to{\mathbb C}[Irr(G_{1})_{\lambda_{1}}]^{st,Aut({\bf G}')}\to{\mathbb C}[Irr(G)],$$
 o\`u la premi\`ere application est la projection naturelle.
 
 \ass{Proposition}{Cette application $transfert$ envoie ${\mathbb C}[Irr(G_{1})^0_{\lambda_{1}}]^{st}$ dans ${\mathbb C}[Irr(G)^0]$.}
 
   Preuve. Soit $\pi'\in {\mathbb C}[Irr(G_{1})^0_{\lambda_{1}}]^{st}$, posons $\pi=transfert(\pi')$. Soit $\epsilon\in G(F)_{p'}$  et soit $X\in \mathfrak{g}_{\epsilon,tn}(F)$. Posons $x=\epsilon exp(X)$ et supposons  $x\in G_{reg}(F)$.  Fixons un ensemble de repr\'esentants $X^{{\bf G}'}(x)$  des classes de conjugaison stable dans $G'(F)$ qui correspondent \`a celle de $x$ dans $G(F)$. Pour tout $x'\in X^{{\bf G}'}(x)$, on fixe un rel\`evement $x'_{1}$ de $x'$ dans $G'_{1}(F)$.  Par d\'efinition, on a l'\'egalit\'e
   $$(1) \qquad \theta_{\pi}(x)D^G(x)^{1/2}=\sum_{x'\in X^{{\bf G}'}(x)}\Delta_{1}(x'_{1},x)\theta_{\pi'}(x'_{1})D^{G'}(x')^{1/2}.$$
   On note $T$ le commutant de $x$ dans $G$ et, pour tout $x'\in X^{{\bf G}'}(x)$, on note $T_{x'}$ celui de $x'$ dans $G'$. On a des isomorphismes $\xi_{T,T_{x'}}:T\to T_{x'}$, cf. \ref{correspondance}. On pose $\epsilon_{x'}=\xi_{T,T_{x'}}(\epsilon)$, $X_{x'}=\xi_{T,T_{x'}}(X)$. Alors $x'=\epsilon_{x'}exp(X_{x'})$ est une $p'$-d\'ecomposition. On fixe un rel\`evement $X_{1,x'}$ de $X_{x'}$ dans $\mathfrak{g}'_{1,tn}(F)$ et un rel\`evement $\epsilon_{1,x'}$ de $\epsilon_{x'}$ dans $G'_{1}(F)$ de sorte que $x'_{1}=\epsilon_{1,x'}exp(X_{1,x'})$. 
   Soit $\lambda\in F^{\times}$ tel que $\lambda X\in \mathfrak{g}_{\epsilon,tn}(F)$, posons $y=\epsilon exp(\lambda X)$.  Alors, pour tout $x'\in X^{{\bf G}'}(x)$, l'\'el\'ement $\epsilon_{x'}exp(\lambda X_{x'})$ correspond \`a $ y$. Montrons que
   
   (2) si $x',x''\in X^{{\bf G}'}(x)$ et $x'\not=x''$, alors $\epsilon_{x'}exp(\lambda X_{x'})$ n'est pas stablement conjugu\'e \`a $\epsilon_{x''}exp(\lambda X_{x''})$.
   
   Supposons qu'ils sont stablement conjugu\'es. Il existe alors  $h\in G'(\bar{F})$ tel que 
   
   \noindent $h^{-1}\epsilon_{x'}exp(\lambda X_{x'})h=\epsilon_{x''}exp(\lambda X_{x''})$. Puisque les \'el\'ements en question sont fortement r\'eguliers, l'op\'erateur de conjugaison par $h$ se restreint en un isomorphisme $Int_{h}:T_{x''}\to T_{x'}$. Le compos\'e $\Xi=\xi_{T,T_{x'}}^{-1}\circ Int_{h}\circ\xi_{T,T_{x''}}$ est un automorphisme de $T$ fixant $y$. De plus, par construction des isomorphismes $\xi_{T,T_{x'}}$ et $\xi_{T,T_{x''}}$, $\Xi$ est forc\'ement de la forme $Int_{h'}$ pour un \'el\'ement $h'\in Norm_{G}(T)(\bar{F})$. Puisqu'il fixe $y$ et que $y$ est fortement r\'egulier, c'est l'identit\'e. Mais alors $Int_{h}(\epsilon_{x''})=\epsilon_{x'}$, $Int_{h}(X_{x''})=X_{x'}$ donc $h^{-1}x'h=x''$, ce qui contredit la d\'efinition de $X^{{\bf G}'}(x)$. Cela d\'emontre (2). 
   
   Cette assertion entra\^{\i}ne que $X^{{\bf G}'}(y)$ a au moins autant d'\'el\'ements que $X^{{\bf G}'}(x)$. La situation \'etant sym\'etrique en $x$ et $y$, ces ensembles ont m\^eme nombre d'\'el\'ements. Alors (2) entra\^{\i}ne que l'on peut choisir pour $X^{{\bf G}'}(y)$ l'ensemble des $\epsilon_{x'}exp(\lambda X_{x'})$ pour $x'\in X^{{\bf G}'}(x)$. On applique cela en rempla\c{c}ant $\lambda$ par $\lambda^2$ pour $\lambda\in \mathfrak{o}_{F}-\{0\}$. L'\'egalit\'e (1) pour $y=\epsilon exp(\lambda^2 X)$ devient
 $$(3)\qquad  \theta_{\pi}(\epsilon exp(\lambda^2X))D^G(\epsilon exp(\lambda^2 X))^{1/2}=\sum_{x'\in X^{{\bf G}'}(x)}\Delta_{1}(\epsilon_{1,x'}exp(\lambda^2X_{1,x'}),\epsilon exp(\lambda^2X))$$
 $$\theta_{\pi'}(\epsilon_{1,x'}exp(\lambda^2X_{1,x'}))D^{G'}(\epsilon_{x'}exp(\lambda^2X_{x'}))^{1/2}.$$ 
 On consid\`ere ces termes comme des fonctions de $\lambda$. Comme on l'a dit plusieurs fois, les discriminants de Weyl contribuent par des puissances de $\vert \lambda\vert _{F}$. D'apr\`es le lemme \ref{facteurdetransfert}, le facteur de transfert est constant. Puisque $\pi'$ est de niveau $0$, le th\'eor\`eme \ref{caracterisation} dit que $\Theta_{\pi'}$ est un quasi-caract\`ere de niveau $0$. Donc les termes  $\theta_{\pi'}(\epsilon_{1,x'}exp(\lambda^2X_{1,x'}))$ appartiennent \`a $E$. Donc $\lambda\mapsto  \theta_{\pi}(\epsilon exp(\lambda^2X))$ appartient \`a $E$, c'est-\`a-dire que $\Theta_{\pi}$ est un quasi-caract\`ere de niveau $0$. D'apr\`es le m\^eme th\'eor\`eme \ref{caracterisation}, $\pi$ appartient \`a   ${\mathbb C}[Irr(G)^0]$.$ \square$
 
 \subsection{Une r\'eciproque partielle\label{unereciproque}}
 Soit ${\bf G}'=(G',{\cal G}',s)\in {\cal E}(G)$.
  L'application $transfert$  se restreint en une application lin\'eaire
    $$transfert_{ell}:{\mathbb C}[Ell(G_{1})_{\lambda_{1}}]^{st}\to {\mathbb C}[Ell(G)].$$
  \ass{Proposition}{Soit $\pi\in {\mathbb C}[Ell(G)^0]$. Supposons que $\pi$ appartienne \`a l'image de l'application $transfert_{ell}$. Alors il existe $\pi'\in   {\mathbb C}[Ell(G_{1})^0_{\lambda_{1}}]^{st}$ telle que $transfert_{ell}(\pi')=\pi$.}
  
  Preuve. Fixons $\sigma'\in{\mathbb C}[Ell(G_{1})_{\lambda_{1}}]^{st}$ tel que $transfert_{ell}(\sigma')=\pi$. D\'ecomposons $\sigma'$ en somme $\sum_{i=1,...,h}\sigma'_{i}$, o\`u $\sigma'_{i}$ est une combinaison lin\'eaire de repr\'esentations elliptiques de $G'_{1}(F)$ dont le caract\`ere central se restreint en un m\^eme caract\`ere $\xi'_{i}$ de $A_{G'_{1}}(F)$. On suppose les $\xi'_{i}$ distincts. On a introduit un caract\`ere $\zeta_{1}$ en \ref{actiondescentres}. Le  caract\`ere $\xi'_{i}(\zeta_{1})_{ \vert A_{G'_{1}}(F)}$ de $A_{G'_{1}}(F)$  se factorise  par la projection $A_{G'_{1}}(F)\to A_{G}(F)$ en un caract\`ere $\xi_{i}$ de ce dernier groupe et les caract\`eres $\xi_{i}$ sont tous distincts. On a l'\'egalit\'e $\pi=\sum_{i=1,...,h}transfert_{ell}(\sigma'_{i})$ et, puisque $\pi$ est de niveau $0$,  on en d\'eduit par interpolation que toutes les composantes $transfert_{ell}(\sigma'_{i})$ sont de niveau $0$. On peut remplacer $\pi$ par chacune de ces composantes et cela ram\`ene le probl\`eme au cas o\`u  $A_{G'_{1}}(F)$ agit sur $\sigma'$ selon un caract\`ere $\xi'$. Ensuite, on peut encore tordre $\sigma'$ par un \'el\'ement de ${\cal A}_{G'_{1}}^*$ et supposer que $\xi'$ est unitaire. 
  
  Pour tout \'el\'ement $\epsilon_{1}\in G'_{1}(F)$ qui est $p'$-compact mod $Z(G'_{1})$, fixons un voisinage $\mathfrak{V}_{\epsilon_{1}}$ de $0$ dans $\mathfrak{g}'_{1,\epsilon_{1}}(F)$ et une famille $(c_{\epsilon_{1},{\cal O}})_{{\cal O}\in Nil(\mathfrak{g}'_{1,\epsilon_{1}})}$ de nombres complexes de sorte que, pour tout $X_{1}\in \mathfrak{V}_{\epsilon_{1}}$ tel que $\epsilon_{1}exp(X_{1})\in G'_{1,reg}(F)$, on ait l'\'egalit\'e
  $$(1) \qquad \theta_{\sigma'}(\epsilon_{1}exp(X_{1}))=\sum_{{\cal O}\in Nil(\mathfrak{g}'_{1,\epsilon_{1}})}c_{\epsilon_{1},{\cal O}}\hat{j}({\cal O},X_{1}).$$
  Le voisinage $\mathfrak{V}_{\epsilon_{1}}$ n'est pas uniquement d\'etermin\'e mais les nombres complexes $c_{\epsilon_{1},{\cal O}}$ le sont. Soit $z_{1}\in Z(G'_{1})_{tu}(F)$. On peut remplacer $\epsilon_{1}$ par $\epsilon_{1}z_{1}$. Cela ne change pas l'ensemble $Nil(\mathfrak{g}'_{1,\epsilon_{1}})$. Montrons que:
  
  (2) pour tout ${\cal O}\in Nil(\mathfrak{g}'_{1,\epsilon_{1}})$, la fonction $z_{1}\mapsto c_{\epsilon_{1}z_{1},{\cal O}}$ est localement constante;
  
  (3) on peut choisir les voisinages $\mathfrak{V}_{\epsilon_{1}}$ de sorte que la fonction $z_{1}\mapsto \mathfrak{V}_{\epsilon_{1}z_{1}}$ soit localement constante.
  
  En effet, choisissons un sous-groupe $\mathfrak{Y}_{\epsilon_{1}}\subset \mathfrak{g}_{\epsilon_{1}}(F)$ de sorte que $\mathfrak{Y}_{\epsilon_{1}}+\mathfrak{V}_{\epsilon_{1}}\subset \mathfrak{V}_{\epsilon_{1}}$.   On v\'erifie que, si $z_{1}\in \mathfrak{Y}_{\epsilon_{1}}$, on peut choisir $\mathfrak{V}_{\epsilon_{1}z_{1}}=\mathfrak{V}_{\epsilon_{1}}$ et que l'on a
   $c_{\epsilon_{1}z_{1},{\cal O}}=c_{\epsilon_{1},{\cal O}}$ pour tout ${\cal O}\in Nil(\mathfrak{g}'_{1,\epsilon_{1}})$. Cela prouve (2) et (3).
   
   Puisque $Z(G'_{1})_{tu}(F)$ est compact, l'assertion (3) se renforce imm\'ediatement en: on peut choisir les voisinages $\mathfrak{V}_{\epsilon_{1}}$ de sorte que $\mathfrak{V}_{\epsilon_{1}z_{1}}=\mathfrak{V}_{\epsilon_{1}}$ pour tout $z_{1}\in Z(G'_{1})_{tu}(F)$. On suppose d\'esormais qu'il en est ainsi. Pour ${\cal O}\in Nil(\mathfrak{g}'_{1,\epsilon_{1}})$, posons
 $$\bar{c}_{\epsilon_{1},{\cal O}}=mes(Z(G'_{1})_{tu}(F))^{-1}\int_{Z(G'_{1})_{tu}(F)}c_{\epsilon_{1}z_{1},{\cal O}}\,dz_{1}.$$
  On a $\bar{c}_{\epsilon_{1}z_{1},{\cal O}}=\bar{c}_{\epsilon_{1},{\cal O}}$ pour tout $z_{1}\in Z(G'_{1})_{tu}(F)$. On d\'efinit une fonction $\tau$  sur $G'_{1,reg}(F)$, \`a support contenu dans l'ensemble $G'_{1,comp}(F)$ des \'el\'ements de $G'_{1}(F)$ qui sont compacts mod $Z(G'_{1})$ de la fa\c{c}on suivante. Pour $x_{1}\in G'_{1,comp}(F)\cap G'_{1,reg}(F)$, on choisit une $p'$-d\'ecomposition $x_{1}=\epsilon_{1}exp(X_{1})$. L'\'el\'ement $\epsilon_{1}$ est $p'$-compact mod $Z(G'_{1})$. On pose
  $$\tau(x_{1})=\sum_{{\cal O}\in Nil(\mathfrak{g}'_{1,\epsilon_{1}})}\bar{c}_{\epsilon_{1},{\cal O}}\hat{j}({\cal O},X_{1}).$$
  Cela ne d\'epend pas de la $p'$-d\'ecomposition choisie car, d'apr\`es \ref{elementsp}(2), toute autre d\'ecomposition est de la forme $x_{1}=(\epsilon_{1}exp(Z_{1}))exp(X_{1}-Z_{1})$, avec $Z_{1}\in \mathfrak{z}(G'_{1})_{tn}(F)$ et on a fait ce qu'il fallait pour que la formule ci-dessus soit insensible \`a un tel changement. Il r\'esulte de cette d\'efinition que, pour tout  \'el\'ement $\epsilon_{1}\in G'_{1}(F)$ qui est $p'$-compact mod $Z(G'_{1})$ et pour tout $X_{1}\in \mathfrak{V}_{\epsilon_{1}}$ tel que $\epsilon_{1}exp(X_{1})\in G_{1,reg}(F)$, on a l'\'egalit\'e
  $$(4)\qquad \tau(\epsilon_{1}exp(X_{1}))=mes(Z(G'_{1})_{tu}(F))^{-1}\int_{Z(G'_{1})_{tu}(F)}\theta_{\sigma'}(\epsilon_{1}z_{1}exp(X_{1}))\,dz_{1}.$$
  Montrons que
  
  (5) la fonction $\tau$ est constante sur les classes de conjugaison stable.
  
  Soient $x_{1},x_{2}\in G_{1,reg}(F)$ que l'on suppose stablement conjugu\'es. On fixe $h\in G_{1,reg}(\bar{F})$ tel que $x_{2}=h^{-1}x_{1}h$. On choisit une $p'$-d\'ecomposition $x_{1}=\epsilon_{1}exp(X_{1})$. Posons $\epsilon_{2}=h^{-1}\epsilon_{1}h$ et $X_{2}=h^{-1}X_{1}h$. Alors $x_{2}=\epsilon_{2}exp(X_{2})$ et cette \'egalit\'e est une $p'$-d\'ecomposition. On veut prouver que $\tau(x_{1})=\tau(x_{2})$. Remarquons que, ou bien  $x_{1}$ et $x_{2}$ sont tous deux compacts mod $Z(G'_{1})$, ou bien aucun d'eux ne l'est. Dans le second cas, on a  $\tau(x_{1})=\tau(x_{2})=0$. Supposons que $x_{1}$ et $x_{2}$ sont  compacts mod $Z(G'_{1})$. Pour $\lambda\in \mathfrak{o}_{F}-\{0\}$, les \'el\'ements $\epsilon_{1}exp(\lambda^2X_{1})$ et $\epsilon_{2}exp(\lambda^2X_{2})$ sont encore stablement conjugu\'es. Si $val_{F}(\lambda)$ est assez grand, on a $\lambda^2X_{1}\in \mathfrak{V}_{\epsilon_{1}}$ et $\lambda^2X_{2}\in \mathfrak{V}_{\epsilon_{2}}$, donc $\tau(\epsilon_{1}exp(\lambda^2X_{1}))$ et $\tau(\epsilon_{2}exp(\lambda^2X_{2}))$ sont calcul\'es par la formule (4). Pour $z_{1}\in Z(G'_{1})_{tu}(F)$, les \'el\'ements $\epsilon_{1}z_{1}exp(\lambda^2X_{1})$ et $\epsilon_{2}z_{1}exp(\lambda^2X_{2})$ sont encore stablement conjugu\'es. Puisque $\sigma'$ est stable, on a donc
  $$\theta_{\sigma'}(\epsilon_{1}z_{1}exp(\lambda^2X_{1}))=\theta_{\sigma'}(\epsilon_{2}z_{1}exp(\lambda^2X_{2})).$$
  La formule (4) entra\^{\i}ne alors l'\'egalit\'e $\tau(\epsilon_{1}exp(\lambda^2X_{1}))=\tau(\epsilon_{2}exp(\lambda^2X_{2}))$. Cela est vrai pour $val_{F}(\lambda)$ assez grande. Mais, comme fonctions de $\lambda$,  les deux termes de cette \'egalit\'e appartiennent \`a $E$. Ils sont donc \'egaux pour tout $\lambda$. En $\lambda=1$, cela entra\^{\i}ne l'\'egalit\'e cherch\'ee $\tau(x_{1})=\tau(x_{2})$ qui prouve (5).

 En particulier, $\tau$   
   est invariante par conjugaison.  Alors $\tau$ est la fonction associ\'ee \`a un quasi-caract\`ere   sur $G'_{1}(F)$ de niveau $0$. Comme $\Theta_{\sigma'}$, ce quasi-caract\`ere se transforme  par $A_{G'_{1}}(F)$ selon le caract\`ere $\xi$ et par $C_{1}(F)$ selon $\lambda_{1}$. D'apr\`es le th\'eor\`eme \ref{undeuxiemetheoreme}, il existe un \'el\'ement $d_{1}\in {\cal D}_{comp}(G'_{1})$ tel que ce quasi-caract\`ere soit \'egal \`a $D^{G'_{1}}[d_{1}]$.  Notons $d_{1,cusp}$ la composante cuspidale de $d_{1}$. Cet \'el\'ement conserve les m\^emes propri\'et\'es que ci-dessus de transformation par $A_{G'_{1}}(F)$ et $C_{1}(F)$. On note $\pi'$ l'\'el\'ement de ${\mathbb C}[Ell(G'_{1})_{\xi'}^0]$ tel que $\Delta_{cusp}(\pi')=d_{1,cusp}$, cf. le (ii) du lemme \ref{surjectivite}.   Par d\'efinition, $\theta_{\pi'}$ co\"{\i}ncide avec $\tau$ sur $G'_{1,ell}(F)$ et cette fonction est invariante par conjugaison stable. Parce que $\pi'$ est elliptique, cela entra\^{\i}ne que $\pi'$ est stable, c'est-\`a-dire $\pi'\in {\mathbb C}[Ell(G'_{1})_{\xi'}^0]^{st}$, cf.  \cite{A3} th\'eor\`eme 6.1. 
    Il est clair que la restriction \`a $C_{1}(F)$ du  caract\`ere central de $\pi'$  est $\lambda_{1}$.

   Pour prouver la proposition, il suffit de prouver que $transfert_{ell}(\pi')=\pi$. Or $\pi$ et $\pi'$ sont elliptiques. D'apr\`es   \cite{A3} th\'eor\`eme 6.2, il suffit de prouver que $transfert_{ell}(\pi')$ et $\pi$ co\"{\i}ncident sur $G_{ell}(F)$. C'est-\`a-dire qu'en posant  $\Pi=transfert_{ell}(\pi')$,
   il suffit de prouver

  (6) pour tout $x\in G_{ell}(F)$, on a l'\'egalit\'e $\theta_{\Pi}(x)=\theta_{\pi}(x)$.
  
  On fixe une $p'$-d\'ecomposition $x=\epsilon exp(X)$. Soit $\lambda\in \mathfrak{o}_{F}-\{0\}$.     
  L'\'egalit\'e (3) de \ref{unpremierresultat} dit que
   $$\theta_{\Pi}( \epsilon exp(\lambda^2X))D^G(\epsilon exp(X))^{1/2}=\sum_{x'\in X^{{\bf G}'}(x)}\Delta_{1}(\epsilon_{1,x'}exp(\lambda^2X_{1,x'}),\epsilon exp(\lambda^2X))$$
 $$\theta_{\pi'}(\epsilon_{1,x'}exp(\lambda^2X_{1,x'}))D^{G'}(\epsilon_{x'}exp(\lambda^2X_{x'}))^{1/2}.$$
 Les \'el\'ements intervenant \`a droite sont elliptiques, on peut remplacer la fonction $\theta_{\pi'}$ par $\tau$. Si $val_{F}$ est assez grand, on utilise (4) et on obtient
    $$\theta_{\Pi}( \epsilon exp(\lambda^2X))D^G(\epsilon exp(X))^{1/2}=mes(Z(G'_{1})_{tu}(F))^{-1}\int_{Z(G'_{1})_{tu}(F)}\sum_{x'\in X^{{\bf G}'}(x)}$$
    $$\Delta_{1}(\epsilon_{1,x'}exp(\lambda^2X_{1,x'}),\epsilon exp(\lambda^2X))\theta_{\sigma'}(\epsilon_{1,x'}z_{1}exp(\lambda^2X_{1,x'}))D^{G'}(\epsilon_{x'}exp(\lambda^2X_{x'}))^{1/2}\,dz_{1}.$$
 Notons $z$ et $z'$ les images de $z_{1}$ dans $Z(G)_{tu}(F)$ et $ Z(G')_{tu}(F)$. D'apr\`es \ref{actiondescentres} (1), on peut remplacer ci-dessus $\Delta_{1}(\epsilon_{1,x'}exp(\lambda^2X_{1,x'}),\epsilon exp(\lambda^2X))$ par $\Delta_{1}(\epsilon_{1,x'}z_{1}exp(\lambda^2X_{1,x'}),\epsilon z exp(\lambda^2X))$. On peut aussi remplacer $D^{G'}(\epsilon_{x'}exp(\lambda^2X_{x'}))$ par $D^{G'}(\epsilon_{x'}z'exp(\lambda^2X_{x'}))$. On obtient
  $$ (7)\qquad \theta_{\Pi}( \epsilon exp(\lambda^2X))D^G(\epsilon exp(X))^{1/2}=mes(Z(G'_{1})_{tu}(F))^{-1}\int_{Z(G'_{1})_{tu}(F)}B(z_{1})\,dz_{1},$$
  o\`u
  $$B(z_{1})=\sum_{x'\in X^{{\bf G}'}(x)}\Delta_{1}(\epsilon_{1,x'}z_{1}exp(\lambda^2X_{1,x'}),\epsilon z exp(\lambda^2X))$$
 $$\theta_{\sigma'}(\epsilon_{1,x'}z_{1}exp(\lambda^2X_{1,x'}))D^{G'}(\epsilon_{x'}z'exp(\lambda^2X_{x'}))^{1/2}.$$
 Il est clair que $B(z_{1})$ est le membre de droite de la formule (3) de \ref{unpremierresultat} appliqu\'ee  \`a la repr\'esentation $\sigma'$, pour le point $\epsilon z exp(\lambda^2X)$ de $G(F)$. Puisque $\pi=transfert_{ell}(\sigma')$, on obtient 
 $$B(z_{1})=\theta_{\pi}(\epsilon z exp(\lambda^2X))D^G(\epsilon z exp(\lambda^2X))^{1/2}.$$
 Le $z$ dispara\^{\i}t du discriminant de Weyl. Il dispara\^{\i}t aussi du premier terme car $\pi$ est de niveau $0$ par hypoth\`ese. D'o\`u
 $$B(z_{1})=\theta_{\pi}(\epsilon  exp(\lambda^2X))D^G(\epsilon exp(\lambda^2X))^{1/2},$$
 puis, d'apr\`es (7),
 $$\theta_{\Pi}( \epsilon exp(\lambda^2X)) =\theta_{\pi}( \epsilon exp(\lambda^2X)).$$
 On a suppos\'e $val_{F}(\lambda)$ assez grande. Mais, puisque $\pi'$ est de niveau $0$,   $\Pi$  l'est aussi d'apr\`es la proposition \ref{unpremierresultat}. Il en est de m\^eme pour $\pi$ par hypoth\`ese. Les deux membres de la formule ci-dessus appartiennent donc \`a $E$. Etant \'egaux pour $val_{F}(\lambda)$ assez grande, ils sont \'egaux pour tout $\lambda$. En $\lambda=1$, cela prouve (6) et la proposition. $\square$
 
  \subsection{Le th\'eor\`eme principal\label{letheoreme}}
  
  Soit ${\bf G}'=(G',{\cal G}',s)\in {\cal E}(G)$.
  \ass{Th\'eor\`eme}{L'application 
  $$transfert:{\mathbb C}[Irr(G_{1})_{\lambda_{1}}]^{st}\to {\mathbb C}[Irr(G)]$$
  v\'erifie l'\'egalit\'e $p^0\circ transfert=transfert\circ p^0$.}
  
  {\bf Remarque.} Les deux $p^0$ de l'\'egalit\'e ne sont pas les m\^emes: le premier vit sur le groupe $G$ et le second sur $G'_{1}$. 
  
  \bigskip
  
  Preuve.  On a l'\'egalit\'e
  $$(1) \qquad {\mathbb C}[Irr(G_{1})_{\lambda_{1}}]^{st}=\oplus_{M'\in \underline{{\cal L}}_{min}^{G'}}Ind_{M'_{1}}^{G'_{1}}({\mathbb C}[Ell(M'_{1})_{{\mathbb R},\lambda_{1}}]^{st,W^{G'}(M')}),$$ 
  cf. \ref{lecasquasideploye} (1), dont on a adapt\'e les notations de fa\c{c}on que l'on esp\`ere compr\'ehensible. En particulier, chaque Levi $M'$ de $G'$ d\'etermine un unique Levi $M'_{1}$ de $G'_{1}$ (son image r\'eciproque dans ce groupe). Notons $\underline{{\cal L}}_{min,G-rel}^{G'}$ l'ensemble des \'el\'ements de $\underline{{\cal L}}_{min}^{G'}$ qui sont "relevants" pour $G$ (dans le cas o\`u $G$ est quasi-d\'eploy\'e, tout \'el\'ement de $\underline{{\cal L}}_{min}^{G'}$ est relevant). Il y a une application naturelle
  de $\underline{{\cal L}}_{min,G-rel}^{G'}$ dans $\underline{{\cal L}} _{min}^G$, cf. \cite{MW} 1.3.4. Si $M'\mapsto M$ par cette application, il se d\'eduit de ${\bf G}'$ une donn\'ee endoscopique ${\bf M}'=(M',{\cal M}',s^M)$ de $M$, qui est elliptique. 
  Soit $M'\in  \underline{{\cal L}}_{min,G-rel}^{G'}$, notons $M\in \underline{{\cal L}} _{min}^G$ le Levi qui lui correspond. On a l'\'egalit\'e ${\cal A}_{M}={\cal A}_{M'}$. L'application $transfert_{ell}: {\mathbb C}[Ell(M'_{1})_{\lambda_{1}}]^{st}\to {\mathbb C}[Ell(M)]$ se prolonge en une application lin\'eaire
  ${\mathbb C}[Ell(M'_{1})_{{\mathbb R},\lambda_{1}}]^{st}\to {\mathbb C}[Ell(M)_{{\mathbb R}}]$ compatible aux op\'erateurs de torsion par des \'el\'ements $\chi\in {\cal A}_{M}^*={\cal A}_{M'}^*$. On la compose avec la projection naturelle sur les invariants par $W^G(M)$ et on obtient par restriction une application lin\'eaire
  $${\mathbb C}[Ell(M'_{1})_{{\mathbb R},\lambda_{1}}]^{st,W^{G'}(M')}\to {\mathbb C}[Ell(M)_{{\mathbb R}}]^{W^G(M)}.$$
  Le transfert commute \`a l'induction et on en d\'eduit une application lin\'eaire de l'espace
   $Ind_{M'_{1}}^{G'_{1}}({\mathbb C}[Ell(M'_{1})_{{\mathbb R},\lambda_{1}}]^{st,W^{G'}(M')})$
   dans le membre de droite de la relation (1) de \ref{representationselliptiques}, c'est-\`a-dire dans  $ {\mathbb C}[Irr(G)]$. Si maintenant $M'\in \underline{{\cal L}}_{min}^{G'}-\underline{{\cal L}}_{min,G-rel}^{G'}$, on envoie  l'espace $Ind_{M'_{1}}^{G'_{1}}({\mathbb C}[Ell(M'_{1})_{{\mathbb R},\lambda_{1}}]^{st,W^{G'}(M')})$ dans $ {\mathbb C}[Irr(G)]$ par l'application nulle. La somme des applications ainsi d\'efinies sur tout $M'\in  \underline{{\cal L}}^{G'}$ est une application lin\'eaire du membre de droite de (1) ci-dessus dans $ {\mathbb C}[Irr(G)]$. C'est
    l'application $transfert$. Il r\'esulte de ces consid\'erations que, pour d\'emontrer le th\'eor\`eme, il suffit de prouver l'\'egalit\'e analogue pour les applications $transfert_{ell}$ associ\'es aux diff\'erents Levi de $G'$ qui sont relevants pour $G$. On ne perd rien \`a consid\'erer seulement l'application $transfert_{ell}$ associ\'ee \`a $G'$ lui-m\^eme. C'est-\`a-dire qu'il suffit de prouver
  
  (2) on a l'\'egalit\'e $p^0\circ transfert_{ell}=transfert_{ell}\circ p^0$.
  
    Comme dans le corollaire \ref{lecasquasideploye}, cela \'equivaut aux deux assertions
  
  (3) si $\pi'\in {\mathbb C}[Ell(G'_{1})_{\lambda_{1}}^0]^{st}$, on a $transfert_{ell}(\pi')\in {\mathbb C}[Ell(G)^0]$;
  
  (4) si $\pi'\in {\mathbb C}[Ell(G'_{1})_{\lambda_{1}}]^{st}$ v\'erifie $p^0(\pi')=0$, alors $p^0\circ transfert_{ell}(\pi')=0$. 
  
  La proposition \ref{unpremierresultat} est justement l'assertion (3). D\'emontrons (4). Comme toujours, on peut d\'ecomposer $\pi'$ en $\sum_{i=1,...,h}\pi'_{i}$ o\`u $\pi'_{i}$ est une combinaison lin\'eaire de repr\'esentations elliptiques de $G'_{1}(F)$ dont le caract\`ere central se restreint en un m\^eme caract\`ere unitaire $\xi'_{i}$ de $A_{G'_{1}}(F)$. Ces repr\'esentations v\'erifient les m\^emes hypoth\`eses que $\pi'$ et il suffit de prouver que $p^0\circ transfert_{ell}(\pi'_{i})=0$ pour tout $i$. En oubliant cette construction, on peut fixer un caract\`ere unitaire $\xi'$ de $A_{G'_{1}}(F)$ et supposer que $\pi'\in {\mathbb C}[Ell(G'_{1})_{\xi',\lambda_{1}}]^{st}$. Le caract\`ere $\xi'$ d\'etermine un caract\`ere $\xi$ de $A_{G}(F)$ comme on l'a vu dans la preuve de \ref{unereciproque} et l'application $transfert_{ell}$ envoie ${\mathbb C}[Ell(G'_{1})_{\xi',\lambda_{1}}]^{st}$ dans ${\mathbb C}[Ell(G)_{\xi}]$. L'action du groupe $Aut({\bf G})$ conserve le premier espace: elle pr\'eserve le caract\`ere $\xi'$ car celui-ci est uniquement d\'etermin\'e par $\xi$. On peut remplacer $\pi'$ par sa projection sur le sous-espace des invariants, cela ne change pas $transfert_{ell}(\pi')$. Cela ne modifie pas non plus la condition $p^0(\pi')=0$ car l'action de $Aut({\bf G}')$ conserve le sous-espace ${\mathbb C}[Ell(G'_{1})^0_{\xi',\lambda_{1}}]^{st}$. Bref, on peut supposer $\pi'\in {\mathbb C}[Ell(G'_{1})_{\xi',\lambda_{1}}]^{st,Aut({\bf G}')}$. On pose $\pi=transfert_{ell}(\pi')$. Il faut maintenant se rappeler la d\'ecomposition (1) de \ref{decomposition}.  L'application $transfert_{ell}$ se restreint en un isomorphisme
  $$ {\mathbb C}[Ell(G'_{1})_{\xi',\lambda_{1}}]^{st,Aut({\bf G}')}\to {\mathbb C}[Ell(G)_{\xi}]_{{\bf G}'}$$
  qui est une similitude pour les produits scalaires elliptiques. On note $c$ le rapport de similitude.  En particulier, $\pi\in {\mathbb C}[Ell(G)_{\xi}]_{{\bf G}'}$. Il r\'esulte de la proposition \ref{decomposition} que $p^0(\pi)$ appartient aussi \`a ${\mathbb C}[Ell(G)_{\xi}]_{{\bf G}'}$, donc appartient \`a l'image de $transfert_{ell}$. D'apr\`es la proposition \ref{unereciproque}, il existe $\pi^{'0}\in {\mathbb C}[Ell(G_{1})^0_{\lambda_{1}}]^{st}$ telle que $transfert_{ell}(\pi^{'0})=p^0(\pi)$. Les m\^emes arguments que ci-dessus permettent de supposer que $\pi^{'0}\in {\mathbb C}[Ell(G'_{1})^0_{\xi',\lambda_{1}}]^{st,Aut({\bf G}')}$. Parce que $p^0$ est une projection orthogonale, on a 
  $<p^0(\pi),p^0(\pi)>_{ell}=<p^0(\pi),\pi>_{ell}$. Donc
  $$<p^0(\pi),p^0(\pi)>_{ell}=<transfert_{ell}(\pi^{'0}),transfert_{ell}(\pi')>_{ell}$$
  $$=c^{-1}<\pi^{'0},\pi'>_{ell}=c^{-1}<\pi^{'0},p^0(\pi')>_{ell},$$
  toujours parce que $p^0$ est une projection orthogonale. Puisque $p^0(\pi')=0$, cela entra\^{\i}ne 
   $<p^0(\pi),p^0(\pi)>_{ell}=0$, d'o\`u $p^0(\pi)=0$. Cela ach\`eve la preuve de (4) et du th\'eor\`eme. $\square$
   
   \section{Quelques cons\'equences} 
   
   \subsection{Le cas quasi-d\'eploy\'e derechef\label{derechef}}
   Supposons $G$ quasi-d\'eploy\'e.
   \ass{Corollaire}{La projection de Bernstein $p^0$ de $I(G)$ se quotiente en une projection $p^0$ de $SI(G)$.}
   
   Preuve. Notons $I^{inst}(G)$ le noyau de la projection $I(G)\to SI(G)$. On doit prouver que $p^0$ conserve ce sous-espace. C'est-\`a-dire, soit $f\in I^{inst}(G)$, on doit prouver que $p^0(f)\in I^{inst}(G)$. Il suffit pour cela de prouver que, pour tout $\pi^{st}\in {\mathbb C}[Irr(G)]^{st}$, on a $\Theta_{\pi^{st}}(p^0(f))=0$. D'apr\`es les propri\'et\'es des projecteurs de Bernstein, on a
   $$\Theta_{\pi^{st}}(p^0(f))=\Theta_{p^0(\pi^{st})}(f).$$
   D'apr\`es le (ii) du corollaire \ref{lecasquasideploye}, on a $p^0(\pi^{st})=p^0\circ p^{st}(\pi^{st})=p^{st}\circ p^0(\pi^{st})$. Donc $p^0(\pi^{st})\in {\mathbb C}[Irr(G)]^{st}$. Alors $\Theta_{p^0(\pi^{st})}(f)=0$ puisque $f\in I^{inst}(G)$. $\square$
   
   \subsection{Transfert de fonctions\label{transfertdefonctions}}
   Soit ${\bf G}'=(G',{\cal G}',s)\in {\cal E}(G)$. 
   \ass{Corollaire}{L'application $transfert:I(G)\to SI_{\lambda_{1}}(G_{1})$ v\'erifie l'\'egalit\'e $p^0\circ transfert=transfert\circ p^0$.}
  Cela se d\'eduit du th\'eor\`eme \ref{letheoreme} comme le corollaire pr\'ec\'edent se d\'eduisait du corollaire  \ref{lecasquasideploye}.
  
  \subsection{Stabilit\'e dans  l'espace ${\cal D}_{cusp}$\label{stabilite}}
  Supposons $G$ quasi-d\'eploy\'e. Disons qu'un \'el\'ement $d\in {\cal D}_{cusp}(G)$ est stable si la distribution $D^G[d]$ l'est  et qu'il est stable sur les elliptiques si la distribution $D^G[d]$ l'est, c'est-\`a-dire si $\theta_{D^G[d]}$ est constante sur les classes de conjugaison stable contenues dans $G_{ell}(F)$. On note ${\cal D}_{cusp}(G)^{st}$ le sous-espace des \'el\'ements stables dans ${\cal D}_{cusp}(G)$.
  
  \ass{Proposition}{Soit $d\in {\cal D}_{cusp}(G)$. Alors $d$ est stable si et seulement si $d$ est stable sur les elliptiques.}
  
  Preuve. Utilisons la d\'ecomposition ${\cal D}_{cusp}(G)=\prod_{\nu\in {\cal N}}{\cal D}_{cusp}(G)^{\nu}$. Il est clair que $d$ est stable, resp. stable sur les elliptiques, si et seulement si chaque composante $d^{\nu}$ est stable, resp. stable sur les elliptiques. On peut donc fixer $\nu\in {\cal N}$ et supposer $d\in {\cal D}_{cusp}(G)^{\nu}$. Notons ${\cal X}$ l'ensemble des restrictions \`a $A_{G}(F)_{c}$ de caract\`eres mod\'er\'ement ramifi\'es de $A_{G}(F)$. Pour $\xi\in {\cal X}$, posons $d_{\xi}=mes(A_{G}(F)_{c})^{-1}\int_{A_{G}(F)_{c}}\xi(a)^{-1}{d^{a}}\, da$. Alors $d=\sum_{\xi\in {\cal X}}d_{\xi}$. De nouveau, $d$ est stable, resp. stable sur les elliptiques, si et seulement si chaque composante $d_{\xi}$ est stable, resp. stable sur les elliptiques. On peut fixer $\xi\in {\cal X}$ et supposer que $d$ se transforme par $A_{G}(F)_{c}$ selon $\xi$. Prolongeons $\xi$ en un caract\`ere   unitaire de $A_{G}(F)$. Alors il existe un unique ${\bf d}\in {\cal D}_{cusp}(G)$ tel que:
  
  ${\bf d}\in {\cal D}_{cusp,\xi}(G)$, c'est-\`a-dire ${\bf d}$ se transforme par $A_{G}(F)$ selon $\xi$;
  
  pour $\nu'\in {\cal N}$, la composante ${\bf d}^{\nu'}$ est nulle si $\nu'\not\in \nu+w_{G}(A_{G}(F))$;
  
  ${\bf d}^{\nu}=d^{\nu}$.
  
  On voit que $d$ est stable, resp. stable sur les elliptiques, si et seulement si ${\bf d}$ l'est. En oubliant cette construction, on peut fixer un caract\`ere unitaire $\xi$ de $A_{G}(F)$ et supposer $d\in {\cal D}_{cusp,\xi}(G)$. D'apr\`es le (ii) du lemme \ref{surjectivite}, il existe un unique $\pi\in {\mathbb C}[Ell(G)_{\xi}^0]$ telle que $\Delta_{cusp}(\pi)=d$. Alors, d'apr\`es \ref{leresultat} (1), $\Theta_{\pi}$ est stable sur les elliptiques. D'apr\`es \cite{A3} th\'eor\`eme 6.1, $\pi$ est stable. Soient $M\in {\cal L}_{min}$ et $P\in {\cal P}(M)$. Puisque $\pi$ est stable, il en est de m\^eme du module de Jacquet $\pi_{P}$. D'apr\`es 
   \ref{leresultat} (1), $\Delta^M_{\pi_{P},cusp}$ est stable sur les elliptiques. A fortiori, sa projection $\Delta^M_{\pi_{M},cusp,G-comp}$ dans ${\cal D}_{cusp,G-comp}(M)$ est stable sur les elliptiques (cette projection consiste \`a restreindre le support \`a un sous-ensemble invariant par conjugaison stable). Supposons que $M$ est un Levi propre. On peut raisonner par r\'ecurrence et supposer notre proposition d\'ej\`a d\'emontr\'ee pour $M$. Donc $\Delta^M_{\pi_{M},cusp,G-comp}$ est stable. C'est-\`a-dire que $D^M(\Delta^M_{\pi_{M},cusp,G-comp})$ est une  distribution stable sur $M(F)$. La stabilit\'e se conserve par induction donc $D^G(\Delta^M_{\pi_{M},cusp,G-comp})$ est une  distribution stable sur $G(F)$. La formule \ref{leresultat} (1) exprime alors $D^G[d]$ comme la diff\'erence entre $\Theta_{\pi}$ qui est stable et la somme sur les $M\not=G$ de distributions stables. Donc $D^G[d]$ est stable et $d$ l'est alors par d\'efinition. $\square$
   
   \subsection{Transfert et module de Jacquet\label{transfertetmoduledeJacquet}}
   Soit ${\bf G}'=(G',{\cal G}',s)\in {\cal E}(G)$.  On a rappel\'e dans la preuve du th\'eor\`eme \ref{letheoreme} l'existence d'une application  
  
  (1) $\underline{{\cal L}}_{min,G-rel}^{G'}\to\underline{{\cal L}}_{min}^G$. 
  
  On note simplement $M'\mapsto M$ cette application. Rappelons quelques-unes de ses propri\'et\'es. Elles sont "bien connues des sp\'ecialistes", bien que nous n'en ayons pas trouv\'e d'\'enonc\'e clair dans la litt\'erature (\cite{MW} est tr\`es incomplet).

  Si $M'\in \underline{{\cal L}}_{min,G-rel}^{G'}$ s'envoie sur $M\in \underline{{\cal L}}_{min}^G$, $M'$ est compl\'et\'e en une donn\'ee endoscopique ${\bf M}'=(M',{\cal M}',s^M)$ de $M$, qui est elliptique.  En particulier, il y a un isomorphisme $A_{M'}\simeq A_{M}$. Les racines de $A_{M'}$ dans $G'$ sont aussi des racines de $A_{M}$ dans $G$. Ainsi, la chambre dans ${\cal A}_{M}^*$ associ\'ee \`a un sous-groupe parabolique $P\in {\cal P}(M)$ est contenue dans la chambre dans ${\cal A}_{M'}^*$ associ\'ee \`a un certain sous-groupe parabolique de composante de Levi $M'$, que l'on note $P'_{P}$.
    
  On sait d\'efinir une correspondance entre classes de conjugaison stable dans $G'_{reg}(F)$ et dans $G_{reg}(F)$. On d\'efinit de m\^eme \`a l'aide de ${\bf M}'$ une correspondance entre classes de conjugaison stable dans $M'_{reg}(F)$ et dans $M_{reg}(F)$. Appelons correspondance pour ${\bf G}'$ la premi\`ere et correspondance pour ${\bf M}'$ la seconde. On a
  
  (2) soient $x'\in M'(F)$ et $x\in M(F)\cap G_{reg}(F)$; supposons que $x$ et $x'$ se correspondent pour ${\bf M}'$; alors ils se correspondent pour ${\bf G}'$.
  
  Signalons que la r\'eciproque est fausse. 
Des donn\'ees auxiliaires que l'on a fix\'ees pour ${\bf G}'$ se d\'eduisent des donn\'ees pour ${\bf M}'$. En particulier $M'_{1}$ est l'image r\'eciproque de $M'$ dans $G'_{1}$.   On a fix\'e un facteur de transfert pour la donn\'ee ${\bf G}'$, notons-le plus pr\'ecis\'ement $\Delta_{1}^{{\bf G}'}$. On peut fixer un facteur de transfert $\Delta_{1}^{{\bf M}'}$ pour la donn\'ee ${\bf M}'$ de sorte que

(3) pour $x$, $x'$ comme en (2) et pour $x'_{1}\in M'_{1}(F)$ au-dessus de $x'$, on ait  
   $\Delta_{1}^{{\bf M}'}(x_{1}',x)=\Delta_{1}^{{\bf G}'}(x'_{1},x)$.
   
   Il existe un plongement $W^{G'}(M')\to W^G(M)$ v\'erifiant la condition suivante. Soient $n'\in Norm_{G'}(M')(F)$ et $n\in Norm_{G}(M)(F)$. Supposons que l'image de $n'$ dans $W^{G'}(M')$ s'envoie par ce plongement sur l'image de $n$ dans $W^G(M)$. On a
   
   (4)  soient $x'\in M'(F)$ et $x\in M(F)\cap G_{reg}(F)$; alors $x$ et $x'$ se correspondent pour ${\bf M}'$ si et seulement si $n'x'n^{'-1}$ et $nxn^{-1}$ se correspondent pour ${\bf M}'$. 
   
   Consid\'erons maintenant un Levi  $M'\in \underline{{\cal L}}_{min}^{G'}-\underline{{\cal L}}_{min,G-rel}^{G'}$. On a
   
   (5) pour $x'\in M'(F)\cap G'_{reg}(F)$, il n'existe pas de $x\in G_{reg}(F)$ correspondant \`a $x'$ pour ${\bf G}'$.

  Soit $\pi'\in {\mathbb C}[ Irr(G')]^{st}$, posons $\pi=transfert(\pi')$. Soient $M\in \underline{{\cal L}}_{min}$ et
  $P\in {\cal P}(M)$. Notons $\pi_{P}$ le module de Jacquet normalis\'e de $\pi$ relatif \`a $P$ et $ \Theta_{\pi_{P}\vert ell}$ la restriction de $\Theta_{\pi_{P}}$ \`a $M_{ell}(F)$.     
    
  \ass{Lemme}{Sous ces hypoth\`eses, on a l'\'egalit\'e
 $$ \Theta_{\pi_{P}\vert ell}=\sum_{M'\in \underline{{\cal L}}_{min}^{G'}; M'\mapsto M}\vert  W^{G'}(M')\vert ^{-1}\sum_{w\in W^{G}(M)}w^{-1}\circ transfert(\Theta_{\pi'_{P'_{w(P)}}\vert ell}).$$}
 
 Preuve. Soit $x\in M_{ell}(F)$.  Notons  $X^{{\bf G}'}(x)$  l'ensemble des classes de conjugaison stable dans $G'(F)$ qui correspondent \`a celle de $x$ pour ${\bf G}'$. Pour $M'\in \underline{{\cal L}}^{G'}_{min,G-rel}$ tel que $M'\mapsto M$,  notons $X^{{\bf M}'}(x)$   l'ensemble des classes de conjugaison stable dans $M'_{reg}(F)$ qui correspondent \`a celle de $x$  pour ${\bf M}'$. Relevons tout \'el\'ement de $W^G(M)$ en un \'el\'ement de $Norm_{G}(M)(F)$ et identifions $W^G(M)$ \`a l'ensemble de ces \'el\'ements.
 On va prouver
 
 (6) toute classe $C^{{\bf M'}}\in \cup_{w\in W^G(M)}X^{{\bf M}'}(wxw^{-1})$ est contenue dans une classe appartenant \`a  $X^{{\bf G}'}(x)$;
 
 on note 
  $\iota^{{\bf M}'}:\cup_{w\in W^G(M)}X^{{\bf M}'}(wxw^{-1})\to X^{{\bf G}'}(x)$ l'application qui, \`a un \'el\'ement $C^{{\bf M'}}$ de l'ensemble de d\'epart,   associe la classe qui la contient; 
 
 (7) la fibre de $\iota^{{\bf M}'}$ au-dessus de tout \'el\'ement de son image a $\vert W^{G'}(M')\vert $ \'el\'ements; 
 
 (8) pour deux \'el\'ements $M'_{1},M'_{2}\in \underline{{\cal L}}^{G'}_{min,G-rel}$ tels que $M'_{1}\not=M'_{2}$ et $M'_{1},M'_{2}\mapsto M$, les images de $\iota^{{\bf M}'_{1}}$ et $\iota^{{\bf M}'_{2}}$ sont disjointes;
 
 (9) $X^{{\bf G}'}(x)$ est r\'eunion des images des $\iota^{{\bf M}'}$ quand $M'$ d\'ecrit les \'el\'ements de $\underline{{\cal L}}_{min,G-rel}$ tels que $M'\mapsto M$.
 
  D\'emontrons (6).  Soient $w\in W^G(M)$, $C^{{\bf M}'}\in X^{{\bf M}'}(wxw^{-1})$ et $x'\in C^{{\bf M}'}$. Alors $x'$ correspond \`a $wxw^{-1}$ pour ${\bf M}'$ donc aussi pour ${\bf G}'$. Puisque $wxw^{-1}$ est conjugu\'e \`a $x$ par un \'el\'ement de $G(F)$, $x'$ correspond aussi \`a $x$ pour ${\bf G}'$. Par d\'efinition, $x'$  appartient donc \`a une unique classe $C^{{\bf G}'} \in X^{{\bf G}'}(x)$ et alors $C^{{\bf M}'}\subset C^{{\bf G}'}$.  
  
 Commen\c{c}ons la preuve de (7). Soient $w\in W^G(M)$, $C^{{\bf M}'}\in  X^{{\bf M}'}(wxw^{-1})$ et $x'\in C^{{\bf M}'}$. On identifie comme ci-dessus $W^{G'}(M')$ \`a un sous-ensemble de $Norm_{G'}(M')(F)$. Pour $w'\in W^{G'}(M')$, $w'x'w^{'-1}$ correspond \`a $w'wx(w'w)^{-1}$ d'apr\`es (4). Donc $w'x'w^{'-1}$ appartient \`a une unique classe $w'(C^{{\bf M}'})\in X^{{\bf M}'}(w'wx(w'w)^{-1})$ et on voit qu'elle  ne d\'epend pas du choix de $x'$.   Par construction, $w'x'w^{'-1}$ et $x'$ sont stablement conjugu\'es dans $G'(F)$ donc $C^{{\bf M}'}$ et $w'(C^{{\bf M}'})$ ont la m\^eme image par  $\iota^{{\bf M}'}$. Montrons que, si $w'\not=1$ (dans $W^{G'}(M')$), 
 les deux classes $C^{{\bf M}'}$ et $w'(C^{{\bf M}'})$ sont distinctes. Sinon, il existe $m'\in M'(\bar{F})$ tel que $ m'w'x'w^{'-1}m^{'-1}=x'$.  Puisque $x'$ est r\'egulier dans $G'(F)$ (parce qu'il correspond \`a $x$ par ${\bf G}'$ et $x$ est r\'egulier dans $G$), cela entra\^{\i}ne $m'w'\in G'_{x'}(\bar{F})\subset M'(\bar{F})$. Alors $w'\in M'(\bar{F})$, d'o\`u $w'=1$ dans $W^{G'}(M')$, contrairement \`a l'hypoth\`ese. Cela d\'emontre que le groupe $W^{G'}(M')$ agit librement sur l'espace de d\'epart de $\iota^{{\bf M}'}$ et que chaque orbite est contenue dans une fibre de cette application. 
 
 On va d\'emontrer la r\'eciproque et en m\^eme temps la relation (8). 
 Pour $i=1,2$, soient $M'_{i}\in \underline{{\cal L}}^{G'}_{min,G-rel}$, $w_{i}\in W^G(M)$, $C_{i}^{{\bf M}_{i}'}\in  X^{{\bf M}_{i}'}(w_{i}xw_{i}^{-1})$ et $x'_{i}\in C_{i}^{{\bf M}_{i}'}$.  Notons $T'_{i}$ le commutant de $x'_{i}$ dans $G'$. Supposons que $\iota^{{\bf M}'_{1}}(C_{1}^{{\bf M}'_{1}})=\iota^{{\bf M}'_{2}}(C_{2}^{{\bf M}'_{2}})$. Alors $x'_{1}$ et $x'_{2}$ sont stablement conjugu\'es dans $G'(F)$. C'est-\`a-dire qu'il existe $h\in G'(\bar{F})$ tel que $x'_{2}=hx'_{1}h^{-1}$. Puisque $x'_{1}$ et $x'_{2}$ appartiennent \`a $G'_{reg}(F)$, on a 
 
 (10) $\sigma(h)\in hT'_{1}(\bar{F})$ pour tout $\sigma\in \Gamma_{F}$. 
 
 Cela implique $hA_{T'_{1}}h^{-1}=A_{T'_{2}}$. 
  Pour $i=1,2$, $x'_{i}$ correspond pour ${\bf M}'_{i}$ \`a $w_{i}xw_{i}^{-1}$ qui est par hypoth\`ese elliptique dans $M(F)$. Donc $x'_{i}$ est elliptique dans $M'_{i}(F)$, ce qui entra\^{\i}ne  $A_{T'_{i}}=A_{M'_{i}}$.   Il en r\'esulte que $hA_{M'_{1}}h^{-1}=A_{M'_{2}}$. Donc la conjugaison par $h$ envoie   le commutant $M'_{1}$ de $A_{M'_{1}}$ sur le commutant $M'_{2}$ de  $A_{M'_{2}}$. Soit $P'_{1}\in {\cal P}(M'_{1})$, posons $P'_{2}=hP'_{1}h^{-1}$. Alors (10) implique que $ P'_{2}$ est d\'efini sur $F$ et $P'_{2}\in {\cal P}(M'_{2})$. De telles paires $(M'_{1},P'_{1})$ et $(M'_{2},P'_{2})$ qui sont conjugu\'es par un \'el\'ement de $G'(\bar{F})$ le sont par un \'el\'ement de $G'(F)$. 
 Par d\'efinition de $\underline{{\cal L}}_{min}$, on a donc $M'_{1}=M'_{2}$. Cela d\'emontre (8). Continuons le calcul en posant simplement $M'=M'_{1}=M'_{2}$. On a maintenant $h\in Norm_{G'}(M')(\bar{F})$. La relation (10) entra\^{\i}ne que l'image de $h$ dans $Norm_{G'}(M')(\bar{F})/M'(\bar{F})$ est fix\'ee par $\Gamma_{F}$. On sait qu'alors cette image est aussi celle d'un \'el\'ement $u\in Norm_{G'}(M')(F)$. 
   Autrement dit, on peut \'ecrire $h=um'$, o\`u $u\in Norm_{G'}(M')(F)$ et $m'\in M'(\bar{F})$. En notant encore $u$ l'image de $u$ dans $W^{G'}(M')$, on montre facilement que $C_{2}^{{\bf M}'}$ est alors l'image de $C_{1}^{{\bf M}'}$ par l'action de $u$. Cela d\'emontre que les fibres de l'application $\iota^{{\bf M}'}$ sont exactement les orbites de l'action du groupe $W^{G'}(M')$ sur l'espace de d\'epart de $\iota^{{\bf M}'}$. D'o\`u (7).
   
   D\'emontrons (9). Soient $C^{{\bf G}'}\in X^{{\bf G}'}(x)$ et $x'\in C^{{\bf G}'}$. Quitte \`a conjuguer $x'$, on peut supposer qu'il existe un Levi $M'\in \underline{{\cal L}}^{G'}_{min}$ tel que $x'\in M'_{ell}(F)$. D'apr\`es (5), $M'$ appartient \`a $\underline{{\cal L}}^{G'}_{min,G-rel}$. Soit $\underline{M}$ l'image de $M'$ dans $\underline{{\cal L}}_{min}$ par l'application (1). Parce que ${\bf M}'$ est une donn\'ee endoscopique elliptique de $\underline{M}$ et que $x'$ est elliptique dans $M'$, on sait que $x'$ correspond pour ${\bf M}'$ \`a un \'el\'ement $y\in \underline{M}_{ell}(F)$. D'apr\`es (2), $x'$ correspond aussi \`a $y$ pour ${\bf G}'$. Cela entra\^{\i}ne que $x$ et $y$ sont stablement conjugu\'es dans $G(F)$. Le m\^eme raisonnement que dans la deuxi\`eme partie de la preuve de (7) montre alors d'une part que $\underline{M}=M$ donc $M'\mapsto M$, d'autre part que la classe de conjugaison stable de $y$ dans $M(F)$ est l'image de celle de $x$ par un \'el\'ement $w\in W^G(M)$. Il en r\'esulte que la classe de conjugaison stable $C^{{\bf M}'}$ de $x'$ appartient \`a $X^{{\bf M}'}(wxw^{-1})$. D'o\`u (10).

 Maintenant, on note pour simplifier $w(x)=wxw^{-1}$ et on identifie
 les ensembles $X^{{\bf G}'}(x)$ et $X^{{\bf M}'}(w(x))$ \`a des ensembles de repr\'esentants des classes en question. 
 En cons\'equence des propri\'et\'es ci-dessus, l'\'egalit\'e \ref{unpremierresultat}(1) peut se r\'ecrire
 
 $$\theta_{\pi}(x)D^G(x)^{1/2}=\sum_{M'\in \underline{{\cal L}}^{G'}_{min,G-rel}; M'\mapsto M} \vert  W^{G'}(M')\vert ^{-1}\sum_{w\in W^{G}(M)}$$
 $$\sum_{x'\in X^{{\bf M}'}(w(x))}\Delta_{1}^{{\bf M}'}(x'_{1},w(x))\theta_{\pi'}(x_{1}')D^{G'}(x')^{1/2}.$$
  Soit $a\in A_{M}(F)$. Pour tout $M'\in \underline{{\cal L}}^{G'}_{min,G-rel}$ tel que $M'\mapsto M$, et tout $w\in W^G(M)$, $w(a)$
   s'identifie \`a un \'el\'ement $a_{w}^{M'}\in A_{M'}(F)$, que l'on rel\`eve en un \'el\'ement $a^{M'}_{w,1}\in A_{M'_{1}}(F)$.  En rempla\c{c}ant $x$ par $ax$ dans la formule ci-dessus, on obtient
$$(11) \qquad \theta_{\pi}(ax)D^G(ax)^{1/2}=\sum_{M'\in \underline{{\cal L}}^{G'}_{min,G-rel}; M'\mapsto M} \vert  W^{G'}(M')\vert ^{-1}\sum_{w\in W^{G}(M)}$$
 $$\sum_{x'\in X^{{\bf M}'}(w(x))}\Delta_{1}^{{\bf M}'}(a^{M'}_{w,1}x_{1}',w(ax))\theta_{\pi'}(a_{w,1}^{M'}x_{1}')D^{G'}(a^{M'}_{w}x')^{1/2}.$$ 
 Supposons $\vert \alpha(a)\vert_{F}$ assez petit pour toute racine $\alpha$ de $A_{M}$ dans l'alg\`ebre de Lie de $U_{P}$. Casselman a prouv\'e qu'alors $\theta_{\pi}(ax)=\delta_{P}(ax)^{1/2}\theta_{\pi_{P}}(ax)$, cf. \cite{C} th\'eor\`eme 5.2. On voit que, sous la m\^eme hypoth\`ese sur $a$, on a $\delta_{P}(ax)^{1/2}D^G(ax)^{1/2}=D^M(ax)^{1/2}$. Le membre de gauche de (11) devient $\theta_{\pi_{P}}(ax)D^M(ax)^{1/2}$. Pour $M'$ et $w$ intervenant ci-dessus, on voit que $a^{M'}_{w}$ v\'erifie  que $\vert \alpha(a^{M'}_{w})\vert _{F}$ est petit pour toute racine $\alpha$ de $A_{M'}$ dans $P'_{w(P)}$. Le m\^eme argument transforme le membre de droite de (11) et on obtient
 
 $$(12)\qquad \theta_{\pi_{P}}(ax)D^M(ax)^{1/2}=\sum_{M'\in \underline{{\cal L}}^{G'}_{min,G-rel}; M'\mapsto M} \vert  W^{G'}(M')\vert ^{-1}\sum_{w\in W^{G}(M)}$$
 $$\sum_{x'\in X^{{\bf M}'}(w(x))}\Delta_{1}^{{\bf M}'}(a^{M'}_{w,1}x'_{1},w(ax))\theta_{\pi'_{P'_{w(P)}}}(a^{M'}_{w,1}x'_{1})D^{M'}(a^{M'}_{w}x')^{1/2}.$$
 Chaque terme est d\'efini pour tout $a\in A_{M}(F)$. En d\'ecomposant $\pi_{P'}$ en repr\'esentations irr\'eductibles, on voit que le membre de gauche est une combinaison lin\'eaire de caract\`eres de $A_{M}(F)$. Il en est de m\^eme du membre de droite (on utilise \ref{actiondescentres} pour les facteurs de transfert). On vient de voir que les deux membres \'etaient \'egaux pour $a$ dans un certain c\^one ouvert. Il en r\'esulte qu'ils sont partout \'egaux. On peut donc remplacer $a$ par (1) dans l'\'egalit\'e pr\'ec\'edente. Pour $M'$ et $w$ intervenant ci-dessus, posons  $\pi_{M',w}=transfert(\pi'_{P'_{w(P)}})$. On reconna\^{\i}t la somme int\'erieure  du membre de droite de (12) (pour $a=1$), c'est $\theta_{\pi_{M',w}}(w(x))D^M(w(x))^{1/2}$. On a $D^M(w(x))^{1/2}=D^M(x)^{1/2}$ et l'\'egalit\'e (12) se transforme en
 $$\theta_{\pi}(x)=\sum_{M'\in \underline{{\cal L}}^{G'}_{min,G-rel}; M'\mapsto M} \vert  W^{G'}(M')\vert ^{-1}\sum_{w\in W^{G}(M)}\theta_{\pi_{M',w}}(w(x)).$$
 Ceci est \'equivalent \`a l'\'egalit\'e de l'\'enonc\'e. $\square$

   \subsection{ L'espace ${\cal D}_{cusp}(G)$ et le transfert\label{finale}}
    
   \ass{Proposition}{Soient $d\in {\cal D}_{cusp}(G)$ et, pour tout ${\bf G}'\in {\cal E}(G)$, soit $d^{{\bf G}'}\in {\cal D}_{cusp,\lambda_{1}}(G'_{1})^{st}$. Supposons que $D^G[d]$ co\"{\i}ncide avec $\sum_{{\bf G}'\in {\cal E}(G)}transfert(D^{G'_{1}}[d^{{\bf G}'}])$    sur $G_{ell}(F)$. Alors 
   $D^G[d]=\sum_{{\bf G}'\in {\cal E}(G)}transfert(D^{G'_{1}}[d^{{\bf G}'}])$.}
   
   Preuve. Soit ${\bf G}'\in {\cal E}(G)$. Des injections naturelles $Z(\hat{G})\to Z(\hat{G}')\to Z(\hat{G}'_{1})$ se d\'eduisent des homomorphismes 
   
   (1) ${\cal N}^{G'_{1}}\to {\cal N}^{G'}\to {\cal N}$. 
   
   Pour $\nu\in {\cal N}$, on d\'efinit la composante $d^{\nu}$ de $d$ comme dans le paragraphe pr\'ec\'edent. On d\'efinit $d^{{\bf G}',\nu}$ comme le produit des composantes $d^{{\bf G}',\nu'}$ sur les $\nu'\in {\cal N}^{G'_{1}}$ qui s'envoient sur $\nu$ par la suite ci-dessus. La famille $((d^{{\bf G}'})_{{\bf G}'\in {\cal E}(G)},d)$ v\'erifie l'hypoth\`ese ou la conclusion de la proposition si et seulement si chaque famille $((d^{{\bf G}',\nu})_{{\bf G}'\in {\cal E}(G)},d^{\nu})$ les v\'erifie. On peut donc fixer $\nu\in {\cal N}$ et supposer $d\in {\cal D}_{cusp}(G)^{\nu}$. Notons ${\cal X}$ l'ensemble des restrictions \`a $A_{G}(F)_{c}$ de caract\`eres mod\'er\'ement ramifi\'es de $A_{G}(F)$. Pour ${\bf G}'\in {\cal E}(G)$, cet ensemble s'identifie \`a celui des restrictions \`a $A_{G'_{1}}(F)_{c}$  de caract\`eres mod\'er\'ement ramifi\'es de $A_{G'_{1}}(F)$ qui co\"{\i}ncident avec $\lambda_{1}$ sur $A_{G'_{1}}(F)\cap C_{1}(F)$: \`a $\xi\in {\cal X}$, on associe d'abord l'image r\'eciproque  de $\xi$ par la surjection $A_{G'_{1}}(F)_{c}\to A_{G}(F)_{c}$, puis le produit $\xi^{{\bf G}'}$ de ce caract\`ere par la restriction de $\zeta_{1}$ \`a $A_{G'_{1}}(F)_{c}$. Ainsi, comme dans le paragraphe pr\'ec\'edent, on associe \`a $d$, resp. $d^{{\bf G}'}$, des composantes $d_{\xi}$, resp. $d^{{\bf G}'}_{\xi^{{\bf G}'}}$, pour $\xi\in {\cal X}$. La famille $((d^{{\bf G}'})_{{\bf G}'\in {\cal E}(G)},d)$ v\'erifie l'hypoth\`ese ou la conclusion de la proposition si et seulement si chaque  famille $((d^{{\bf G}'}_{\xi^{{\bf G}'}})_{{\bf G}'\in {\cal E}(G)},d_{\xi})$ les v\'erifie. On peut fixer $\xi\in {\cal X}$ et supposer que $d=d_{\xi}$ et $d^{{\bf G}'}=d^{{\bf G}'}_{\xi^{{\bf G}'}}$ pour tout ${\bf G}'$.  Prolongeons $\xi$ en un caract\`ere   unitaire de $A_{G}(F)$, qui s'identifie comme ci-dessus pour tout ${\bf G}'$ \`a un caract\`ere  $\xi^{{\bf G}'}$ de $A_{G'_{1}}(F)$ dont la restriction \`a $A_{G'_{1}}(F)\cap C_{1}(F)$ co\"{\i}ncide avec la restriction de $\lambda_{1}$. On construit un \'el\'ement ${\bf d}\in {\cal D}_{cusp,\xi}(G)$ comme dans le paragraphe pr\'ec\'edent et de fa\c{c}on similaire des \'el\'ements ${\bf d}^{{\bf G}'}\in {\cal D}_{cusp,\lambda_{1},\xi^{{\bf G}'}}(G'_{1})$. La famille $((d^{{\bf G}'})_{{\bf G}'\in {\cal E}(G)},d)$ v\'erifie l'hypoth\`ese ou la conclusion de la proposition si et seulement si la famille $(({\bf d}^{{\bf G}'})_{{\bf G}'\in {\cal E}(G)},{\bf d})$ les v\'erifient. En r\'esum\'e, on peut fixer un caract\`ere unitaire $\xi$ de $A_{G}(F)$, qui d\'etermine pour tout ${\bf G}'$ un tel caract\`ere $\xi^{{\bf G}'}$ de $A_{G'_{1}}(F)$, et supposer $d\in {\cal D}_{cusp,\xi}(G)$ et $d^{{\bf G}'}\in {\cal D}_{cusp,\lambda_{1},\xi^{{\bf G}'}}(G'_{1})$ pour tout ${\bf G}'$.

   A ce point, on peut simplifier le probl\`eme en se ramenant au cas o\`u une unique donn\'ee ${\bf G}'$ intervient. Pour cela, on introduit la repr\'esentation $\pi\in {\mathbb C}[Ell(G)_{\xi}^0]$ telle que $\Delta_{cusp}(\pi)=d$. On utilise la d\'ecomposition (1) de \ref{decomposition} et on \'ecrit conform\'ement $\pi=\sum_{{\bf G}'\in {\cal E}(G)}\pi_{{\bf G}'}$. On sait que $\pi_{{\bf G}'}$ v\'erifie les m\^emes propri\'et\'es que $\pi$ d'apr\`es la proposition \ref{decomposition}. On note $d_{{\bf G}'}=\Delta_{cusp}(\pi_{{\bf G}'})$. On a $d=\sum_{{\bf G}'}d_{{\bf G}'}$. On voit que si $((d^{{\bf G}'})_{{\bf G}'\in {\cal E}(G)},d)$ v\'erifie l'hypoth\`ese  de l'\'enonc\'e, alors, pour tout ${\bf G}'\in {\cal E}(G)$, $D^G[d_{{\bf G}'}]$ co\"{\i}ncide sur $G_{ell}(F)$ avec $transfert(D^{G'_{1}}[d^{{\bf G}'}])$. Inversement, si $D^G[d_{{\bf G}'}]=transfert(D^{G'_{1}}[d^{{\bf G}'}])$ pour tout ${\bf G}'$, la conclusion de l'\'enonc\'e est v\'erifi\'ee. Cela nous ram\`ene au cas o\`u la famille $(d^{{\bf G}'})_{{\bf G}'\in {\cal E}(G)}$ a au plus une composante non nulle. On note d\'esormais ${\bf G}'$ l'indice de cette composante et on note $(d',d)\in {\cal D}_{cusp,\lambda_{1};\xi'}(G'_{1})^{st}\times {\cal D}_{cusp,\xi}(G)$  le couple auquel se r\'eduisent les donn\'ees de d\'epart.

   On introduit la repr\'esentation $\pi\in {\mathbb C}[Ell(G)_{\xi}^0]$ telle que $\Delta_{\pi,cusp}=d$ et  la repr\'esentation $\pi'\in {\mathbb C}[Ell(G'_{1})^0_{\lambda_{1},\xi'}]^{st}$ telle que $\Delta^{G'_{1}}_{\pi',cusp}=d'$. Alors,  d'apr\`es \ref{leresultat} (1), $\Theta_{\pi}$ co\"{\i}ncide avec $\Theta_{transfert(\pi')}$ sur les elliptiques.  D'apr\`es \cite{A3} th\'eor\`eme 6.2,  $\pi=transfert(\pi')$. Soit $P\in {\cal P}(M)$.  D'apr\`es \ref{leresultat} (1), le lemme \ref{transfertetmoduledeJacquet} entra\^{\i}ne que $D^M[\Delta^M_{\pi_{P},cusp}]$ co\"{\i}ncide sur $M_{ell}(F)$ avec 
 $$(3)\qquad  \sum_{M'\in \underline{{\cal L}}_{min,G-rel}^{G'}; M'\mapsto M}\vert  W^{G'}(M')\vert ^{-1}\sum_{w\in W^{G}(M)}w^{-1}\circ transfert( D^{M'_{1}}[\Delta^{M'_{1}}_{\pi'_{P'_{w(P)}},cusp}]).$$
 Remarquons que, pour $M'$ et $w$ intervenant ci-dessus, on peut transformer la donn\'ee endoscopique ${\bf M}'$ de $M$ en une donn\'ee $w^{-1}({\bf M}')$. Le terme que l'on somme peut se r\'ecrire $transfert(D^{w^{-1}(M'_{1})}[d^{w^{-1}(M'_{1})}])$ pour un \'el\'ement convenable  $d^{w^{-1}(M'_{1})}\in {\cal D}_{cusp,w^{-1}(\lambda_{1})}(w^{-1}(M'_{1}))$.  C'est-\`a-dire que la somme (3) est de la m\^eme forme que celle qui intervient dans l'\'enonc\'e de la proposition. 
 Supposons $M\not=G$. En raisonnant par r\'ecurrence, on suppose notre proposition prouv\'ee pour $M$. Alors  $D^M[\Delta^M_{\pi_{P},cusp}]$ est \'egal \`a (3). Pour un couple $(x'_{1},x)\in G'_{1,reg}(F)\times G_{reg}(F)$ d'\'el\'ements qui se correspondent, $x$ est compact mod $Z(G)$ si et seulement si $x'_{1}$ est compact mod $Z(G'_{1})$: cela parce
  que $Z(G)(F)$ est un sous-groupe de $Z(G')(F)$ d'indice fini.  Restreindre  $D^M[\Delta^M_{\pi_{P},cusp}]$ \`a $M(F)\cap G_{comp}(F)$ \'equivaut \`a restreindre chaque terme  de (3) index\'e par $M'$ \`a $M'(F)\cap G'_{1,comp}(F)$. Apr\`es restriction \`a $M(F)\cap G_{comp}(F)$, on sait que le terme $\Delta^M_{\pi_{P},cusp}$ ne d\'epend plus de $P$. Il en est de m\^eme des termes intervenant dans (3). On obtient que $D^M[\Delta^M_{\pi_{M},cusp,G-comp}]$ est \'egal \`a
  $$   \sum_{M'\in \underline{{\cal L}}_{min,G-rel}^{G'}; M'\mapsto M}\vert  W^{G'}(M')\vert ^{-1}\sum_{w\in W^{G}(M)} w^{-1}\circ transfert(D^{M'_{1}}[\Delta^{M'_{1}}_{\pi'_{M'},cusp,G'_{1}-comp}]).$$
  Induisons \`a $G$. Evidemment, l'induction est insensible \`a la torsion par un \'el\'ement de $W^G$. Les $w$ disparaissent de la formule et la  somme en $w$ est remplac\'ee par la multiplication par $\vert W^G(M)\vert $. D'autre part, l'induction commute au transfert. On obtient
  $$(4) \qquad D^G[\Delta^M_{\pi_{M},cusp,G-comp}]=   \sum_{M'\in \underline{{\cal L}}_{min,G-rel}^{G'}; M'\mapsto M}\vert  W^{G'}(M')\vert ^{-1}\vert W^G(M)\vert 
 $$
 $$ transfert(D^{G'_{1}}[\Delta^{M'_{1}}_{\pi'_{M'},cusp,G'_{1}-comp}]).$$
  On utilise l'\'egalit\'e \ref{leresultat} (2) qui 
   peut se r\'ecrire
  $$ D^G[d]= \Theta_{\pi\vert G_{comp}(F)}-\sum_{M\in \underline{{\cal L}}_{min},M\not=G}\vert W^G(M)\vert ^{-1}D^G[\Delta^M_{\pi_{M},cusp,G-comp}],$$
  o\`u $ \Theta_{\pi\vert G_{comp}(F)}$ est la restriction de $\Theta_{\pi}$ \`a $G_{comp}(F)$. 
  En utilisant (4), on obtient
  $$ D^G[d]= \Theta_{\pi\vert G_{comp}(F)}-\sum_{M\in \underline{{\cal L}}_{min},M\not=G}\sum_{M'\in \underline{{\cal L}}_{min,G-rel}^{G'}; M'\mapsto M}\vert  W^{G'}(M')\vert ^{-1} transfert(D^{G'_{1}}[\Delta^{M'_{1}}_{\pi'_{M'},cusp,G'_{1}-comp}])$$
  $$=\Theta_{\pi\vert G_{comp}(F)}-\sum_{M'\in \underline{{\cal L}}_{min,G-rel}^{G'}; M'\not=G'}\vert  W^{G'}(M')\vert ^{-1} transfert(D^{G'_{1}}[\Delta^{M'_{1}}_{\pi'_{M'},cusp,G'_{1}-comp}]).$$
  On peut remplacer l'ensemble de sommation $ \underline{{\cal L}}_{min,G-rel}^{G'}$ par $ \underline{{\cal L}}_{min}^{G'}$: pour un Levi non relevant, le terme que l'on somme est nul.
  En utilisant de nouveau \ref{leresultat} (2) cette fois dans $G'_{1}$, on obtient
  $$ D^G[d]= \Theta_{\pi\vert G_{comp}(F)}+transfert(D^{G'_{1}}[d']-\Theta_{\pi'\vert G'_{1,comp}(F)}).$$
  On a d\'ej\`a dit que, gr\^ace \`a Arthur, on avait $\pi=transfert(\pi')$ donc les termes $\Theta_{\pi\vert G_{comp}(F)}$ et $transfert(\Theta_{\pi'\vert G'_{1,comp}(F)})$ disparaissent. D'o\`u
  $$ D^G[d]= transfert(D^{G'_{1}}[d']),$$
  ce qu'il fallait d\'emontrer. $\square$

jean-loup.waldspurger@imj-prg.fr

CNRS- Institut de Math\'ematiques de Jussieu- Paris rive gauche

4 place Jussieu

Bo\^{\i}te courrier 247

75252 Paris Cedex 05


\begin{thebibliography}{r}

\bibitem{A1} J. Arthur: {\it The local behaviour of weighted orbital integrals}, Duke Math. J. 56 (1988), p. 223-293

\bibitem{A2} J. Arthur: {\it On elliptic tempered characters}, Acta Math. 171 (1993), p. 73-138



\bibitem{A3} J. Arthur: {\it On local character relations}, Selecta Math. 2 (1996), p. 501-579


\bibitem{C}  W. Casselman: {\it Characters and Jacquet modules}, Math. Ann. 230 (1977), p. 101-105

\bibitem{D} S. DeBacker: {\it Homogeneity results for invariant distributions of a reductive $p$-adic groups}, Ann. Sc. Ec. Norm. Sup. 35 (2002), p. 391-422

\bibitem{DR} S. DeBacker, M.  Reeder: {\it Depth-zero supercuspidal $L$-packets and their stability}, Annals of Mah. 169 (2009), p. 795-901

\bibitem {F} A. Ferrari: {\it Th\'eor\`eme de l'indice et formule des traces}, manuscripta math. 124 (2007), p. 363-390


\bibitem{HR} T. Haines, M. Rapoport: {\it On parahoric subgroups}, appendice \`a G.Pappas, M. Rapoport {\it Twisted loop groups and their affine flag varieties}, Adv. in Math. 219 (2008), p. 118-198

\bibitem{HV} G. Henniart, M.-F. Vign\'eras: {\it A Satake isomorphism for representations mod $p$ of reductive groups over local fields}, Journal f\" ur die r. und ang. Math. 701 (2015), p. 33-75

\bibitem{Lanard1} T. Lanard: {\it Sur les $\ell$-blocs de niveau z\'ero des groupes $p$-adiques}, Compositio Math. 154 (2018), p. 1473-1507

\bibitem{Lanard2} T. Lanard: {\it Sur les $\ell$-blocs de niveau z\'ero des groupes $p$-adiques II}, arXiv RT 1806.09543


\bibitem{L} R. P. Langlands: {\it Stable conjugacy: definitions and lemmas}, Can. J. Math. XXXI (1979), p. 700-725

\bibitem{LS} R. P. Langlands, D. Shelstad: {\it On the definition of transfer factors}, Math. Ann. 278 (1987), p. 219-271



\bibitem{MW} C. Moeglin, J.-L. Waldspurger: {\it Stabilisation de la formule des traces tordue}, volume 1, Progress in Math. 316, Birkh\" auser 2016

\bibitem {Oi} M. Oi: {\it Depth preserving property of the local Langlands correspondence for quasi-split classical groups in a large residual characteristic}, arXiv NT 1804.10901
 
 \bibitem{P} G. Prasad: {\it Finite group actions on reductive groups and buildings and tamely-ramified descent in Bruhat-Tits theory}, arXiv RT 1705.02906
 
\bibitem{PY} G. Prasad, J.-K. Yu: {\it On finite groups actions on reductive groups and buildings}, Invent. Math. 147 (2002), p. 545-560

\bibitem{W1} J.-L. Waldspurger: {\it Une formule int\'egrale reli\'ee \`a la conjecture locale de Gross-Prasad}, Compositio Math. 146 (2010), p. 1180-1290

 


\bibitem{W3} J.-L. Waldspurger: {\it Une formule int\'egrale reli\'ee \`a la conjecture locale de Gross-Prasad, $2^{e}$ partie: extension aux repr\'esentations temp\'er\'ees}, Ast\'erisque 346 (2012), p. 171-312

\bibitem{W2} J.-L. Waldspurger: {\it Caract\`eres de repr\'esentations de niveau $0$}, \`a para\^{\i}tre aux Ann. de la Fac. Sc. de Toulouse

\bibitem{W4}  J.-L. Waldspurger: {\it  Le lemme fondamental implique le transfert}, Compositio Math. 105 (1997), p. 153-236

\end{thebibliography}
 \end{document}